\documentclass{book}
\usepackage[centertags]{amsmath}
\usepackage{amsfonts}
\usepackage{amssymb}
\usepackage{amsthm}
\usepackage{array}

\newlength{\singlebaselineskip}
\newfont{\fib}{cmfib8}

\def\@begintheorem#1#2{
  \par
  \bgroup{
    \bf #1\ #2.
  }
  \it
  \ignorespaces
}
\def\@opargbegintheorem#1#2#3{
  \par
  \bgroup{
    \bf #1\ #2\ (#3).
  }
  \it
  \ignorespaces
}
\def\@endtheorem{\egroup}

\newtheorem{thm}{Theorem}[section]
\newtheorem{cor}[thm]{Corollary}
\newtheorem{lem}[thm]{Lemma}
\newtheorem{prop}[thm]{Proposition}
\newtheorem{defn}{Definition}[section]

\newcommand{\s}{\scriptstyle}
\newcommand{\g}{\mathfrak{g}}
\newcommand{\h}{\mathfrak{h}}
\newcommand{\gl}{\mathfrak{gl}}
\newcommand{\si}{\mathfrak{sp}}
\newcommand{\so}{\mathfrak{so}}
\newcommand{\sll}{\mathfrak{sl}}
\newcommand{\Hom}{\operatorname{Hom}}
\newcommand{\R}{\mathbb{R}}

\newcommand{\C}{\mathbb{C}}
\newcommand{\N}{\mathbb{N}}
\newcommand{\F}{\mathbb{F}}
\newcommand{\pnth}[1]{\left( #1 \right)}
\newcommand{\ad}[1]{\operatorname{ad}\left( #1 \right)}
\newcommand{\add}{\operatorname{ad}}
\newcommand{\e}{\mathrm{e}}
\newcommand{\f}{\mathrm{f}}
\newcommand{\X}{\mathrm{X}}
\newcommand{\Y}{\mathrm{Y}}
\newcommand{\A}{\mathbf{A}}
\newcommand{\D}{\mathbf{D}}

\newcommand{\abs}[1]{\left\vert #1 \right\vert}
\newcommand{\To}{\rightarrow}
\newcommand{\la}{\lambda}

\newcommand{\be}{\beta}
\newcommand{\ga}{\gamma}
\newcommand{\om}{\omega}
\newcommand{\vep}{\varepsilon}
\newcommand{\del}{\delta}

\newcommand{\brkt}[1]{\left[ #1 \right]}

\newcolumntype{L}{>{$\s} l <{$}}
\newcolumntype{K}{>{\scriptsize} l   }
\newcolumntype{C}{>{$\s } c <{$}}

\begin{document}

\newpage
\begin{titlepage}

  \addtocounter{page}{1}
  \mbox{}

  \begin{center}
     \begin{tabular}{c}

  %

        A CLASSIFICATION OF REAL INDECOMPOSABLE SOLVABLE LIE ALGEBRAS \\

  \ \\    
          OF SMALL DIMENSION WITH CODIMENSION ONE NILRADICALS \\[.25in]


                                  by \\[.25in]

                               Alan R. Parry \\[.5in]



                 A thesis submitted in partial fulfillment \\
                 of the requirements for the degree \\[.125in]

                                  of \\[.125in]

                          MASTER OF SCIENCE \\[.125in]

                                   in \\[.125in]


                             Mathematics \\[.5in]

      \end{tabular}
   \end{center}

  \begin{tabbing}
  \mbox{\hspace{95mm}}\= \kill

    Approved: \\[.25in]

    \rule{60mm}{.1mm}            \>        \rule{60mm}{.1mm}\\
    Ian Anderson               \>          Charles Torre \\
    Major Professor              \>        Committee Member \\[.2in]

    \rule{60mm}{.1mm}            \>        \rule{60mm}{.1mm}\\
    David Brown                   \>       Byron Burnham \\
    Committee Member             \>        Dean of Graduate Studies\\[.5in]

  \end{tabbing}

  \begin{center}
    \begin{tabular}{c}
                         UTAH STATE UNIVERSITY \\
                              Logan, Utah \\[.125in]
                                  2007
    \end{tabular}
  \end{center}

\vfill

\end{titlepage}


\newpage
\pagenumbering{roman}
\setcounter{page}{2}
%
%
%

\begin{center}

  \null \vskip 260pt
        Copyright \copyright \hspace{10pt} Alan R. Parry \hspace{5pt} 2007 \\
  \vskip 14pt All Rights Reserved\\

\end{center}

\vfill\eject


\newpage
%
%
%

\begin{center}
  ABSTRACT
\end{center}

\addcontentsline{toc}{chapter}{ABSTRACT}{

  \setlength{
    \baselineskip}{8mm}

    \begin{center}
      \begin{tabular}{c}
        \\[.05mm]

  %
  %

    {A Classification of Real Indecomposable Solvable Lie Algebras}
    \\

  \ \\    

       {of Small Dimension with Codimension One Nilradicals}


     \\[8mm]  

                              by
     \\[8mm]  

  %
                 Alan R. Parry, Master of Science

     \\[4mm]  

                  Utah State University, 2007

     \\[6mm]  

      \end{tabular}
    \end{center}

    \noindent
    Major Professor: Dr. Ian M. Anderson \\[-4mm]
    Department: Mathematics and Statistics

    \vskip 11pt

%
%
\indent
This thesis was concerned with classifying the real indecomposable
solvable Lie algebras with codimension one nilradicals of dimensions
two through seven.  This thesis was organized into three chapters.

In the first, we described the necessary concepts and definitions
about Lie algebras as well as a few helpful theorems that are
necessary to understand the project.  We also reviewed many concepts
from linear algebra that are essential to the research.

The second chapter was occupied with a description of how we went about classifying the Lie algebras.  In
particular, it outlined the basic premise of the classification: that we can use the automorphisms of the
nilradical of the Lie algebra to find a basis with the simplest structure equations possible.  In addition, it
outlined a few other methods that also helped find this basis.  Finally, this chapter included a discussion of the
canonical forms of certain types of matrices that arose in the project.

The third chapter presented a sample of the classification of the
seven dimensional Lie algebras.  In it, we proceeded step-by-step
through the classification of the Lie algebras whose nilradical was
one of four specifically chosen because they were representative of
the different types that arose during the project.

In the appendices, we presented our results in a list of the
multiplication tables of the isomorphism classes found.



\par

\setlength{\baselineskip}{\singlebaselineskip}}




\newpage
%
%

\begin{center}
   ACKNOWLEDGMENTS
\end{center}

\addcontentsline{toc}{chapter}{ACKNOWLEDGMENTS}

{ \setlength{\baselineskip}{8mm} \sloppy


I would like to thank Dr. Ian M. Anderson, my major professor,
without whose immense help and expertise this thesis would certainly
have been impossible, not to mention my success as a student at USU.

I would also like to thank Dr. Mark Fels, who has always been a source of encouragement no matter how many times
he had to have me in his classes.

In addition, I give my thanks to Dr. David Brown, who was always good for a laugh and was of great support to me
in my collegiate career; to him, I offer the universal greeting, "Bah weep grawna weep mini bom!"

I also want to give a special thanks to Dr. Charles Torre, who was
willing to be a member of my advisory committee on such short
notice.

Most especially, I thank my beautiful wife Alesha and my daughter Madison, whose support was unwavering and love
unconditional.  You were always there when I didn't think I'd make it.

Financial support from Ian Anderson's NSF grant DMF 0410373 is gratefully acknowledged.

\hfill Alan R. Parry

\par
\setlength{\baselineskip}{\singlebaselineskip}
}




\tableofcontents








\newpage
\pagenumbering{arabic}
%
%
%
%
%
%
%
\setlength{\baselineskip}{22pt}

%
%

\chapter{A REVIEW OF FREQUENTLY USED CONCEPTS}\label{c1}

\thispagestyle{empty}

This chapter will deal with reviewing those concepts and definitions
that are most helpful in discussing this classification of solvable
Lie algebras.  We will begin our discussion with a few basic
definitions and concepts of Lie theory.  This will be followed by
the major linear algebra concepts needed for our classification.

\section{Lie Algebras: Definitions and Concepts}

We first define a Lie algebra and the various types of Lie algebra
homomorphisms.
\begin{defn} A {\em Lie algebra} is a vector space, $\g$, over a field, $\F$, coupled
with a mapping $$[,]:\g \times \g \To \g$$ that is
\begin{enumerate}
  \item bilinear $\pnth{[cX+Y,Z]=c[X,Z]+[Y,Z]; \ \ [X,cY+Z]=c[X,Y]+[X,Z]}$,
  \item skew-symmetric $\pnth{[X,Y]=-[Y,X]}$,
  \item and satisfies the Jacobi property $$\brkt{X,[Y,Z]}+[Y,[Z,X]]+[Z,[X,Y]]=0$$
\end{enumerate}
for all $c \in \F$ and $X,Y,Z \in \g$.
\end{defn}

The map $[,]$ is called the Lie bracket.

\begin{defn} If $\g$ and $\h$ are Lie algebras, then a linear
transformation $T: \g \To \h$ is a {\em Lie algebra homomorphism} if
for all $X,Y \in \g$ we have that $$T([X,Y]) = [T(X),T(Y)].$$ If $\g
= \h$, then $T$ is a {\em Lie algebra endomorphism}.  If $T$ is
bijective, then $T$ is a {\em Lie algebra isomorphism}.  If $T$ is
both a Lie algebra endomorphism and an isomorphism, we call $T$ a
{\em Lie algebra automorphism}. \end{defn}

There are a number of different properties of Lie algebras.  In
fact, the classification of Lie algebras into isomorphism classes is
done by finding canonical forms for algebras with certain properties
that are preserved by isomorphism.  As such, it is necessary to
discuss a few of these properties here.

\begin{defn} The {\em derived algebra}, denoted $D\g$, is the set
$$\{X \in \g \ \left| \ \mathrm{for \ some} \ Y,Z \in \g, \ X=[Y,Z] \right. \}.$$
\end{defn}
An alternative way to describe the derived algebra is
$$D\g=[\g,\g],$$ where $[\g,\g]$ denotes all possible Lie brackets
between vectors in $\g$. This gives rise to another concept called
the {\em derived series}, given by a series of $D^{i}\g$, for all $i
\in \N$, where each $D^{i}\g$ is defined inductively by
$$D^{i}\g=[D^{i-1}\g, D^{i-1}\g] \ \ \mathrm{with} \ \
D^{1}\g=D\g=[\g,\g].$$

Another useful series is the {\em lower central series}, given by a
series of $D_{i}\g$, for all $i \in \N$, where each $D_{i}\g$ is
defined by $$D_{i}\g=[D_{i-1}\g, \g] \ \ \mathrm{with} \ \
D_{1}\g=D\g=[\g,\g].$$ Note that for all $i \in \N$, $D^{i}\g
\subseteq D_{i}\g$.

The dimensions of each $D^{i}\g$ and $D_{i}\g$ is invariant under a
Lie algebra isomorphism.

These facts give us enough information to discuss the difference
between two special types of Lie algebras.  So we give two
definitions.  When in context, let $0$ denote the zero vector.

\begin{defn} A Lie algebra is {\em solvable} if for some $i \in \N$,
$D^{i}\g = \{ 0 \}$. \end{defn}

\begin{defn} A Lie algebra is {\em nilpotent} if for some $i \in \N$,
$D_{i}\g = \{ 0 \}$. \end{defn}

It is clear then that because $D^{i}\g \subseteq D_{i}\g$ for all $i
\in \N$, all nilpotent Lie algebras are also solvable.  In addition,
as the dimensions of $D^{i}\g$ and $D_{i}\g$ are invariant under a
Lie algebra isomorphism, it follows that whether an algebra is
solvable or nilpotent is invariant as well.

There is another property that is important to note as it plays a
large role in the type of algebras we classify in this paper.  So we
offer the following definition.

\begin{defn} A Lie algebra, $\g$, is said to be {\em decomposable} if
there exist lower dimensional algebras $\g_{1}$ and $\g_{2}$ such
that $$\g=\g_{1} \oplus \g_{2}$$ where $\oplus$ denotes a Lie
algebra direct sum.  A Lie algebra is said to be {\em
indecomposable} if it is not decomposable. \end{defn}

Another useful tool is the concept of the nilradical of a solvable
Lie algebra.

\begin{defn} The {\em nilradical} of a Lie algebra, $\g$, is its
maximal nilpotent ideal.  We'll denote it $NR(\g)$.  The {\em
codimension} of the nilradical of a Lie algebra is the difference
between the dimension of the entire Lie algebra and the dimension of
the nilradical, that is $\dim \g - \dim NR(\g)$.
\end{defn}

The nilradical is similar to the radical of a Lie algebra, which is
defined to be the maximal solvable ideal.  The radical is useful in
classifying non-solvable Lie algebras, as the Levi decomposition
states that every Lie algebra can be written as the semi-direct
product of its radical and a semisimple Lie algebra. However, in a
solvable Lie algebra, the radical is obviously the entire algebra.
So we use the nilradical instead. The nilradical is unique in a Lie
algebra and if two algebras are isomorphic, then their nilradicals
are isomorphic as well. Thus it becomes the perfect object by which
to classify the algebra.

Moreover, there is a useful theorem that states that the derived
algebra, $D\g$, of a solvable Lie algebra, $\g$, is contained in the
nilradical of $\g$.  That is, $D\g \subseteq NR(\g)$
\cite{varadarajan}.

In this classification, we classify the indecomposable solvable Lie
algebras over $\R$ of dimensions two through seven with codimension
one nilradicals. So we have most of the necessary definitions and
concepts. However, there are a few more useful tools that we'll use
in the classification that we'll describe here.

The first is a derivation.

\begin{defn} A linear transformation $T: \g \To \g$ is a {\em
derivation} if it satisfies the Leibniz rule.  That is for all $X,Y
\in \g$, we have $$T([X,Y]) = [T(X),Y] + [X,T(Y)].$$\end{defn}

The space of derivations form a Lie algebra with the Lie bracket
given by the commutator $[T,U] = T \circ U - U \circ T$.  The
exponential of a derivation, $e^{T}$ is an automorphism of $\g$.

The next tool is a representation of a Lie algebra.

\begin{defn} A {\em Lie algebra representation} is a Lie algebra homomorphism
$\rho: \g \to \gl(V)$, where $\gl(V)$ denotes the Lie algebra of all
linear transformations of a vector space $V$ whose Lie bracket is
given by the commutator. \end{defn}

A special Lie algebra representation is the adjoint or $\add$
representation of a Lie algebra $\g$.

\begin{defn} The {\em $\add$ representation} of a Lie algebra, $\g$,
is a Lie algebra homomorphism $\add: \g \To \gl(\g)$ such that for
all $X,Y \in \g$, $$\ad{X}(Y)=[X,Y].$$
\end{defn}

We move on now to review a few Linear algebra concepts.

\section{Linear Algebra: Change of Basis Matrices and Real
Jordan Canonical Form}

Here we give a discussion about invertible Linear transformations on
a Vector space and how this idea extends to Lie Algebra
isomorphisms.  Most importantly though, we'll talk about how these
invertible operators can be used to ``move'' a matrix into Jordan
Canonical Form.  We begin by defining the matrix representation of a
transformation.

\begin{defn}
    Let $V$ and $W$ be finite dimensional vector spaces over a
    field $\F$ with bases $\beta=\{v_{1},\ldots,v_{n}\}$ and
    $\gamma=\{w_{1},\ldots,w_{m}\}$ respectively.  Let $T:V \to W$
    be a linear transformation, then we could express $T$ evaluated
    on every basis vector $v_{i}$ in the following way
    $$T(v_{i})=\sum_{j=1}^{m} c_{i}^{j} w_{j}, \ \ \ \mathrm{for} \
    1\leq i \leq n,$$ where $c_{i}^{j} \in \F$.  Then the {\em matrix
    representation} of $T$ denoted $[T]_{\beta}^{\gamma}$ is the
    $m \times n$ matrix given by
    $$\pnth{[T]_{\beta}^{\gamma}}_{i,j}=c_{i}^{j}.$$
\end{defn}

We also note that if $v \in V$, then $v=\sum_{i=1}^{n} a_{i} v_{i}$
and so we could define the $m \times 1$ matrix $[v]_{\beta}$ by
$\pnth{[v]_{\beta}}_{i}=a_{i}.$  This will give us the property that
$$[T]_{\beta}^{\gamma}[v]_{\beta}=[T(v)]_{\gamma}.$$  This is how
the matrix representation is useful.  It simplifies evaluation of a
linear transformation on a vector to matrix multiplication.

Now we can define an invertible linear transformation and see how
this relates to its matrix representation.

\begin{defn}
    Let $T:V \To W$ and $U:W \To V$ be linear transformations.  Let
    $I_{V}:V \To V$ and $I_{W}:W \To W$ be the identity
    transformations on $V$ and $W$ respectively.  If $T \circ
    U=I_{V}$ and $U \circ T=I_{W}$, then $T$ is {\em invertible} and
    $T^{-1}=U$.
\end{defn}
As we've seen that evaluation of a transformation $T$ can be reduced
to left multiplication of $[T]_{\gamma}^{\beta}$, then we can see
that $[T \circ U]_{\beta}^{\beta}=[T]_{\gamma}^{\beta}
[U]_{\beta}^{\gamma}$. Because of this, we have that $T$ is an
invertible linear transformation if and only if
$[T]_{\gamma}^{\beta}$ is an invertible matrix.  Which reduces to
$\det[T]_{\gamma}^{\beta} \neq 0$.  This shows quite easily then
that if $T:V \to W$ invertible, then $\dim V= \dim W$. Of course a
linear transformation is invertible if and only if it is bijective.
Hence any Lie algebra isomorphism is an invertible linear
transformation.  It should be noted that an invertible $n \times n$
matrix is an element of the group $GL(n,\F)$ where $\F$ is the field
of scalars.  We will usually denote such a matrix in this manner.

Next we define a change of basis matrix.  These are important because the matrix representation of any Lie algebra
automorphism can be interpreted as a change of basis matrix.

\begin{defn}
    Let $V$ be a finite-dimensional vector space and let $\beta$ and
    $\beta'$ be two bases for $V$.  Let $Q=
    [I_{V}]_{\beta'}^{\beta}$, then we call $Q$ the {\em change
    of basis matrix} that changes $\beta'$-coordinates into
    $\beta$-coordinates.
\end{defn}

Let $T:V \To V$ be a linear transformation.  Then if $Q$ is the change of basis matrix that changes
$\beta'$-coordinates into $\beta$-coordinates, we have that $$[T]_{\beta'}^{\beta'} = Q^{-1} [T]_{\beta}^{\beta}
Q.$$  This is because $$Q[T]_{\beta'}^{\beta'} = [I_{V}]_{\beta'}^{\beta} [T]_{\beta'}^{\beta'} =
[I_{V}T]_{\beta'}^{\beta} =[T I_{V}]_{\beta'}^{\beta} = [T]_{\beta}^{\beta}[I_{V}]_{\beta'}^{\beta} =
[T]_{\beta}^{\beta} Q.$$  Hence, we can see how a changing the basis on a vector space affects the matrix
representation of a linear transformation. Analogously, we also can see how changing the basis of a Lie algebra
affects the matrix representation of a Lie algebra endomorphism, which is a concept we use almost constantly
throughout the computations in this paper.  It also plays a key role in the discussion of Jordan canonical form,
which we'll review next.

The reader is cited to any standard Linear Algebra text for a more
in depth study of Jordan canonical form including the proofs of the
following statements.  The real Jordan canonical form is explained
in detail in the book by Hirsch and Smale \cite{hirsch}.

Let $T:V \To V$ be a linear transformation on a vector space $V$ over a field $\F$.

\begin{defn}
    An {\em eigenvector} of $T$ is a vector, $v \in V$, such that $(T-\la I)(v) = 0$ for some $\la \in \F$.  The
    value $\la$ is called the {\em eigenvalue} corresponding to the eigenvector $v$. A {\em generalized
    eigenvector}, corresponding to the eigenvalue $\la \in \F$, of $T$ is a vector, $v' \in V$, such that
    $(T-\la I)^{p}(v) = 0$ for some positive integer $p$.
\end{defn}

\begin{defn}
    The subspace of $V$ $$\{v \in V \ | \ (T-\la I)^{p}(v) = 0 \ \mathrm{for \ some \ positive \ integer} \ p\}$$
    is called the {\em generalized eigenspace} of $T$ corresponding to $\la$ and is denoted $K_{\la}$.  If
    $K_{\la}$ consists only of eigenvectors of $T$, then $K_{\la}$ is simply called the {\em eigenspace} of $T$
    corresponding to $\la$.
\end{defn}

\begin{defn}
    The polynomial $q(\la) = \det(T-\la I)$ is called the {\em characteristic polynomial} of $T$ and has degree
    equal to the dimension of $V$.
\end{defn}
The roots of the characteristic polynomial are the eigenvalues of $T$.  In addition, the multiplicity of any root
of the characteristic polynomial is equal to the dimension of its corresponding generalized eigenspace.

\begin{defn}
    Let $v$ be a generalized eigenvector of $T$.  If $p$ is the smallest integer such that $(T-\la I)^{p}(v) = 0$,
    then the set $\{(T-\la I)^{p-1}(v), (T-\la I)^{p-2}(v), \ldots, (T-\la I)(v), v\}$ is called a {\em cycle of
    generalized eigenvectors} corresponding to $\la$ of length $p$.  The vectors $(T-\la I)^{p-1}(v)$ and $v$ are
    called the {\em initial vector} and {\em end vector} of the cycle respectively.
\end{defn}
Clearly the initial vector of any cycle is an eigenvector of $T$.

\begin{thm}
    If $\be = \{v_{1},\ldots,v_{p}\}$ is a cycle of generalized eigenvectors corresponding to $\la \in \F$ where
    $v_{i}=(T-\la I)^{p-i}(v_{p})$.  Then $T(v_{1})=\la v_{1}$ and for all $i$ such that $2 \leq i \leq p$,
    $$T(v_{i}) = \la v_{i} + v_{i-1}.$$
\end{thm}

{\em Proof.} We know that $v_{1}$ is an eigenvector and so $(T-\la I)(v_{1}) = 0$ or equivalently $T(v_{1}) = \la
v_{1}$.  Also $$(T-\la I)(v_{i}) = (T-\la I)^{p-(i-1)}(v_{p}) = v_{i-1}.$$  The result follows immediately.
$\blacksquare$ \vspace{.1 in}

Every cycle of generalized eigenvectors is a linearly independent
set.  Moreover, if the initial vectors of the cycles are linearly
independent, then the union of any number of cycles is linearly
independent. Furthermore, for any $K_{\la}$, there exists a basis
for $K_{\la}$ consisting of disjoint cycles of generalized
eigenvectors.

Thus if the characteristic polynomial splits over $\F$, then $V$ is
the direct sum of generalized eigenspaces and has a basis, $\be$,
consisting of disjoint cycles of generalized eigenvectors.  In this
basis, $T$ is in what is called {\em Jordan canonical form}.  The
Jordan canonical form of a transformation $T$ is unique up to the
ordering of its eigenvalues. If $A$ is the matrix representation of
$T$ in a given basis, $\ga$, then the change of basis matrix, $Q$,
that changes $\ga$ coordinates into $\be$ coordinates is such that
$$Q^{-1} A Q = [T]_{\be}^{\be}.$$  Thus the Jordan canonical form of the matrix representation, $A$, of
a linear transformation, $T$, if it exists, is a matrix of the form
$$\pnth{\begin{array}{ccccccccc} \la_{1} & a_{1,1} & & & & & \\ &
\la_{1} &a_{1,2} &  & & & \\ & & \ddots
& & & \\ & & & \la_{1} & & & \\ & & & & \la_{2} & a_{2,1} & \\ & & & & & \ddots & \\
& & & & & & \la_{n} \end{array}}$$ where the $a_{i,j} \in \{0,1\}$. In addition, the $\la_{i}$ are the eigenvalues
of $A$.  The square block belonging to each $\la_{i}$ corresponds to a generalized eigenspace of $A$ and the block
is called a Jordan block.

For example, let $A$ be a $6 \times 6$ matrix whose characteristic
polynomial splits over $\R$.  If, for some invertible $Q$,
$$Q^{-1} A Q = \pnth{ \begin{array}{cccccc} 2&0&0&0&0&0 \\
0&2&0&0&0&0 \\ 0&0&3&0&0&0 \\ 0&0&0&5&1&0 \\ 0&0&0&0&5&1 \\
0&0&0&0&0&5 \end{array}},$$ then this is a Jordan canonical form of
$A$.  The eigenvalues of $A$ are 2, 3, and 5.  The Jordan blocks are
$$\pnth{ \begin{array}{cc} 2&0 \\ 0&2 \end{array}}, \ \ \
\pnth{\begin{array}{c} 3 \end{array}}, \ \ \
\pnth{\begin{array}{ccc} 5&1&0 \\ 0&5&1 \\ 0&0&5 \end{array}},$$ and
each of these blocks correspond to a generalized eigenspace of $A$.

We are considering real Lie algebras which are vector spaces over
$\R$. However, not all polynomials split over $\R$. It is true,
though, that as $\C$ is algebraically closed, every polynomial with
real coefficients splits over $\C$.  Thus the complexification of
$V$, $V_{c}$, can always be written as the direct sum of generalized
eigenspaces and the complex Jordan canonical form always exists for
such a matrix. But it is important to note that as any real linear
transformation $T$ will have a characteristic polynomial with real
coefficients, any complex root will come in a conjugate pair. Hence
if $z$ is a complex eigenvalue of a linear transformation $T$, then
$\overline{z}$ is as well. In fact, if $v$ is an eigenvector of $A$
corresponding to $z$, then $\overline{v}$ is an eigenvector
corresponding to $\overline{z}$.

Let $z=a+bi$ be an eigenvalue of $T$ corresponding to an eigenvector
$v_{1}$.  Then $v_{2}=\overline{v_{1}}$ is an eigenvector
corresponding to $a-bi$ and a submatrix of the Jordan canonical form
of the matrix representation of $T$ is $$\pnth{\begin{array}{cc}
a+bi &0 \\ 0&a-bi
\end{array}}.$$  Let $v_{1}'=v_{1}+v_{2}$ and
$v_{2}'=-i(v_{1}-v_{2})$.  Then $v_{1}'$ and $v_{2}'$ are real vectors that form a basis for the space spanned by
the eigenspaces of $v_{1}$ and $v_{2}$.  In this basis the submatrix block above becomes $$\pnth{\begin{array}{cc}
a&b \\ -b&a \end{array}}.$$ If a cycle of generalized eigenvectors of length two corresponds to a complex
eigenvalue, then we have a Jordan block of the form $$\pnth{\begin{array}{cccc} a+bi&1&0&0 \\ 0&a+bi&0&0 \\ 0&0&a-bi&1 \\
0&0&0&a-bi \end{array}}.$$ Then using a similar trick as described above but on both the eigenvectors and the
generalized eigenvectors, this block will become $$\pnth{ \begin{array}{cccc} a&b&1&0 \\ -b&a&0&1 \\ 0&0&a&b \\
0&0&-b&a \end{array}}.$$  This can be easily extended to a generalized eigenspace of any dimension.

If we use this to deal with a characteristic polynomial with complex roots, then the resulting form of $T$ is
called the {\em real Jordan canonical form} and the basis for $V_{c}$ that puts $T$ into real Jordan canonical
form is also a basis for $V$. Thus we can still consider $V$ as a vector space over $\R$ and every linear
transformation $T:V \To V$ has a real Jordan canonical form.

This completes our review of frequently used concepts.  The reader
should have the necessary information to understand the mechanics
behind the paper.  In the next chapter, we'll discuss the methods we
use in classifying Lie algebras and also the canonical forms of a
few special types of matrices.

%
\setlength{\baselineskip}{11pt}
%

%
\setlength{\baselineskip}{22pt}

%
%

\chapter{A DISCUSSION OF METHODS USED AND CANONICAL FORMS}\label{c2}

This chapter is devoted to discussing in detail methods that will be
utilized frequently throughout the classification.  As they are
discussed here, when they are used in the text the reader will be
referred back to this chapter for a more detailed explanation.  In
this chapter, we will also be building off the topics reviewed in
the previous chapter.

\section{Classifying Solvable Lie Algebras with Codimension One Abelian Nilradicals}

The first of these methods that we will discuss is that of how we'll
classify a solvable Lie algebra with a codimension one abelian
nilradical.

To begin, let $\g$ be an $n$-dimensional solvable Lie algebra with
such a nilradical, $NR(\g)$.  Now choose a basis for $\g$,
$\beta=\{\e_{1},\e_{2},\ldots,\e_{n}\}$, such that the set
$\gamma=\{\e_{1},\e_{2},\ldots, \e_{n-1}\}$ forms a basis for the
nilradical.  As $\g$ is solvable, we know that $D\g \subseteq
NR(\g)$, and this yields $$[\e_{i},\e_{j}]=\sum_{k=1}^{n-1}
A_{i,j}^{k} \e_{k}$$ for all $1 \leq i,j \leq n$. Also as $NR(\g)$
is abelian we have that $[\e_{i},\e_{j}]=0$ for all $1 \leq i,j \leq
n-1$. Combining these two facts, we have that the structure
equations of $\g$ are simply
$$[\e_{i},\e_{n}]=\sum_{k=1}^{n-1} A_{i,n}^{k} \e_{k}$$
for all $1 \leq i \leq n-1$. And so it would only make sense to
direct our discussion towards $\ad{\e_{n}}$.  From the structure
equations we have that $$\ad{\e_{n}}_{\beta}^{\beta} =
\begin{pmatrix}
 -A_{1,n}^{1}&-A_{2,n}^{1}& \cdots & -A_{n-1,n}^{1} & 0 \\
-A_{1,n}^{2}& -A_{2,n}^{2}& \cdots &  -A_{n-1,n}^{2} & 0 \\ \vdots & \vdots & \ddots & \vdots & \vdots \\
-A_{1,n}^{n-1} & -A_{2,n}^{n-1} & \cdots & -A_{n-1,n}^{n-1} & 0 \\ 0
&0 & \cdots & 0 &0 \end{pmatrix}$$  It is clear that we lose no
information of the structure constants if we restrict this
transformation to act only on the nilradical.  Thus we will consider
$\ad{\e_{n}}$ in the following way
$$\left.\ad{\e_{n}}\right|_{NR(\g)} = \begin{pmatrix}
-A_{1,n}^{1}&-A_{2,n}^{1}& \cdots & -A_{n-1,n}^{1} \\
-A_{1,n}^{2}& -A_{2,n}^{2}& \cdots &  -A_{n-1,n}^{2} \\ \vdots &
\vdots & \ddots & \vdots \\ -A_{1,n}^{n-1} & -A_{2,n}^{n-1} & \cdots
& -A_{n-1,n}^{n-1} \end{pmatrix}$$ but we'll refer to it as
$\ad{\e_{n}}$.

Now we note that as $NR(\g)$ is abelian ($[a,b]=0$ for all $a,b \in
NR(\g)$), any change of basis strictly on the vectors in the
nilradical will not change the structure constants of the
nilradical. To compute the matrix representation of $\ad{\e_{n}}$
after we apply such a change of basis, we first construct a change
of basis matrix $M$ and conjugate $\ad{\e_{n}}$ by it in the
following way $$M^{-1} \ad{\e_{n}} M.$$  We can pick $M \in
GL(NR(\g))$ arbitrarily so this expression can move $\ad{\e_{n}}$
into real Jordan canonical form.

This can always be done when the nilradical is abelian and codimension one, and so this method is used quite often
in our classification.  It would also be helpful in our classification discussion to be familiar with the possible
real Jordan canonical forms for each dimension through dimension six.  We will discuss this a little later.

We should also note that we can make one more change that won't affect the structure equations of the nilradical.
This basis change is simply exchanging $\f_{6}$ by a constant multiple of itself.  This change won't affect the
nilradical basis at all, but can change values in the $\ad{\f_{6}}$ matrix.  So this is usually the type of basis
change we'll make at the end of each classification.  However, for some real Jordan canonical forms, we add into
the basis change a scaling of one or more of the nilradical basis vectors.  This is done to keep certain values
the same while scaling $\f_{6}$.  For instance, when we have a Jordan block that has any number of ones on the
super diagonal.  These added scalings are simply reapplying one or more of the automorphisms of the nilradical.
This doesn't mean much in this case where any invertible linear transformation is an automorphism of the
nilradical, but it will become more relevant as we move further through our classification. So, to keep from being
confusing in our use of the automorphisms, we'll simply make the change when making our final change of basis.

\section{Classifying Solvable Lie Algebras with Codimension One Non-abelian Nilradicals}

This section is a natural extension of the preceding one. We will
see many of the same ideas used, as well as extensions to them.

We are given a Lie algebra $\g$ with a codimension one non-abelian
nilradical, where the first $n-1$ vectors in a basis for $\g$
constitute a basis for $NR(\g)$.  We will always apply the necessary
isomorphism to $NR(\g)$ to move it into a canonical form.  We
utilize known texts for a classification of nilpotent Lie algebras
of the required dimensions \cite{winternitz-1,gong}.

As $D\g \subseteq NR(\g)$, we can still view $\ad{\e_{n}}$
restricted to the nilradical and not lose any information on the
structure equations involving $\e_{n}$. Arbitrarily, this is given
by $$\ad{\e_{n}}=
\pnth{\begin{array}{cccc} A^{1}_{1} &A^{1}_{2}&\cdots&A^{1}_{n-1} \\
A^{2}_{1} &A^{2}_{2} &\cdots& A^{2}_{n-1} \\ \vdots & \vdots& \ddots
& \vdots \\ A^{n-1}_{1} & A^{n-1}_{2} & \cdots & A^{n-1}_{n-1}
\end{array}}.$$  This is the matrix that we would like to simplify.
And we have quite a few tools to do it.

\subsection{The Jacobi Property}

The next step is to require the algebra to satisfy the Jacobi property.  This will eliminate several of the
$A_{j}^{i}$ in the $\ad{\e_{n}}$ matrix simultaneously.  In fact, as $D\g \subseteq NR(\g)$, $\ad{\e_{n}}$ maps
$NR(\g)$ to $NR(\g)$, and by the Jacobi property, we have, for $i,j < n$, that
\begin{align}
    \ad{\e_{n}}\pnth{[\e_{i},\e_{j}]} &= \brkt{\e_{n},[\e_{i},\e_{j}]} =
    - \brkt{\e_{i},\brkt{\e_{j},\e_{n}}} -
    \brkt{\e_{j},\brkt{\e_{n},\e_{i}}} \notag \\
    &= \brkt{\e_{i},\brkt{\e_{n},\e_{j}}} +
    \brkt{\brkt{\e_{n},\e_{i}},\e_{j}} =
    \brkt{\e_{i},\ad{\e_{n}}(\e_{j})} +
    \brkt{\ad{\e_{n}}(\e_{i}),\e_{j}}. \notag
\end{align}
Then $\ad{\e_{n}}$ is a derivation of the nilradical.  This implies
that in general, $\ad{\e_{n}}$ is an arbitrary derivation of the
nilradical.  This significantly simplifies the possible general form
of $\ad{\e_{n}}$.

Note that we didn't check the Jacobi property when the codimension one nilradical was abelian.  If the nilradical
is abelian then $[u,v]=0$ for all $u,v \in NR(\g)$.  Then, as $D\g \subseteq NR(\g)$, we see that $[\e_{n},\e_{i}]
\in NR(\g)$ for all $i <n$.  This yields that for all $i,j <n$, we have that $$\brkt{\e_{n},\brkt{\e_{i},\e_{j}}}
+ \brkt{\e_{i},\brkt{\e_{j},\e_{n}}} + \brkt{\e_{j},\brkt{\e_{n},\e_{i}}} =0.$$  Hence the Jacobi property was
already satisfied.

\subsection{Perturbing $\e_{n}$}

We begin our discussion here by making an observation. Consider an
$n$-dimensional Lie algebra with a codimension one {\em abelian}
nilradical and basis $\beta=\{\e_{1},\e_{2},\ldots,\e_{n}\}$, where
the first $n-1$ vectors constitute a basis for the nilradical.

Pick a new basis
\begin{align}
    \notag \e_{1}&=\e_{1} \\
    \notag \e_{2} &= \e_{2} \\
    \notag &\ \ \vdots \\
    \notag \e_{n-1}&= \e_{n-1} \\
    \notag \tilde{\e_{n}} &= \e_{n}+ \sum_{k=1}^{n-1} \la_{k} \e_{k}
\end{align}
and consider the structure equations.  We have that, for $i,j \leq
n-1$, $[\e_{i},\e_{j}]=0$ simply because the first $n-1$ vectors
still form a basis for the abelian nilradical.  Then the only
brackets left to consider are those of the form
$$[\tilde{\e_{n}},\e_{i}].$$  Thus we have that for $i \leq n-1$ (as of course $[\tilde{\e_{n}},
\tilde{\e_{n}}]=0$)
\begin{align}
    \notag [\tilde{\e_{n}},\e_{i}] &= \brkt{\e_{n}+
    \sum_{k=1}^{n-1} \la_{k}\e_{k},\e_{i}} \\
    \notag &= [\e_{n},\e_{i}]+\sum_{k=1}^{n-1} \la_{k}
    [\e_{k},\e_{i}] \\
    \notag &= [\e_{n},\e_{i}]
\end{align}
Thus we can see that the structure equations do not change for this
kind of change of basis.  That is, for an algebra with an abelian
nilradical, perturbing a vector outside of the nilradical by a
linear combination of vectors inside the nilradical will not change
the structure equations of the algebra. However, a scaling on
$\e_{n}$ would change the structure equations, something we used
quite often in classifying algebras of this type.

The important implication of this observation is that if the
nilradical is not abelian, then perturbing a vector outside the
nilradical by a linear combination of vectors inside the nilradical
{\em will} change the structure equations of the nilradical.  But it
will only change those structure equations that have a nonzero
bracket involving the vector we perturbed by.

We will use this idea to our advantage and perturb $\e_{n}$ by a
particular linear combination of vectors in the nilradical in order
to simplify the structure equations (or rather to zero out a few of
the $A_{j}^{i}$ terms in the $\ad{\e_{n}}$ matrix). We will also
remember that we can still change the structure equations by scaling
$\e_{n}$.

One more simple way to see how exactly we can change the structure
equations involving $\e_{n}$ is that the idea explained above
reduces to perturbing $\ad{\e_{n}}$ by linear combinations of the
$\add$ matrices of the nilradical basis vectors.  This is because
$$\ad{\e_{n}+ \sum_{k=1}^{n-1} \la_{k} \e_{k}} = \ad{\e_{n}} +
\sum_{k=1}^{n-1} \la_{k} \ad{\e_{k}}.$$ This gives us a better
viewpoint to see exactly what $A^{i}_{j}$'s we can change or
eliminate.

\subsection{The Automorphisms of $NR(\g)$}

Recall that when $NR(\g)$ was abelian, we could apply any change of
basis matrix to $\ad{\e_{n}}$, which is what allowed us to put
$\ad{\e_{n}}$ into real Jordan canonical form.  We could do this
because we knew that it wouldn't change the structure equations of
the nilradical.

We can than permute the basis of the nilradical in any way that
doesn't change its structure equations.  In other words, we can
conjugate $\ad{\e_{n}}$ by any change of basis matrix representing a
transformation that doesn't change the structure equations of the
nilradical. The obvious and only choice of transformations are the
automorphisms of the nilradical. This is consistent with our
previous work with abelian nilradicals because the automorphisms of
an abelian algebra consist of all invertible linear transformations
on that algebra.

These automorphisms will not affect the structure equations of the
nilradical by definition, but they do affect the makeup of
$\ad{\e_{n}}$.  We can then use these to our advantage to simplify
$\ad{\e_{n}}$.

We stated last chapter that if $T$ is a derivation of a Lie algebra,
then $e^{T}$ is an automorphism of that Lie algebra.  As such, the
automorphisms of $NR(\g)$ can be found by first computing the
derivations of $NR(\g)$ and then exponentiating them.

If the Lie algebra of derivations has a semisimple part via its Levi
decomposition, then the group corresponding to that semisimple
subalgebra is a subgroup of the automorphism group.  As
$\ad{\e_{n}}$ is an arbitrary derivation, we can then use this
subgroup of the automorphism group to move the semisimple part of
$\ad{\e_{n}}$ into a canonical form.  In the next section, we
discuss the semisimple algebras that arose as we classified from one
nilradical to the next. Specifically, we find the canonical forms of
the matrices representing these semisimple subalgebras when only
conjugation by its corresponding group is allowed.

\section{Canonical Forms}

In this section, we will discuss the possible canonical forms of
different sets of matrices.  In each case, we only allow conjugation
by a particular group of matrices.  We will first consider the set
of all $n \times n$ matrices, which is the Lie algebra $\gl(n,\R)$,
and allow conjugation by elements of the group $GL(n,\R)$; this, of
course, results in the real Jordan canonical form. The second set of
matrices considered are those in the symplectic Lie algebra
$\si(4,\R)$ and we allow conjugation by the group $SP(4,\R)$.  The
final set of matrices considered form a representation of the Lie
algebra $\so(3,1,\R)$, but an non-equivalent representation to the
usual one. In that case, we will allow conjugation by the group of
matrices corresponding to that representation.

\subsection{Real Jordan Canonical Forms of Matrices of Small Dimension}

Let $a \in \gl(n,\R)$ and $A \in GL(n\R)$.  The canonical forms of
$a$, of course, is the real Jordan canonical form, which has already
been described in detail in Chapter \ref{c1}. In this section, we
simply enumerate explicitly the real Jordan canonical forms of
matrices of dimension six or less. We can describe these possible
forms by taking a look at how the generalized eigenspaces of a
transformation split the vector space, $V$, into a direct sum of
smaller vector spaces over $\R$.

Dimension one is trivial as the one dimensional $\R$ cannot be
split. Thus the only real Jordan canonical form for a one
dimensional linear transformation is $$\begin{pmatrix} \lambda
\end{pmatrix}.$$

For dimension two, $\R^{2}$ can be split in the following ways
\begin{align}
    \notag \R^{2} &= \R \oplus \R \\
    \notag \R^{2} &= \R^{2}
\end{align}
The case where $\R^{2}=\R \oplus \R$ is easy as we can build it from
all the combinations of one dimensional real Jordan canonical forms.
This will give us the case where $a$ is diagonalizable over $\R$.
When $\R^{2}=\R^{2}$, we have two cases: either the transformation
is diagonalizable over $\C$ but not $\R$, or it can be put into
general real Jordan canonical form over $\R$.  Thus we have three
possibilities for the two-dimensional case.
\begin{align}
    \notag &\begin{pmatrix}\lambda_{1} & 0 \\ 0 & \lambda_{2}
    \end{pmatrix},& & \begin{pmatrix} \lambda_{1} & \lambda_{2}
    \\ -\lambda_{2} & \lambda_{1} \end{pmatrix}, & & \begin{pmatrix}
    \lambda & 1 \\ 0 & \lambda \end{pmatrix}.
\end{align}
We will order the following lists in a similar manner to the pattern
above.

Next we consider linear transformations on three-dimensional vector spaces.  We can split $\R^{3}$ in the
following ways
\begin{align}
    \notag \R^{3} &=  \R \oplus \R \oplus \R \displaybreak[0] \\
    \notag \R^{3} &= \R^{2} \oplus \R \displaybreak[0] \\
    \notag \R^{3} &= \R^{3}
\end{align}
As will always be the case, the splitting of $\R^{n}$ into
one-dimensional eigenspaces, here that is $\R^{3}=\R \oplus \R
\oplus \R$, yields the situation that $a$ is diagonalizable over the
reals.  The case where $\R^{3}=\R^{2} \oplus \R$ does something of
note.  The $\R^{2}$ will follow the two cases for $\R^{2}$ given
above, while the $\R$ piece simply yields a one-dimensional
eigenspace.  Actually it will follow that however $\R^{n}$ breaks up
into eigenspaces, the dimension of each eigenspace will determine
what real Jordan canonical forms that block can have.  As such, we
will refer quite often to information found previously when
determining the general real Jordan canonical forms for a particular
dimension. Finally, $\R^{3}$ only yields one case of a
three-dimensional generalized eigenspace.  The complex eigenspaces
do not come into play here as they must appear in conjugate pairs
(as we're working with characteristic polynomials with real
coefficients) and we have an odd-dimensional space. Thus dimension
three only yields these four cases, up to an ordering of eigenvalues
of course.
\begin{align}
    \notag & \begin{pmatrix} \lambda_{1} &0&0 \\ 0&\lambda_{2} &0 \\
    0&0& \lambda_{3} \end{pmatrix}, & & \begin{pmatrix} \lambda_{1} & \lambda_{2} &0 \\
    -\lambda_{2}& \lambda_{1} &0 \\ 0&0& \lambda_{3} \end{pmatrix}, \displaybreak[0]
    \\
    \notag & \begin{pmatrix} \lambda_{1} &
    1 &0 \\ 0 & \lambda_{1}&0 \\ 0&0& \lambda_{2} \end{pmatrix},& & \begin{pmatrix} \lambda&1&0 \\
    0& \lambda&1 \\ 0&0& \lambda
    \end{pmatrix}.
\end{align}

We now consider dimension four.  We have the splitting of $\R^{4}$, thus
\begin{align}
    \notag \R^{4}&=\R \oplus \R \oplus \R \oplus \R \displaybreak[0] \\
    \notag \R^{4}&= \R^{2} \oplus \R \oplus \R \displaybreak[0] \\
    \notag \R^{4}&= \R^{2} \oplus \R^{2} \displaybreak[0] \\
    \notag \R^{4}&= \R^{3} \oplus \R \displaybreak[0] \\
    \notag \R^{4}&= \R^{4}
\end{align}
Using the method described above, we can look at the previous
dimensions for the possible forms of each eigenspace of that
dimension and combining them together according to the direct sum
listed.  This will cover the first four cases.  As for
$\R^{4}=\R^{4}$, we have a few possibilities.  We could simply have
a four-dimensional real generalized eigenspace, or we could have two
two-dimensional complex generalized eigenspaces, in which case, we
would move them to the real case using $2 \times 2$ blocks with a
two-dimensional identity block in the strictly upper triangular
piece.  This is described in more detail in the Chapter \ref{c1}.
For dimension four then, we have the following nine cases (always up
to an ordering of the eigenvalues).
\begin{align}
    \notag & \begin{pmatrix} \lambda_{1} &0&0&0 \\
    0&\lambda_{2} &0&0 \\ 0&0& \lambda_{3} &0 \\ 0&0&0& \lambda_{4}
    \end{pmatrix}, & & \begin{pmatrix}
    \lambda_{1} &\lambda_{2}&0&0 \\ -\lambda_{2}&\lambda_{1}& 0&0 \\
    0&0& \lambda_{3} &0 \\ 0&0&0& \lambda_{4} \end{pmatrix}, & &
    \begin{pmatrix} \lambda_{1}&1&0&0 \\ 0& \lambda_{1}&0&0 \\
    0&0&\lambda_{2}&0 \\ 0&0&0& \lambda_{3} \end{pmatrix}, \displaybreak[0] \\
    \notag & \begin{pmatrix} \lambda_{1} & \lambda_{2} &0&0 \\ -\lambda_{2}&
    \lambda_{1} &0&0 \\ 0&0& \lambda_{3} & \lambda_{4} \\ 0&0&
    -\lambda_{4}& \lambda_{3} \end{pmatrix},& & \begin{pmatrix} \lambda_{1} &1 &0&0 \\
    0& \lambda_{1}&0&0 \\ 0&0& \lambda_{2}&\lambda_{3} \\ 0&0&
    -\lambda_{3}& \lambda_{2} \end{pmatrix},& & \begin{pmatrix}
    \lambda_{1}&1&0&0 \\ 0& \lambda_{1}&0&0 \\ 0&0& \lambda_{2}&1 \\
    0&0&0& \lambda_{2} \end{pmatrix}, \displaybreak[0] \\
    \notag & \begin{pmatrix} \lambda_{1}&1&0&0 \\ 0&\lambda_{1}&1&0 \\ 0&0&
    \lambda_{1}&0 \\ 0&0&0& \lambda_{2} \end{pmatrix},& & \begin{pmatrix} \lambda_{1} & \lambda_{2} &1&0 \\
    -\lambda_{2}& \lambda_{1} &0&1 \\ 0&0& \lambda_{1} & \lambda_{2} \\
    0&0& -\lambda_{2}& \lambda_{1} \end{pmatrix}, & &
    \begin{pmatrix} \lambda &1&0&0 \\0& \lambda&1&0 \\ 0&0& \lambda&1 \\
    0&0&0& \lambda \end{pmatrix}.
\end{align}

For dimension five, we split $\R^{5}$ in the following ways
\begin{align}
    \notag \R^{5}&=\R \oplus \R \oplus \R \oplus \R \oplus \R \displaybreak[0] \\
    \notag \R^{5}&=\R^{2} \oplus \R \oplus \R \oplus \R \displaybreak[0] \\
    \notag \R^{5}&=\R^{2} \oplus \R^{2} \oplus \R \displaybreak[0] \\
    \notag \R^{5}&=\R^{3} \oplus \R \oplus \R \displaybreak[0] \\
    \notag \R^{5}&=\R^{3} \oplus \R^{2} \displaybreak[0] \\
    \notag \R^{5}&=\R^{4} \oplus \R \displaybreak[0] \\
    \notag \R^{5}&=\R^{5}
\end{align}
We again use the same technique already described and also note that the only case that $\R^{5}=\R^{5}$ generates
is a five-dimensional generalized eigenspace.  Thus we have the following 12 cases for dimension five.
\begin{align}
    \notag &\begin{pmatrix} \lambda_{1} &0&0&0&0 \\
    0&\lambda_{2} &0&0&0 \\ 0&0& \lambda_{3} &0&0 \\ 0&0&0&
    \lambda_{4}&0 \\ 0&0&0&0& \lambda_{5} \end{pmatrix}, & &
    \begin{pmatrix} \lambda_{1} &\lambda_{2}&0&0&0 \\
    -\lambda_{2}&\lambda_{1}& 0&0&0 \\ 0&0& \lambda_{3} &0&0 \\ 0&0&0& \lambda_{4}&0 \\ 0&0&0&0& \lambda_{5}
    \end{pmatrix}, \displaybreak[0] \\
    \notag &\begin{pmatrix}\lambda_{1}&1&0&0&0 \\0& \lambda_{1}&0&0&0\\
    0&0&\lambda_{2}&0&0 \\ 0&0&0& \lambda_{3}&0 \\ 0&0&0&0& \lambda_{4}
    \end{pmatrix},& & \begin{pmatrix} \lambda_{1} & \lambda_{2} &0&0&0 \\ -\lambda_{2}&
    \lambda_{1} &0&0&0 \\ 0&0& \lambda_{3} &
    \lambda_{4}&0 \\ 0&0& -\lambda_{4}& \lambda_{3}&0 \\ 0&0&0&0&
    \lambda_{5} \end{pmatrix}, \displaybreak[0] \\
    \notag &\begin{pmatrix} \lambda_{1} &1 &0&0&0
    \\ 0& \lambda_{1}&0&0&0 \\ 0&0& \lambda_{2}& \lambda_{3}&0 \\ 0&0&
    -\lambda_{3}& \lambda_{2}&0 \\ 0&0&0&0& \lambda_{4} \end{pmatrix}, & &
    \begin{pmatrix} \lambda_{1}&1&0&0&0 \\ 0& \lambda_{1}&0&0&0 \\
    0&0& \lambda_{2}&1&0 \\ 0&0&0& \lambda_{2}&0
    \\ 0&0&0&0& \lambda_{3} \end{pmatrix}, \displaybreak[0] \\
    \notag &\begin{pmatrix} \lambda_{1}&1&0&0&0
    \\ 0&\lambda_{1}&1&0&0 \\ 0&0&\lambda_{1} &0&0 \\ 0&0&0& \lambda_{2}&0 \\
    0&0&0&0& \lambda_{3} \end{pmatrix},& & \begin{pmatrix} \lambda_{1}&1&0&0&0 \\
    0&\lambda_{1}&1&0&0 \\ 0&0& \lambda_{1} &0&0 \\ 0&0&0& \lambda_{2}& \lambda_{3} \\ 0&0&0& -\lambda_{3}&
    \lambda_{2} \end{pmatrix}, \displaybreak[0] \\
    \notag &\begin{pmatrix} \lambda_{1}&1&0&0&0 \\
    0&\lambda_{1}&1&0&0 \\ 0&0& \lambda_{1} &0&0
    \\ 0&0&0& \lambda_{2}&1 \\ 0&0&0&0& \lambda_{2} \end{pmatrix},& &
    \notag \begin{pmatrix} \lambda_{1} & \lambda_{2} &1&0&0 \\ -\lambda_{2}& \lambda_{1} &0&1&0 \\
    0&0& \lambda_{1} &\lambda_{2}&0 \\ 0&0& -\lambda_{2}& \lambda_{1}&0 \\ 0&0&0&0&
    \lambda_{3} \end{pmatrix}, \displaybreak[0] \\
    \notag &\begin{pmatrix} \lambda_{1} &1&0&0&0
    \\0& \lambda_{1}&1&0&0 \\ 0&0& \lambda_{1}&1&0
    \\ 0&0&0& \lambda_{1}&0 \\0&0&0&0& \lambda_{2} \end{pmatrix},& &
    \begin{pmatrix} \lambda  &1&0&0&0 \\0& \lambda&1&0&0 \\ 0&0&
    \lambda&1&0 \\ 0&0&0& \lambda&1 \\0&0&0&0& \lambda
    \end{pmatrix}.
\end{align}

Finally, for dimension six, we can split $\R^{6}$ up into eigenspaces thus
\begin{align}
    \notag \R^{6} &= \R \oplus \R \oplus \R \oplus \R \oplus \R \oplus \R \displaybreak[0] \\
    \notag \R^{6} &= \R^{2} \oplus \R \oplus \R \oplus \R \oplus \R \displaybreak[0] \\
    \notag \R^{6} &= \R^{2} \oplus \R^{2} \oplus \R \oplus \R \displaybreak[0] \\
    \notag \R^{6} &= \R^{2} \oplus \R^{2} \oplus \R^{2} \displaybreak[0] \\
    \notag \R^{6} &= \R^{3} \oplus \R \oplus \R \oplus \R \displaybreak[0] \\
    \notag \R^{6} &= \R^{3} \oplus \R^{2} \oplus \R \displaybreak[0] \\
    \notag \R^{6} &= \R^{3} \oplus \R^{3} \displaybreak[0] \\
    \notag \R^{6} &= \R^{4} \oplus \R \oplus \R \displaybreak[0] \\
    \notag \R^{6} &= \R^{4} \oplus \R^{2} \displaybreak[0] \\
    \notag \R^{6} &= \R^{5} \oplus \R \displaybreak[0] \\
    \notag \R^{6} &= \R^{6}
\end{align}
Following the same method we can find all the forms up to the last case.  For $\R^{6}=\R^{6}$, we can have either
a six-dimensional real generalized eigenspace, or a six-dimensional complex generalized eigenspace, again with
identity matrices on the strictly upper triangular part.  Thus for dimension six, we have the following 23
possible real Jordan canonical forms
\begin{align}
    \notag & \begin{pmatrix} \lambda_{1} &0&0&0&0&0 \\ 0&\lambda_{2} &0&0&0&0 \\
    0&0& \lambda_{3}
    &0&0&0 \\ 0&0&0& \lambda_{4}&0&0 \\ 0&0&0&0& \lambda_{5}&0 \\
    0&0&0&0&0& \lambda_{6} \end{pmatrix},& & \begin{pmatrix} \lambda_{1} & \lambda_{2}&0&0&0&0 \\
    -\lambda_{2}& \lambda_{1} &0&0&0&0 \\ 0&0& \lambda_{3} &0&0&0 \\
    0&0&0& \lambda_{4}&0&0 \\ 0&0&0&0& \lambda_{5}&0
    \\ 0&0&0&0&0& \lambda_{6} \end{pmatrix}, \displaybreak[0] \\
    \notag & \begin{pmatrix}
    \lambda_{1} &1&0&0&0&0\\0&\lambda_{1} &0&0&0&0\\0&0& \lambda_{2}
    &0&0&0 \\ 0&0&0& \lambda_{3}&0&0 \\ 0&0&0&0& \lambda_{4}&0 \\
    0&0&0&0&0& \lambda_{5} \end{pmatrix},& & \begin{pmatrix} \lambda_{1} & \lambda_{2}&0&0&0&0 \\
    -\lambda_{2}&\lambda_{1} &0&0&0&0 \\ 0&0& \lambda_{3} & \lambda_{4} &0&0 \\ 0&0&-\lambda_{4}& \lambda_{3}&0&0 \\
    0&0&0&0& \lambda_{5}&0 \\ 0&0&0&0&0& \lambda_{6} \end{pmatrix}, \displaybreak[0] \\
    \notag & \begin{pmatrix} \lambda_{1} &1&0&0&0&0 \\
    0&\lambda_{1} &0&0&0&0 \\ 0&0& \lambda_{2} & \lambda_{3} &0&0 \\
    0&0&-\lambda_{3}& \lambda_{2}&0&0 \\ 0&0&0&0& \lambda_{4}&0 \\
    0&0&0&0&0&
    \lambda_{5} \end{pmatrix},& & \begin{pmatrix} \lambda_{1} &1&0&0&0&0 \\
    0&\lambda_{1} &0&0&0&0 \\ 0&0& \lambda_{2} &1&0&0 \\ 0&0&0& \lambda_{2}&0&0 \\ 0&0&0&0& \lambda_{3}&0 \\
    0&0&0&0&0& \lambda_{4} \end{pmatrix}, \displaybreak[0] \\
    \notag & \begin{pmatrix} \lambda_{1}
    &\lambda_{2} &0&0&0&0 \\ -\lambda_{2}&\lambda_{1} &0&0&0&0 \\ 0&0& \lambda_{3} &\lambda_{4}&0&0 \\
    0&0&-\lambda_{4}&  \lambda_{3}&0&0 \\ 0&0&0&0& \lambda_{5}&\lambda_{6} \\ 0&0&0&0&-\lambda_{6}& \lambda_{5}
    \end{pmatrix},& & \begin{pmatrix} \lambda_{1} &1&0&0&0&0 \\ 0&\lambda_{1} &0&0&0&0 \\ 0&0&
    \lambda_{2} &\lambda_{3}&0&0 \\ 0&0&-\lambda_{3}& \lambda_{2}&0&0 \\ 0&0&0&0& \lambda_{4}&\lambda_{5} \\
    0&0&0&0&-\lambda_{5}& \lambda_{4} \end{pmatrix}, \displaybreak[0] \\
    \notag & \begin{pmatrix} \lambda_{1} &1&0&0&0&0 \\ 0&\lambda_{1}
    &0&0&0&0 \\ 0&0& \lambda_{2} &1&0&0 \\ 0&0&0& \lambda_{2}&0&0 \\
    0&0&0&0& \lambda_{3}& \lambda_{4} \\ 0&0&0&0&-\lambda_{4}&
    \lambda_{3} \end{pmatrix},& & \begin{pmatrix} \lambda_{1}
    &1&0&0&0&0 \\ 0&\lambda_{1}
    &0&0&0&0 \\ 0&0& \lambda_{2} &1&0&0 \\ 0&0&0& \lambda_{2}&0&0 \\
    0&0&0&0& \lambda_{3}&1 \\ 0&0&0&0&0& \lambda_{3} \end{pmatrix}, \displaybreak[0] \\
    \notag & \begin{pmatrix} \lambda_{1} &1&0&0&0&0 \\
    0&\lambda_{1} &1&0&0&0 \\ 0&0& \lambda_{1} &0&0&0 \\ 0&0&0& \lambda_{2}&0&0 \\ 0&0&0&0& \lambda_{3}&0 \\
    0&0&0&0&0& \lambda_{4}
    \end{pmatrix},& & \begin{pmatrix} \lambda_{1} &1&0&0&0&0 \\ 0&\lambda_{1}
    &1&0&0&0 \\ 0&0& \lambda_{1} &0&0&0 \\ 0&0&0&
    \lambda_{2}&\lambda_{3}&0 \\ 0&0&0&-\lambda_{3}& \lambda_{2}&0 \\
    0&0&0&0&0& \lambda_{4} \end{pmatrix}, \displaybreak[0] \\
    \notag & \begin{pmatrix}
    \lambda_{1} &1&0&0&0&0 \\ 0&\lambda_{1}
    &1&0&0&0 \\ 0&0& \lambda_{1} &0&0&0 \\ 0&0&0& \lambda_{2}&1&0 \\
    0&0&0&0& \lambda_{2}&0 \\ 0&0&0&0&0& \lambda_{3} \end{pmatrix},& & \begin{pmatrix} \lambda_{1} &1&0&0&0&0 \\
    0&\lambda_{1}&1&0&0&0 \\ 0&0& \lambda_{1} &0&0&0 \\ 0&0&0& \lambda_{2}&1&0 \\
    0&0&0&0& \lambda_{2}&1 \\ 0&0&0&0&0& \lambda_{2} \end{pmatrix}, \displaybreak[0]
    \\
    \notag & \begin{pmatrix} \lambda_{1}
    &\lambda_{2}&1&0&0&0 \\ -\lambda_{2}&\lambda_{1} &0&1&0&0\\
    0&0& \lambda_{1} &\lambda_{2}&0&0 \\ 0&0&-\lambda_{2}& \lambda_{1}&0&0 \\ 0&0&0&0& \lambda_{3}&0 \\ 0&0&0&0&0&
    \lambda_{4}
    \end{pmatrix},& & \begin{pmatrix} \lambda_{1}
    &\lambda_{2}&1&0&0&0 \\  -\lambda_{2}&\lambda_{1} &0&1&0&0 \\
    0&0& \lambda_{1} &\lambda_{2}&0&0 \\ 0&0&-\lambda_{2}& \lambda_{1}&0&0 \\ 0&0&0&0& \lambda_{3}&\lambda_{4} \\
    0&0&0&0&-\lambda_{4}& \lambda_{3} \end{pmatrix}, \displaybreak[0] \\
    \notag & \begin{pmatrix} \lambda_{1} & \lambda_{2}&1&0&0&0 \\
    -\lambda_{2}&\lambda_{1}
    &0&1&0&0 \\ 0&0& \lambda_{1} &\lambda_{2}&0&0 \\ 0&0&-\lambda_{2}& \lambda_{1}&0&0 \\ 0&0&0&0& \lambda_{3}&1 \\
    0&0&0&0&0& \lambda_{3} \end{pmatrix},& &
    \begin{pmatrix} \lambda_{1} &1&0&0&0&0 \\ 0&\lambda_{1}
    &1&0&0&0 \\ 0&0& \lambda_{1} &1&0&0 \\ 0&0&0& \lambda_{1}&0&0 \\
    0&0&0&0& \lambda_{2}&0 \\ 0&0&0&0&0& \lambda_{3} \end{pmatrix}, \displaybreak[0]
    \\
    \notag & \begin{pmatrix} \lambda_{1} &1&0&0&0&0 \\ 0&\lambda_{1}
    &1&0&0&0 \\ 0&0& \lambda_{1} &1&0&0 \\ 0&0&0& \lambda_{1}&0&0 \\
    0&0&0&0& \lambda_{2}& \lambda_{3} \\ 0&0&0&0&-\lambda_{3}&
    \lambda_{2} \end{pmatrix},& & \begin{pmatrix} \lambda_{1}
    &1&0&0&0&0 \\ 0&\lambda_{1}
    &1&0&0&0 \\ 0&0& \lambda_{1} &1&0&0 \\ 0&0&0& \lambda_{1}&0&0 \\
    0&0&0&0& \lambda_{2}&1 \\ 0&0&0&0&0& \lambda_{2} \end{pmatrix}, \displaybreak[0] \\
    \notag & \begin{pmatrix} \lambda_{1} &1&0&0&0&0 \\
    0&\lambda_{1} &1&0&0&0 \\ 0&0& \lambda_{1} &1&0&0 \\ 0&0&0& \lambda_{1}&1&0 \\ 0&0&0&0& \lambda_{1}&0 \\
    0&0&0&0&0& \lambda_{2}
    \end{pmatrix},& & \begin{pmatrix} \lambda_{1} &\lambda_{2}&1&0&0&0 \\ -\lambda_{2}&\lambda_{1} &0&1&0&0 \\
    0&0& \lambda_{1} &\lambda_{2}&1&0 \\ 0&0&-\lambda_{2}& \lambda_{1}&0&1 \\ 0&0&0&0& \lambda_{1}&\lambda_{2} \\
    0&0&0&0&-\lambda_{2}& \lambda_{1} \end{pmatrix}, \displaybreak[0] \\
    \notag & \begin{pmatrix} \lambda &1&0&0&0&0 \\ 0&\lambda &1&0&0&0 \\
    0&0& \lambda &1&0&0 \\ 0&0&0& \lambda&1&0 \\ 0&0&0&0& \lambda&1 \\
    0&0&0&0&0& \lambda \end{pmatrix}.
\end{align}

This completes the list as far as we require.

\subsection{Canonical Forms of Matrices in $\si(4,\R)$}

Let $V=\R^{2n}$.  Let $J \in GL(V)$ be given as $$J= \begin{pmatrix} 0&I_{n}
\\ -I_{n}&0 \end{pmatrix}$$ where $I_{n}$ is the $n \times n$ identity.  Note that
$J$ is skew-symmetric and has the property that $J^{2} = I_{2n}$. The symplectic Lie algebra, $\si(2n,\R)$, is
defined as
\begin{equation}\label{eq:A}\si(2n,\R)=\{a \in \Hom(V,V) \ | \ a^{t}J + Ja = 0 \}. \
\cite{sattinger} \notag \end{equation}

\begin{lem}\label{A}
    Let $a \in \si(2n,\R)$.  Then $a^{t} \in \si(2n,\R)$.
\end{lem}

{\em Proof.} If $a \in \si(2n,\R)$, then $a^{t}J+Ja=0$.  Multiply on the left by $J$ to obtain $Ja^{t}J-a=0$. Then
multiply on the right: $-Ja^{t}-aJ=0$ or equivalently $(a^{t})^{t}J + J a^{t} = 0$. $\blacksquare$ \vspace{.1 in}

Define $\om:V \times V \To \R$ by $\omega(x,y)=x^{t}Jy$ for all $x,y
\in V$.  This is called a {\em symplectic form}. It can also be
viewed as a map $V_{c} \times V_{c} \To \C$ with the same rule of
assignment. We say that $\om$ is {\em non-degenerate} if whenever $z
\in V$ is such that $\om(z,y)=0$ for all $y \in V$, then $z=0$.

\begin{prop}\label{B}
    The symplectic form, $\om:V \times V \To \R$ ($\om:V_{c} \times V_{c} \To \C$),
    is a non-degenerate skew-symmetric bilinear form on $V$ ($V_{c}$).
\end{prop}

{\em Proof.} From its construction, it is clearly bilinear over $\R$ or $\C$. To prove skew-symmetry, we note that
as $\omega(x,y)$ can be viewed as a $1\times 1$ matrix, we have that $(\om(x,y))^{t}=\om(x,y)$. Then, as
$J^{t}=-J$, this yields that
$$\om(x,y)=\pnth{\om(x,y)}^{t} = (x^{t}Jy)^{t} = y^{t}J^{t}x = -y^{t}Jx = -\om(y,x).$$  Finally to prove
non-degeneracy, let $z \in V$ be such that $\om(z,y)=z^{t}Jy=0$ for all $y \in V$.  As $J$ is invertible, it has
zero kernel, thus $z=0$. $\blacksquare$ \vspace{.1 in}

We will always denote the symplectic form by $\om$ as it will be clear from the context whether we mean the
complex or real form.

Now we define the symplectic group.  The symplectic group, $Sp(2n,\R)$ is defined as $$Sp(2n,\R)= \{ A \in GL(V) \
| \ \om(Ax,Ay) = \om(x,y) \}.$$  That is, it is the group that preserves the symplectic form.

\begin{lem}\label{C}
    $a \in \si(2n,\R)$ if and only if $\om(ax,y)=-\om(x,ay)$ for all $x,y \in V$. $A \in Sp(2n,\R)$ if an only if
    $A^{t}JA = J$.
\end{lem}

{\em Proof.} As $a \in \si(2n,\R)$, we have that $a^{t}J + Ja=0$.  Let $x,y \in V$.  This gives us
$$\om(ax,y)=(ax)^{t}Jy = x^{t} a^{t}Jy = -x^{t}J ay = -\om(x,ay).$$

Next assume that $\om(ax,y)=-\om(x,ay)$ for all $x,y \in V$.  Then
\begin{align}
    \om(ax,y) &= -\om(x,ay) \notag\\
    (ax)^{t} J y &= -x^{t} J a y \notag\\
    x^{t} a^{t} J y &= -x^{t} J a y \notag
\end{align}
As this is true for all $x,y \in V$, this implies that $a^{t}J = -Ja$ or equivalently that $a^{t}J +Ja = 0$. Hence
$a \in \si(2n,\R)$.

To prove the statement about $Sp(2n,\R)$, first assume that $A \in Sp(2n,\R)$.  Then
\begin{align}
    \om(Ax,Ay) &= \om(x,y) \notag \\
    (Ax)^{t} J Ay &= x^{t}Jy \notag\\
    x^{t}A^{t}JAy &= x^{t}Jy \notag
\end{align}
As this is true for all $x,y \in V$, we see that $A^{t}JA = J$.

If we assume first that $A^{t}JA = J$, then $A$ clearly preserves the symplectic form and hence $A \in Sp(2n,\R)$.
$\blacksquare$ \vspace{.1 in}

We now wish to conjugate $a \in \si(2n,\R)$ by an arbitrary element $A \in Sp(2n,\R)$, that is, $A^{-1}aA$.

\begin{lem} \label{D}
    If $a \in \si(2n,\R)$ and $A \in Sp(2n,\R)$, then $A^{-1}aA \in \si(2n,\R)$.
\end{lem}

{\em Proof.} If $A \in Sp(2n,\R)$, then $\om(Ax,y) = \om(A^{-1}Ax,A^{-1}y) = \om(x,A^{-1}y)$.  This yields
$$\om(A^{-1}aA x, y) = \om(aA x, Ay) = -\om(Ax,aAy) = -\om(x,A^{-1}aAy).$$  Then by Lemma \ref{C}, $A^{-1}aA \in
\si(2n,\R)$. $\blacksquare$ \vspace{.1 in}

If $a_{1},a_{2} \in \si(2n,\R)$ are such that $A^{-1}a_{1}A=a_{2}$ for some $A \in Sp(2n,\R)$, we say that $a_{1}$
and $a_{2}$ are {\em symplectically similar}.

This naturally brings up the question: what kind of canonical forms could $a \in \si(2n,\R)$ have if this were the
only kind of change of basis allowed? In four dimensions, the result is as follows

\begin{thm} \label{E}
    Let $a \in \si(4,\R)$, then $a$ is symplectically similar to one
    of the following ten matrices.  We call this the {\em real
    symplectic canonical form} of the matrix.
    \begin{align}
        (1) \ &\begin{pmatrix} \la &0&0&0 \\ 0&\mu&0&0 \\ 0&0&-\la&0 \\
        0&0&0&-\mu \end{pmatrix}, \  \la,\mu \in \R, & (2) \ & \begin{pmatrix}
        \la &0&0&0 \\ 0&0&0&\vep \\ 0&0&-\la&0 \\ 0&0&0&0
        \end{pmatrix}, \ \begin{array}{l} \la \in \R, \\ \vep^{2}=1, \end{array}
        \notag \displaybreak[0]\\
        (3) \ & \begin{pmatrix} \la &1 &0&0 \\ 0&\la&0&0 \\ 0&0&-\la&0 \\
        0&0&-1&-\la \end{pmatrix}, \  \la \in \R, & (4) \ & \begin{pmatrix} 0&0&\vep&0 \\
        0&0&0&\vep \\ 0&0&0&0 \\ 0&0&0&0 \end{pmatrix}, \  \vep^{2} =
        1, \notag \displaybreak[0]\\
        (5) \ &\begin{pmatrix} 0&1&0&0 \\ 0&0&0&\vep
        \\ 0&0&0&0 \\ 0&0&-1&0 \end{pmatrix}, \  \vep^{2} = 1, & (6) \ & \begin{pmatrix} \la &0&0&0 \\
        0&0&0&\vep \mu \\ 0&0&-\la&0 \\ 0&-\vep\mu &0&0
        \end{pmatrix}, \ \begin{array}{l} \la \in \R, \\ \mu > 0, \\ \vep^{2} =1, \end{array} \notag \displaybreak[0]\\
        (7) \ &\begin{pmatrix} 0&0&\vep&0 \\ 0&0&0&\del \mu \\ 0&0&0&0 \\ 0&-\del\mu &0&0 \end{pmatrix},
        \ \begin{array}{l} \mu > 0, \\ \del^{2}=\vep^{2}=1, \end{array} & (8) \ &\begin{pmatrix} \la&\mu&0&0 \\
        -\mu&\la&0&0 \\ 0&0&-\la&\mu \\ 0&0&-\mu&-\la \end{pmatrix}, \ \begin{array}{l} \la \in \R, \\ \mu \neq 0,
        \end{array} \notag \displaybreak[0]\\
        (9) \ &\begin{pmatrix} 0&0&\vep \mu&0 \\ 0&0&0&\vep \eta \\ -\vep
        \mu&0&0&0 \\ 0&-\vep \eta&0&0 \end{pmatrix}, \ \begin{array}{l} \eta \neq -\mu, \\
        \eta,\mu \neq 0, \\ \vep^{2} = 1, \end{array} & (10) \ &
        \begin{pmatrix} 0&\mu&\vep&0 \\ -\mu&0&0&\vep \\ 0&0&0&\mu
        \\ 0&0&-\mu&0 \end{pmatrix}, \ \begin{array}{l} \mu \neq 0, \\ \vep^{2} = 1. \end{array}
        \notag \displaybreak[0]
    \end{align}
\end{thm}

The remainder of this section is the proof of this theorem.  First
we present some general theory and then move on to the particulars
of the four dimensional case.

\begin{lem}\label{F}
    Let $a \in \si(2n,\R)$ and let $\la$ be an eigenvalue of $a$.  Then $-\la$ is also an eigenvalue of $a$.
\end{lem}

{\em Proof.} As $a \in \si(2n,\R)$, we have that $a^{t}J+Ja=0$ or equivalently $JaJ^{-1}=-a^{t}$.  If $\la$ is an
eigenvalue of $a$, then $\det(a-\la I_{2n}) = 0$.  This yields the following
\begin{align}
     0 &= \det(a-\la I_{2n}) = \det(J(a-\la I_{2n})J^{-1}) = \det(JaJ^{-1} - \la JI_{2n}J^{-1}) \notag \\
     &= \det(-a^{t}-\la I_{2n}) = \det(-(a+\la I_{2n})^{t}) = (-1)^{2n}\det((a+\la I_{2n})^{t}) =
     \det(a+\la I_{2n}) \notag
\end{align}
which, of course, implies that $-\la$ is also an eigenvalue of $a$. $\blacksquare$ \vspace{.1 in}

We know that the characteristic polynomial of $a$ always splits over $\C$, thus we have that the complexification
of $V$, $V_{c}$, is the direct sum of the generalized eigenspaces of $a$.  This allows us to write $V_{c}$ in the
following manner
$$V_{c} = K_{0} \oplus K_{\la_{1}} \oplus K_{-\la_{1}} \oplus K_{\la_{2}} \oplus K_{-\la_{2}} \oplus \cdots \oplus
K_{\la_{k}} \oplus K_{-\la_{k}}$$ where $1 \leq k \leq n$, $K_{\mu}$ is the generalized eigenspace corresponding
to the eigenvalue $\mu$, $\la_{i} \neq 0$ for all $i$, and the $\la_{i}$ are distinct and such that $\la_{i} \neq
-\la_{j}$ for any pair $i,j$. Note that if all the eigenvalues of $a$ are real, then its characteristic polynomial
splits over $\R$ and $V$ decomposes in the manner above where all the summands are subspaces over $\R$.

\begin{prop}\label{G}
    If $\mu \neq -\eta$, then $K_{\eta}$ and $K_{\mu}$ are symplectically orthogonal to each other.
    That is, $\om(K_{\eta},K_{\mu})=0$.
\end{prop}

{\em Proof.} Let $v \in K_{\eta}$ be an eigenvector of $a$ (i.e. $av = \eta v$) and let $w \in K_{\mu}$ be an
eigenvector of $a$.  Then we have that$$\eta \om(v,w) = \om(\eta v,w) = \om(a v,w)= -\om(v, a w)= -\om(v, \mu w)=
-\mu \om(v,w).$$ This implies that $(\eta+\mu) \om(v,w)=0$.  As $\mu \neq -\eta$, this implies that $\om(v,w)=0$.

Now let $\{w,w'\} \subseteq K_{\mu}$ be a cycle of generalized
eigenvectors. Then we have that
\begin{align}
    \notag \eta \om(v,w') &= \om(\eta v,w') = \om(a v,w')= -\om(v,a
    w') \\
    \notag &= -\om(v, \mu w' + w) = -\mu\om(v,w') - \om(v,w)= -\mu\om(v,w').
\end{align}
This implies that $(\eta+\mu)\om(v,w')=0$.  And as $\mu \neq -\eta$,
this gives us that $\om(v,w')=0$.

Now we have that any eigenvector in $K_{\eta}$ is symplectically orthogonal to any eigenvector and ``one-step''
generalized eigenvector in $K_{\mu}$.  Proceed inductively on the step of the generalized eigenvector in $K_{\mu}$
and this will imply that any eigenvector in $K_{\eta}$ is symplectically orthogonal to any vector in $K_{\mu}$.

Next let $\{v,v'\} \subseteq K_{\eta}$ be a cycle of generalized eigenvectors.  Then for any eigenvector $w \in
K_{\mu}$, we have that $\om(v',w)=0$ by the above argument and the skew-symmetry of $\om$.  Let $\{w,w'\}$ be a
cycle of generalized eigenvectors and recall that $\om(v,w')=0$. This will give us
\begin{align}
     \eta \om(v',w') &= \om(v,w') + \eta \om(v',w') = \om(\eta v' + v, w') = \om(a v', w') \notag \\
     &= -\om(v',a w') = -\om(v',\mu w' + w) = -\mu\om(v',w') - \om(v',w) = -\mu\om(v',w') \notag
\end{align}
which implies that $(\eta+\mu)\om(v',w')=0$.  And as $\mu \neq -\eta$, this yields that $\om(v',w')=0$. Therefore
any one-step generalized eigenvector in $K_{\eta}$ is symplectically orthogonal to any eigenvector and any
one-step generalized eigenvector in $K_{\mu}$.  Again proceed inductively on the step of the generalized
eigenvector in $K_{\mu}$ and we'll obtain that any one-step generalized eigenvector in $K_{\eta}$ is
symplectically orthogonal to any vector in $K_{\mu}$.

Finally, continue this argument inductively on the step of the generalized eigenvector in $K_{\eta}$ and this will
yield that every vector in $K_{\eta}$ is symplectically orthogonal to every vector in $K_{\mu}$. Therefore
whenever $\mu \neq -\eta$, $K_{\eta}$ is symplectically orthogonal to $K_{\mu}$ as proposed. $\blacksquare$
\vspace{.1 in}

Note that in the above proposition, we make no assumption as to the value of $\eta$.  It is, in fact, true if
$\eta=0$.  However, if $\eta \neq 0$, then this leads immediately to a useful corollary.

\begin{cor}\label{H}
    If $\eta \neq 0$, then $K_{\eta}$ is symplectically orthogonal to itself.
\end{cor}

\begin{prop}\label{I}
    $\om$ is non-degenerate on $K_{\mu} \oplus K_{-\mu}$ for all $\mu \neq 0$.  In addition, $\om$ is
    non-degenerate on $K_{0}$.
\end{prop}

{\em Proof.} First, consider the case when $\mu \neq 0$.  Let $z \in K_{\mu} \oplus K_{-\mu}$ be such that
$\om(z,y) =0$ for all $y \in K_{\mu} \oplus K_{-\mu}$.  Then as $K_{\mu}$ and $K_{-\mu}$ are both symplectically
orthogonal to the rest of the space, we have that $\om(z,y)=0$ for all $y \in V$, which implies that $z=0$ because
$\om$ is non-degenerate on $V$.

For an analogous reason, $\om$ is non-degenerate on $K_{0}$. $\blacksquare$ \vspace{.1 in}

\begin{cor}\label{J}
    If $v \in K_{\mu}$ is nonzero, then there exists a nonzero $w \in K_{-\mu}$ such that $\om(v,w) \neq 0$.
\end{cor}

{\em Proof:} If $\mu = 0$, then this is a direct consequence of the proposition.  If $\mu \neq 0$, then by
Corollary \ref{H}, $K_{\mu}$ is symplectically orthogonal to itself.  Thus $\om$ is totally degenerate on
$K_{\mu}$. However, as $\om$ is non-degenerate on $K_{\mu} \oplus K_{-\mu}$, we have that for any nonzero $v \in
K_{\mu}$, there exists a $z \in K_{\mu} \oplus K_{-\mu}$ such that $\om(v,z) \neq 0$.  But we can write $z=u+w$,
where $u \in K_{\mu}$ and $w \in K_{-\mu}$.  Then as $\om(v,u)=0$, it must be that $\om(v,w) \neq 0$.
$\blacksquare$ \vspace{.1 in}

We can say a little more when $v \in K_{\mu}$ is an eigenvector.

\begin{prop}\label{K}
    If $v \in K_{\mu}$ is an eigenvector and $w \in K_{-\mu}$ such that $\om(v,w) \neq 0$, then $w$ is the end
    vector in any cycle of generalized eigenvectors to which it belongs.
\end{prop}

{\em Proof.} Assume to the contrary that there is a cycle of generalized eigenvectors such that $w$ is not the end
vector.  This implies that there is a generalized eigenvector $w'$ in the cycle such that $aw' = -\mu w'+w$. Then
we have that
\begin{align}
     \mu \om(v,w') &= \om(\mu v, w') = \om(a v, w') = -\om(v, a w') = -\om(v, -\mu w' + w) = \mu \om(v, w') -
     \om(v, w) \notag
\end{align}
But this implies that $-\om(v,w)=0$, which is a contradiction to the fact that $\om(v,w) \neq 0$.  Thus $w$ is the
end vector in any cycle of generalized eigenvectors to which it belongs. $\blacksquare$ \vspace{.1 in}

We now investigate what happens when two eigenvectors $v \in K_{\mu}$ and $w \in K_{-\mu}$ are such that $\om(v,w)
\neq 0$.  We already know by the previous proposition that any cycle of generalized eigenvectors to which they
belong must necessarily have length 1.  But there is considerably more that we can say as well.  For convenience
in the following proposition, let $(\cdot)$ denote the span of the vectors inside and $(\cdot)^{\perp}$ denote the
symplectically orthogonal complement.

\begin{prop}\label{L}
    If $v \in K_{\mu}$ and $w \in K_{-\mu}$ are eigenvectors such that $\om(v,w) \neq 0$.  Then
    \begin{enumerate}
        \item $(v,w)^{\perp} \oplus (v,w) = V.$
        \item $(v,w)^{\perp}$ is $a$-invariant.
        \item $\om$ is non-degenerate on $(v,w)^{\perp}$.
    \end{enumerate}
\end{prop}

{\em Proof.  1.} First, we show that $(v,w)^{\perp} \bigcap (v,w) =
0$.  Let $z \in (v,w)^{\perp} \bigcap (v,w)$, then $z = c_{1}v +
c_{2} w$ for some $c_{1},c_{2} \in \R$.  This yields that $z -
c_{1}v - c_{2}w = 0$, which implies that $\om(z-c_{1}v-c_{2}w,y)=0$
for all $y \in V$.  In particular, as $\om(z,w) = \om(z,v)=0$, we
have
\begin{align}
    0 &= \om(z-c_{1}v-c_{2}w,w) = \om(z,w) - c_{1}\om(v,w) - c_{2} \om(w,w) = -c_{1}\om(v,w), \notag \\
    0 &= \om(z-c_{1}v-c_{2}w, v) = \om(z,v) - c_{1}\om(v,v) - c_{2} \om(w,v) = -c_{2}\om(w,v). \notag
\end{align}
As $\om(v,w) \neq 0$, this implies that $c_{1}=c_{2}=0$ and hence $z = 0$.

Next, we show that $(v,w)^{\perp} + (v,w) = V$.  Let $z \in V$ be arbitrary.  As $\om(v,w) \neq 0$, we can let
$\tilde{z} = z - \frac{\om(z,w)}{\om(v,w)} v - \frac{\om(z,v)}{\om(w,v)} w$.  Then
\begin{align}
    \om(\tilde{z},v) &= \om\pnth{z - \frac{\om(z,w)}{\om(v,w)} v - \frac{\om(z,v)}{\om(w,v)} w,v} = \om(z,v) -
    \frac{\om(z,w)}{\om(v,w)} \om(v,v) - \frac{\om(z,v)}{\om(w,v)} \om(w,v) = 0, \notag \\
\intertext{and}
    \om(\tilde{z},w) &= \om\pnth{z - \frac{\om(z,w)}{\om(v,w)} v - \frac{\om(z,v)}{\om(w,v)} w, w} =
    \om(z,w) - \frac{\om(z,w)}{\om(v,w)} \om(v,w) - \frac{\om(z,v)}{\om(w,v)} \om(w,w) = 0. \notag
\end{align}
Thus $\tilde{z} \in (v,w)^{\perp}$ and we can write $z$ as $$z = \tilde{z} + \frac{\om(z,w)}{\om(v,w)} v +
\frac{\om(z,v)}{\om(w,v)} w.$$  Therefore $(v,w)^{\perp} + (v,w) = V$. \vspace{.1 in}

{\em 2.} Let $z \in (v,w)^{\perp}$.  Then
\begin{align}
    \om(a z, v) &= -\om(z, a v) = -\om(z, \mu v) = -\mu \om(z,v) = 0, \notag \\
    \om(a z, w) &= -\om(z, a w) = -\om(z, -\mu w) = \mu \om(z, w) = 0. \notag
\end{align}
Thus $az \in (v,w)^{\perp}$. \vspace{.1 in}

{\em 3.} Let $z \in  (v,w)^{\perp}$ be such that $\om(z,y)=0$ for all $y \in (v,w)^{\perp}$.  Then as
$\om(z,v)=\om(z,w)=0$, we have that $\om(z,y)=0$ for all $y \in (v,w)$.  But $V= (v,w)^{\perp} \oplus (v,w)$, so
$\om(z,y) = 0$ for all $y \in V$.  Thus $z =0$ by the non-degeneracy of $\om$ on $V$. $\blacksquare$ \vspace{.1
in}

This leads to a useful corollary.

\begin{cor}\label{M}
    If $v,w$ is a cycle of generalized eigenvectors in $K_{0}$ such that $\om(v,w) \neq 0$, then the conclusion of
    Proposition \ref{L} holds.
\end{cor}

{\em Proof.} Parts {\em 1.} and {\em 3.} hold for exactly the same reasons.  It suffices to prove that {\em 2.}
holds in this case. Note that if $v,w$ is a cycle of generalized eigenvectors such that $\om(v,w) \neq 0$, then
$v$ and $w$ belong to a generalized eigenspace corresponding to the eigenvalue 0.

Let $z \in (v,w)^{\perp}$.  Then $$\om(a z, v) = -\om(z, a v) = -\om(z, 0) = 0,$$ and $$\om(a z, w) = -\om(z, a w)
= -\om(z, v) = 0.$$  Thus $az \in (v,w)^{\perp}$ and we see that {\em 2.} holds. $\blacksquare$ \vspace{.1 in}

Note that Parts {\em 1.} and {\em 3.} of Proposition \ref{L} actually hold for any pair of vectors, $v,w \in V$
such that $\om(v,w) \neq 0$.  We can then prove quite easily the following proposition.

\begin{prop}\label{P}
    The dimension of $K_{0}$ is even and $\dim K_{\mu} = \dim
    K_{-\mu}$ for all $\mu \neq 0$.
\end{prop}

{\em Proof.} Let $V=K_{0}$.  If $\dim K_{0}=0$, then the result is
trivial. Assume then that $\dim K_{0} >0$.  We also have that $\dim
K_{0}\neq 1$, otherwise $\om$ is degenerate on $K_{0}$.  Thus $\dim
K_{0} \geq 2$.

Let $v_{1} \in K_{0}$ be nonzero, then as $\om$ is non-degenerate on $K_{0}$, there exists a nonzero $w_{1} \in
K_{0}$ such that $\om(v_{1},w_{1}) \neq 0$. Then as parts {\em 1.} and {\em 3.} of Proposition \ref{L} hold for
$v_{1}$ and $w_{1}$, we have that
$$K_{0}=(v_{1},w_{1}) \oplus (v_{1},w_{1})^{\perp}$$ and that $\om$ is non-degenerate on
$(v_{1},w_{1})^{\perp}$.  Also note that \linebreak[0] $\dim
(v_{1},w_{1})^{\perp} = \dim K_{0} - 2$.  If $\dim
(v_{1},w_{1})^{\perp} \linebreak[0] =0$ then we're done. If
$\dim(v_{1},w_{1})^{\perp} \neq 0$, then it must be that
$\dim(v_{1},w_{1})^{\perp} \geq 2$, otherwise $\om$ is degenerate on
$K_{0}$.

Let $v_{2} \in (v_{1},w_{1})^{\perp}$ be nonzero, then as $\om$ is non-degenerate on $(v_{1},w_{1})^{\perp}$,
there exists a nonzero $w_{2} \in (v_{1},w_{1})^{\perp}$ such that $\om(v_{2},w_{2}) \neq 0$.  Then, because of
parts {\em 1.} and {\em 3.} of Proposition \ref{L}, we can write
$$K_{0}=(v_{1},w_{1}) \oplus (v_{2},w_{2}) \oplus
(v_{1},w_{1},v_{2},w_{2})^{\perp}.$$  We continue this process until it terminates.  It is guaranteed to terminate
because $\dim K_{0} < \infty$. Thus $K_{0}$ is even dimensional.

The proof that $\dim K_{\mu} = \dim K_{-\mu}$ for all $\mu \neq 0$ is almost identical except that the $v_{i}$
vectors come from $K_{\mu}$ and the $w_{i}$ vectors come from $K_{-\mu}$. Then $K_{\mu} \oplus K_{-\mu}$
decomposes into a direct sum of subspaces each spanned by two vectors, one from $K_{\mu}$ and one from $K_{-\mu}$.
This puts a basis from $K_{\mu}$ into one-to-one correspondence with a basis for $K_{-\mu}$ and yields the result.
$\blacksquare$ \vspace{.1 in}

A final lemma that will prove useful, particularly in finding the canonical forms of $a \in \si(4,\R)$, is the
following.

\begin{lem}\label{N}
    Any matrix in the following sets of matrices in $\si(2,\R)$ or
    $\si(4,\R)$ is {\em not} symplectically similar to the other matrix
    in the set to which it belongs.
    \begin{align}
        (1) & \left\{ a_{1}= \begin{pmatrix} 0&1 \\ 0&0 \end{pmatrix}, \ \ a_{2}= \begin{pmatrix} 0&-1
        \\ 0&0 \end{pmatrix} \right\} \notag \\
        (2) & \left\{ a_{3}= \begin{pmatrix} 0&\mu \\ -\mu&0 \end{pmatrix}, \ \ a_{4}=\begin{pmatrix}
        0&-\mu \\ \mu&0 \end{pmatrix} \right\} \notag \displaybreak[0] \\
        (3) & \left\{ a_{5}=\begin{pmatrix} 0&1&0&0 \\ 0&0&0&1 \\
        0&0&0&0 \\ 0&0&-1&0 \end{pmatrix}, \ \ a_{6}=\begin{pmatrix} 0&1&0&0 \\ 0&0&0&-1 \\ 0&0&0&0 \\
        0&0&-1&0 \end{pmatrix} \right\} \notag \\
        (4) & \left\{ a_{7}=\begin{pmatrix} 0&\mu&1&0 \\ -\mu&0&0&1 \\ 0&0&0&\mu \\ 0&0&-\mu&0 \end{pmatrix}, \
        \ a_{8}=\begin{pmatrix} 0&\mu&-1&0 \\ -\mu&0&0&-1 \\ 0&0&0&\mu \\
        0&0&-\mu&0 \end{pmatrix} \right\} \notag
    \end{align}
\end{lem}

{\em Proof.} First we prove that $\si(2,\R) = \sll(2,\R)$.  In two dimensions, $J=\begin{pmatrix} 0&1 \\
-1&0 \end{pmatrix}$.  Let $b = \begin{pmatrix} c_{1}&c_{2}
\\ c_{3}&c_{4} \end{pmatrix}$ be an arbitrary element of $\gl(2,\R)$.
Then forcing $b$ to satisfy the defining symplectic condition would yield the equation
$$\begin{pmatrix} 0&0 \\ 0&0 \end{pmatrix}= \begin{pmatrix} c_{1}&c_{3} \\ c_{2}&c_{4} \end{pmatrix}
\begin{pmatrix} 0&1 \\ -1&0 \end{pmatrix} + \begin{pmatrix} 0&1 \\
-1&0 \end{pmatrix}\begin{pmatrix} c_{1}&c_{2} \\ c_{3}&c_{4}
\end{pmatrix} = \begin{pmatrix} 0&c_{1}+c_{4} \\
-c_{1}-c_{4}& 0 \end{pmatrix}$$  Solving this equation, we find that
$c_{4}=-c_{1}$, that is, $b$ must be a trace free matrix. Hence
$\si(2,\R) = \sll(2,\R)$ and consequently, $Sp(2,\R) = SL(2,\R)$.

To prove that the two matrices in set (1) are not symplectically similar to each other first we'll call the two
matrices $a_{1}$ and $a_{2}$ respectively.  Then let $A \in GL(2,\R)$ be arbitrary. Solving the equation
$A^{-1}a_{1}A=a_{2}$ is equivalent to solving $a_{1}A
= Aa_{2}$ as long as we make the restriction that $\det A \neq 0$.  If $A=\begin{pmatrix} d_{1} & d_{2} \\
d_{3}& d_{4} \end{pmatrix}$, then this equation becomes
\begin{align}
     \begin{pmatrix} 0&1 \\ 0&0 \end{pmatrix} \begin{pmatrix} d_{1} & d_{2} \\ d_{3}& d_{4}
    \end{pmatrix} &= \begin{pmatrix} d_{1} & d_{2} \\ d_{3}& d_{4} \end{pmatrix} \begin{pmatrix} 0&-1
    \\ 0&0 \end{pmatrix} \notag \\
     \begin{pmatrix} d_{3}&d_{4} \\ 0&0 \end{pmatrix} &= \begin{pmatrix} 0&-d_{1} \\
    0&d_{3} \end{pmatrix} \notag
\end{align}
Thus in order to satisfy the equation $A^{-1}a_{1}A=a_{2}$, $A$
would have to have the form $\begin{pmatrix} d_{1}&d_{2} \\
0&-d_{1} \end{pmatrix}$ with $d_{1} \neq 0$, which always has negative determinant. As every element in
$Sp(2,\R)=SL(2,\R)$ has determinant 1, no symplectic group element can make the required move. Thus $a_{1}$ and
$a_{2}$ are not symplectically similar.

Next we prove that the matrices in set (2) are symplectically dissimilar.  Let $A \in GL(2,\R)$ be arbitrary just
as in the previous part. Again, the equation $A^{-1}a_{3}A = a_{4}$ is equivalent to solving $a_{3}A=Aa_{4}$ if we
require $\det A \neq 0$. Solving this equation, we obtain that $A$ must have the form $\begin{pmatrix} d_{1}&d_{2}
\\ d_{2}&-d_{1}
\end{pmatrix},$ which has the property that $\det A =
-\pnth{(d_{1})^2+(d_{2})^{2}}$.  This is a nonpositive value, which implies that $A \not \in Sp(2,\R)$ because
$\det A = 1$ whenever $A \in Sp(2,\R)$. Thus there is no symplectic change of basis that will, by conjugation,
move $a_{3}$ to $a_{4}$.

In order to prove, that the matrices in set (3) are symplectically dissimilar, we let $A \in GL(4,\R)$ be given by
the matrix
$$A=\begin{pmatrix} d_{1}&d_{2}&d_{3}&d_{4} \\
d_{5}&d_{6}&d_{7}&d_{8} \\ d_{9}&d_{10}&d_{11}&d_{12} \\
d_{13}&d_{14}&d_{15}&d_{16} \end{pmatrix}$$ with $\det A \neq 0$. We consider again the equation $a_{5}A=Aa_{6}$.
Solving this equation, we obtain the following
\begin{align}
    \begin{pmatrix} 0&1&0&0 \\ 0&0&0&1 \\
    0&0&0&0 \\ 0&0&-1&0 \end{pmatrix} \begin{pmatrix} d_{1}&d_{2}&d_{3}&d_{4} \\
    d_{5}&d_{6}&d_{7}&d_{8} \\ d_{9}&d_{10}&d_{11}&d_{12} \\ d_{13}&d_{14}&d_{15}&d_{16} \end{pmatrix} &=
    \begin{pmatrix} d_{1}&d_{2}&d_{3}&d_{4} \\ d_{5}&d_{6}&d_{7}&d_{8} \\ d_{9}&d_{10}&d_{11}&d_{12} \\
    d_{13}&d_{14}&d_{15}&d_{16} \end{pmatrix} \begin{pmatrix} 0&1&0&0 \\ 0&0&0&-1 \\ 0&0&0&0 \\ 0&0&-1&0
    \end{pmatrix} \notag \\
    \begin{pmatrix} d_{5}&d_{6}&d_{7}&d_{8} \\ d_{13}&d_{14}&d_{15}&d_{16} \\ 0&0&0&0 \\
    -d_{9}&-d_{10}&-d_{11}&-d_{12} \end{pmatrix} &= \begin{pmatrix} 0&d_{1}&-d_{4}&-d_{2} \\ 0&
    d_{5}&-d_{8}&-d_{6} \\ 0&d_{9}&-d_{12}&-d_{10} \\ 0&d_{13}&-d_{16}&-d_{14} \end{pmatrix} \notag
\end{align}
This implies that any invertible linear transformation $A$ that by conjugation will move $a_{5}$ to $a_{6}$ is of
the form $$\begin{pmatrix} d_{1}&d_{2}&d_{3}&d_{4} \\ 0&d_{1}&-d_{4}&-d_{2} \\ 0&0&-d_{1}&0 \\
0&0&d_{2}&-d_{1} \end{pmatrix}.$$  If we require this matrix to satisfy $A^{t}JA = J$ so that it is also in the
symplectic group, then we get, aside from other things, that $(d_{1})^{2} = -1$ and hence there is no real
symplectic matrix that will, by conjugation, move $a_{5}$ to $a_{6}$.

Finally, to prove that the matrices in set (4) are symplectically dissimilar, we let $A \in GL(4,\R)$ be arbitrary
as above. We consider again the equation $a_{7}A=Aa_{8}$.  Solving this equation, we have that any invertible
matrix $A$ that by conjugation will move $a_{7}$ to $a_{8}$ is of the form
$$\begin{pmatrix} d_{1}&d_{2}&d_{3}&d_{4} \\
-d_{2}&d_{1}&-d_{4}&d_{3}
\\ 0&0&-d_{1}&-d_{2} \\ 0&0&d_{2}&-d_{1} \end{pmatrix}.$$  Requiring this to satisfy $A^{t}JA = J$, will require,
amongst other things, that $(d_{1})^{2} + (d_{2})^{2} = -1$, which is impossible as $d_{1},d_{2} \in \R$.
Therefore the two matrices are symplectically dissimilar. $\blacksquare$ \vspace{.1 in}

We now have the tools necessary to compute the canonical forms of $a \in \si(4,\R)$.  We will consider first those
$a$ that have real eigenvalues and then those where one or more of the eigenvalues are complex.  We'll classify
according to the decomposition of $V$ into generalized eigenspaces of $a$ when the eigenvalues are real and
$V_{c}$ when some or all of them are complex.  The cases will follow these nine decompositions of $V$ or $V_{c}$.
\begin{enumerate}
    \item $V=K_{\la} \oplus K_{-\la} \oplus K_{\mu} \oplus K_{-\mu}$
    where $\la,\mu \in \R$ are both nonzero and $\mu \neq \pm \la$.
    \item $V=K_{\la} \oplus K_{-\la}$ where $\la \in \R$ is nonzero.
    \item $V=K_{\la} \oplus K_{-\la} \oplus K_{0}$ where $\la \in \R$ is nonzero.
    \item $V=K_{0}$.
    \item $V_{c}=K_{\la} \oplus K_{-\la} \oplus K_{\mu i} \oplus
    K_{-\mu i}$ where $\la,\mu \in \R$ are both nonzero.
    \item $V_{c}=K_{0} \oplus K_{\mu i} \oplus K_{-\mu i}$ where
    $\mu \in \R$ is nonzero.
    \item $V_{c}=K_{z} \oplus K_{-z} \oplus K_{\overline{z}} \oplus
    K_{-\overline{z}}$ where $z=\la + \mu i$ for some nonzero $\la, \mu \in \R$.
    \item $V_{c}=K_{\mu i} \oplus K_{-\mu i} \oplus K_{\eta i}
    \oplus K_{-\eta i}$ where $\mu, \eta \in \R$ are both nonzero
    and $\eta \neq \pm \mu$.
    \item $V_{c}=K_{\mu i} \oplus K_{-\mu i}$ where $\mu \in \R$ is
    nonzero.
\end{enumerate}
This enumerates every possible decomposition of $V$ or $V_{c}$ into generalized eigenspaces of $a \in \si(4,\R)$.

\subsubsection{Case 1: $V=K_{\la} \oplus K_{-\la} \oplus K_{\mu} \oplus
K_{-\mu}$}

We also assume in this case that $\la,\mu \in \R$ are both nonzero and $\mu \neq \pm \la$.  As $\dim V = 4$ it
must be that each of these generalized eigenspaces are one dimensional. This implies, of course, that $a$ is
diagonalizable.

Then pick eigenvectors $v_{1} \in K_{\la}$ and $w_{1} \in K_{-\la}$ such that $\om(v_{1},w_{1}) \neq 0$ and
decompose them out as in Proposition \ref{L}.  Let $v_{2} \in K_{\mu}$ and $w_{2} \in K_{-\mu}$.  Then
$\{v_{2},w_{2}\}$ is a basis for $(v_{1},w_{1})^{\perp}$ and $\om(v_{2},w_{2}) \neq 0$ by the non-degeneracy of
$\om$ on $(v_{1},w_{1})^{\perp}$.

Make the additional change that if $\om(v_{1},w_{1}) <0$, replace
$w_{1}$ with $-w_{1}$ and relabel it as $w_{1}$.  In this way, we
can ensure that $\om(v_{1},w_{1}) >0$.  Do the same for
$\om(v_{2},w_{2})$. Then we pick a basis, $\be$, for $V$ consisting
of eigenvectors in the following way, $$\be=\left\{
\frac{1}{\sqrt{\om(v_{1},w_{1})}} v_{1}, \
\frac{1}{\sqrt{\om(v_{2},w_{2})}} v_{2}, \
\frac{1}{\sqrt{\om(v_{1},w_{1})}} w_{1}, \
\frac{1}{\sqrt{\om(v_{2},w_{2})}} w_{2}\right\}$$ and relabel the
vectors $v_{1}'$, $v_{2}'$, $w_{1}'$, and $w_{2}'$. Then $\be$ has
the properties that $\om(v_{i}',v_{j}')=\om(w_{i}',w_{j}')=0$ and
$\om(v_{i}',w_{j}') = \del_{ij}$. Thus the matrix $$W=
\begin{pmatrix} v_{1}' & v_{2}' & w_{1}' & w_{2}'
\end{pmatrix}$$ is in $Sp(4,\R)$ as $W^{t} J W = J$.  In this basis,
$a$ is of form (1) presented in Theorem \ref{E} with $\la\neq \pm \mu$

\subsubsection{Case 2: $V=K_{\la} \oplus K_{-\la}$}

In this case, we also assume that $\la \in \R$ is nonzero.  Note that as $\dim K_{\la} = \dim K_{-\la}$, both of
these generalized eigenspaces are of dimension two.

If $a$ is diagonalizable, then we proceed exactly as in Case 1 with
$\mu = \pm \la$ depending on how the eigenvectors are to be ordered.
Then the real symplectic canonical form of $a$ is of form (1)
presented in Theorem \ref{E} with $\mu = \pm \la$.

On the other hand, if $a$ is not diagonalizable, then either $K_{\la}$ contains a cycle of generalized
eigenvectors of length two or $K_{-\la}$ does.

\begin{lem}\label{Q}
    $K_{\la}$ contains a cycle of generalized eigenvectors of
    length two if and only if $K_{-\la}$ contains one.
\end{lem}

{\em Proof:} First, let $v,v' \in K_{\la}$ be a cycle of length two. Let $w \in K_{-\la}$ be an eigenvector, then
by Proposition \ref{K}, $\om(w,v)=0$.  Thus if $K_{-\la}$ contains only eigenvectors, then $\om$ is degenerate as
$v$ is symplectically orthogonal to the entire space.  As $\dim K_{-\la}=2$, it must then contain a cycle of
generalized eigenvectors of length two.  The converse is analogous. $\blacksquare$ \vspace{.1 in}

As $K_{\la}$ and $K_{-\la}$ both contain a cycle of length two, let $\{v_{1},v_{2}\}$ be a cycle of generalized
eigenvectors in $K_{\la}$ and $\{w_{1},w_{2}\}$ be one in $K_{-\la}$.  Then each these sets forms a basis for its
respective space.  In addition, as $v_{1}$ is an eigenvector in $K_{\la}$ and $w_{1}$ is not the end vector in its
cycle, we have by Proposition \ref{K} that $\om(v_{1},w_{1})=0$, and by Corollary \ref{H},
$\om(v_{1},v_{2})=\om(w_{1},w_{2})=0$.  This implies, as $\om$ is non-degenerate on $K_{\la} \oplus K_{-\la}$,
that $\om(v_{1},w_{2}) \neq 0$ and $\om(v_{2},w_{1}) \neq 0$.  It is, however, possible that $\om(v_{2},w_{2})
\neq 0$.  In that case, consider the cycle of generalized eigenvectors $\left\{w_{1}, w_{2} -
\frac{\om(v_{2},w_{2})}{\om(v_{2},w_{1})} w_{1}\right\}$.  If we relabel the elements of this cycle as $w_{1}'$
and $w_{2}'$, then it has the properties that $\om(v_{2},w_{2}')=0$ and $\om(v_{1},w_{2}')=\om(v_{1},w_{2})$.  In
addition, we have the following
\begin{align}
     \la \om(v_{2},w_{2}') + \om(v_{1},w_{2}') &= \om(\la v_{2} + v_{1}, w_{2}') = \om(a v_{2}, w_{2}')
     = - \om(v_{2}, a w_{2}') \notag \\
     &= - \om(v_{2}, -\la w_{2}' + w_{1}') = \la \om(v_{2},w_{2}') - \om(v_{2},w_{1}') \notag
\end{align}
which implies that $\om(v_{1},w_{2}') = -\om(v_{2},w_{1}')$.  Thus
if $\om(v_{1},w_{2}') <0$, then $\om(v_{2},w_{1}') >0$.  If this is
the case, then replace $w_{2}'$ with $-w_{2}'$ and relabel to
guarantee that $\om(v_{1},w_{2}')>0$ as well.  Note now that $a
w_{2}' =-\la w_{2}' -w_{1}'$.  If $\om(v_{1},w_{2}') >0$ then
$\om(v_{2},w_{1}')<0$.  In this case, replace $w_{1}'$ with
$-w_{1}'$ and relabel to guarantee that $\om(v_{2},w_{1}') >0$.
Again, this will make $a w_{2}' = -\la w_{2}' - w_{1}'$. After
making either of these changes, we pick a new basis, $\be$, for $V$
in the following way $$\be =\left\{
\frac{1}{\sqrt{\om(v_{1},w_{2}')}} v_{1},
\frac{1}{\sqrt{\om(v_{2},w_{1}')}} v_{2},
\frac{1}{\sqrt{\om(v_{1},w_{2}')}} w_{2}',
\frac{1}{\sqrt{\om(v_{2},w_{1}')}} w_{1}'\right\}$$ and relabel the
vectors $v_{1}''$, $v_{2}''$, $w_{1}''$, and $w_{2}''$ respectively.
Then the matrix $$W=\begin{pmatrix} v_{1}'' & v_{2}'' & w_{1}'' &
w_{2}'' \end{pmatrix}$$ is such that $W^{t} J W = J$ and is in
$Sp(4,\R)$.  In addition, in this basis, $a$ is of form (3)
presented in Theorem \ref{E} with $\la \neq 0$.

\subsubsection{Case 3: $V=K_{\la} \oplus K_{-\la} \oplus K_{0}$}

We also assume that $\la \in \R$ is nonzero.  By Proposition \ref{P}, $\dim K_{0}$ is even. Then we have here that
$\dim K_{0}=2$ and $\dim K_{\la} = \dim K_{-\la} = 1$.

If $a$ is diagonalizable, then these generalized eigenspaces consist
only of eigenvectors and we proceed as in Case 1.  In this case, the
real symplectic canonical form of $a$ is form (1) from Theorem
\ref{E} with $\mu = 0$ and $\la \neq 0$.

If $a$ is not diagonalizable, then $K_{0}$ contains a cycle of
generalized eigenvectors of length two.  Let $v_{1}$ and $w_{1}$ be
eigenvectors in $K_{\la}$ and $K_{-\la}$ respectively.  Then
$\om(v_{1},w_{1}) \neq 0$ and we can split them off as in
Proposition \ref{L}. Pick a basis for $(v_{1},w_{1})^{\perp}=K_{0}$,
$\{v_{2}',w_{2}'\}$, such that the two vectors form a cycle of
generalized eigenvectors of length two.  We know that
$\om(v_{2}',w_{2}') \neq 0$ because $\om$ is non-degenerate.  As
with the previous case if $\om(v_{1},w_{1}) < 0$, then replace
$w_{1}$ with $-w_{1}$ and relabel.  However, if $\om(v_{2}',w_{2}')
<0$, we can't replace $w_{2}'$ so trivially.  This is because if we
change the sign of $w_{2}'$, then $\om(v_{2}',w_{2}') >0$ but $a
w_{2}' = -v_{2}'$ which by Lemma \ref{N} is symplectically
dissimilar to the case when $\om(v_{2}',w_{2}') >0$ in the first
place.  However, we still make the change so that
$\om(v_{2}',w_{2}')>0$.

Then we can pick a basis $\be$ given by
$$\be= \left\{ \frac{1}{\sqrt{\om(v_{1},w_{1})}} v_{1}, \ \frac{1}{\sqrt{\om(v_{2}',w_{2}')}} v_{2}', \
\frac{1}{\sqrt{\om(v_{1},w_{1})}} w_{1}, \
\frac{1}{\sqrt{\om(v_{2}',w_{2}')}} w_{2}' \right\}$$ and relabel
the vectors $v_{1}''$, $v_{2}''$, $w_{1}''$, and $w_{2}''$
respectively. Then $\be$ has the properties that
$\om(v_{i}'',v_{j}'')=\om(w_{i}'',w_{j}'')=0$ and
$\om(v_{i}'',w_{j}'') = \del_{ij}$. Then the matrix $$W=
\begin{pmatrix} v_{1}'' & v_{2}'' & w_{1}'' & w_{2}'' \end{pmatrix}$$
is in $Sp(4,\R)$ as $W^{t} J W = J$.  In this basis, $a$ is of form
(2) presented in Theorem \ref{E} with $\la \neq 0$.

\subsubsection{Case 4: $V=K_{0}$}

In this case, $K_{0}$ is obviously of dimension four.

If $a$ is diagonalizable, then, as its only eigenvalue is 0, $a$ is
the zero transformation and is already in real symplectic canonical
form, which is form (1) in Theorem \ref{E} with $\la=\mu=0$.

If $a$ is not diagonalizable, then one of the following statements
is true.
\begin{enumerate}
    \item There is a basis for $K_{0}$ consisting of three cycles
    of generalized eigenvectors, two of length one and one of length
    two.
    \item There is a basis for $K_{0}$ consisting of two cycles of
    generalized eigenvectors of length two.
    \item There is a basis for $K_{0}$ consisting of one cycle of
    generalized eigenvectors of length four.
\end{enumerate}
A basis for $K_{0}$ consisting of a cycle of length three and a
cycle of length one is impossible.  To show this, let $\{v,v',v''\}$
be a cycle of length three and let $\{w\}$ be a cycle of length one
such that $\{v,v',v'',w\}$ is a basis for $K_{0}$. Then, by
Proposition \ref{K} and the non-degeneracy of $\om$, $\om(v'',w)
\neq 0$, and $\om(v,v'') \neq 0$. Moreover, $w$ and $v$ are
symplectically orthogonal to everything else in the space and the
nonzero vector $v - \frac{\om(v,v'')}{\om(w,v'')} w$ is
symplectically orthogonal to the entire space, which is impossible.

If the first possibility is true, let a basis for $V$ be given by
$\{v_{1}, w_{1}, v_{2}, w_{2}\}$ where $v_{1}$ and $w_{1}$ are
cycles of generalized eigenvectors of length one and $v_{2}, w_{2}$
is a cycle of length two. Then $\om(v_{1},v_{2})=0$ and
$\om(w_{1},v_{2})=0$ because $v_{2}$ is not the end vector in its
cycle.  This implies that if $\om(v_{1},w_{1}) = 0$, then
$\om(v_{1},w_{2}) \neq 0$ because $\om$ is non-degenerate.  But this
tells us that the nonzero vector $u=v_{2} -
\frac{\om(v_{2},w_{2})}{\om(v_{1}, w_{2})} v_{1}$ is symplectically
orthogonal to the entire space.  This, of course, can't happen and
hence $\om(v_{1},w_{1}) \neq 0$.  Then we follow the second argument
in Case 3 to obtain that the real symplectic canonical form of $a$
is form (2) in Theorem \ref{E} with $\la = 0$.

For the second possibility above, let $\{v_{1},v_{2},w_{1},w_{2}\}$ be a basis for $K_{0}$ such that $v_{1},v_{2}$
and $w_{1},w_{2}$ are cycles of generalized eigenvectors.  We know that $\om(v_{1},w_{1}) = 0$ by Proposition
\ref{K}. If $\om(v_{1},v_{2})=\om(w_{1},w_{2}) = 0$, then the second argument presented in Case 2 holds here as
well. This allows $\la = 0$ in form (3) in Theorem \ref{E}.

If one of these values is nonzero, then order the basis so that $\om(v_{1},v_{2}) \neq 0$.  Then $v_{1}$ and
$v_{2}$ decompose off as in Corollary \ref{M}. Let $\{w_{1}',w_{2}'\}$ be a cycle of generalized eigenvectors in
$(v_{1},v_{2})^{\perp}$. This implies that $\om(w_{1}',w_{2}') \neq 0$ because of the non-degeneracy of $\om$ on
$(v_{1},v_{2})^{\perp}$ and the fact that $\{w_{1}',w_{2}'\}$ form a basis for $(v_{1},v_{2})^{\perp}$.  This
gives us three possibilities,
\begin{itemize}
    \item $\om(v_{1},v_{2})>0$ and $\om(w_{1}',w_{2}')>0$.
    \item $\om(v_{1},v_{2})>0$ and $\om(w_{1}',w_{2}')<0$.  In which case, we replace $w_{2}'$ with $-w_{2}'$ and
    relabel it so that $\om(w_{1}',w_{2}')>0$.  Making this change will, however, result in $a w_{2}'=-w_{1}'$.
    This possibility includes the instance where $\om(v_{1},v_{2})<0$ and $\om(w_{1}',w_{2}')>0$ because we can
    just relabel the vectors appropriately to get the desired condition.
    \item $\om(v_{1},v_{2})<0$ and $\om(w_{1}',w_{2}')<0$.  In which case, we replace $v_{2}$ and $w_{2}'$ with
    $-v_{2}$ and $-w_{2}'$ respectively and relabel them so that $\om(v_{1},v_{2})>0$ and $\om(w_{1}',w_{2}')>0$.
    However, this will result in $a v_{2} = -v_{1}$ and $a w_{2}' = -w_{1}'$.
\end{itemize}
By Lemma \ref{N}, none of these three possibilities are symplectically similar to any other one.  However, the
second situation is symplectically similar to the canonical form presented in the last argument of Case 3 with
$\la = 0$.  This is because $a$,
in the basis the second case would produce, is of the form $$\begin{pmatrix} 0&0&1&0 \\ 0&0&0&-1 \\
0&0&0&0 \\ 0&0&0&0 \end{pmatrix}.$$ And if we can conjugate by the
matrix $$\begin{pmatrix} \frac{1}{2} &-\frac{1}{2}&0&\frac{1}{2} \\ 0&0&\frac{1}{2}&\frac{1}{2} \\0&0&1&-1 \\
-1&-1&0&1 \end{pmatrix},$$ which is in $Sp(2n,\R)$ because $A^{t}JA = J$, we obtain $$\begin{pmatrix} 0&1&0&0 \\
0&0&0&0 \\ 0&0&0&0 \\ 0&0&-1&0 \end{pmatrix}.$$  However, as the first and third possibilities are not
symplectically similar to the second one, they are not similar to the canonical form presented in Case 3 either.

Thus in each remaining case, after making its outlined correction, we pick a basis, $\be$ for $K_{0}$ in the
following way $$\be = \left\{ \frac{1}{\sqrt{\om(v_{1},v_{2})}} v_{1}, \ \frac{1}{\sqrt{\om(w_{1}',w_{2}')}}
w_{1}', \ \frac{1}{\sqrt{\om(v_{1},v_{2})}} v_{2}, \ \frac{1}{\sqrt{\om(w_{1}',w_{2}')}} w_{2}' \right\}$$ and
relabel the vectors $v_{1}'',w_{1}'',v_{2}'',w_{2}''$.  Then the matrix
$$W=\begin{pmatrix} v_{1}'' &w_{1}'' &v_{2}''& w_{2}''
\end{pmatrix}$$ is in $Sp(4,\R)$ because $W^{t}JW=J$.  In this
basis, $a$ is of form (4) in Theorem \ref{E}.

Finally, in the third instance, where $K_{0}$ has a basis consisting of a cycle of generalized eigenvector of
length four, we let $\{v_{1},v_{2},v_{3},v_{4}\}$ be a cycle of generalized eigenvectors.  Then this set also
defines a basis for $K_{0}$. We have, by Proposition \ref{K}, that $\om(v_{1},v_{2})=\om(v_{1},v_{3})=0$ and that
$\om(v_{1},v_{4}) \neq 0$.  We also note that
$$\om(v_{1},v_{4}) = \om(a v_{2}, v_{4})= - \om(v_{2},a v_{4}) = -
\om(v_{2}, v_{3})$$ implying as well that $\om(v_{2},v_{3}) \neq 0$ and $$\om(v_{2},v_{4}) = \om(a v_{3},v_{4}) =
-\om(v_{3}, a v_{4}) = -\om(v_{3},v_{3}) = 0.$$  Next, we pick a new basis in the following manner $$\{v_{1}, \
v_{2}, \ v_{3}- \frac{\om(v_{3},v_{4})}{2\om(v_{1},v_{4})}v_{1}, \ v_{4}- \frac{\om(v_{3}, v_{4})}{2\om(v_{1},
v_{4})} v_{2} \}$$ and relabel the vectors $v_{1}'$, $v_{2}'$, $v_{3}'$, and $v_{4}'$.  This new basis is still a
cycle of generalized eigenvectors with all the properties above, except that
\begin{align}
     \om(v_{3}',v_{4}') &= \om\left(v_{3}-\frac{\om(v_{3},v_{4})}{2\om(v_{1},v_{4})}v_{1}, v_{4}-
    \frac{\om(v_{3}, v_{4})}{2\om(v_{1}, v_{4})} v_{2} \right) \notag \\
     &= \om(v_{3},v_{4}) - \frac{\om(v_{3},v_{4})}{2\om(v_{1},v_{4})} \om(v_{3},v_{2}) -
    \frac{\om(v_{3},v_{4})}{2\om(v_{1},v_{4})} \om(v_{1},v_{4}) + \pnth{\frac{\om(v_{3},v_{4})}{2
    \om(v_{1},v_{4})}}^{2}\om(v_{1},v_{2}) \notag \\
     &= \om(v_{3},v_{4}) - \frac{\om(v_{3},v_{4})}{2\om(v_{3},v_{2})} \om(v_{3},v_{2}) -
    \frac{1}{2}\om(v_{3},v_{4}) \notag \\
     &= 0 \notag
\end{align}

If $\om(v_{1}',v_{4}')>0$, then $\om(v_{2}',v_{3}')<0$ and we replace $v_{3}'$ with $-v_{3}'$ and relabel.  This
will yield that $\om(v_{1}',v_{4}')=\om(v_{2}',v_{3}') > 0$, but that $av_{4}' = -v_{3}'$ and $av_{3}' = -v_{2}'$.
On the other hand, if $\om(v_{1}',v_{4}')<0$, then $\om(v_{2}',v_{3}')>0$ and we replace $v_{4}'$ with $-v_{4}'$
and relabel. This gives us that $\om(v_{1}',v_{4}')=\om(v_{2}',v_{3}') > 0$, but that $av_{4}' = -v_{3}'$.  We
know by Lemma \ref{N} that these two cases are not symplectically similar.  However, we make whatever change is
necessary to guarantee that $\om(v_{1}',v_{4}')=\om(v_{2}',v_{3})>0$.

In either case, we pick a basis, $\be$, for $V$ in the following way
$$\be=\left\{ \frac{1}{\sqrt{\om(v_{1}',v_{4}')}} v_{1}', \
\frac{1}{\sqrt{\om(v_{1}',v_{4}')}}v_{2}', \ \frac{1}{\sqrt{\om(v_{1}',v_{4}')}} v_{4}', \
\frac{1}{\sqrt{\om(v_{1}',v_{4}')}} v_{3}' \right\}$$ and relabel the vectors $v_{1}''$, $v_{2}''$, $v_{3}''$, and
$v_{4}''$.  Then the matrix $$W= \begin{pmatrix} v_{1}''& v_{2}''& v_{3}''& v_{4}''
\end{pmatrix}$$ is such that $W^{t}JW = J$ and hence is in
$Sp(4,\R)$.  In this basis, $a$ is of form (5) in Theorem \ref{E}.

\subsubsection{Case 5: $V_{c}=K_{\la} \oplus K_{-\la} \oplus K_{\mu i} \oplus
K_{-\mu i}$}

We also assume here that $\la,\mu \in \R$ are both nonzero.  In this case, we note that $a$ is diagonalizable over
$\C$ and that each of the generalized eigenspaces in the decomposition of $V_{c}$ is one dimensional.

Let $v_{1} \in K_{\la}$ and $w_{1} \in K_{-\la}$. Then
$\om(v_{1},w_{1}) \neq 0$, because $\om$ is non-degenerate on
$K_{\la} \oplus K_{-\la}$, and $v_{1}$ and $w_{1}$ split off as in
Proposition \ref{L}.  Pick $v_{2} \in K_{\mu i}$ and $w_{2} \in
K_{-\mu i}$ such that $w_{2}=\overline{v_{2}}$. Then
$\om(v_{2},w_{2}) \neq 0$ again because of the non-degeneracy of
$\om$ on $(v_{1},w_{1})^{\perp} = K_{\mu i} \oplus K_{-\mu i}$. We
pick two new vectors $v_{2}'=v_{2}+w_{2}$ and $w_{2}' =
-i(v_{2}-w_{2})$. The set $\{v_{2}',w_{2}'\}$ is still a basis for
$K_{\mu i} \oplus K_{-\mu i}$ but the vectors are real.  Moreover
the set $\{v_{1}, w_{1}, v_{2}', w_{2}'\}$ forms a basis for $V$
when considered over $\R$.  In addition,
\begin{align}
    \om(v_{2}',w_{2}') &= \om(v_{2}+w_{2},-i(v_{2}-w_{2})) \notag \\
    &= -i \om(v_{2},v_{2}) + i \om(v_{2},w_{2}) - i
    \om(w_{2},v_{2}) + i \om(w_{2},w_{2})
    = 2i \om(v_{2},w_{2})
    \neq 0 \notag
\end{align}
and $\om(v_{2}',w_{2}') \in \R$ because $v_{2}'$ and $w_{2}'$ are real vectors.  We also have that $$a v_{2}' =
a(v_{2}+w_{2}) = a v_{2} +a w_{2} = \mu i v_{2} - \mu i w_{2} = \mu i (v_{2}-w_{2}) = -\mu w_{2}'$$ and $$a w_{2}'
= -i (a v_{2} - a w_{2}) = -i (\mu i v_{2} + \mu i w_{2}) = -\mu i^{2} (v_{2}+w_{2}) = \mu (v_{2}+w_{2}) = \mu
v_{2}'.$$

If $\om(v_{1},w_{1}) <0$, then replace $w_{1}$ with $-w_{1}$ and relabel just as we've done before.  Also if
$\om(v_{2}',w_{2}')<0$, then we'll replace $w_{2}'$ with $-w_{2}'$ and relabel.  While this will make
$\om(v_{2}',w_{2}')>0$, it will also have the effect that $a v_{2}' = \mu w_{2}'$ and $a w_{2}' = -\mu v_{2}'$. We
know by Lemma \ref{N} that these two cases a symplectically dissimilar, but we nonetheless make the change so that
$\om(v_{2}',w_{2}')>0$.

In either instance, after ensuring that $\om(v_{1},w_{1})$ and $\om(v_{2}',w_{2}')$ are both positive, we pick a
new basis, $\be$, for $V$ in the following way $$\be = \left\{ \frac{1}{\sqrt{\om(v_{1},w_{1})}} v_{1}, \
\frac{1}{\sqrt{\om(v_{2}',w_{2}')}} v_{2}', \ \frac{1}{\sqrt{\om(v_{1},w_{1})}} w_{1}, \
\frac{1}{\sqrt{\om(v_{2}',w_{2}')}} w_{2}' \right\}$$ and relabel the vectors $v_{1}''$, $v_{2}''$, $w_{1}''$, and
$w_{2}''$.  Then the matrix $$W= \begin{pmatrix} v_{1}'' & v_{2}'' & v_{3}'' &v_{4}''
\end{pmatrix}$$ is such that $W^{t}JW=J$ and hence $W \in Sp(4,\R)$.
In this basis, $a$ is of form (6) in Theorem \ref{E} with $\la \neq 0$. We also assume $\mu >0$ to make a
distinction between the dissimilar cases.

\subsubsection{Case 6: $V_{c}=K_{0} \oplus K_{\mu i} \oplus K_{-\mu
i}$}

Also assume in this section that $\mu \in \R$ is nonzero.  Note that $\dim K_{0}=2$ and $\dim K_{\mu i} = \dim
K_{-\mu i}=1$ as well.

If $a$ is diagonalizable over $\C$, then we proceed as in Case 5,
except that $v_{1}$ and $w_{1}$ are eigenvectors in $K_{0}$ such
that $\om(v_{1},w_{1}) \neq 0$.  This will lead to form (6) in
Theorem \ref{E} with $\la =0$.

On the other hand, if $a$ is not diagonalizable over $\C$, then
there exists a basis for $K_{0}$ consisting of a cycle of
generalized eigenvectors of length two.  Let $\{v_{1},w_{1}\}$ be a
cycle of generalized eigenvectors in $K_{0}$.  Then this set also
forms a basis for $K_{0}$ and it must be that $\om(v_{1},w_{1}) \neq
0$ or $\om$ is degenerate on $K_{0}$. Thus the span of these two
vectors, $K_{0}$, decomposes out as in Corollary \ref{M} and note
that $(v_{1},w_{1})^{\perp}=K_{\mu i} \oplus K_{-\mu i}$. Next let
$v_{2} \in K_{\mu i}$ and let $w_{2} \in K_{-\mu i}$ be such that
$w_{2}=\overline{v_{2}}$.  Then $\om(v_{2},w_{2}) \neq 0$. Let
$v_{2}'=v_{2}+w_{2}$ and $w_{2}'=-i(v_{2}+w_{2})$.  Clearly $v_{2}'$
and $w_{2}'$ span $K_{\mu i} \oplus K_{-\mu i}$ and have the
property that $\om(v_{2}',w_{2}') \neq 0$.  Note that $$a v_{2}' =
-\mu w_{2}', \ \ \ \ \ a w_{2}' = \mu v_{2}'.$$  The vectors
$\{v_{1},w_{1},v_{2}',w_{2}'\}$ now form a basis for $V$ when
considered as a vector space over $\R$.

If $\om(v_{1},w_{1})<0$, then we replace $w_{1}$ with $-w_{1}$ and
relabel so that this value is guaranteed to be positive.  However,
this gives $w_{1}$ the property that $a w_{1} = -v_{1}$.  By Lemma
\ref{N}, we know that this is not symplectically similar to the case
where $\om(v_{1},w_{1})>0$ in the first place.  Similarly, if
$\om(v_{2}',w_{2}') <0$, then we replace $w_{2}'$ with $-w_{2}'$ and
relabel so that this values is positive. Making this change will
yield the equations $a v_{2}' = \mu w_{2}'$ and $a w_{2}' = -\mu
v_{2}'$.  And again by Lemma \ref{N}, this is not symplectically
similar to the case where $\om(v_{2}',w_{2}')>0$ in the first place.
After making the necessary changes, we pick a new basis, $\be$, for
$V$, as follows
$$\be = \left\{ \frac{1}{\sqrt{\om(v_{1},w_{1})}} v_{1}, \
\frac{1}{\sqrt{\om(v_{2}',w_{2}')}} v_{2}', \ \frac{1}{\sqrt{\om(v_{1},w_{1})}} w_{1}, \
\frac{1}{\sqrt{\om(v_{2}',w_{2}')}} w_{2}' \right\}$$ and relabel the vectors $v_{1}''$, $v_{2}''$, $w_{1}''$, and
$w_{2}''$.  Then the matrix $$W= \begin{pmatrix} v_{1}'' & v_{2}'' & v_{3}'' &v_{4}''
\end{pmatrix}$$ is such that $W^{t}JW=J$ and hence $W \in Sp(4,\R)$.
In this basis, $a$ is of form (7) in Theorem \ref{E}.  We also assume that $\mu >0$ to make a distinction between
the dissimilar cases.

\subsubsection{Case 7: $V_{c}=K_{z} \oplus K_{-z} \oplus K_{\overline{z}} \oplus
K_{-\overline{z}}$}

In this section, we also assume that $z=\la + \mu i$ for some $\la,
\mu \in \R$ such that $\la,\mu \neq 0$.  Also it is clear that $a$
is diagonalizable and that every generalized eigenspace of $a$ in
the decomposition of $V_{c}$ is one-dimensional.  Before we proceed
further, we present a lemma that will be helpful.

\begin{lem}\label{O}
    $\overline{\om(v,w)} = \om(\overline{v},\overline{w})$ for all $v,w \in V$.
\end{lem}

{\em Proof.} Let $v,w \in V$.  Then $\overline{\om(v,w)} =
\overline{v^{t}Jw} = \overline{v}^{t} J \overline{w} =
\om(\overline{v},\overline{w})$. $\blacksquare$ \vspace{.1 in}

Let $v_{1} \in K_{z}$, $v_{2} \in K_{\overline{z}}$, $w_{1} \in
K_{-z}$, and $w_{2} \in K_{-\overline{z}}$ be eigenvectors such that
$v_{2}=\overline{v_{1}}$ and $w_{2}=\overline{w_{1}}$. Then
$\{v_{1},v_{2},w_{1},w_{2}\}$ is a basis for $V_{c}$. In addition,
$\om(v_{1},w_{1}) \neq 0$ and $\om(v_{2},w_{2}) \neq 0$ while any
other pair in the basis is symplectically orthogonal. From this
basis for $V_{c}$, we construct a basis for $V$ in the following way
$$v_{1}' = v_{1}+v_{2}, \ \ v_{2}'=-i(v_{1}-v_{2}), \ \ w_{1}' = w_{1}+w_{2}, \ \ w_{2}'=-i(w_{1}-w_{2}).$$
Note that $$\overline{\om(v_{1},w_{1})}=
\om(\overline{v_{1}},\overline{w_{1}}) = \om(v_{2},w_{2}).$$ Let
$\Re(x)$ and $\Im(x)$ denote the real and imaginary parts of $x$
respectively. It follows then that
\begin{align}
     \om(v_{1}',w_{1}') &= \om(v_{1}+v_{2},w_{1}+w_{2}) = \om(v_{1},w_{1}) +\om(v_{1},w_{2}) +\om(v_{2},w_{1})
     +\om(v_{2},w_{2}) \notag \\
     &= \om(v_{1},w_{1}) + \om(v_{2},w_{2}) = 2 \Re(\om(v_{1},w_{1}))
     \notag \displaybreak[0]\\
\intertext{and}
     \om(v_{2}',w_{2}') &= \om(-i(v_{1}-v_{2}),-i(w_{1}-w_{2}))
     = -\om(v_{1}-v_{2},w_{1}-w_{2}) \notag \\
     &=
     -\om(v_{1},w_{1})+\om(v_{1},w_{2})+\om(v_{2},w_{1})-\om(v_{2},w_{2})
     \notag \\
     &= -(\om(v_{1},w_{1})+\om(v_{2},w_{2}))
     = -\om(v_{1}',w_{1}') \notag \displaybreak[0]\\
\intertext{Also, we have that}
     \om(v_{1}',w_{2}') &= \om(v_{1}+v_{2},-i(w_{1}-w_{2}))
     = -i\om(v_{1},w_{1})+i\om(v_{1},w_{2})-i\om(v_{2},w_{1})+i\om(v_{2},w_{2}) \notag \\
     &= -i \pnth{\om(v_{1},w_{1}) - \om(v_{2},w_{2})}
     = 2 \Im(\om(v_{1},w_{1})) \notag \displaybreak[0] \\
     \om(v_{2}',w_{1}') &= \om(-i(v_{1}-v_{2}),w_{1}+w_{2})
     = -i\om(v_{1},w_{1})-i\om(v_{1},w_{2})+i\om(v_{2},w_{1})+i\om(v_{2},w_{2}) \notag \\
     &= -i \pnth{ \om(v_{1},w_{1}) - \om(v_{2},w_{2})}
     = \om(v_{1}',w_{2}') \notag \displaybreak[0] \\
\intertext{and}
     \om(v_{1}',v_{2}') &= \om(v_{1}+v_{2},-i(v_{1}-v_{2})) = 0 \notag \\
     \om(w_{1}',w_{2}') &= \om(w_{1}+w_{2},-i(w_{1}-w_{2})) = 0
     \notag
\end{align}
Furthermore, we have that
\begin{align}
    a v_{1}' &= a(v_{1}+v_{2}) = (\la+\mu i)v_{1} + (\la -\mu i)v_{2}
    =\la(v_{1}+v_{2}) + \mu i (v_{1}-v_{2}) = \la v_{1}' -\mu
    v_{2}', \notag \\
    a v_{2}' &= a( -i(v_{1}-v_{2})) = -i(\la +\mu i) v_{1} + i(\la -\mu
    i)v_{2} = \mu(v_{1}+v_{2}) - \la i(v_{1}-v_{2}) = \mu v_{1}' +\la
    v_{2}', \notag
\end{align}
and similarly $$a w_{1}' = -\la w_{1}' +\mu w_{2}', \ \ \ \ \ a w_{2}' = -\mu w_{1}' -\la w_{2}'.$$

As $\om(v_{1},w_{1}) \neq 0$, we know that either
$\Re(\om(v_{1},w_{1})) \neq 0$, or $\Re(\om(v_{1},w_{1}))=0$ and
$\Im(\om(v_{1},w_{1})) \neq 0$.

If $\Re(\om(v_{1},w_{1})) \neq 0$, then let $$v_{1}'' = v_{1}' + \frac{\om(v_{1}',w_{2}')}{\om(v_{1}',w_{1}')}
v_{2}', \ \ v_{2}'' = v_{2} - \frac{\om(v_{1}',w_{2}')}{\om(v_{1}',w_{1}')} v_{1}', \ \ w_{1}'' = w_{1}', \ \
w_{2}''=w_{2}'.$$  The set $\{v_{1}'',v_{2}'',w_{1}'',w_{2}''\}$ is still a basis for $V$ and has the following
properties
\begin{align}
     \om(v_{1}'',w_{2}'') &= \om\pnth{v_{1}' + \frac{\om(v_{1}',w_{2}')}{\om(v_{1}',w_{1}')} v_{2}',
     w_{2}'} = \om(v_{1}',w_{2}') + \frac{\om(v_{1}',w_{2}')}{\om(v_{1}',w_{1}')} \om(v_{2}',w_{2}') \notag \\
     &= \om(v_{1}',w_{2}') - \frac{\om(v_{1}',w_{2}')}{\om(v_{1}',w_{1}')} \om(v_{1}',w_{1}') = 0 \notag
     \displaybreak[0] \\
     \om(v_{2}'',w_{1}'') &= \om\pnth{v_{2}' - \frac{\om(v_{1}',w_{2}')}{\om(v_{1}',w_{1}')} v_{1}',
     w_{1}'} = \om(v_{2}',w_{1}') - \frac{\om(v_{1}',w_{2}')}{\om(v_{1}',w_{1}')} \om(v_{1}',w_{2}') \notag \\
     &= \om(v_{2}',w_{1}') - \frac{\om(v_{2}',w_{1}')}{\om(v_{1}',w_{1}')} \om(v_{1}',w_{1}') = 0 \notag
     \displaybreak[0] \\
     \om(v_{1}'',w_{1}'') &= \om\pnth{v_{1}' + \frac{\om(v_{1}',w_{2}')}{\om(v_{1}',w_{1}')} v_{2}',
     w_{1}'} = \om(v_{1}',w_{1}') + \frac{\om(v_{1}',w_{2}')}{\om(v_{1}',w_{1}')} \om(v_{2}',w_{1}') \notag \\
     \om(v_{2}'',w_{2}'') &= \om\pnth{v_{2}' - \frac{\om(v_{1}',w_{2}')}{\om(v_{1}',w_{1}')} v_{1}',
     w_{2}'} = \om(v_{2}',w_{2}') - \frac{\om(v_{1}',w_{2}')}{\om(v_{1}',w_{1}')} \om(v_{1}',w_{2}') \notag \\
     &= -\om(v_{1},w_{1}') - \frac{\om(v_{1}',w_{2}')}{\om(v_{1}',w_{1}')} \om(v_{2}',w_{1}') =
     - \om(v_{1}'',w_{1}'') \notag
\end{align}
and $\om(v_{1}'',v_{2}'')=\om(w_{1}'',w_{2}'') = 0$. Clearly, as
$\om$ is non-degenerate, $\om(v_{1}'',w_{1}'') =
-\om(v_{2}'',w_{2}'') \neq 0$.  In addition,
\begin{align}
    a v_{1}'' &= a \pnth{v_{1}' + \frac{\om(v_{1}',w_{2}')}{\om(v_{1}',w_{1}')} v_{2}'} = \la v_{1}' -
    \mu v_{2}' + \frac{\om(v_{1}',w_{2}')}{\om(v_{1}',w_{1}')} (\mu v_{1}' + \la v_{2}') = \la v_{1}''
    - \mu v_{2}'', \notag \\
    a v_{2}'' &= a \pnth{v_{2}' - \frac{\om(v_{1}',w_{2}')}{\om(v_{1}',w_{1}')} v_{1}'}
    = \mu v_{1}' + \la v_{2}' - \frac{\om(v_{1}',w_{2}')}{\om(v_{1}',w_{1}')} (\la v_{1}' - \mu
    v_{2}) = \mu v_{1}'' + \la v_{2}'' \notag
\end{align}
and as $w_{1}'' = w_{1}'$ and $w_{2}'' = w_{1}'$, their relationship remains unchanged.

If $\Re(\om(v_{1},w_{1}))=0$, then $\Im(\om(v_{1},w_{1})) \neq 0$. This implies that
$\om(v_{1}',w_{1}')=\om(v_{2}',w_{2}') =0$ and $\om(v_{1}',w_{2}')=\om(v_{2}',w_{1}') \neq 0$. Let
$$v_{1}''=v_{1}'+v_{2}', \ \ v_{2}'' = v_{2}'-v_{1}', \ \ w_{1}''=w_{1}'+w_{2}', \ \ w_{2}''=w_{2}'-w_{1}'.$$  The
set $\{v_{1}'',v_{2}'',w_{1}'',w_{2}''\}$ is still a basis for $V$ and has the properties that
\begin{align}
     \om(v_{1}'',w_{2}'') &= \om(v_{1}'+v_{2}', w_{2}' - w_{1}') = \om(v_{1}',w_{2}') -\om(v_{1}',w_{1}') +
     \om(v_{2}',w_{2}') - \om(v_{2}',w_{1}') = 0 \notag \\
     \om(v_{2}'',w_{1}'') &= \om(v_{2}'-v_{1}',w_{1}'+w_{2}') = \om(v_{2}',w_{1}') +\om(v_{2}'+w_{2}') -
     \om(v_{1}',w_{1}') - \om(v_{1}',w_{2}') = 0 \notag \\
     \om(v_{1}'',w_{1}'') &= \om(v_{1}'+v_{2}',w_{1}'+w_{2}') = \om(v_{1}',w_{1}') + \om(v_{1}',w_{2}') +
     \om(v_{2}',w_{1}') + \om(v_{2}',w_{2}') = 2 \om(v_{1}',w_{2}') \notag \\
     \om(v_{2}'',w_{2}'') &= \om(v_{2}'-v_{1}',w_{2}'-w_{1}') = \om(v_{2}',w_{2}') - \om(v_{2}',w_{1}') -
     \om(v_{1}',w_{2}') + \om(v_{1}',w_{1}') = -\om(v_{1}'',w_{1}'') \notag
\end{align}
and $\om(v_{1}'',v_{2}'')=\om(w_{1}'',w_{2}'')=0$.  Then $\om(v_{1}'',w_{1}'')=-\om(v_{2}'',w_{2}'') \neq 0$.  In
addition,
\begin{align}
    av_{1}'' &= a (v_{1}'+v_{2}') = \la v_{1}' - \mu v_{2}' +\mu v_{1}'
    +\la v_{2}' = \la v_{1}'' - \mu v_{2}'', \notag \\
    a v_{2}'' &= a(v_{2}'-v_{1}') = \mu v_{1}' +\la v_{2}' - \la v_{1}' +\mu v_{2}' =
    \mu v_{1}'' +\la v_{2}'', \notag
\end{align}
and similarly $$a w_{1}'' = -\la w_{1}'' + \mu w_{2}'', \ \ \ \ \ a w_{2}'' = -\mu w_{1}'' -\la w_{2}''.$$

In either case, we know that $\om(v_{1}'',w_{1}'') = -\om(v_{2}'',w_{2}'') \neq 0$. If $\om(v_{1}'',w_{1}'') <0$,
then $\om(v_{2}'',w_{2}'')>0$ and we replace $w_{1}''$ with $-w_{1}''$ and relabel.  This will make both values
positive, but will yield that $a w_{1}'' = -\la w_{1}'' -\mu w_{2}''$ and $a w_{2}'' = \mu w_{1}''-\la w_{2}''$.
On the other hand, if $\om(v_{1}'',w_{1}'')>0$, then $\om(v_{2}'',w_{2}'') <0$ and we replace $w_{2}''$ with
$-w_{2}''$ and relabel. This will again make both values positive but will yield that $a w_{1}'' = -\la w_{1}''
-\mu w_{2}''$ and $a w_{2}'' = \mu w_{1}'' -\la w_{2}''$.

After this final change, we pick a new basis, $\be$, for $V$ as
follows $$\be= \left\{ \frac{1}{\sqrt{\om(v_{1}'',w_{1}'')}}
v_{1}'', \ \frac{1}{\sqrt{\om(v_{2}'',w_{2}'')}} v_{2}'', \
\frac{1}{\sqrt{\om(v_{1}'',w_{1}'')}} w_{1}'', \
\frac{1}{\sqrt{\om(v_{2}'',w_{2}'')}} w_{2}'' \right\}$$ and relabel
the vectors $v_{1}'''$, $v_{2}'''$, $w_{1}'''$, and $w_{2}'''$
respectively.  This implies that the matrix
$$W= \begin{pmatrix} v_{1}''' & v_{2}''' & v_{3}''' & v_{4}''' \end{pmatrix}$$ is in $Sp(4,\R)$ as $W^{t} J W
= J$.  In this basis, $a$ is of form (8) in Theorem \ref{E} with $\la \neq 0$.

\subsubsection{Case 8: $V_{c}=K_{\mu i} \oplus K_{-\mu i} \oplus K_{\eta i}
\oplus K_{-\eta i}$}

Also we assume that $\mu, \eta \in \R$ are both nonzero and $\eta
\neq \pm \mu$. Let $v_{1} \in K_{\mu i}$, then
$v_{2}=\overline{v_{1}} \in K_{-\mu i}$. We are guaranteed that
$\om(v_{1},v_{2}) \neq 0$. This pair then decomposes out according
to Proposition \ref{L}. However, note that $(v_{1},v_{2})^{\perp} =
K_{\eta i} \oplus K_{-\eta i}$. Let $w_{1} \in K_{\eta i}$, then
$w_{2}= \overline{w_{1}} \in K_{-\eta i}$ and we are again
guaranteed that $\om(w_{1},w_{2}) \neq 0$.  In addition, every pair
of vectors from the set $\{v_{1},v_{2},w_{1},w_{2}\}$ other than
$v_{1},v_{2}$ and $w_{1},w_{2}$ are symplectically orthogonal.
Moreover, this set forms a basis for $V_{c}$.

Let $$v_{1}'=v_{1}+v_{2}, \ \ \ v_{2}' = -i(v_{1}-v_{2}), \ \ \
w_{1}' = w_{1}+w_{2}, \ \ \ w_{2}' = -i(w_{1}-w_{2}),$$ then the set
$\{v_{1}',v_{2}',w_{1}',w_{2}'\}$ forms a basis for $V$ over $\R$
such that $$av_{1}' = -\mu v_{2}', \ \ \ a v_{2}' = \mu v_{1}', \ \
\ a w_{1}' = -\eta w_{2}', \ \ \ a w_{2}' = \eta w_{2}'$$ and
$\om(v_{1}',v_{2}') \neq 0$ and $\om(w_{1}',w_{2}') \neq 0$.

If $\om(v_{1}',v_{2}')<0$, then replace $v_{2}'$ with $-v_{2}'$ and
relabel so that $\om(v_{1}',v_{2}')>0$.  Doing this will, however,
yield that $av_{1}' = \mu v_{2}$ and $a v_{2}' = -\mu v_{2}'$. By
Lemma \ref{N}, this is not symplectically similar to the case where
$\om(v_{1}',v_{2}')>0$ in the first place.  The same can be said for
$w_{1}'$ and $w_{2}'$.

After making the necessary changes, we pick a new basis, $\be$, for $V$ as follows $$\be = \left\{
\frac{1}{\sqrt{\om(v_{1}',v_{2}')}} v_{1}', \frac{1}{\sqrt{\om(w_{1}',w_{2}')}} w_{1}',
\frac{1}{\sqrt{\om(v_{1}',v_{2}')}} v_{2}', \frac{1}{\sqrt{\om(w_{1}',w_{2}')}} w_{2}' \right\}$$ and relabel the
vectors $\{v_{1}'',w_{1}'',v_{2}'',w_{2}''\}$.  Then the matrix
$$W= \begin{pmatrix} v_{1}'' & w_{1}'' & v_{2}'' &w_{2}''
\end{pmatrix}$$ is in $Sp(4,\R)$ as $W^{t}JW = J$.  In this base $a$
is of form (9) in Theorem \ref{E} with $\eta \neq \pm \mu$, $\eta,\mu \neq 0$.

\subsubsection{Case 9: $V_{c}=K_{\mu i} \oplus K_{-\mu i}$}

Also assume that $\mu \in \R$ is nonzero.  Note that $\dim K_{\mu i} = \dim K_{-\mu i}$.

If $a$ is diagonalizable, then let $v_{1} \in K_{\mu i}$, then
$v_{2}=\overline{v_{1}} \in K_{-\mu i}$.  If $\om(v_{1},v_{2})= 0$,
then we let $w_{1} \in K_{-\mu i}$ and pick $w_{2}=\overline{w_{1}}
\in K_{\mu i}$ such that $w_{1} \neq v_{2}$, and we're guaranteed
that $\om(v_{1},w_{1}) \neq 0$ and $\om(v_{2},w_{2}) \neq 0$.  We
also know that $\om(v_{1},w_{2}) =\om(w_{1},v_{2}) = 0$ because
$v_{1},w_{2} \in K_{\mu i}$ and $v_{2},w_{1} \in K_{-\mu i}$.
Moreover,
\begin{align}
    \overline{\om(v_{1},w_{1})} &= \om(\overline{v_{1}},\overline{w_{1}})
    = \om(v_{2},w_{2}), \notag \\
    \overline{\om(w_{1},w_{2})} &= \om(\overline{w_{1}},\overline{w_{2}}) = \om(w_{2},w_{1})= -
    \om(w_{1},w_{2}). \notag
\end{align}
Let $$w_{1}' = w_{1}-\frac{\om(w_{1},w_{2})}{2\om(v_{2},w_{2})} v_{2}, \ \ \ w_{2}' = w_{2} +
\frac{\om(w_{1},w_{2})}{2\om(v_{2},w_{2})}v_{1}.$$  Then we have that
\begin{align}
    \om(w_{1}',w_{2}') &= \om\pnth{w_{1}-\frac{\om(w_{1},w_{2})}{2\om(v_{2},w_{2})} v_{2},w_{2} +
    \frac{\om(w_{1},w_{2})}{2\om(v_{2},w_{2})} v_{1}} \notag \\
    &= \om(w_{1},w_{2}) + \frac{\om(w_{1},w_{2})}{2 \om(v_{1},w_{2})} \om(w_{1},v_{1}) -
    \frac{\om(w_{1},w_{2})}{2\om(v_{2},w_{2})} \om(v_{2},w_{2}) = 0 \notag
\end{align}
Furthermore, $$\overline{w_{1}'} = \overline{w_{1}-\frac{\om(w_{1},w_{2})}{2\om(v_{2},w_{2})}
v_{2}}=\overline{w_{1}} - \frac{\overline{\om(w_{1},w_{2})}}{2\overline{\om(v_{2},w_{2})}} \overline{v_{2}}= w_{2}
+ \frac{\om(w_{1},w_{2})}{2\om(v_{1},w_{1})} v_{1}= w_{2}'$$ Relabel $w_{1}'$ and $w_{2}'$ as $w_{1}$ and $w_{2}$
and then proceed precisely as in Case 7.  This allows $\la=0$ in form (8) of Theorem \ref{E}.

Now if $\om(v_{1},v_{2})\neq 0$, then these two decompose out
according to Proposition \ref{L} and we pick $w_{1} \in K_{\mu i}$
and $w_{2}=\overline{w_{1}} \in K_{-\mu i}$ such that $w_{1} \neq
v_{1}$.  Then $\om(w_{1},w_{2}) \neq 0$ and we proceed precisely as
in Case 8.  This allows $\eta = \mu$ in form (9) of Theorem \ref{E}.
However, we must still require that
$\eta \neq -\mu$ because the following three matrices in $\si(4,\R)$, $$a_{1} = \begin{pmatrix} 0&0&-\mu&0 \\ 0&0&0&\mu \\
\mu&0&0&0 \\ 0&-\mu&0&0 \end{pmatrix}, \ \ a_{2} = \begin{pmatrix} 0&0&\mu&0 \\ 0&0&0&-\mu \\ -\mu&0&0&0
\\ 0&\mu&0&0 \end{pmatrix}, \ \ a_{3} = \begin{pmatrix} 0&\mu&0&0 \\ -\mu&0&0&0 \\ 0&0&0&\mu \\
0&0&-\mu&0 \end{pmatrix},$$ are symplectically similar conjugating by the matrices in $Sp(4,\R)$ $$A_{1}=
\begin{pmatrix} 0&1&0&0 \\ 1&0&0&0 \\ 0&0&0&1 \\ 0&0&1&0 \end{pmatrix}, \ \
A_{2}=\begin{pmatrix} 1&1&\frac{1}{2}&0 \\ 1&1&0&\frac{1}{2} \\ -1&1&0& \frac{1}{2} \\
1&-1&\frac{1}{2}&0
\end{pmatrix}.$$  That is, $(A_{1})^{-1}a_{1} A_{1} = a_{2}$ and $(A_{2})^{-1} a_{2} A_{2} =
a_{3}$.

If $a$ is not diagonalizable, then either $K_{\mu i}$ contains a cycle of generalized eigenvectors of length two
or $K_{-\mu i}$ does.  However, Lemma \ref{Q} shows, in fact, that they both do.

Let $\{v_{1},v_{2}\}$ be a cycle of generalized eigenvectors in $K_{\mu i}$, then $\{w_{1}=\overline{v_{1}},
w_{2}=\overline{v_{2}}\}$ is a cycle of generalized eigenvectors in $K_{-\mu i}$.

By Proposition 1.10 and the non-degeneracy of $\om$, we know that
$\om(v_{1},w_{2}) \neq 0$, $\om(v_{2},w_{1}) \neq 0$, and
$\om(v_{1},w_{1})=0$. In fact, we have
\begin{align}
    \mu i \om(v_{2},w_{2}) + \om(v_{1},w_{2}) &= \om(\mu i v_{2}+v_{1},w_{2}) =\om(a
    v_{2},w_{2})=-\om(v_{2},a w_{2}) \notag \\
    &= -\om(v_{2},-\mu i w_{2} + w_{1}) = \mu i \om(v_{2},w_{2}) -
    \om(v_{2},w_{1}), \notag
\end{align}
which implies that $\om(v_{1},w_{2}) = - \om(v_{2},w_{1})$.  With this information, we pick a new basis for
$V_{c}$ as follows
$$v_{1}' = v_{1}+w_{1}, \ \ v_{2}'=v_{2}+w_{2}, \ \ w_{1}' =
-i(v_{1} - w_{1}), \ \ w_{2}' = -i(v_{2}-w_{2}).$$  As these four vectors are real, they actually form a basis for
$V$.  In addition, we have that $$a v_{1}' = -\mu w_{1}', \ \ a w_{1}' = \mu v_{1}', \ \ a v_{2}' = -\mu w_{2}' +
v_{1}', \ \ a w_{2}' = \mu v_{2}' + w_{1}'.$$  We also have the following
\begin{align}
     \om(v_{1}',w_{2}') &= \om(v_{1}+w_{1},-i(v_{2}-w_{2}))
     = -i \om(v_{1},v_{2}) + i \om(v_{1},w_{2}) - i \om(w_{1},v_{2}) + i \om(w_{1},w_{2})
     = 0, \notag \\
     \om(v_{2}',w_{1}') &= \om(v_{2}+w_{2}, -i(v_{1} - w_{1}))
     = -i \om(v_{2},v_{1}) + i \om(v_{2},w_{1}) - i \om(w_{2},v_{1}) +i \om(w_{2},w_{1})
     = 0, \notag \displaybreak[0]\\
     \om(v_{1}',w_{1}') &= \om(v_{1}+w_{1},-i(v_{1}-w_{1})) = 0, \notag \displaybreak[0]\\
     \om(v_{1}',v_{2}') &= \om(v_{1}+w_{1},v_{2}+w_{2})
     = \om(v_{1},v_{2}) + \om(v_{1},w_{2}) + \om(w_{1},v_{2}) + \om(w_{1},w_{2})
     = 2 \om(v_{1},w_{2}), \notag \\
     \om(w_{1}',w_{2}') &= \om(-i(v_{1}-w_{1}),-i(v_{2}-w_{2}))
     = -\om(v_{1}-w_{1},v_{2} - w_{2}) \notag \\
     &= -\om(v_{1},v_{2}) + \om(v_{1},w_{2}) + \om(w_{1},v_{2}) - \om(w_{1},w_{2})
     = \om(v_{1}',v_{2}'). \notag
\end{align}
However, it is possible that $\om(v_{2}',w_{2}')$ is nonzero.  To
fix this potential problem, we pick another new basis for $V$ in the
following way $$v_{1}'' = v_{1}', \ \ v_{2}'' = v_{2}' +
\frac{\om(v_{2}',w_{2}')}{2\om(v_{2}',v_{1}')} w_{1}', \ \ w_{1}'' =
w_{1}', \ \ w_{2}'' = w_{2}' -
\frac{\om(v_{2}',w_{2}')}{2\om(v_{2}',v_{1}')} v_{1}'.$$  The set
$\{v_{1}'',v_{2}'',w_{1}'',w_{2}''\}$ still forms a basis for $V$,
but now
\begin{align}
    \om(v_{2}'',w_{2}'') &= \om\pnth{v_{2}' + \frac{\om(v_{2}',w_{2}')}{2\om(v_{2}',v_{1}')}
    w_{1}', w_{2}' - \frac{\om(v_{2}',w_{2}')}{2\om(v_{2}',v_{1}')}
    v_{1}'} \notag \\
    &= \om(v_{2}',w_{2}') - \frac{\om(v_{2}',w_{2}')}{2\om(v_{2}',v_{1}')}\om(v_{2}',v_{1}') +
    \frac{\om(v_{2}',w_{2}')}{2\om(v_{2}',v_{1}')} \om(w_{1}',w_{2}') \notag \\
    &= \om(v_{2}',w_{2}') - \frac{\om(v_{2}',w_{2}')}{2\om(v_{2}',v_{1}')}\om(v_{2}',v_{1}') -
    \frac{\om(v_{2}',w_{2}')}{2\om(w_{1}',w_{2}')} \om(w_{1}',w_{2}') \notag \\
    &= \om(v_{2}',w_{2}') -\frac{\om(v_{2}',w_{2}')}{2} - \frac{\om(v_{2}',w_{2}')}{2} \notag \\
    &= 0. \notag
\end{align}
In addition, $$a v_{2}'' = a \pnth{v_{2}' + \frac{\om(v_{2}',w_{2}')}{2\om(v_{2}',v_{1}')} w_{1}'} = a v_{2}' +
\frac{\om(v_{2}',w_{2}')}{2 \om(v_{2}',v_{1}'} a w_{1}' = -\mu w_{2}' + v_{1}' + \mu \frac{\om(v_{2}',w_{2}')}{2
\om(v_{2}',v_{1}')} v_{1}' = -\mu w_{2}'' + v_{1}''$$ and $$a w_{2}'' = a \pnth{w_{2}' -
\frac{\om(v_{2}',w_{2}')}{2\om(v_{2}',v_{1}')} v_{1}'} = a w_{2}' - \frac{\om(v_{2}',w_{2}')}{2
\om(v_{2}',v_{1}')} a v_{1}' = \mu v_{2}' + w_{1}' + \mu \frac{\om(v_{2}',w_{2}')}{2 \om(v_{2}',v_{1}')} w_{1}' =
\mu v_{2}'' + w_{1}''.$$  Finally, $\om(v_{1}'',w_{2}'')=\om(v_{2}'',w_{1}'')=0$ and $\om(v_{1}'',v_{2}'') =
\om(v_{1}',v_{2}') = \om(w_{1}',w_{2}') = \om(w_{1}'',w_{2}'')$.

If $\om(v_{1}'',v_{2}'')=\om(w_{1}'',w_{2}'') < 0$, then replace the
vectors $v_{2}''$ and $w_{2}''$ with $-v_{2}''$, and $-w_{2}''$
respectively and relabel.  This will make the values in question
positive, but will result in $$a v_{2}'' = -\mu w_{2}'' - v_{1}'', \
\ \ \ \ a w_{2}'' = \mu v_{2}'' - w_{1}''.$$  Lemma \ref{N} shows,
however, that this case is not symplectically similar to the case
where $\om(v_{1}'',v_{2}'') \linebreak[0] =\om(w_{1}'',w_{2}'') >0$
in the first place.  But we still make the change to guarantee that
$\om(v_{1}'',v_{2}'')=\om(w_{1}'',w_{2}'')>0$.

After making the appropriate changes, we can pick a final basis for $V$, $\be$, as follows $$\be=\left\{
\frac{1}{\sqrt{\om(v_{1}'',v_{2}'')}} v_{1}'', \frac{1}{\sqrt{\om(v_{1}'',v_{2}'')}} w_{1}'',
\frac{1}{\sqrt{\om(v_{1}'',v_{2}'')}} v_{2}'', \frac{1}{\sqrt{\om(v_{1}'',v_{2}'')}} w_{2}'' \right\}$$ and
relabel the vectors $v_{1}'''$, $w_{1}'''$, $v_{2}'''$, and $w_{2}'''$. Then the matrix
$$W=\begin{pmatrix} v_{1}''' & w_{1}''' & v_{2}''' & w_{2}''' \end{pmatrix}$$ is in $Sp(4,\R)$ because
$W^{t}JW=J$.  In this basis, $a$ is of form (10) in Theorem \ref{E}.

This covers all of the possible cases of $a \in \si(4,\R)$ and
completes the proof of Theorem \ref{E}.

\subsection{Canonical Forms of Matrices in a Nonstandard Representation of $\so(3,1,\R)$}

Let $V=\R^{4}$.  Let $J_{1},J_{2} \in GL(V)$ be given as $$J_{1}= \begin{pmatrix} 0&0&1&0 \\ 0&0&0&1 \\ -1&0&0&0 \\
0&-1&0&0 \end{pmatrix}, \ \ \ \ J_{2}=\begin{pmatrix} 0&0&0&1 \\
0&0&-1&0 \\ 0&1&0&0 \\ -1&0&0&0 \end{pmatrix}.$$ Note that $J_{1}$
and $J_{2}$ are skew-symmetric and have the property that $J_{i}^{2}
= -I_{4}$, where $I_{4}$ denotes the $4 \times 4$ identity. Let
$\h(J_{2})$ be defined by
\begin{equation}\h(J_{2})=\{a \in \Hom(V,V) \ | \ a^{t}J_{i} + J_{i}a = 0 \ \text{for all} \ i \in \{1,2\} \}.
\notag \end{equation} Clearly $\h(J_{2})$ is a subalgebra of
$\si(4,\R)$ as $J_{1}$ is the $J$ used in the section on
$\si(4,\R)$. This is why we use the notation $\h(J_{2})$ to denote
it; that is, it is the subalgebra, $\h$, of $\si(4,\R)$ whose
elements are also skew-symmetric about the matrix $J_{2}$. As such,
many of our proofs will refer to those in that section.

However, before we get started, we should note a few interesting
properties of $\h(J_{2})$.  Let $D_{1}$ be an arbitrary $4 \times 4$
matrix.  Then we can use Maple to solve the equation
$D_{1}^{t}J_{i}+J_{i}D_{1}=0$ for all $i \in \{1,2\}$.  Doing so
will force $D_{1}$ to be of the form $$D_{1}=\begin{pmatrix} A&B \\
C&-A^{t}
\end{pmatrix},$$ with $A,B,C$ all $2 \times2$ matrices such that
$A=\bigl{(}\begin{smallmatrix} p&-q \\ q&p
\end{smallmatrix}\bigr{)}$ and $B$ and $C$ are trace-free symmetric.
This implies that a basis for $\h(J_{2})$ is
\begin{align}
    \notag \begin{pmatrix} 1&0&0&0
    \\ 0&1&0&0 \\ 0&0&-1&0 \\ 0&0&0&-1 \end{pmatrix},
    & & \begin{pmatrix} 0&-1&0&0
    \\ 1&0&0&0 \\ 0&0&0&-1 \\ 0&0&1&0 \end{pmatrix}, &
    & \begin{pmatrix} 0&0&0&0
    \\ 0&0&0&0 \\ 1&0&0&0 \\ 0&-1&0&0  \end{pmatrix}, \\
    \notag \begin{pmatrix} 0&0&0&0 \\
    0&0&0&0 \\ 0&1&0&0 \\ 1&0&0&0 \end{pmatrix}, & &
    \begin{pmatrix} 0&0&1&0
    \\ 0&0&0&-1 \\ 0&0&0&0 \\ 0&0&0&0 \end{pmatrix}, & &
    \begin{pmatrix} 0&0&0&1 \\
    0&0&1&0 \\ 0&0&0&0 \\ 0&0&0&0 \end{pmatrix},
\end{align}
Label these elements $\X_{1}, \ldots, \X_{6}$. As the Lie bracket in
a matrix Lie algebra is given by the commutator, then using this
basis, the Lie algebra $\h(J_{2})$ has the multiplication table
$$\begin{array}{cccccccc}
    & \vline & \X_{1}&\X_{2}&\X_{3}&\X_{4}&\X_{5}&\X_{6} \\
    \hline
    \X_{1} & \vline & 0&0&-2\X_{3}&-2\X_{4}&2\X_{5}&2\X_{6} \\
    \X_{2} & \vline & 0&0&2\X_{4}&-2\X_{3}&2\X_{6}&-2\X_{5} \\
    \X_{3} & \vline & 2\X_{3}&-2\X_{4}&0&0&-\X_{1}&-\X_{2} \\
    \X_{4} & \vline & 2\X_{4}&2\X_{3}&0&0&\X_{2}&-\X_{1} \\
    \X_{5} & \vline & -2\X_{5}&-2\X_{6}&\X_{1}&-\X_{2}&0&0 \\
    \X_{6} & \vline & -2\X_{6}&2\X_{5}&\X_{2}&\X_{1}&0&0
\end{array}$$
We make the change of basis
\begin{align}
    \notag \Y_{1}&=\frac{1}{2}\X_{1}, & \Y_{2}&=
    \frac{\sqrt{2}}{4}(\X_{3}+\X_{4}+\X_{5}+\X_{6}), & \Y_{3}&=
    \frac{\sqrt{2}}{4}(\X_{3}-\X_{4}+\X_{5}-\X_{6}), \\
    \notag \Y_{4}&=\frac{\sqrt{2}}{4}(\X_{3}+\X_{4}-\X_{5}-\X_{6}),
    & \Y_{5}&=\frac{\sqrt{2}}{4}(\X_{3}-\X_{4}-\X_{5}+\X_{6}), &
    \Y_{6}&=-\frac{1}{2}\X_{2}.
\end{align}
In this basis, the multiplication table becomes
$$\begin{array}{cccccccc}
    & \vline & \Y_{1}&\Y_{2}&\Y_{3}&\Y_{4}&\Y_{5}&\Y_{6} \\
    \hline
    \Y_{1} & \vline & 0&-\Y_{4}&-\Y_{5}&-\Y_{2}&-\Y_{3}&0 \\
    \Y_{2} & \vline & \Y_{4}&0&-\Y_{6}&\Y_{1}&0&-\Y_{3} \\
    \Y_{3} & \vline & \Y_{5}&\Y_{6}&0&0&\Y_{1}&\Y_{2} \\
    \Y_{4} & \vline & \Y_{2}&-\Y_{1}&0&0&\Y_{6}&-\Y_{5} \\
    \Y_{5} & \vline & \Y_{3}&0&-\Y_{1}&-\Y_{6}&0&\Y_{4} \\
    \Y_{6} & \vline & 0&\Y_{3}&-\Y_{2}&\Y_{5}&-\Y_{4}&0
\end{array}$$
We will show that this is the multiplication table for $\so(3,1,\R)$
and hence $\h(J_{2})$ is a representation of $\so(3,1,\R)$. We'll
call this representation $\rho_{1}$.  We'll also show that
$\rho_{1}$ is not equivalent to the standard representation of
$\so(3,1,\R)$.

Let $\rho_{2}$ denote the standard representation of $\so(3,1,\R)$,
which consists of those matrices that are skew-symmetric about the
matrix $$M=\begin{pmatrix} -1&0&0&0 \\
0&1&0&0 \\ 0&0&1&0 \\ 0&0&0&1 \end{pmatrix},$$ that is, it consists
of those matrices, $a$, such that $a^{t}M + Ma = 0$. Let $D_{2}$ be
an arbitrary $4 \times 4$ matrix. Then we use Maple to solve the
equation $D_{2}^{t}M+MD_{2}=0$. We find then that to satisfy the
equation, $D_{2}$ must be of
the form $$D_{2}= \begin{pmatrix} 0&d_{1}&d_{2}&d_{3} \\
d_{1}&0&-d_{4}&-d_{5} \\ d_{2}&d_{4}&0&-d_{6}
\\ d_{3}&d_{5}&d_{6}&0 \end{pmatrix}.$$  This implies that a basis
for the standard representation of $\so(3,1,\R)$ is given by the
following six matrices
\begin{align}
    \notag \begin{pmatrix} 0&1&0&0 \\ 1&0&0&0 \\ 0&0&0&0 \\ 0&0&0&0
    \end{pmatrix}, & & \begin{pmatrix} 0&0&1&0 \\ 0&0&0&0 \\ 1&0&0&0
    \\ 0&0&0&0 \end{pmatrix}, & & \begin{pmatrix} 0&0&0&1 \\ 0&0&0&0 \\ 0&0&0&0 \\ 1&0&0&0
    \end{pmatrix},\\
    \notag \begin{pmatrix} 0&0&0&0 \\ 0&0&-1&0 \\
    0&1&0&0 \\ 0&0&0&0 \end{pmatrix},  & & \begin{pmatrix} 0&0&0&0
    \\ 0&0&0&-1 \\ 0&0&0&0 \\ 0&1&0&0 \end{pmatrix}, & &
    \begin{pmatrix} 0&0&0&0 \\ 0&0&0&0 \\ 0&0&0&-1 \\ 0&0&1&0
    \end{pmatrix}.
\end{align}
Call these elements $\ga_{1}, \ldots, \ga_{6}$ respectively.  If we
compute the multiplication table of this Lie algebra in this basis,
we'll obtain the second multiplication table given in this section.
Thus, as abstract Lie algebras these two representations are
isomorphic. This shows that $\rho_{1}$ is a representation of
$\so(3,1,\R)$. But $\rho_{1}$ and $\rho_{2}$ are not equivalent as
there exists no $T \in GL(\R^{4})$ such that $T \circ \rho_{1} =
\rho_{2} \circ T$. We show this by proving the equivalent statement
that there exists no invertible $4 \times 4$ matrix $T$ such that
$T\Y_{i} = \ga_{i}T$ for $1 \leq i \leq 6$. Let $T$ be an arbitrary
$4 \times 4$ matrix. Then we use Maple to solve the system $T\Y_{i}
= \ga_{i}T$ for $1 \leq i \leq 6$, which yields that $T=0$. Hence
$\rho_{1}$ and $\rho_{2}$ are not equivalent.

As the two representations are not equivalent, we cannot rely on the
classification of the canonical forms of matrices in $\rho_{2}$ to
discover those canonical forms of the matrices in the representation
of $\rho_{1}$, that is, the canonical forms of matrices in
$\h(J_{2})$. Instead we find these canonical forms directly. We
begin by proving a few facts about $\h(J_{2})$.

\begin{lem}\label{AA}
    Let $a \in \h(J_{2})$.  Then $a^{t} \in \h(J_{2})$.
\end{lem}

{\em Proof.} If $a \in \h(J_{2})$, then
$a^{t}J_{1}+J_{1}a=a^{t}J_{2}+J_{2}a=0$.  For either case, multiply
on the left by $J_{i}$ to obtain $J_{i}a^{t}J_{i}-a=0$. Then
multiply on the right by $J_{i}$: $-J_{i}a^{t}-aJ_{i}=0$ or
equivalently $(a^{t})^{t}J_{1} + J_{1} a^{t} =
(a^{t})^{t}J_{2}+J_{2}a^{t}= 0$. $\blacksquare$ \vspace{.1 in}

Define $\om_{i}:V \times V \To \R$ by $\omega_{i}(x,y)=x^{t}J_{i}y$
for all $x,y \in V$ and $i \in \{1,2\}$. Then the $\om_{i}$ are a
pair of symplectic forms. They can also be viewed as maps $V_{c}
\times V_{c} \To \C$ with the same rules of assignment.

\begin{prop}\label{BB}
    These forms, $\om_{i}:V \times V \To \R$ ($\om_{i}:V_{c} \times V_{c} \To \C$) for $i \in \{1,2\}$,
    are non-degenerate skew-symmetric bilinear forms on $V$ ($V_{c}$).
\end{prop}

{\em Proof.} From their construction, they are clearly bilinear over $\R$ or $\C$. To prove skew-symmetry, we note
that as $\omega_{i}(x,y)$ can be viewed as a $1\times 1$ matrix, we have that $(\om_{i}(x,y))^{t}=\om_{i}(x,y)$.
Then, as $J_{i}^{t}=-J_{i}$, this yields that for $i \in \{1,2\}$
$$\om_{i}(x,y)=\pnth{\om_{i}(x,y)}^{t} = (x^{t}J_{i}y)^{t} = y^{t}J_{i}^{t}x = -y^{t}J_{i}x = -\om_{i}(y,x).$$
Finally, to prove non-degeneracy, let $z \in V$ be such that
$\om_{i}(z,y)=z^{t}J_{i}y=0$ for all $y \in V$ for some $i \in
\{1,2\}$. As both $J_{i}$ are invertible, each has zero kernel, thus
$z=0$. $\blacksquare$ \vspace{.1 in}

We will always denote these symplectic forms by simply $\om_{1}$ or
$\om_{2}$ as it will be clear from the context whether we mean the
complex or real form.

\begin{prop}\label{RR}
    The forms $\om_{1}$ and $\om_{2}$ have the following properties for all $x,y \in V$.
    \begin{enumerate}
        \item $\om_{1}(J_{2}x,y)=\om_{1}(x,J_{2}y)=-\om_{2}(J_{1}x,y)=-\om_{2}(x,J_{1}y)$.
        \item $\om_{i}(J_{i}x,y)=-\om_{i}(x,J_{i}y)$ for all $i \in \{1,2\}$.
        \item $\om_{i}(J_{1}J_{2}x,y)=\om_{i}(x,J_{1}J_{2}y) = -\om_{i}(J_{2}J_{1}x,y)
        = -\om_{i}(x,J_{2}J_{1}y)$ for all $i \in \{1,2\}$.
    \end{enumerate}
\end{prop}

{\em Proof.}  First note that $J_{1}J_{2}=-J_{2}J_{1}$.  Let $x,y \in V$.  Then we have that
\begin{align}
    \om_{1}(x,J_{2}y)&=x^{t}J_{1}J_{2}y = -x^{t}J_{2}J_{1}y= (J_{2}x)^{t}J_{1}y = \om_{1}(J_{2}x,y), \notag \\
    \om_{2}(J_{1}x,y)&=(J_{1}x)^{t}J_{2}y = -x^{t}J_{1}J_{2}y = x^{t}J_{2}J_{1}y = \om_{2}(x,J_{1}y), \notag \\
    \om_{1}(x,J_{2}y)&=x^{t}J_{1}J_{2}y = -(J_{1}x)^{t}J_{2}y = -\om_{2}(J_{1}x,y). \notag
\intertext{Also we have, for all $i \in \{1,2\}$, that}
    \om_{i}(J_{i}x,y)&=(J_{i}x)^{t}J_{i}y = -x^{t}J_{i}J_{i}y = -\om_{i}(x,J_{i}y). \notag
\end{align}
The third statement is a direct consequence of the previous two and
the fact that $J_{1}J_{2}=-J_{2}J_{1}$. $\blacksquare$ \vspace{.05
in}

Now we define the group preserving $J_{1}J_{2}$.  This group, which
we'll call $H(J_{2})$ is defined as
$$H(J_{2})= \{ A \in GL(V) \ | \ \om_{i}(Ax,Ay) = \om_{i}(x,y) \ \text{for all} \ i \in \{1,2\} \}.$$  That is,
it is the group that preserves the symplectic forms.

\begin{lem}\label{CC}
    $a \in \h(J_{2})$ if and only if $\om_{i}(ax,y)=-\om_{i}(x,ay)$ for all $x,y \in V$ and $i \in \{1,2\}$.
    $A \in H(J_{2})$ if an only if $A^{t}J_{i}A = J_{i}$ for all $i \in \{1,2\}$.
\end{lem}

{\em Proof.} As $a \in \h(J_{2})$, we have that $a^{t}J_{i} +
J_{i}a=0$ for all $i \in \{1,2\}$.  Let $x,y \in V$. This gives us
that for all $i \in \{1,2\}$ $$\om_{i}(ax,y)=(ax)^{t}J_{i}y = x^{t}
a^{t}J_{i}y = -x^{t}J_{i} ay = -\om_{i}(x,ay).$$

Next assume that $\om_{i}(ax,y)=-\om_{i}(x,ay)$ for all $x,y \in V$ and $i \in \{1,2\}$.  Then
\begin{align}
    \om_{i}(ax,y) &= -\om_{i}(x,ay) \notag\\
    (ax)^{t} J_{i} y &= -x^{t} J_{i} a y \notag\\
    x^{t} a^{t} J_{i} y &= -x^{t} J_{i} a y \notag
\end{align}
for all $i \in \{1, 2\}$.  As this is true for all $x,y \in V$, this
implies that $a^{t}J_{i} = -J_{i}a$ or equivalently that $a^{t}J_{i}
+J_{i}a = 0$. Hence $a \in \h(J_{2})$.

To prove the statement about $H(J_{2})$, first assume that $A \in
H(J_{2})$.  Then for all $i \in \{1,2\}$
\begin{align}
    \om_{i}(Ax,Ay) &= \om_{i}(x,y) \notag \\
    (Ax)^{t} J_{i} Ay &= x^{t}J_{i}y \notag\\
    x^{t}A^{t}J_{i}Ay &= x^{t}J_{i}y \notag
\end{align}
As this is true for all $x,y \in V$, we see that $A^{t}J_{i}A = J_{i}$.

If we assume first that $A^{t}J_{i}A = J_{i}$, then $A$ clearly
preserves the symplectic forms and hence $A \in H(J_{2})$.
$\blacksquare$ \vspace{.1 in}

We now wish to conjugate $a \in \si(2n,\R)$ by an arbitrary element $A \in Sp(2n,\R)$, that is, $A^{-1}aA$.

\begin{lem} \label{DD}
    If $a \in \h(J_{2})$ and $A \in H(J_{2})$, then $A^{-1}aA \in \h(J_{2})$.
\end{lem}

{\em Proof.} If $A \in H(J_{2})$, then $\om_{i}(Ax,y) =
\om_{i}(A^{-1}Ax,A^{-1}y) = \om_{i}(x,A^{-1}y)$.  This yields
$$\om_{i}(A^{-1}aA x, y) = \om_{i}(aA x, Ay) = -\om_{i}(Ax,aAy) =
-\om_{i}(x,A^{-1}aAy).$$  Then by Lemma \ref{CC}, $A^{-1}aA \in
\h(J_{2})$. $\blacksquare$ \vspace{.1 in}

If $a_{1},a_{2} \in \h(J_{2})$ are such that $A^{-1}a_{1}A=a_{2}$
for some $A \in H(J_{2})$, we say that $a_{1}$ and $a_{2}$ are {\em
$\h$-symplectically similar}.

This naturally brings up the question: what kind of canonical forms
could $a \in \h(J_{2})$ have if this were the only kind of change of
basis allowed? The result is as follows

\begin{thm} \label{EE}
    Let $a \in \h(J_{2})$, then $a$ is $\h$-symplectically similar to one
    of the following three matrices.  We call this the {\em real $\h$-symplectic canonical form} of the matrix.
    \begin{align}
        (1) \ &\begin{pmatrix} \la &0&0&0 \\ 0&\la&0&0 \\ 0&0&-\la&0 \\
        0&0&0&-\la \end{pmatrix}, \  \la \in \R, & (2) \ & \begin{pmatrix} 0&0&1&0 \\ 0&0&0&-1  \\ 0&0&0&0 \\
        0&0&0&0 \end{pmatrix}, \notag \\
        (3) \ & \begin{pmatrix} \la&\vep\mu&0&0 \\ -\vep\mu&\la&0&0 \\
        0&0&-\la&\vep\mu \\ 0&0&-\vep\mu&-\la \end{pmatrix}, \ \begin{array}{l} \la \geq0, \\ \mu > 0, \\
        \vep^{2}=1. \end{array} \notag
    \end{align}
\end{thm}

The remainder of this section is the proof of this theorem.  We
present first some preliminary facts that will allow us to do so.

\begin{lem}\label{FF}
    Let $a \in \h(J_{2})$ and let $\la$ be an eigenvalue of $a$.  Then $-\la$ is also an eigenvalue of $a$.
\end{lem}

{\em Proof.} As $\h(J_{2}) \subseteq \si(4,\R)$, this result follows
by Lemma \ref{F}. $\blacksquare$ \vspace{.1 in}

\begin{lem}\label{QQ}
    If $v \in V$ is an eigenvector of $a \in \h(J_{2})$ corresponding to the eigenvalue $\la$,
    then $J_{1}J_{2}v$ is as well. Furthermore, if $v$ and $J_{1}J_{2}v$ are linearly dependent,
    then $v = \pm i J_{1}J_{2}v$.  This
    implies that if $v$ is real, then $v$ and $J_{1}J_{2}v$ are linearly independent.  Finally, if $v= \pm
    i J_{1}J_{2}v$ and we write $v= u + i w$, then $w = \pm J_{1}J_{2} u$.
\end{lem}

{\em Proof.} If $v \in V$ is an eigenvector of $a \in \h(J_{2})$
corresponding to the eigenvalue $\la \in \C$, then $av = \la v$.
Also by Lemma \ref{AA}, $a^{t} \in \h(J_{2})$. This yields that
$$a(J_{1}J_{2}v) = -J_{1}a^{t}J_{2}v = J_{1}J_{2}av = \la(J_{1}J_{2}
v).$$

Before we go on, note that $(J_{1}J_{2})^{2}=-J_{2}J_{1}J_{1}J_{2} =
-I_{4}$.  Assume that $v$ and $J_{1}J_{2}v$ are linearly dependent.
Then consider the equation
\begin{align}
    c_{1}v + c_{2}J_{1}J_{2}v &= 0 \notag \\
\intertext{Multiply on the left by $J_{1}J_{2}$.}
    c_{1}J_{1}J_{2}v - c_{2} v &=0 \notag \\
\intertext{Then we see that $c_{1}J_{1}J_{2}v = c_{2}v$.  Now multiply the original equation by $c_{1}$ and make
this substitution.}
    (c_{1})^{2} v + c_{2} c_{1} J_{1}J_{2}v &=0 \notag \\
    (c_{1})^{2}v + (c_{2})^{2} v &=0 \notag \\
    \pnth{(c_{1})^{2} + (c_{2})^{2}} v &=0 \notag
\end{align}
As $v \neq 0$, this implies that $(c_{1})^{2} + (c_{2})^{2} =0$ or equivalently that $c_{2}=\pm c_{1} i$.  Thus if
$v$ and $J_{1}J_{2}v$ are linearly dependent, then $v= \pm i J_{1}J_{2} v$ (simply let $c_{1}=1$).  If $v$ is
real, then $c_{1},c_{2} \in \R$ and we see that $c_{1}=c_{2}=0$.

If $v$ is complex, then we can write $v=u+iw$ where $u,w \in V$.  If $c_{2}=\pm c_{1} i$, then let $c_{1}=1$ and
we obtain that $\pm iJ_{1}J_{2}v=v$. This yields the following
\begin{align}
    \pm iJ_{1}J_{2}v &= v \notag \\
    \pm i(J_{1}J_{2}u + i J_{1}J_{2}w) &= u+ i w \notag \\
    \mp J_{1}J_{2}w \pm iJ_{1}J_{2}u &= u+iw \notag
\end{align}
This implies that $w= \pm J_{1}J_{2}u$.  Then $w=\pm J_{1}J_{2}u$ and we can write $v=u \pm iJ_{1}J_{2}u$.
$\blacksquare$ \vspace{.1 in}

\begin{lem}\label{TT}
    If $\{v_{1},\ldots, v_{k}\}$ is a cycle of generalized eigenvectors corresponding to $\la \in \C$, then so
    is $\{J_{1}J_{2}v_{1},\ldots, J_{1}J_{2}v_{k}\}$.
\end{lem}

{\em Proof.} First, note that if $\{v_{1},\ldots,v_{k}\}$ is a cycle
of generalized eigenvectors corresponding to $\la \in \C$, then, by
definition, $v_{i}=(a-\la I_{n})^{k-i} v_{k}$ for all $1 \leq i \leq
k$.

We now proceed by induction if $k=1$, then by Lemma \ref{QQ}, if
$v_{1}$ is an eigenvector then so is $J_{1}J_{2}v_{1}$ and the
statement holds.

Assume that the statement is true for $k=m$.

Finally let $k=m+1$, that is $\{v_{1},\ldots,v_{m+1}\}$ is a cycle
of generalized eigenvectors.  Then $\{v_{1},\ldots,v_{m}\}$ is a
cycle of generalized eigenvectors of length $m$.  Thus, by
assumption, $\{J_{1}J_{2}v_{1}, \ldots, \linebreak[0]
J_{1}J_{2}v_{m}\}$ is a cycle of generalized eigenvectors of length
$m$ as well.  This implies that $J_{1}J_{2}v_{i} = (a-\la
I_{n})^{m-i} J_{1}J_{2} v_{m}$.

As $\{v_{1},\ldots,v_{m+1}\}$ is a cycle, then we have that
$v_{m}=(a-\la I)v_{m+1}$. This implies
\begin{align}
    \notag (a-\la I_{n}) J_{1}J_{2} v_{m+1} &= (aJ_{1}J_{2}-\la I_{n}
    J_{1}J_{2}) v_{m+1} = \pnth{-J_{1}a^{t}J_{2} - J_{1}J_{2}(\la I_{n})} v_{m+1} \\
    \notag &= \pnth{J_{1}J_{2}a - J_{1}J_{2}(\la I_{n})} v_{m+1} =
    J_{1}J_{2} \pnth{(a-\la I_{n}) v_{m+1}} = J_{1}J_{2} v_{m}. \\
\intertext{This yields that}
    \notag J_{1}J_{2}v_{i} &= (a- \la I_{n})^{m-i} J_{1}J_{2} v_{m}
    = (a - \la I_{n})^{m-i} (a - \la I_{n}) J_{1}J_{2} v_{m+1} \\
    &= (a-
    \la I_{n})^{m+1-i} J_{1}J_{2} v_{m+1}.
\end{align}
Hence $\{J_{1}J_{2}v_{1}, \ldots, J_{1}J_{2} v_{m+1}\}$ is a cycle
of generalized eigenvectors.  Therefore, by induction, the statement
holds for all $k \in \N$.  Note that if $v_{1}$ and
$J_{1}J_{2}v_{1}$ are linearly independent, then these cycles are
linearly independent of each other. $\blacksquare$ \vspace{.1 in}

We know that the characteristic polynomial of $a$ always splits over $\C$, thus we have that the complexification
of $V$, $V_{c}$, is the direct sum of the generalized eigenspaces of $a$.  This allows us to write $V_{c}$ in the
following manner
$$V_{c} = K_{0} \oplus K_{\la_{1}} \oplus K_{-\la_{1}} \oplus K_{\la_{2}} \oplus K_{-\la_{2}} \oplus \cdots \oplus
K_{\la_{k}} \oplus K_{-\la_{k}}$$ where $1 \leq k \leq n$, $K_{\mu}$ is the generalized eigenspace corresponding
to the eigenvalue $\mu$, $\la_{i} \neq 0$ for all $i$, and the $\la_{i}$ are distinct and such that $\la_{i} \neq
-\la_{j}$ for any pair $i,j$. Note that if all the eigenvalues of $a$ are real, then its characteristic polynomial
splits over $\R$ and $V$ decomposes in the manner above where all the summands are subspaces over $\R$.

\begin{prop}\label{GG}
    If $\mu \neq -\eta$, then $K_{\eta}$ and $K_{\mu}$ are $\h$-symplectically orthogonal to each other.
    That is, $\om_{i}(K_{\eta},K_{\mu})=0$.
\end{prop}

{\em Proof.} $K_{\eta}$ and $K_{\mu}$ are orthogonal with respect to
$\om_{1}$ by Lemma \ref{G}.  The proof that they are orthogonal with
respect to $\om_{2}$ is analogous to the proof of that lemma.
$\blacksquare$ \vspace{.1 in}

Note that the above proposition is true if $\eta=0$.  However, if
$\eta \neq 0$, then this leads immediately to a useful corollary.

\begin{cor}\label{HH}
    If $\eta \neq 0$, then $K_{\eta}$ is $\h$-symplectically orthogonal to itself.
\end{cor}

\begin{prop}\label{II}
    For $i \in \{1,2\}$, $\om_{i}$ is non-degenerate on $K_{\mu} \oplus K_{-\mu}$ for all $\mu \neq 0$.  In addition,
    $\om_{i}$ is non-degenerate on $K_{0}$.
\end{prop}

{\em Proof.} As $\h(J_{2}) \subseteq \si(4,\R)$, then by Proposition
\ref{I}, this result is clear for $\om_{1}$. The proof for $\om_{2}$
is analogous to the proof of that proposition. $\blacksquare$
\vspace{.1 in}

\begin{cor}\label{JJ}
    If $v \in K_{\mu}$ is nonzero, then there exists a nonzero $w_{1} \in K_{-\mu}$ such that $\om_{1}(v,w_{1}) \neq 0$
    and a nonzero $w_{2} \in K_{-\mu}$ such that $\om_{2}(v,w_{2}) \neq 0$.
\end{cor}

{\em Proof:} The existence of $w_{1}$ is guaranteed by Corollary
\ref{J}. The proof of the existence of $w_{2}$ is analogous to the
proof of that corollary. Note that it is possible that
$w_{2}=w_{1}$. $\blacksquare$ \vspace{.1 in}

We can say a little more when $v \in K_{\mu}$ is an eigenvector.

\begin{prop}\label{KK}
    If $v \in K_{\mu}$ is an eigenvector and $w \in K_{-\mu}$ such that $\om_{1}(v,w) \neq 0$ or $\om_{2}(v,w) \neq 0$,
    then $w$ is the end vector in any cycle of generalized eigenvectors to which it belongs.
\end{prop}

{\em Proof.} As $\h(J_{2}) \subseteq \si(4,\R)$, this is true for
$\om_{1}$ by Proposition \ref{K}.  The proof for $\om_{2}$ is
analogous to the proof of that proposition. $\blacksquare$
\vspace{.1 in}

\begin{prop}\label{PP}
    The dimension of $K_{0}$ is even and $\dim K_{\mu} = \dim
    K_{-\mu}$ for all nonzero $\mu \in \C$.  Moreover, if $\mu \in \R$, then $\dim K_{\mu}$ is even.
\end{prop}

{\em Proof.} As $\h(J_{2}) \subseteq \si(4,\R)$, the first two
properties follow immediately.

Let $\mu \in \R$.  If $\dim K_{\mu}=0$, then we're done.  Assume
then that $\dim K_{\mu} \neq 0$.  Let $\{v_{1},\ldots, v_{k}\}
\subseteq K_{\mu}$ be a cycle of generalized eigenvectors.  Then by
Lemma \ref{TT}, $\{J_{1}J_{2}v_{1},\ldots, J_{1}J_{2}v_{k}\}$ is
also a cycle of generalized eigenvectors.  Furthermore, As the
$v_{i}$ are real vectors, then by Lemma \ref{QQ}, $v_{1}$ and
$J_{1}J_{2}v_{1}$ are linearly independent and hence their entire
cycles are linearly independent of each other.  Therefore every
cycle of generalized eigenvectors is paired with another cycle of
generalized eigenvectors of the same length to which it is linearly
independent. As there exists a basis for $K_{\mu}$ consisting of
disjoint cycles of generalized eigenvectors of $a$, we have that
$\dim K_{\mu}$ must be even. $\blacksquare$ \vspace{.1 in}

\begin{lem}\label{SS}
    Let $\mu i$ be an eigenvalue of $a \in \h(J_{2})$.  Then $\dim K_{\mu i} =2$ and $a$ is diagonalizable.
\end{lem}

{\em Proof.} As $\dim K_{-\mu i} = \dim K_{\mu i}$ and $\dim V=4$, then clearly $\dim K_{\mu i} \leq 2$.

Now assume to the contrary the $\dim K_{\mu i} = 1$.  Let $v_{1} \in
K_{\mu i}$ be nonzero.  By Lemma \ref{QQ} and as $\dim K_{\mu i}=1$,
we have that $v_{1}$ and $J_{1}J_{2}v_{1}$ must be linearly
dependent and consequently $v_{1}=\pm i J_{1}J_{2} v_{1}$.  Assume
first that $v_{1} = i J_{1}J_{2}v_{1}$. Then by the same lemma, we
can write $v_{1}=u + i J_{1}J_{2} u$ for some real vector $u$. This
implies that $v_{2}=u - i J_{1}J_{2} u \in K_{-\mu i}$.  In
addition, as $\dim K_{-\mu i} = \dim K_{\mu i}$, we have that
$\{v_{1},v_{2}\}$ is a basis for $K_{\mu i} \oplus K_{-\mu i}$.
However, for $j \in \{1,2\}$, we have by Lemma \ref{RR} that
\begin{align}
    \om_{j}(v_{1},v_{2}) &= \om_{j}(u+iJ_{1}J_{2}u, u - i J_{1}J_{2} u) \notag \\
    &= \om_{j}(u,u) - i \om_{j}(u,J_{1}J_{2}u) + i \om_{j}(J_{1}J_{2}u,u) +
    \om_{j}(J_{1}J_{2}u,J_{1}J_{2}u) \notag \\
    &= -i \om_{j}(u,J_{1}J_{2}u) + i \om_{j}(u,J_{1}J_{2}u) = 0. \notag
\end{align}
This implies that $\om_{j}$ is degenerate on $K_{\mu i} \oplus K_{-\mu i}$ contradicting Proposition \ref{II}. By
a similar argument the same result can be shown if $v_{1}=-iJ_{1}J_{2}v_{1}$.  Thus $\dim K_{\mu i} \geq 2$ and
consequently $\dim K_{\mu i} = 2$.

Next we show that $a$ is diagonalizable.  As $\dim K_{\mu i} = \dim
K_{-\mu i} = 2$, we have that $V= K_{\mu i} \oplus K_{-\mu i}$. Thus
if $a$ is not diagonalizable, then there exists a basis for $K_{\mu
i}$ consisting of a cycle of generalized eigenvectors of length two.
Let $\{v_{1},v_{2}\}$ be such a cycle. Then, by Lemma \ref{TT},
$\{J_{1}J_{2}v_{1}, J_{1}J_{2}v_{2}\}$ is also a cycle of
generalized eigenvectors in $K_{\mu i}$. However, as $\dim K_{\mu i}
= 2$, this implies that $v_{1}$ and $J_{1}J_{2}v_{1}$ are linearly
dependent and hence $v_{2}$ and $J_{1}J_{2}v_{2}$ are as well.  By
Lemma \ref{QQ}, we have then that $v_{1} = \pm i J_{1}J_{2} v_{1}$.

Assume first that $v_{1} = i J_{1}J_{2} v_{1}$.  Then $v_{2} = i J_{1}J_{2} v_{2}$ as well or $\{J_{1}J_{2}v_{1},
J_{1}J_{2}v_{2}\}$ is not a cycle of generalized eigenvectors.  This is because if $v_{2}= - i J_{1}J_{2} v_{2}$,
then
\begin{align}
    -i a(J_{1}J_{2} v_{2}) &= a (-i J_{1}J_{2} v_{2}) = a v_{2} = \mu i v_{2} + v_{1} = \mu i (-i J_{1}J_{2}
    v_{2}) + i J_{1}J_{2} v_{1} = - i ( \mu i J_{1} J_{2}v_{2} - J_{1}J_{2} v_{1}) \notag
\intertext{which implies that}
    a(J_{1}J_{2} v_{2}) &= \mu i J_{1}J_{2}v_{2} - J_{1}J_{2}v_{1} \notag
\end{align}
Now as $v_{1}=i J_{1}J_{2}v_{1}$ and $v_{2} = i J_{1}J_{2} v_{2}$,
then, by Lemma \ref{QQ} again, we have that $v_{1}=u_{1} + i
J_{1}J_{2} u_{1}$ and $v_{2}= u_{2} + i J_{1}J_{2} u_{2}$ for real
vectors $u_{1},u_{2}$.

Let $w_{1}=\overline{v_{1}}=u_{1} - i J_{1}J_{2} u_{1}$ and
$w_{2}=\overline{v_{2}}=u_{2} - i J_{1}J_{2} u_{2}$, then
$\{w_{1},w_{2}\}$ is a cycle of generalized eigenvectors in $K_{-\mu
i}$.  For $j \in \{1,2\}$, we have by Corollary \ref{HH} that
$\om_{j}(v_{1},v_{2})=0$ and by Proposition \ref{KK} that
$\om_{j}(v_{1},w_{1}) = 0$.  However, we also have by Lemma \ref{RR}
that
\begin{align}
    \om_{j}(v_{1},w_{2}) &= \om_{j}(u_{1} + i J_{1}J_{2} u_{1}, u_{2} - i J_{1}J_{2}u_{2}) \notag \\
    &= \om_{j}(u_{1},u_{2}) - i \om_{j}(u_{1},J_{1}J_{2}u_{2}) + i \om_{j}(J_{1}J_{2}u_{1},u_{2}) +
    \om_{j}(J_{1}J_{2}u_{1},J_{1}J_{2}u_{2}) \notag \\
    &= \om_{j}(u_{1},u_{2}) - i \om_{j}(u_{1},J_{1}J_{2}u_{2}) + i \om_{j}(u_{1},J_{1}J_{2}u_{2}) -
    \om_{j}(u_{1},u_{2}) = 0 \notag
\end{align}
As $v_{1} \neq 0$, this implies that $\om_{j}$ is degenerate on $V$,
which is impossible.  By a similar argument, we get the same result
if $v_{1}=-i J_{1}J_{2} v_{1}$.  Therefore $a$ is diagonalizable.
$\blacksquare$ \vspace{.1 in}

\begin{lem}\label{NN}
    The following two matrices in $\h(J_{2})$ are not $\h$-symplectically similar for all $\la,\mu \neq 0$.  $$a_{1} =
    \begin{pmatrix} \la & \mu &0&0 \\ -\mu &\la &0&0 \\ 0&0&-\la&\mu
    \\ 0&0&-\mu&-\la \end{pmatrix}, \ \ \ \ \ a_{2} =
    \begin{pmatrix} \la &-\mu&0&0 \\ \mu &\la &0&0 \\ 0&0&-\la&-\mu
    \\ 0&0&\mu &-\la \end{pmatrix}.$$
\end{lem}

{\em Proof:} We prove this by showing that the equation $a_{1} A = A
a_{2}$ has no solution in $H(J_{2})$.  Let $A \in GL(V)$ be given
by $$A=\begin{pmatrix} d_{1}&d_{2}&d_{3}&d_{4} \\
d_{5}&d_{6}&d_{7}&d_{8} \\ d_{9}&d_{10}&d_{11}&d_{12} \\
d_{13}&d_{14}&d_{15}&d_{16} \end{pmatrix}.$$  By requiring $A$ to
satisfy the equation above, we find that $A$ must have the
form $$\begin{pmatrix} d_{1}&d_{2}&0&0 \\ d_{2}&-d_{1}&0&0 \\
0&0&d_{11}&d_{12} \\ 0&0&d_{12}&-d_{11} \end{pmatrix}.$$  However,
if we also require $A$ to be an element in $H(J_{2})$, that is,
require it to satisfy the additional equations $A^{t}J_{i}A = J_{i}$
for all $i \in \{1,2\}$, then we get, amongst other things, that
$d_{1}d_{11} +d_{2}d_{12}$ is equal to both 1 and -1, which is, of
course, impossible.  Thus there is no solution to $a_{1}A=Aa_{2}$ in
$H(J_{2})$ and consequently $a_{1}$ and $a_{2}$ are not
$\h$-symplectically similar. $\blacksquare$ \vspace{.1 in}

We now have the tools necessary to compute the canonical forms of $a
\in \h(J_{2})$.  We will consider first those $a$ that have real
eigenvalues and then those where one or more of the eigenvalues are
complex.  We'll classify according to the decomposition of $V$ into
generalized eigenspaces of $a$ when the eigenvalues are real and
$V_{c}$ when some or all of them are complex.  The cases will follow
these four decompositions of $V$ or $V_{c}$.
\begin{enumerate}
    \item $V=K_{\la} \oplus K_{-\la}$ where $\la \in \R$ is nonzero.
    \item $V=K_{0}$.
    \item $V_{c}=K_{z} \oplus K_{-z} \oplus K_{\overline{z}} \oplus
    K_{-\overline{z}}$ where $z=\la + \mu i$ for some $\la, \mu \in \R$
    such that $\la,\mu > 0$.
    \item $V_{c}=K_{\mu i} \oplus K_{-\mu i}$ where $\mu \in \R$ is
    positive.
\end{enumerate}
By Proposition \ref{PP} and Lemma \ref{SS}, this enumerates every
possible decomposition of $V$ or $V_{c}$ into generalized
eigenspaces of $a \in \h(J_{2})$.

\subsubsection{Case 1: $V=K_{\la} \oplus K_{-\la}$}

Also assume that $\la \in \R$ is nonzero.  We know, by Proposition
\ref{PP}, that $\dim K_{\la} = \dim K_{-\la} = 2$ and, by Lemma
\ref{QQ}, $K_{\la}$ and $K_{-\la}$ must contain at least two
dimensions of eigenvectors. Thus$K_{\la}$ and $K_{-\la}$ contain
only eigenvectors and hence $a$ is diagonalizable.

Let $v_{1} \in K_{\la}$, then by Corollary \ref{JJ}, there exists a
$w_{1} \in K_{-\la}$ such that $\om_{1}(v_{1},w_{1}) \neq 0$.  Let
$v_{2}=J_{1}J_{2}v_{1}$ and $w_{2}=J_{1}J_{2}w_{1}$.  Then by Lemma
\ref{QQ}, $v_{2}$ and $w_{2}$ are eigenvectors of $a$ corresponding
to $\la$ and $-\la$ respectively and are linearly independent of
$v_{1}$ and $w_{1}$.  Then the set $\{v_{1},v_{2},w_{1},w_{2}\}$ is
a basis for $V$ and has the following properties. Because
$v_{1},v_{2} \in K_{\la}$ and $w_{1},w_{2} \in K_{-\la}$, we have
that $\om_{i}(v_{1},v_{2})=\om_{i}(w_{1},w_{2})=0$.  In addition,
though, we have that
\begin{align}
    \om_{1}(v_{1},w_{2}) &= \om_{1}(v_{1},J_{1}J_{2}w_{1}) = -\om_{1}(J_{1}v_{1},J_{2}w_{1}) =
    \om_{2}((J_{1})^{2} v_{1}, w_{1}) = -\om_{2}(v_{1},w_{1}), \notag \\
    \om_{2}(v_{1},w_{2}) &= \om_{2}(v_{1},J_{1}J_{2}w_{1}) = -\om_{1}(v_{1},(J_{2})^{2}w_{1}) =
    \om_{1}(v_{1},w_{1}). \notag \displaybreak[0] \\
\intertext{Also}
    \om_{1}(v_{2},w_{1}) &= \om_{1}(J_{1}J_{2}v_{1},w_{1}) = -\om_{1}(w_{1},J_{1}J_{2}v_{1}) =
    \om_{2}(w_{1},v_{1}) = -\om_{2}(v_{1},w_{1}), \notag \\
    \om_{2}(v_{2},w_{1}) &= \om_{2}(J_{1}J_{2}v_{1},w_{1}) = -\om_{2}(w_{1},J_{1}J_{2}v_{1}) =
    -\om_{1}(w_{1},v_{1}) = \om_{1}(v_{1},w_{1}), \notag \displaybreak[0] \\
\intertext{and}
    \om_{1}(v_{2},w_{2}) &= -\om_{1}(J_{1}J_{2}v_{1},J_{1}J_{2}w_{1}) = -\om_{1}((J_{1})^{2}J_{2} v_{1},
    J_{2}w_{1}) \notag \\
    &= \om_{1}(J_{2}v_{1},J_{2}w_{1}) =\om_{1}(v_{1},(J_{2})^{2}w_{1}) = -\om_{1}(v_{1},w_{1}), \notag \\
    \om_{2}(v_{2},w_{2}) &= \om_{2}(J_{1}J_{2}v_{1},J_{1}J_{2}w_{1}) = \om_{2}((J_{1})^{2}J_{2}v_{1},J_{2}w_{1}) \notag \\
    &= -\om_{2}(J_{2}v_{1},J_{2}w_{1}) = \om_{2}((J_{2})^{2}v_{1},w_{1}) = -\om_{2}(v_{1},w_{1}). \notag
\end{align}

At this point, we pick a new basis in the following way
\begin{align}
    v_{1}' &= v_{1}, &v_{2}' &= v_{2} \notag \\
    w_{1}' &= w_{1} - \frac{\om_{1}(v_{2},w_{1})}{\om_{1}(v_{2},w_{2})} w_{2}, & w_{2}' &= J_{1}J_{2}w_{1}' \notag
\end{align}
Note that as $\om_{1}(v_{2},w_{2})=-\om_{1}(v_{1},w_{1})$, we have that $\om_{1}(v_{2},w_{2}) \neq 0$ and so this
basis change makes sense.  Note that $w_{1}',w_{2}' \in K_{-\la}$ and so are still eigenvectors.  This implies
that $\om_{i}(v_{1}',v_{2}') = \om_{i}(w_{1}',w_{2}') = 0$.  Furthermore as $w_{2}' = J_{1}J_{2}w_{1}'$, the
properties computed above still hold. In addition,
\begin{align}
    \om_{1}(v_{2}',w_{1}') &= \om_{1}\pnth{v_{2},w_{1} - \frac{\om_{1}(v_{2},w_{1})}{\om_{1}(v_{2},w_{2})} w_{2}} =
    \om_{1}(v_{2},w_{1}) - \frac{\om_{1}(v_{2},w_{1})}{\om_{1}(v_{2},w_{2})} \om_{1}(v_{2},w_{2}) = 0, \notag \\
\intertext{which implies that}
    \om_{2}(v_{2}',w_{2}') &= \om_{1}(v_{1}',w_{2}') = -\om_{2}(v_{1}',w_{1}') = \om_{1}(v_{2}',w_{1}') = 0. \notag
    \displaybreak[0] \\
\intertext{Moreover,}
    \om_{1}(v_{1}',w_{1}') &= \om_{1}\pnth{v_{1},w_{1}-\frac{\om_{1}(v_{2},w_{1})}{\om_{1}(v_{2},w_{2})} w_{2}} =
    \om_{1}(v_{1},w_{1}) -\frac{\om_{1}(v_{2},w_{1})}{\om_{1}(v_{2},w_{2})} \om_{1}(v_{1},w_{2}) \notag \\
    &= \om_{1}(v_{1},w_{1}) + \frac{\om_{1}(v_{1},w_{1})}{\om_{1}(v_{1},w_{2})} \om_{1}(v_{1},w_{2}) = 2
    \om_{1}(v_{1},w_{1}) \notag \\
\intertext{and hence}
    \om_{2}(v_{1}',w_{2}') &= \om_{2}(v_{2}',w_{1}') = -\om_{1}(v_{2}',w_{2}') = \om_{1}(v_{1}',w_{1}') = 2
    \om_{1}(v_{1},w_{1}) \neq 0. \notag
\end{align}

Finally, if $\om_{1}(v_{1}',w_{1}')>0$, then replace $v_{1}'$, $w_{1}'$, and $w_{2}'$ with $-v_{1}'$, $-w_{1}'$,
and $-w_{2}'$ and relabel them $v_{1}'$, $w_{1}'$, and $w_{2}'$ respectively.  If $\om_{1}(v_{1}',w_{1}')<0$, then
replace $v_{2}'$, $w_{1}'$, and $w_{2}'$ with $-v_{2}'$, $-w_{1}'$, and $-w_{2}'$ and relabel them $v_{2}'$,
$w_{1}'$, and $w_{2}'$ respectively.  This will ensure the following relationship $$\om_{1}(v_{1}',w_{1}') =
\om_{1}(v_{2}',w_{2}') = \om_{2}(v_{1}',w_{2}') = -\om_{2}(v_{2}',w_{1}') >0.$$  Then we pick a new basis for $V$,
call it $\be$, as follows $$\be = \left\{ \frac{1}{\sqrt{\om_{1}(v_{1}'w_{1}')}} v_{1}',
\frac{1}{\sqrt{\om_{1}(v_{1}'w_{1}')}} v_{2}', \frac{1}{\sqrt{\om_{1}(v_{1}'w_{1}')}} w_{1}',
\frac{1}{\sqrt{\om_{1}(v_{1}'w_{1}')}} w_{2}' \right\},$$ and relabel the vectors $v_{1}''$, $v_{2}''$, $w_{1}''$,
and $w_{2}''$ respectively.  Then the matrix $$W=\begin{pmatrix} v_{1}'' & v_{2}'' & w_{1}'' & w_{2}''
\end{pmatrix}$$ has the properties that $W^{t}J_{1}W=J_{1}$ and $W^{t} J_{2} W = J_{2}$ and is thus in
$H(J_{2})$.  Finally, in this basis, $a$ is of form (1) in Theorem
\ref{EE} with $\la \neq 0$.

\subsubsection{Case 2: $V=K_{0}$}

If $a$ is diagonalizable, then, as its only eigenvalue is 0, $a$ is
the zero transformation and is already in form (1) in Theorem
\ref{EE} with $\la = 0$.

If $a$ is not diagonalizable, then there exists a cycle of
generalized eigenvectors of length two in $K_{0}$.  Let
$\{v_{1},v_{2}\}$ be such a cycle.  Let $w_{1}=J_{1}J_{2}v_{1}$ and
$w_{2}=J_{1}J_{2}v_{2}$. Then by Lemma \ref{TT}, $\{w_{1},w_{2}\}$
is also a cycle of generalized eigenvectors and linearly independent
to $\{v_{1},v_{2}\}$.  Then $\{v_{1},v_{2},w_{1},w_{2}\}$ is a basis
for $K_{0}$.  By Proposition \ref{KK}, we see that
$\om_{i}(v_{1},w_{1}) = 0$.  Furthermore
\begin{align}
    \om_{1}(v_{2},w_{2}) &= \om_{1}(v_{2},J_{1}J_{2}v_{2}) = -\om_{1}(J_{1}v_{2},J_{2}v_{2}) = \om_{2}(J_{1}
    v_{2},J_{1}v_{2}) = 0 \notag \\
    \om_{2}(v_{2},w_{2}) &= \om_{2}(v_{2},J_{1}J_{2}v_{2}) = -\om_{1}(J_{2}v_{2},J_{2}v_{2}) = 0 \notag
\intertext{In addition, we have that}
    \om_{1}(v_{1},w_{2}) &= \om_{1}(v_{1},J_{1}J_{2}v_{2}) = -\om_{1}(J_{1}v_{1},J_{2}v_{2}) =
    \om_{2}((J_{1})^{2} v_{1}, v_{2}) = -\om_{2}(v_{1},v_{2}), \notag \\
    \om_{2}(v_{1},w_{2}) &= \om_{2}(v_{1},J_{1}J_{2}v_{2}) = -\om_{1}(v_{1},(J_{2})^{2}v_{2}) =
    \om_{1}(v_{1},v_{2}). \notag \displaybreak[0] \\
\intertext{Also}
    \om_{1}(w_{1},v_{2}) &= \om_{1}(J_{1}J_{2}v_{1},v_{2}) = -\om_{1}(v_{2},J_{1}J_{2}v_{1}) =
    \om_{2}(v_{2},v_{1}) = -\om_{2}(v_{1},v_{2}), \notag \\
    \om_{2}(w_{1},v_{2}) &= \om_{2}(J_{1}J_{2}v_{1},v_{2}) = -\om_{2}(v_{2},J_{1}J_{2}v_{1}) =
    -\om_{1}(v_{2},v_{1}) = \om_{1}(v_{1},v_{2}), \notag \displaybreak[0] \\
\intertext{and}
    \om_{1}(w_{1},w_{2}) &= -\om_{1}(J_{1}J_{2}v_{1},J_{1}J_{2}v_{2}) = -\om_{1}((J_{1})^{2}J_{2} v_{1},
    J_{2}v_{2}) \notag \\
    &= \om_{1}(J_{2}v_{1},J_{2}v_{2}) =\om_{1}(v_{1},(J_{2})^{2}v_{2}) = -\om_{1}(v_{1},v_{2}), \notag \\
    \om_{2}(w_{1},w_{2}) &= \om_{2}(J_{1}J_{2}v_{1},J_{1}J_{2}v_{2}) = \om_{2}((J_{1})^{2}J_{2}v_{1},J_{2}v_{2}) \notag \\
    &= -\om_{2}(J_{2}v_{1},J_{2}v_{2}) = \om_{2}((J_{2})^{2}v_{1},v_{2}) = -\om_{2}(v_{1},v_{2}). \notag
\end{align}
As $\om$ is non-degenerate on $V$, we have that $\om_{1}(v_{1},v_{2})$ and $\om_{2}(v_{1},v_{2})$ are not both
zero.

First assume that $\om_{2}(v_{1},v_{2}) \neq 0$.  Then we solve the equation $$c_{1}x^{2} + 2c_{2}x -c_{1} = 0$$
where $c_{1}=\om_{2}(v_{1},v_{2})$ and $c_{2}=\om_{1}(v_{1},v_{2})$.  The discriminant of this equation is
$4(c_{2})^{2} + 4(c_{1})^{2}$.  As $c_{1},c_{2} \in \R$ and $c_{1} \neq 0$, we see that the discriminant is
strictly positive and this quadratic equation has two real solutions.  Let $x_{0} \in \R$ be a solution to the
equation. Then pick a new basis as follows
\begin{align}
    v_{1}' &= v_{1} - x_{0} w_{1}, & v_{2}' &= v_{2} - x_{0} w_{2},
    \notag \\
    w_{1}' &= J_{1}J_{2} v_{1}', & w_{2}' &= J_{1}J_{2} v_{2}'.
    \notag
\end{align}
Note that $\{v_{1}',v_{2}'\}$ is still a cycle of generalized eigenvectors.  Then by Lemma \ref{TT},
$\{w_{1}',w_{2}'\}$ is one as well.  Also note that
\begin{align}
    \om_{i}(v_{1}',w_{1}') &= \om_{i}(v_{2}',w_{2}')=0 \ \text{for all $i
    \in \{1,2\}$,} \notag \displaybreak[0] \\
\intertext{and by a similar argument to that above, we have that}
    \om_{1}(v_{1}',v_{2}') &= \om_{2}(v_{1}',w_{2}') = \om_{2}(w_{1}',v_{2}') = -\om_{1}(w_{1}',w_{2}'), \notag \\
    \om_{2}(v_{1}',v_{2}') &= -\om_{1}(v_{1}',w_{2}') = -\om_{1}(w_{1}',v_{2}') = -\om_{2}(w_{1}',w_{2}'). \notag
\end{align}
However, in this basis we also have that
\begin{align}
    \om_{2}(v_{1}',v_{2}') &= \om_{2}(v_{1}-x_{0}w_{1},v_{2}-x_{0} w_{2}) = \om_{2}(v_{1},v_{2}) - x_{0}
    \om_{2}(v_{1},w_{2}) -x_{0}\om_{2}(w_{1},v_{2}) + (x_{0})^{2} \om_{2}(w_{1},w_{2}) \notag \\
    &= \om_{2}(v_{1},v_{2}) - 2 x_{0} \om_{1}(v_{1},v_{2}) - (x_{0})^{2} \om_{2}(v_{1},v_{2}) = -(c_{1}(x_{0})^{2}
    +2c x_{0} - c_{1}) = 0. \notag
\end{align}
Then as $\om$ is non-degenerate, we must have that $$\om_{1}(v_{1}',v_{2}') = \om_{2}(v_{1}',w_{2}') =
\om_{2}(w_{1}',v_{2}') = -\om_{1}(w_{1}',w_{2}') \neq 0.$$

Now assume that $\om_{2}(v_{1},v_{2})=0$ in the first place, then $\om_{1}(v_{1},v_{2})\neq 0$.  In this case,
simply relabel $v_{1}$, $v_{2}$, $w_{1}$, and $w_{2}$ as $v_{1}'$, $v_{2}'$, $w_{1}'$, and $w_{2}'$ respectively
and we have the exact situation as described above.

If $\om_{1}(v_{1}',v_{2}')>0$, then let
\begin{align}
    v_{1}'' &= v_{1}' & v_{2}'' &= v_{2}' \notag \\
    w_{1}'' &= -w_{1}' & w_{2}'' &= w_{2}' \notag \\
\intertext{On the other hand, if $\om_{1}(v_{1}',v_{2}') <0$, then let}
    v_{1}'' &= -w_{1}' & v_{2}'' &= - w_{2}' \notag \\
    w_{1}'' &= -v_{1}' & w_{2}'' &= v_{2}' \notag
\end{align}
Either case will yield that $$\om_{1}(v_{1}'',v_{2}'') = \om_{1}(w_{1}'',w_{2}'') = \om_{2}(v_{1}'',w_{2}'') =
-\om_{2}(w_{1}'',v_{2}'')$$ and any other pair on either $\om_{i}$ products to 0.  However, it also yields that
$$a v_{2}'' = v_{1}'', \ \ \text{but} \ \ a w_{2}'' = -w_{1}''.$$  Then we pick a final basis for $V$, $$\left\{
\frac{1}{\sqrt{\om_{1}(v_{1}'',v_{2}'')}} v_{1}'', \ \frac{1}{\sqrt{\om_{1}(v_{1}'',v_{2}'')}} w_{1}'', \
\frac{1}{\sqrt{\om_{1}(v_{1}'',v_{2}'')}} v_{2}'', \ \frac{1}{\sqrt{\om_{1}(v_{1}'',v_{2}'')}} w_{2}'' \right\}$$
and relabel the vectors $v_{1}'''$, $w_{1}'''$, $v_{2}'''$, and $w_{2}'''$ respectively.  Then the matrix
$$W=\begin{pmatrix} v_{1}''' & w_{1}''' & v_{2}''' & w_{2}''' \end{pmatrix}$$ is in $H(J_{2})$ because
$W^{t}J_{i} W = J_{i}$ for all $i \in \{1,2\}$.  Finally, in this basis, $a$ is of form (2) in Theorem \ref{EE}.

\subsubsection{Case 3: $V_{c}=K_{z} \oplus K_{-z} \oplus K_{\overline{z}} \oplus K_{-\overline{z}}$}

We also assume that $z=\la + \mu i$ such that $\la, \mu
>0$. Before we begin, we extend a lemma proved in the section on the
symplectic Lie algebra.

\begin{lem}\label{OO}
    $\overline{\om_{i}(v,w)} = \om_{i}(\overline{v},\overline{w})$ for all $v,w \in V$ and $i \in \{1,2\}$.
\end{lem}

{\em Proof.}  Let $v,w \in V$ and $i \in \{1,2\}$.  Then
$\overline{\om_{i}(v,w)} = \overline{v^{t}J_{i}w} = \overline{v}^{t}
J_{i} \overline{w} = \om_{i}(\overline{v},\overline{w})$.
$\blacksquare$ \vspace{.1 in}

Clearly in this case, each eigenspace must be of dimension one and
consequently contains only eigenvectors.  Hence $a$ is
diagonalizable.  Let $v_{1} \in K_{z}$ and $w_{1} \in K_{-z}$, and
let $v_{2}=\overline{v_{1}}$ and $w_{2}=\overline{w_{1}}$.  Then
$v_{2} \in K_{\overline{z}}$ and $w_{2} \in K_{-\overline{z}}$.  By
the non-degeneracy of $\om_{i}$, we know that $\om_{i}(v_{1},w_{1})
\neq 0$ and $\om_{i}(v_{2},w_{2}) \neq 0$ for all $i \in \{1,2\}$
and any other pair products to 0 using either form.  In fact,
$$\overline{\om_{i}(v_{1},w_{1})} =
\om_{i}(\overline{v_{1}},\overline{w_{1}})= \om_{i}(v_{2},w_{2}).$$
Now we make the change of basis
\begin{align}
    v_{1}'&=v_{1}+v_{2} & v_{2}' &= J_{1}J_{2}v_{1}' \notag \\
    w_{1}'&=w_{1}+w_{2} & w_{2}' &= J_{1}J_{2}w_{1}' \notag
    \displaybreak[0] \\
\intertext{As $J_{1}J_{2}v_{1} \in K_{z}$ and $K_{z}$ is one dimensional, we see that $v_{1}$ and
$J_{1}J_{2}v_{1}$ are linearly dependent. Then by Lemma \ref{QQ}, we have that $\Im(v_{1}) = \pm
J_{1}J_{2}\Re(v_{1})$.  The same can be said of $w_{1}$ and $J_{1}J_{2}w_{1}$.  This implies that we have}
    v_{1}'&= 2\Re(v_{1}) & v_{2}' &= \pm 2 \Im(v_{1})=\mp i (v_{1}-v_{2}) \notag \\
    w_{1}'&= 2\Re(w_{1}) & w_{2}' &= \pm 2 \Im(w_{2})=\mp i (w_{1}-w_{2}) \notag
    \displaybreak[0]
\end{align}
Clearly $\om_{i}(v_{1}',v_{2}')=\om_{i}(w_{1}',w_{2}')=0$, but in
addition, this implies, as before, that
\begin{align}
    \om_{1}(v_{1}',w_{1}')&=\om_{2}(v_{1}',w_{2}')=\om_{2}(v_{2}',w_{1}')
    = -\om_{1}(v_{2}',w_{2}') \notag \\
    \om_{2}(v_{1}',w_{1}')&=-\om_{1}(v_{1}',w_{2}') =
    -\om_{2}(v_{2}',w_{1}') = -\om_{2}(v_{2}',w_{2}') \notag
    \displaybreak[0] \\
\intertext{Furthermore, we have that}
    \om_{1}(v_{1}',w_{1}') &= \om_{1}(v_{1}+v_{2},w_{1}+w_{2}) = \om_{1}(v_{1},w_{1}) +\om_{1}(v_{1},w_{2}) +\om_{1}(v_{2},w_{1})
    +\om_{1}(v_{2},w_{2}) \notag \\
    &= \om_{1}(v_{1},w_{1}) + \om_{1}(v_{2},w_{2}) = 2 \Re(\om_{1}(v_{1},w_{1}))
    \notag \displaybreak[0]\\
\intertext{and}
    \om_{1}(v_{1}',w_{2}') &= \om_{1}(v_{1}+v_{2},\mp
    i(w_{1}-w_{2})) \notag \\
    &= \mp i\om_{1}(v_{1},w_{1})\pm i\om_{1}(v_{1},w_{2})\mp i\om_{1}(v_{2},w_{1})\pm i\om_{1}(v_{2},w_{2}) \notag \\
    &= \mp i \pnth{\om_{1}(v_{1},w_{1}) - \om_{1}(v_{2},w_{2})}
    = \pm 2 \Im(\om_{1}(v_{1},w_{1})) \notag \displaybreak[0] \\
\intertext{Together, these imply that}
    2\Re(\om_{1}(v_{1},w_{1})) &= \om_{1}(v_{1}',w_{1}')=\om_{2}(v_{1}',w_{2}')=\om_{2}(v_{2}',w_{1}')
    = -\om_{1}(v_{2}',w_{2}') \notag \\
    \pm 2 \Im(\om_{1}(v_{1},w_{1})) &= -\om_{2}(v_{1}',w_{1}')=\om_{1}(v_{1}',w_{2}') =
    \om_{1}(v_{2}',w_{1}') = \om_{2}(v_{2}',w_{2}'). \notag
\end{align}
Finally, we also have that
\begin{align}
    a v_{1}' &= \la v_{1}' \mp \mu v_{2}', & a v_{2}' &= \pm \mu v_{1}' +
    \la v_{2}', \notag \\
    a w_{1}' &= -\la v_{1}' \pm  \mu v_{2}', & a w_{2}' &= \mp \mu v_{1}'
    -\la v_{1}'. \notag
\end{align}

As $\om_{1}(v_{1},w_{1}) \neq 0$, we have that either
$\Re(\om_{1}(v_{1},w_{1})) \neq 0$, or $\Re(\om_{1}(v_{1},w_{1}))=0$
and $\Im(\om_{1}(v_{1},w_{1})) \linebreak[0] \neq 0$.

If $\Re(\om_{1}(v_{1},w_{1})) \neq 0$, then let
\begin{align}
    v_{1}'' &= v_{1}' + \frac{\om_{1}(v_{1}',w_{2}')}{\om_{1}(v_{1}',w_{1}')} v_{2}', & v_{2}'' &=
    J_{1}J_{2}v_{1}'', \notag \\
    w_{1}'' &= w_{1}', & w_{2}''&=w_{2}'. \notag
\end{align}
The set $\{v_{1}'',v_{2}'',w_{1}'',w_{2}''\}$ is still a basis for $V$.  We know that
\begin{align}
    v_{2}'' &= J_{1}J_{2} v_{1}'' = J_{1}J_{2} \pnth{v_{1}' +
    \frac{\om_{1}(v_{1}',w_{2}')}{\om_{1}(v_{1}',w_{1}')} v_{2}'} \notag \\
    &= v_{2}' + \frac{\om_{1}(v_{1}',w_{2}')}{\om_{1}(v_{1}',w_{1}')} (J_{1}J_{2})^{2}v_{1}' = v_{2}' -
    \frac{\om_{1}(v_{1}',w_{2}')}{\om_{1}(v_{1}',w_{1}')} v_{1}'. \notag \displaybreak[0] \\
\intertext{Then we have the following properties}
     \om_{1}(v_{1}'',w_{2}'') &= \om_{1}\pnth{v_{1}' + \frac{\om_{1}(v_{1}',w_{2}')}{\om_{1}(v_{1}',w_{1}')} v_{2}',
     w_{2}'} = \om_{1}(v_{1}',w_{2}') + \frac{\om_{1}(v_{1}',w_{2}')}{\om_{1}(v_{1}',w_{1}')} \om_{1}(v_{2}',w_{2}') \notag \\
     &= \om_{1}(v_{1}',w_{2}') - \frac{\om_{1}(v_{1}',w_{2}')}{\om_{1}(v_{1}',w_{1}')} \om_{1}(v_{1}',w_{1}') = 0, \notag
     \displaybreak[0] \\
     \om_{1}(v_{1}'',w_{1}'') &= \om_{1}\pnth{v_{1}' + \frac{\om_{1}(v_{1}',w_{2}')}{\om_{1}(v_{1}',w_{1}')} v_{2}',
     w_{1}'} = \om_{1}(v_{1}',w_{1}') + \frac{\om_{1}(v_{1}',w_{2}')}{\om_{1}(v_{1}',w_{1}')} \om_{1}(v_{2}',w_{1}'),
     \notag \\
\intertext{and $\om_{1}(v_{1}'',v_{2}'')=\om_{1}(w_{1}'',w_{2}'') = 0$. Clearly, as $\om_{1}$ and $\om_{2}$ are
non-degenerate, we have that $\om_{1}(v_{1}'',w_{1}'') \neq 0$.  Thus we have as well that}
    \om_{1}(v_{1}'',w_{1}'') &=\om_{2}(v_{1},w_{2}) = \om_{2}(v_{2},w_{1}) = - \om_{1}(v_{2},w_{2}) \neq 0 \notag \\
    \om_{2}(v_{1}'',w_{1}'') &= - \om_{1}(v_{1}'',w_{2}'') = -\om_{1}(v_{2}'',w_{1}'') = -\om_{2}(v_{2},w_{2}) = 0
    \notag
\end{align}
In addition,
\begin{align}
    a v_{1}'' &= a \pnth{v_{1}' + \frac{\om_{1}(v_{1}',w_{2}')}{\om_{1}(v_{1}',w_{1}')} v_{2}'} = \la v_{1}' \mp
    \mu v_{2}' + \frac{\om_{1}(v_{1}',w_{2}')}{\om_{1}(v_{1}',w_{1}')} (\pm \mu v_{1}' + \la v_{2}') = \la v_{1}''
    \mp \mu v_{2}'', \notag \\
    a v_{2}'' &= a \pnth{v_{2}' - \frac{\om_{1}(v_{1}',w_{2}')}{\om_{1}(v_{1}',w_{1}')} v_{1}'}
    = \pm \mu v_{1}' + \la v_{2}' - \frac{\om_{1}(v_{1}',w_{2}')}{\om_{1}(v_{1}',w_{1}')} (\la v_{1}' \mp \mu
    v_{2}) = \pm \mu v_{1}'' + \la v_{2}'' \notag
\end{align}
and as $w_{1}'' = w_{1}'$ and $w_{2}'' = w_{1}'$, their relationship remains unchanged.

Now we consider the other possibility.  If $\Re(\om_{1}(v_{1},w_{1}))=0$, then $\Im(\om_{1}(v_{1},w_{1})) \neq 0$.
Then let
\begin{align}
    v_{1}''&= v_{1}'+v_{2}', & v_{2}'' &= J_{1}J_{2}(v_{1}''),  \notag \\
    w_{1}''&=w_{1}'+w_{2}', & w_{2}''&= J_{1}J_{2}(w_{1}''). \notag
\end{align}
The set $\{v_{1}'',v_{2}'',w_{1}'',w_{2}''\}$ is still a basis for $V$.  We also have that
\begin{align}
    v_{2}'' &= J_{1}J_{2}v_{1}'' = J_{1}J_{2}(v_{1}'+v_{2}') = v_{2}' + (J_{1}J_{2})^{2} v_{1}' = v_{2}'-v_{1}',
    \notag \\
    w_{2}'' &= J_{1}J_{2}w_{1}'' = J_{1}J_{2}(w_{1}'+w_{2}') = w_{2}' + (J_{1}J_{2})^{2} w_{1}' = w_{2}'-w_{1}'.
    \notag \displaybreak[0] \\
\intertext{Then we have the following properties}
     \om_{1}(v_{1}'',w_{2}'') &= \om_{1}(v_{1}'+v_{2}', w_{2}' - w_{1}') = \om_{1}(v_{1}',w_{2}') -\om_{1}(v_{1}',w_{1}') +
     \om_{1}(v_{2}',w_{2}') - \om_{1}(v_{2}',w_{1}') = 0 \notag \\
     \om_{1}(v_{1}'',w_{1}'') &= \om_{1}(v_{1}'+v_{2}',w_{1}'+w_{2}') = \om_{1}(v_{1}',w_{1}') + \om_{1}(v_{1}',w_{2}') +
     \om_{1}(v_{2}',w_{1}') + \om_{1}(v_{2}',w_{2}') \notag \\
     &= 2 \om_{1}(v_{1}',w_{2}') \notag \displaybreak[0] \\
\intertext{and
$\om_{1}(v_{1}'',v_{2}'')=\om_{1}(w_{1}'',w_{2}'')=0$.  Then
$\om_{1}(v_{1}'',w_{1}'')\neq 0$ and we have that}
    \om_{1}(v_{1}'',w_{1}'') &=\om_{2}(v_{1},w_{2}) = \om_{2}(v_{2},w_{1}) = - \om_{1}(v_{2},w_{2}) \neq 0 \notag \\
    \om_{2}(v_{1}'',w_{1}'') &= - \om_{1}(v_{1}'',w_{2}'') = -\om_{1}(v_{2}'',w_{1}'') = -\om_{2}(v_{2},w_{2}) = 0
    \notag
\end{align}
In addition,
\begin{align}
    av_{1}'' &= a (v_{1}'+v_{2}') = \la v_{1}' \mp \mu v_{2}' \pm \mu v_{1}'
    +\la v_{2}' = \la v_{1}'' \mp \mu v_{2}'', \notag \\
    a v_{2}'' &= a(v_{2}'-v_{1}') = \pm \mu v_{1}' +\la v_{2}' - \la v_{1}' \pm \mu v_{2}' =
    \pm \mu v_{1}'' +\la v_{2}'', \notag
\end{align}
and similarly $$a w_{1}'' = -\la w_{1}'' \pm \mu w_{2}'', \ \ \ \ \ a w_{2}'' = \mp \mu w_{1}'' -\la w_{2}''.$$

In either case, we know that $$\om_{1}(v_{1}'',w_{1}'')
=\om_{2}(v_{1},w_{2}) = \om_{2}(v_{2},w_{1}) =
-\om_{1}(v_{2}'',w_{2}'') \neq 0.$$ If $\om_{1}(v_{1}'',w_{1}'')
<0$, then replace $v_{1}''$ with $-v_{1}''$ and relabel.  This will
make $$\om_{1}(v_{1}'',w_{1}'') =\om_{2}(v_{1},w_{2}) =
-\om_{2}(v_{2},w_{1}) = \om_{1}(v_{2}'',w_{2}'') > 0,$$ but will
yield that $a v_{1}'' = \la v_{1}'' \pm \mu v_{2}''$ and $a v_{2}''
= \mp \mu v_{1}''+\la v_{2}''$. On the other hand, if
$\om_{1}(v_{1}'',w_{1}'')>0$, then replace $v_{2}''$ with $-v_{2}''$
and relabel. This will again give us the same situation as above.

Finally, we pick a basis, $\be$, for $V$ in the following way $$\be=
\left\{ \frac{1}{\sqrt{\om_{1}(v_{1}'',w_{1}'')}} v_{1}'', \
\frac{1}{\sqrt{\om_{1}(v_{1}'',w_{1}'')}} v_{2}'', \
\frac{1}{\sqrt{\om_{1}(v_{1}'',w_{1}'')}} w_{1}'', \
\frac{1}{\sqrt{\om_{1}(v_{1}'',w_{1}'')}} w_{2}'' \right\}$$ and
relabel the vectors $v_{1}'''$, $v_{2}'''$, $w_{1}'''$, and
$w_{2}'''$ respectively.  Then the matrix
$$W=\begin{pmatrix} v_{1}''' & v_{2}''' & w_{1}''' & w_{2}''' \end{pmatrix}$$ is in $H(J_{2})$ because
$W^{t}J_{i}W=J_{i}$ for all $i \in \{1,2\}$.  In this basis, $a$ is
of form (3) in Theorem \ref{EE} with $\la \neq 0$.  The $\vep$ is
present because by Lemma \ref{NN}, the two cases are not
$\h$-symplectically similar.

\subsubsection{Case 4: $V_{c}=K_{\mu i} \oplus K_{-\mu i}$}

We also assume that $\mu \in \R$ such that $\mu >0$.  In addition, by Lemma \ref{SS}, we have that $a$ is
diagonalizable.

Now let $v_{1} \in K_{\mu i}$ be an eigenvector.  As $v_{1}$ is complex, it is not clear that $J_{1}J_{2}v_{1}$ is
linearly independent of $v_{1}$.  If $v_{1}$ and $J_{1}J_{2}v_{1}$ are linearly independent, then let
$v_{1}'=v_{1}+iJ_{1}J_{2}v_{1}$ and $v_{2}'=v_{1}-iJ_{1}J_{2}v_{1}$.  Then $v_{1}'$ and $v_{2}'$ are linearly
independent and such that $$iJ_{1}J_{2}v_{1}'=v_{1}', \ \ \text{and} \ \ -iJ_{1}J_{2}v_{2}'=v_{2}'.$$  Let
$w_{1}'=\overline{v_{1}'}$ and $w_{2}=\overline{v_{2}'}$.  Then $w_{1}'$ and $w_{2}'$ are eigenvectors in $K_{-\mu
i}$.  Moreover, by Proposition \ref{GG}, we know that $\om_{i}(v_{1}',v_{2}')=\om_{i}(w_{1}',w_{2}')=0$.
Furthermore
\begin{align}
    w_{1}'&=\overline{v_{1}'}=\overline{v_{1}+iJ_{1}J_{2}v_{1}} =\overline{v_{1}}+\overline{iJ_{1}J_{2}v_{1}} =
    \overline{v_{1}}-iJ_{1}J_{2}\overline{v_{1}} \notag \\
    w_{2}'&=\overline{v_{2}'}=\overline{v_{1}-iJ_{1}J_{2}v_{1}} =\overline{v_{1}}-\overline{iJ_{1}J_{2}v_{1}} =
    \overline{v_{1}}+iJ_{1}J_{2}\overline{v_{1}} \notag \\
\intertext{This implies by Proposition \ref{RR} that for all $i \in \{1,2\}$}
    \om_{i}(v_{1}',w_{1}') &= \om_{i}(v_{1}+iJ_{1}J_{2}v_{1}, \overline{v_{1}}-iJ_{1}J_{2}\overline{v_{1}}) \notag \\
    &= \om_{i}(v_{1},\overline{v_{1}}) -i \om_{i}(v_{1},J_{1}J_{2}\overline{v_{1}})+i
    \om_{i}(J_{1}J_{2}v_{1},\overline{v_{1}}) + \om_{i}(J_{1}J_{2}v_{1},J_{1}J_{2}\overline{v_{1}}) \notag \\
    &= \om_{i}(v_{1},\overline{v_{1}}) - \om_{i}(v_{1},\overline{v_{1}}) = 0. \notag \\
    \om_{i}(v_{2}',w_{2}') &= \om_{i}(v_{1}-iJ_{1}J_{2}v_{1}, \overline{v_{1}}+iJ_{1}J_{2}\overline{v_{1}}) \notag \\
    &= \om_{i}(v_{1},\overline{v_{1}}) +i \om_{i}(v_{1},J_{1}J_{2}\overline{v_{1}})-i
    \om_{i}(J_{1}J_{2}v_{1},\overline{v_{1}}) + \om_{i}(J_{1}J_{2}v_{1},J_{1}J_{2}\overline{v_{1}}) \notag \\
    &= \om_{i}(v_{1},\overline{v_{1}}) - \om_{i}(v_{1},\overline{v_{1}}) = 0. \notag
\end{align}
Then $\om_{i}(v_{1}',w_{2}') \neq 0$ by the non-degeneracy of
$\om_{i}$.  We relabel our vectors in the following manner
\begin{align}
    v_{1}'' &= v_{1}', & v_{2}'' &= w_{1}', & w_{1}'' &= w_{2}', & w_{2}'' &= v_{2}'. \notag
\end{align}
Then we have that $\om_{i}(v_{1}'',v_{2}'')=0$,
$\om_{i}(w_{1}'',w_{2}'')=0$, $\om_{i}(v_{1}'',w_{2}'') = 0$, and
$\om_{i}(v_{1}'',w_{1}'') \neq 0$ and $v_{2}'' = \overline{v_{1}''}$
and $w_{2}''=\overline{w_{1}''}$.  From this point, we follow the
same argument given in the previous case with relabeling $v_{i}''$
as $v_{i}$ and $w_{i}''$ as $w_{i}$ and recalling the fact that
$J_{1}J_{2}v_{1}''$ and $v_{1}''$ are linearly dependent, as are
$J_{1}J_{2}w_{1}''$ and $w_{1}''$. Then the real $\h$-symplectic
canonical form of $a$ is of form (3) in Theorem \ref{EE} with $\la =
0$. \vspace{.1 in}

Now consider the situation where $v_{1}$ and $J_{1}J_{2}v_{1}$ are
linearly dependent, then $\pm i J_{1}J_{2}v_{1}=v_{1}$.  Let $v_{2}
\in K_{\mu i}$ be another eigenvector linearly independent of
$v_{1}$. If $J_{1}J_{2}v_{2}$ and $v_{2}$ are linearly independent
then we proceed as above using $v_{2}$ as $v_{1}$. If
$J_{1}J_{2}v_{2}$ and $v_{2}$ are linearly dependent, then we know
that $\pm i J_{1}J_{2}v_{2}=v_{2}$.

If the sign on $i$ is different between the two above equations,
then the two vectors $$v_{1}+v_{2} \ \ \text{and} \ \
J_{1}J_{2}(v_{1}+v_{2})$$ are linearly independent.  To see this,
write $\pm i J_{1}J_{2}v_{1}=v_{1}$, then $\mp i
J_{1}J_{2}v_{2}=v_{2}$.  Then $$J_{1}J_{2}(v_{1}+v_{2}) =
J_{1}J_{2}v_{1} + J_{1}J_{2}v_{2} = \mp i v_{1} \pm i v_{2} = \mp i
(v_{1} - v_{2})$$ which is clearly independent of $v_{1}+v_{2}$.  If
this is the case, we relabel $v_{1}+v_{2}$ as simply $v_{1}$ and
proceed as in the situation where $v_{1}$ and $J_{1}J_{2}v_{1}$ were
linearly independent to begin with.

If the sign on $i$ is the same in the two equations, then let
$v_{1}'=v_{1} + i v_{2}$ and $v_{2}'=v_{1}-iv_{2}$. Then $v_{1}'$
and $v_{2}'$ are linearly independent and $\pm i J_{1}J_{2} v_{1}' =
v_{1}'$ and $\mp i J_{1}J_{2} v_{2}' = v_{2}'$.  If $-iJ_{1}J_{2}
v_{1}' =v_{1}'$, then relabel $v_{1}'$ and $v_{2}'$ as $v_{2}'$ and
$v_{1}'$ respectively.  After making this change if necessary, this
ensures that $$i J_{1}J_{2}v_{1}' = v_{1}', \ \ \text{and} \ \
-iJ_{1}J_{2}v_{2}' = v_{2}'.$$  Then we proceed as in the argument
above after it reached this relation.

This covers every possible decomposition of $V$ or $V_{c}$ into generalized eigenspaces of $a$ and so completes
the proof of Theorem \ref{EE}.

%
\setlength{\baselineskip}{11pt}

%
\setlength{\baselineskip}{22pt}

%
%

\chapter{A SAMPLE OF THE CLASSIFICATION OF DIMENSION SEVEN}

As has been stated before, the method we use to find all the
possible solvable indecomposable Lie algebras with codimension one
nilradicals focuses on the nilradical of the algebra.  The possible
nilradicals of a Lie algebra of dimension $n$ are all the nilpotent
algebras of dimension $n-1$.  We reference the classifications of
Winternitz and Gong for a list of possible nilpotent algebras of
dimensions one through six and give the list in Appendix \ref{a1}
\cite{winternitz-1,gong}.  Since much of the classification of these
algebras is the same from one nilradical to the next, we offer in
this section the step by step classification of the seven
dimensional algebras that stem from four nilradicals of dimension
six.  The derivations of the first nilradical form a solvable Lie
algebra.  The semi-simple part of the Lie algebra formed by the
derivations of the second nilradical is isomorphic to $\sll(2,\R)$.
The semi-simple part of the Lie algebra formed by the derivations of
the third nilradical is isomorphic to $\si(4,\R)$.  And the
semi-simple part of the Lie algebra formed by the derivations of the
fourth nilradical is isomorphic to the representation of
$\so(3,1,\R)$ discussed in Chapter \ref{c2}. These four nilradicals
are representative of the different situations we have to deal with
when classifying these algebras through dimension seven, and as
such, the classification of the algebras from the other nilradicals
are similar. In Appendix \ref{a2}, we give a complete table of all
solvable indecomposable Lie algebras with codimension one
nilradicals from dimension two through dimension seven.

\section{A Solvable Derivation Algebra}

In this section, we classify those seven-dimensional algebras, $\g$,
whose nilradical, $NR(\g)$, is isomorphic to one with structure
equations
\begin{align}
    \notag [\f_{2},\f_{5}]&=\f_{1}, & [\f_{3},\f_{4}]&=-\f_{1}, &
    [\f_{3},\f_{6}]&=\f_{2}, & [\f_{4},\f_{6}]&=\f_{3}, &
    [\f_{5},\f_{6}]&=\f_{4}.
\end{align}
This is Nilradical 16 listed under the six dimensional nilradicals
in Appendix \ref{a1}.  Let $\{\f_{1},\ldots,\f_{6}\}$ be a basis for
the nilradical and pick the vectors so that they have the above
structure equations.  Complete this to a basis for $\g$ by including
a vector $\f_{7} \not \in NR(\g)$.  Then as $D\g \in NR(\g)$, we can
view $\ad{\f_{7}}$ as a transformation on the nilradical. After
requiring the algebra to satisfy the Jacobi property, we have that
$\ad{\f_{7}}$ is of the form $$\ad{\f_{7}}=\begin{pmatrix}
2x_{5}+3y_{6}& b_{1}& c_{1}&
d_{1}& x_{1}& y_{1} \\ 0& x_{5}+3y_{6}& c_{2}& 0& x_{2}& -d_{1} \\
0& 0& x_{5}+2y_{6}& c_{2}& 0& c_{1} \\ 0& 0& 0& x_{5}+y_{6}& c_{2}&
-b_{1} \\ 0& 0& 0& 0& x_{5}& 0 \\ 0& 0& 0& 0& 0& y_{6}
\end{pmatrix}.$$

The other nonzero $\add$ matrices are
\begin{align}
    \notag \ad{\f_{2}} &= \begin{pmatrix} 0&0&0&0&1&0 \\ 0&0&0&0&0&0 \\ 0&0&0&0&0&0
    \\ 0&0&0&0&0&0 \\ 0&0&0&0&0&0 \\ 0&0&0&0&0&0 \end{pmatrix}, &
    \ad{\f_{3}} &= \begin{pmatrix} 0&0&0&-1&0&0 \\ 0&0&0&0&0&1 \\ 0&0&0&0&0&0
    \\ 0&0&0&0&0&0 \\ 0&0&0&0&0&0 \\ 0&0&0&0&0&0 \end{pmatrix},
    \displaybreak[0] \\
    \notag \ad{\f_{4}} &= \begin{pmatrix} 0&0&1&0&0&0 \\ 0&0&0&0&0&0 \\
    0&0&0&0&0&1 \\ 0&0&0&0&0&0 \\ 0&0&0&0&0&0 \\ 0&0&0&0&0&0
    \end{pmatrix}, & \ad{\f_{5}} &= \begin{pmatrix} 0&-1&0&0&0&0 \\ 0&0&0&0&0&0
    \\ 0&0&0&0&0&0 \\ 0&0&0&0&0&1 \\ 0&0&0&0&0&0 \\ 0&0&0&0&0&0
    \end{pmatrix}, \displaybreak[0] \\
    \notag \ad{\f_{6}} &= \begin{pmatrix} 0&0&0&0&0&0 \\ 0&0&-1&0&0&0 \\
    0&0&0&-1&0&0 \\ 0&0&0&0&-1&0 \\ 0&0&0&0&0&0 \\ 0&0&0&0&0&0
    \end{pmatrix}.
\end{align}
This will allow us to zero out the $x_{1}$, $d_{1}$, $c_{1}$,
$b_{1}$, and $c_{2}$ positions in one simple perturbing of $\f_{7}$.

Then the only entries in $\ad{\f_{7}}$ that we don't have under
control are $x_{5}$, $y_{6}$, $x_{2}$ and $y_{1}$. We'll find that
if $x_{5}$ and $y_{6}$ were simultaneously 0, then algebra becomes
nilpotent, so clearly that cannot happen.  However, to deal with the
other two entries, we now turn our attention to computing the
automorphisms of the nilradical.

By using Maple to compute the derivations and exponentiate them, we
find that the automorphism group of the nilradical is generated by
the following nine one-parameter groups of transformations.
\begin{align}
    \notag \A_{1} &= \begin{pmatrix} (s_{1})^{3}&0&0&0&0&0 \\ 0&1&0&0&0&0 \\
    0&0&s_{1}&0&0&0 \\ 0&0&0&(s_{1})^{2}&0&0 \\ 0&0&0&0&(s_{1})^{3}&0 \\
    0&0&0&0&0&\frac{1}{s_{1}} \end{pmatrix}, & \A_{2} &= \begin{pmatrix} 1&s_{2}&0&0&0&0 \\ 0&1&0&0&0&0 \\
    0&0&1&0&0&0 \\ 0&0&0&1&0&-s_{2} \\ 0&0&0&0&1&0 \\ 0&0&0&0&0&1
    \end{pmatrix}, \displaybreak[0] \\
    \notag \A_{3} &= \begin{pmatrix} 1&0&0&0&0&0 \\ 0&(s_{3})^{3}&0&0&0&0 \\
    0&0&s_{3}&0&0&0 \\ 0&0&0&\frac{1}{s_{3}}&0&0 \\ 0&0&0&0&\pnth{\frac{1}{s_{3}}}^{3}&0 \\
    0&0&0&0&0&(s_{3})^{2} \end{pmatrix}, & \A_{4} &= \begin{pmatrix} 1&0&s_{4}&0&0&\frac{(s_{4})^{2}}{2} \\
    0&1&0&0&0&0 \\ 0&0&1&0&0&s_{4} \\ 0&0&0&1&0&0 \\ 0&0&0&0&1&0 \\ 0&0&0&0&0&1
    \end{pmatrix}, \displaybreak[0] \\
    \notag \A_{5} &= \begin{pmatrix} 1&0&0&0&0&0 \\ 0&1&s_{5}&\frac{(s_{5})^{2}}{2}& \frac{(s_{5})^{3}}{6}&0 \\
    0&0&1&s_{5}&\frac{(s_{5})^{2}}{2}&0 \\ 0&0&0&1&s_{5}&0 \\ 0&0&0&0&1&0 \\ 0&0&0&0&0&1
    \end{pmatrix}, & \A_{6} &= \begin{pmatrix} 1&0&0&s_{6}&0&0 \\ 0&1&0&0&0&0 \\
    0&0&1&0&0&-s_{6} \\ 0&0&0&1&0&0 \\ 0&0&0&0&1&0 \\ 0&0&0&0&0&1
    \end{pmatrix}, \displaybreak[0] \\
    \notag \A_{7} &= \begin{pmatrix} 1&0&0&0&s_{7}&0 \\ 0&1&0&0&0&0 \\
    0&0&1&0&0&0 \\ 0&0&0&1&0&0 \\ 0&0&0&0&1&0 \\ 0&0&0&0&0&1
    \end{pmatrix}, & \A_{8} &= \begin{pmatrix} 1&0&0&0&0&0 \\ 0&1&0&0&s_{8}&0 \\
    0&0&1&0&0&0 \\ 0&0&0&1&0&0 \\ 0&0&0&0&1&0 \\ 0&0&0&0&0&1
    \end{pmatrix}, \displaybreak[0] \\
    \notag \A_{9} &= \begin{pmatrix} 1&0&0&0&0&s_{9} \\ 0&1&0&0&0&0 \\
    0&0&1&0&0&0 \\ 0&0&0&1&0&0 \\ 0&0&0&0&1&0 \\ 0&0&0&0&0&1
    \end{pmatrix}.
\end{align}
However, conjugation by $\A_{2}$, $\A_{4}$, $\A_{5}$, $\A_{6}$, and
$\A_{7}$ will only affect the positions that we'll zero out by a
basis change.  As such, they become less useful in simplifying the
$\ad{\f_{7}}$ matrix.

As there is no semisimple part of the derivation algebra, we use
only a single parent case.

\subsection{Parent Case 1:}

We start with the $\ad{\f_{7}}$ matrix of the form
$$\ad{\f_{7}}=\begin{pmatrix} 2x_{5}+3y_{6}& b_{1}& c_{1}&
d_{1}& x_{1}& y_{1} \\ 0& x_{5}+3y_{6}& c_{2}& 0& x_{2}& -d_{1} \\
0& 0& x_{5}+2y_{6}& c_{2}& 0& c_{1} \\ 0& 0& 0& x_{5}+y_{6}& c_{2}&
-b_{1} \\ 0& 0& 0& 0& x_{5}& 0 \\ 0& 0& 0& 0& 0& y_{6}
\end{pmatrix}.$$ Conjugate $\ad{\f_{7}}$ by $\A_{8}$ and note that if $y_{6} \neq 0$,
then we could let $s_{8}=-\frac{x_{2}}{3y_{6}}$ which would zero out
the $x_{2}$ position.  That is, if we apply the automorphism
$\A_{8}$ with $s_{8}=-\frac{x_{2}}{3y_{6}}$, the resulting
$\ad{\f_{7}}$ matrix would have a 0 in the $x_{2}$ position.  It
also happens to change the $x_{1}$ position, but we will simply
relabel what ends up in that position as $x_{1}$.  We'll usually use
this same general method when conjugating by an automorphism.  At
any rate, this yields two possible cases.
\begin{enumerate}
    \item In this first case, either $y_{6} \neq 0$ and we moved the $x_{2}$ position to
    0, or $y_{6}=0$ and $x_{2}=0$ already.
    \item In this second case, $y_{6}=0$, but $x_{2} \neq 0$.
\end{enumerate}

\subsubsection{Subcase 1:}

The resulting $\ad{\f_{7}}$ matrix is as follows
$$\ad{\f_{7}}=\begin{pmatrix} 2x_{5}+3y_{6}& b_{1}& c_{1}&
d_{1}& x_{1}& y_{1} \\ 0& x_{5}+3y_{6}& c_{2}& 0& 0& -d_{1} \\
0& 0& x_{5}+2y_{6}& c_{2}& 0& c_{1} \\ 0& 0& 0& x_{5}+y_{6}& c_{2}&
-b_{1} \\ 0& 0& 0& 0& x_{5}& 0 \\ 0& 0& 0& 0& 0& y_{6}
\end{pmatrix}.$$  Next conjugate by $\A_{9}$. In similar manner
as above, if $x_{5} \neq -y_{6}$, then we can pick
$s_{9}=-\frac{y_{1}}{2(y_{6}+x_{5})}$ and this conjugation will
result in moving $y_{1}$ to 0.  Thus we have another two cases
\begin{enumerate}
    \item $x_{5} \neq -y_{6}$ and we can move $y_{1}$ to 0, or
    $x_{5}=-y_{6}$ and $y_{1}=0$ already.
    \item $x_{5}=-y_{6}$, but $y_{1} \neq 0$.
\end{enumerate}

\subsubsection{Subcase 1.1:}

We have $\ad{\f_{7}}$ as follows $$\ad{\f_{7}}=\begin{pmatrix}
2x_{5}+3y_{6}& b_{1}& c_{1}&
d_{1}& x_{1}& 0 \\ 0& x_{5}+3y_{6}& c_{2}& 0& 0& -d_{1} \\
0& 0& x_{5}+2y_{6}& c_{2}& 0& c_{1} \\ 0& 0& 0& x_{5}+y_{6}& c_{2}&
-b_{1} \\ 0& 0& 0& 0& x_{5}& 0 \\ 0& 0& 0& 0& 0& y_{6}
\end{pmatrix}.$$  We already know that $x_{5}$ and $y_{6}$ are not
simultaneously 0.  This allows us to bifurcate on this as well.
\begin{enumerate}
    \item $y_{6} \neq 0$.
    \item $y_{6}=0$, which implies that $x_{5} \neq 0$ or the
    algebra is nilpotent.
\end{enumerate}

\subsubsection{Subcase 1.1.1:}

In this section, simply make the basis change $$\e_{i}=\f_{i} \ \
\text{for} \ \ 1 \leq i \leq 6, \ \ \text{and} \ \
\e_{7}=-\frac{1}{y_{6}}(\f_{7}-x_{1}\f_{2}+d_{1}\f_{3}-c_{1}\f_{4}
+b_{1}\f_{5}+c_{2}\f_{6})$$ and let $a = \frac{x_{5}}{y_{6}}$.
Notice how we perturbed $\f_{7}$.  This makes use of the other
$\add$ matrices and zeros out the $x_{1}$, $d_{1}$, $c_{1}$,
$b_{1}$, and $c_{2}$ positions.  This yields the structure equations
\begin{align}
    \notag [\e_{2},\e_{5}]&=\e_{1}, & [\e_{3},\e_{6}]&=\e_{2}, &
    [\e_{4},\e_{6}]&=\e_{3}, &[\e_{4},\e_{6}]&=\e_{3}, \\
    \notag [\e_{1},\e_{7}]&=(2a+3)\e_{1}, & [\e_{2},\e_{7}]&=
    (a+3)\e_{2}, &[\e_{3},\e_{7}]&=(a+2)\e_{3}, & [\e_{4},\e_{7}]&=
    (a+1)\e_{4}, \\
    \notag [\e_{5},\e_{7}]&=a\e_{5}, & [\e_{6},\e_{7}]&=\e_{6},
\end{align}
with $a \in \R$.  In the table in Appendix \ref{a2}, this is
[7,[6,16],1,1].

\subsubsection{Subcase 1.1.2:}

Here we assumed that $y_{6}=0$ and so we have the $\ad{\f_{7}}$
matrix $$\ad{\f_{7}}=\begin{pmatrix} 2x_{5}& b_{1}&
c_{1}& d_{1}& x_{1}& 0 \\ 0& x_{5}& c_{2}& 0& 0& -d_{1} \\
0& 0& x_{5}& c_{2}& 0& c_{1} \\ 0& 0& 0& x_{5}& c_{2}& -b_{1}
\\ 0& 0& 0& 0& x_{5}& 0 \\ 0& 0& 0& 0& 0& 0
\end{pmatrix}.$$  Then we make the basis change $$\e_{i}=\f_{i} \ \
\text{for} \ \ 1 \leq i \leq 6, \ \ \text{and} \ \
\e_{7}=-\frac{1}{x_{5}}(\f_{7}-x_{1}\f_{2}+d_{1}\f_{3}-c_{1}\f_{4}
+b_{1}\f_{5}+c_{2}\f_{6}),$$  which yields the structure equations
\begin{align}
    \notag [\e_{2},\e_{5}]&=\e_{1}, & [\e_{3},\e_{6}]&=\e_{2}, &
    [\e_{4},\e_{6}]&=\e_{3}, &[\e_{4},\e_{6}]&=\e_{3}, &[\e_{1},\e_{6}]&=2\e_{1}, \\
    \notag [\e_{2},\e_{7}]&=\e_{2}, &[\e_{3},\e_{7}]&=\e_{3}, &
    [\e_{4},\e_{7}]&=\e_{4}, &[\e_{5},\e_{7}]&=\e_{5}.
\end{align}
This is [7,[6,16],1,2] in the table.

\subsubsection{Subcase 1.2:}

In this section, we assumed that $x_{5}=-y_{6}$ and $y_{1} \neq 0$.
This will yield the $\ad{\f_{7}}$ matrix
$$\ad{\f_{7}}=\begin{pmatrix} y_{6}&
b_{1}& c_{1}& d_{1}& x_{1}& y_{1} \\ 0& 2y_{6}& c_{2}& 0& 0& -d_{1} \\
0& 0& y_{6}& c_{2}& 0& c_{1} \\ 0& 0& 0& 0& c_{2}& -b_{1}
\\ 0& 0& 0& 0& -y_{6}& 0 \\ 0& 0& 0& 0& 0& y_{6}
\end{pmatrix},$$  and implies that $y_{6} \neq 0$ or the algebra is
nilpotent. At this point, we conjugate by $\A_{3}$. As $y_{6}$ and
$y_{1}$ are both nonzero, we can pick
$s_{3}=\frac{\sqrt{\abs{y_{1}y_{6}}}}{y_{1}}$, which will scale the
$y_{1}$ position to $\pm y_{6}$.  This also requires us to relabel
some of the entries of $\ad{\f_{7}}$ as we did previously.  This is
how we'll usually deal with an automorphism that scales the entries
of $\ad{\f_{7}}$.  Then we make the basis change
$$\e_{i}=\f_{i} \ \ \text{for} \ \ 1 \leq i \leq 6, \ \ \text{and} \
\ \e_{7}=-\frac{1}{y_{6}}(\f_{7}-x_{1}\f_{2}+d_{1}\f_{3}-c_{1}\f_{4}
+b_{1}\f_{5}+c_{2}\f_{6})$$ and let $\vep = \pm 1$.  This will yield
the structure equations
\begin{align}
    \notag [\e_{2},\e_{5}]&=\e_{1}, & [\e_{3},\e_{6}]&=\e_{2}, &
    [\e_{4},\e_{6}]&=\e_{3}, &[\e_{4},\e_{6}]&=\e_{3}, &[\e_{1},\e_{6}]&=\e_{1}, \\
    \notag [\e_{2},\e_{7}]&=2\e_{2}, &[\e_{3},\e_{7}]&=\e_{3}, &
    [\e_{5},\e_{7}]&=-\e_{5}, &[\e_{6},\e_{7}]&=\vep \e_{1}+\e_{6}.
\end{align}
with $\vep^{2}=1$.  In the table, this is [7,[6,16],1,3].

\subsubsection{Subcase 2:}

In this section, we assumed that $y_{6}=0$ and $x_{2} \neq 0$.  This
gives us the $\ad{\f_{7}}$ matrix $$\ad{\f_{7}}=\begin{pmatrix}
2x_{5}& b_{1}& c_{1}&
d_{1}& x_{1}& y_{1} \\ 0& x_{5}& c_{2}& 0& x_{2}& -d_{1} \\
0& 0& x_{5}& c_{2}& 0& c_{1} \\ 0& 0& 0& x_{5}& c_{2}& -b_{1}
\\ 0& 0& 0& 0& x_{5}& 0 \\ 0& 0& 0& 0& 0& 0
\end{pmatrix},$$  and implies that $x_{5} \neq 0$ or the algebra is
nilpotent.  We now conjugate by $\A_{9}$.  As $x_{5} \neq 0$, we
pick $s_{9}=-\frac{y_{1}}{2x_{5}}$, which will zero out the $y_{1}$
position.  Then we conjugate by $\A_{1}$, and as $x_{5}$ and $x_{2}$
are both nonzero, we can let $s_{1}=\frac{\sqrt[3]{
x_{5}(x_{2})^{2}}}{x_{2}}$, which will scale the $x_{2}$ position to
$x_{5}$.  Finally, we make the basis change $$\e_{i}=\f_{i} \ \
\text{for} \ \ 1 \leq i \leq 6, \ \ \text{and} \ \
\e_{7}=-\frac{1}{x_{5}}(\f_{7}-x_{1}\f_{2}+d_{1}\f_{3}-c_{1}\f_{4}
+b_{1}\f_{5}+c_{2}\f_{6}).$$  This gives us the structure equations
\begin{align}
    \notag [\e_{2},\e_{5}]&=\e_{1}, & [\e_{3},\e_{6}]&=\e_{2}, &
    [\e_{4},\e_{6}]&=\e_{3}, &[\e_{4},\e_{6}]&=\e_{3}, &[\e_{1},\e_{6}]&=2\e_{1}, \\
    \notag [\e_{2},\e_{7}]&=\e_{2}, &[\e_{3},\e_{7}]&=\e_{3}, &
    [\e_{4},\e_{7}]&=\e_{4}, &[\e_{5},\e_{7}]&=\e_{2}+\e_{5}.
\end{align}
This is [7,[6,16],1,4] in the table and completes the classification
of seven dimensional algebras with this nilradical. \vspace{.1 in}

Now that the reader has the idea of how we use the automorphisms to
simplify the $\ad{\f_{7}}$ matrix, we will, to conserve space,
refrain from giving the explicit value of the automorphism parameter
that will make the desired change, but simply state its existence
and the parameter it pivots on. We now move on to the classification
of an algebra with a nilradical whose derivation algebra has a
semisimple part isomorphic to $\sll(2,\R)$.

\section{A Derivation Algebra with Semisimple Part Isomorphic to
$\sll(2,\R)$}

In this section, we classify those algebras, $\g$, with nilradical,
$NR(\g)$, isomorphic to the six dimensional nilpotent algebra with
structure equations
\begin{align}
    \notag [\f_{4},\f_{5}]&=\f_{2}, & [\f_{4},\f_{6}]&=\f_{3}, &
    [\f_{5},\f_{6}]&=\f_{4}.
\end{align}
This is Nilradical 4 from the list of six dimensional nilradicals in
Appendix \ref{a1}. We will also assume that we have a basis for $\g$
such that $\{\f_{1},\ldots,\f_{6}\}$ forms a basis for $NR(\g)$ and
has the structure equations given above. Let the last vector in
$\g$, $\f_{7}$, be an arbitrary vector not contained in $NR(\g)$. By
the Jacobi property, the $\add$ matrix of $\f_{7}$ must be of the
form
$$\ad{\f_{7}}=\begin{pmatrix} a_{1}&0&0&0&x_{1}&y_{1} \\
a_{2}&2x_{5}+y_{6}&y_{5}&-y_{4}&x_{2}&y_{2} \\
a_{3}&x_{6}&x_{5}+2y_{6}&x_{4}&x_{3}&y_{3} \\
0&0&0&x_{5}+y_{6}&x_{4}&y_{4} \\ 0&0&0&0&x_{5}&y_{5} \\
0&0&0&0&x_{6}&y_{6} \end{pmatrix}.$$

First, we look at the other nonzero $\add$ matrices.  They are
\begin{align}
    \notag \ad{\f_{4}} &= \begin{pmatrix} 0&0&0&0&0&0 \\ 0&0&0&0&1&0
    \\ 0&0&0&0&0&1 \\ 0&0&0&0&0&0 \\ 0&0&0&0&0&0 \\ 0&0&0&0&0&0
    \end{pmatrix} & \ad{\f_{5}} &= \begin{pmatrix} 0&0&0&0&0&0 \\
    0&0&0&-1&0&0 \\ 0&0&0&0&0&0 \\ 0&0&0&0&0&1 \\ 0&0&0&0&0&0 \\
    0&0&0&0&0&0 \end{pmatrix}, \\
    \notag \ad{\f_{6}}&= \begin{pmatrix}0&0&0&0&0&0 \\
    0&0&0&0&0&0 \\ 0&0&0&-1&0&0 \\ 0&0&0&0&-1&0 \\ 0&0&0&0&0&0 \\
    0&0&0&0&0&0 \end{pmatrix}.
\end{align}
This will allow us, by perturbing $\f_{7}$, to annihilate, $y_{4}$,
$x_{4}$, and, if we can move the $x_{2}$ value to the $y_{3}$ value,
both of these as well, by perturbing $\f_{7}$.

We now take a look at the automorphisms.  Again we use Maple to
compute a basis for the derivation algebra and to compute the Levi
decomposition of it.  Then we find that a basis for the semisimple
part is
\begin{align}
    \notag  \begin{pmatrix} 0&0&0&0&0&0 \\ 0&0&1&0&0&0 \\
    0&0&0&0&0&0 \\ 0&0&0&0&0&0 \\ 0&0&0&0&0&1 \\ 0&0&0&0&0&0
    \end{pmatrix}, & & \begin{pmatrix} 0&0&0&0&0&0 \\ 0&1&0&0&0&0 \\
    0&0&-1&0&0&0 \\ 0&0&0&0&0&0 \\ 0&0&0&0&1&0 \\ 0&0&0&0&0&-1
    \end{pmatrix}, & & \begin{pmatrix} 0&0&0&0&0&0 \\ 0&0&0&0&0&0 \\
    0&1&0&0&0&0 \\ 0&0&0&0&0&0 \\ 0&0&0&0&0&0 \\ 0&0&0&0&1&0
    \end{pmatrix}.
\end{align}
Label these as $\D_{1},\D_{2}, \D_{3}$ respectively.  These three
vectors form a basis for a subalgebra isomorphic to $\sll(2,\R)$.
This is clear by matching these basis vectors with the corresponding
standard basis vectors of $\sll(2,\R)$,
\begin{align}
    \notag  \D_{1} &\mapsto \begin{pmatrix} 0&1 \\ 0&0 \end{pmatrix},
    & \D_{2} &\mapsto \begin{pmatrix} 1&0 \\ 0&-1 \end{pmatrix}, &
    \D_{3}&\mapsto \begin{pmatrix} 0&0 \\ 1&0 \end{pmatrix}.
\end{align}
However, we can use another basis vector, this time from the radical
of the derivation algebra, namely, $$\begin{pmatrix} 0&0&0&0&0&0 \\ 0&3&0&0&0&0 \\
0&0&3&0&0&0 \\ 0&0&0&2&0&0 \\ 0&0&0&0&1&0 \\ 0&0&0&0&0&1
\end{pmatrix},$$ which we'll call $\D_{4}$, to do the following.
Note that
\begin{align}
    \notag \frac{1}{2}(D_{2}+D_{4}) &= \begin{pmatrix} 0&0&0&0&0&0
    \\ 0&2&0&0&0&0 \\ 0&0&1&0&0&0 \\ 0&0&0&1&0&0 \\ 0&0&0&0&1&0 \\
    0&0&0&0&0&0 \end{pmatrix}, & -\frac{1}{2}(D_{2}-D_{4}) &=
    \begin{pmatrix} 0&0&0&0&0&0 \\ 0&1&0&0&0&0 \\ 0&0&2&0&0&0 \\
    0&0&0&1&0&0 \\ 00&0&0&0&0&0 \\ 0&0&0&0&0&1 \end{pmatrix}.
\end{align}
Then the subalgebra with basis $\left\{\frac{1}{2}(D_{2}+D_{4}),
D_{1},D_{3}, -\frac{1}{2}(D_{2}-D_{4}) \right\}$ is isomorphic to
$\gl(2,\R)$ via the map
\begin{align}
    \notag \frac{1}{2}(D_{2}+D_{4})&\mapsto \begin{pmatrix} 1&0 \\ 0&0
    \end{pmatrix}, & D_{1} &\mapsto \begin{pmatrix} 0&1 \\ 0&0
    \end{pmatrix}, & D_{3} &\mapsto \begin{pmatrix} 0&0 \\ 1&0
    \end{pmatrix}, & -\frac{1}{2}(D_{2}-D_{4}) &\mapsto
    \begin{pmatrix} 0&0 \\ 0&1 \end{pmatrix},
\end{align}
which is easily checked, as the structure equations of the two
algebras are the same.

Note that the lower right hand $2 \times 2$ blocks of the derivation
vectors are identical to the basis vectors of $\gl(2,\R)$.  This
leads us to suspect that if we exponentiate the derivation vectors,
then we can find a nilradical automorphism that has an arbitrary
$GL(2,\R)$ matrix in the lower right hand $2 \times 2$ block.  And
sure enough, we use Maple to check that the matrix
$$\A_{0}=\begin{pmatrix} 1&0&0&0&0&0 \\ 0&a(ad-bc)&b(ad-bc)&0&0&0 \\
0&c(ad-bc)&d(ad-bc)&0&0&0 \\ 0&0&0&ad-bc&0&0 \\ 0&0&0&0&a&b \\
0&0&0&0&c&d \end{pmatrix}$$ is an automorphism of the nilradical.

We next use Maple to compute a complete basis for the derivation
algebra and exponentiate it.  Then we find that the following
one-parameter groups of transformations generate the automorphism
group of the nilradical.
\begin{align}
    \notag \A_{1}&=\begin{pmatrix} s_{1}&0&0&0&0&0 \\ 0&1&0&0&0&0 \\
    0&0&1&0&0&0 \\ 0&0&0&1&0&0 \\ 0&0&0&0&1&0 \\ 0&0&0&0&0&1
    \end{pmatrix}, &\A_{2}&=\begin{pmatrix} 1&0&0&0&0&0 \\ s_{2}&1&0&0&0&0 \\
    0&0&1&0&0&0 \\ 0&0&0&1&0&0 \\ 0&0&0&0&1&0 \\ 0&0&0&0&0&1
    \end{pmatrix}, \displaybreak[0] \\
    \notag \A_{3} &= \begin{pmatrix} 1&0&0&0&0&0 \\ 0&1&0&0&0&0 \\
    s_{3}&0&1&0&0&0 \\ 0&0&0&1&0&0 \\ 0&0&0&0&1&0 \\ 0&0&0&0&0&1
    \end{pmatrix}, & \A_{4} &= \begin{pmatrix} 1&0&0&0&0&0 \\ 0&s_{4}&0&0&0&0 \\
    0&0&(s_{4})^{2}&0&0&0 \\ 0&0&0&s_{4}&0&0 \\ 0&0&0&0&1&0 \\
    0&0&0&0&0&s_{4} \end{pmatrix} \displaybreak[0] \\
    \notag \A_{5} &= \begin{pmatrix} 1&0&0&0&0&0 \\ 0&1&0&0&0&0 \\
    0&s_{5}&1&0&0&0 \\ 0&0&0&1&0&0 \\ 0&0&0&0&1&0 \\ 0&0&0&0&s_{5}&1
    \end{pmatrix}, & \A_{6} &= \begin{pmatrix} 1&0&0&0&0&0 \\ 0&1&s_{6}&0&0&0 \\
    0&0&1&0&0&0 \\ 0&0&0&1&0&0 \\ 0&0&0&0&1&s_{6} \\ 0&0&0&0&0&1
    \end{pmatrix}, \displaybreak[0] \\
    \notag \A_{7} &= \begin{pmatrix} 1&0&0&0&0&0 \\ 0&(s_{7})^{2}&0&0&0&0 \\
    0&0&s_{7}&0&0&0 \\ 0&0&0&s_{7}&0&0 \\ 0&0&0&0&s_{7}&0 \\ 0&0&0&0&0&1
    \end{pmatrix}, & \A_{8} &= \begin{pmatrix} 1&0&0&0&0&0 \\ 0&1&0&s_{8}&0&-\frac{(s_{8})^{2}}{2} \\
    0&0&1&0&0&0 \\ 0&0&0&1&0&-s_{8} \\ 0&0&0&0&1&0 \\ 0&0&0&0&0&1
    \end{pmatrix}, \displaybreak[0] \\
    \notag \A_{9} &= \begin{pmatrix} 1&0&0&0&0&0 \\ 0&1&0&0&0&0 \\
    0&0&1&s_{9}&\frac{(s_{9})^{2}}{2}&0 \\ 0&0&0&1&s_{9}&0 \\ 0&0&0&0&1&0 \\ 0&0&0&0&0&1
    \end{pmatrix}, & \A_{10} &= \begin{pmatrix} 1&0&0&0&s_{10}&0 \\ 0&1&0&0&0&0 \\
    0&0&1&0&0&0 \\ 0&0&0&1&0&0 \\ 0&0&0&0&1&0 \\ 0&0&0&0&0&1
    \end{pmatrix}, \displaybreak[0] \\
    \notag \A_{11} &= \begin{pmatrix} 1&0&0&0&0&0 \\ 0&1&0&0&s_{11}&0 \\
    0&0&1&0&0&0 \\ 0&0&0&1&0&0 \\ 0&0&0&0&1&0 \\ 0&0&0&0&0&1
    \end{pmatrix}, & \A_{12} &= \begin{pmatrix} 1&0&0&0&0&0 \\ 0&1&0&0&0&0 \\
    0&0&1&0&s_{12}&0 \\ 0&0&0&1&0&0 \\ 0&0&0&0&1&0 \\ 0&0&0&0&0&1
    \end{pmatrix}, \displaybreak[0] \\
    \notag \A_{13} &= \begin{pmatrix} 1&0&0&0&0&s_{13} \\ 0&1&0&0&0&0 \\
    0&0&1&0&0&0 \\ 0&0&0&1&0&0 \\ 0&0&0&0&1&0 \\ 0&0&0&0&0&1
    \end{pmatrix}, &\A_{14} &= \begin{pmatrix} 1&0&0&0&0&0 \\ 0&1&0&0&0&s_{14} \\
    0&0&1&0&0&0 \\ 0&0&0&1&0&0 \\ 0&0&0&0&1&0 \\ 0&0&0&0&0&1
    \end{pmatrix}, \displaybreak[0] \\
    \notag \A_{15} &= \begin{pmatrix} 1&0&0&0&0&0 \\ 0&1&0&0&0&0 \\
    0&0&1&0&0&s_{15} \\ 0&0&0&1&0&0 \\ 0&0&0&0&1&0 \\ 0&0&0&0&0&1
    \end{pmatrix}.
\end{align}
Note that $\A_{8}$ and $\A_{9}$ could be used to zero out the
$y_{4}$ and $x_{4}$ positions.  As we can do this with a simple
perturbing of $\f_{7}$, these automorphisms will be less useful than
the others in simplifying $\ad{\f_{7}}$.

Now, we start classifying the possible forms of $\ad{\f_{7}}$.
First, we note that by block multiplication, we can conjugate
$\ad{\f_{7}}$ by $\A_{0}$ and pick $a,b,c,d$ to put the lower right
hand $2 \times 2$ block of $\ad{\f_{7}}$ into real Jordan form. This
yields three parent cases depending on the real Jordan form of this
$2\times2$ block.  They are
\begin{align}
    \notag 1. \ \ \ad{\f_{7}}&=\begin{pmatrix} a_{1}&0&0&0&x_{1}&y_{1} \\
    a_{2}&2\la_{1}+\la_{2}&0&-y_{4}&x_{2}&y_{2} \\
    a_{3}&0&\la_{1}+2\la_{2}&x_{4}&x_{3}&y_{3} \\
    0&0&0&\la_{1}+\la_{2}&x_{4}&y_{4} \\ 0&0&0&0&\la_{1}&0 \\
    0&0&0&0&0&\la_{2} \end{pmatrix}, \displaybreak[0]\\
    \notag 2. \ \ \ad{\f_{7}}&=\begin{pmatrix} a_{1}&0&0&0&x_{1}&y_{1} \\
    a_{2}&3\la_{1}&\la_{2}&-y_{4}&x_{2}&y_{2} \\
    a_{3}&-\la_{2}&3\la_{1}&x_{4}&x_{3}&y_{3} \\
    0&0&0&2\la_{1}&x_{4}&y_{4} \\ 0&0&0&0&\la_{1}&\la_{2} \\
    0&0&0&0&-\la_{2}&\la_{1} \end{pmatrix}, \ \begin{array}{l} \text{with $\la_{2}
    \neq 0$ and the eigenvalues ordered} \\ \text{so that
    $\frac{\la_{1}}{\la_{2}} \geq 0$.} \end{array} \displaybreak[0]\\
    \notag 3. \ \ \ad{\f_{7}}&=\begin{pmatrix} a_{1}&0&0&0&x_{1}&y_{1} \\
    a_{2}&3\la&1&-y_{4}&x_{2}&y_{2} \\
    a_{3}&0&3\la&x_{4}&x_{3}&y_{3} \\
    0&0&0&2\la&x_{4}&y_{4} \\ 0&0&0&0&\la&1 \\
    0&0&0&0&0&\la \end{pmatrix}.
\end{align}

\subsection{Parent Case 1:}

We consider here, the $\ad{\f_{7}}$ matrix
$$\ad{\f_{7}}=\begin{pmatrix} a_{1}&0&0&0&x_{1}&y_{1} \\
a_{2}&2\la_{1}+\la_{2}&0&-y_{4}&x_{2}&y_{2} \\
a_{3}&0&\la_{1}+2\la_{2}&x_{4}&x_{3}&y_{3} \\
0&0&0&\la_{1}+\la_{2}&x_{4} &y_{4} \\ 0&0&0&0&\la_{1}&0 \\
0&0&0&0&0&\la_{2} \end{pmatrix}.$$  We note that if $a_{1}$,
$\la_{1}$, and $\la_{2}$ are all simultaneously 0, then the algebra
is nilpotent.

We start by conjugating by $\A_{2}$, which will allow us to move the
$a_{2}$ position to 0 if $a_{1} \neq 2\la_{1}+\la_{2}$.  We have two
cases then.
\begin{enumerate}
    \item $a_{1} \neq 2 \la_{1}+\la_{2}$ and we move $a_{2}$ to
    0, or $a_{1}=2\la_{1}+\la_{2}$ and $a_{2}=0$ already.
    \item $a_{1}=2\la_{1}+\la_{2}$, but $a_{2} \neq 0$.
\end{enumerate}

\subsubsection{Subcase 1.1:}

In this section, the $a_{2}$ position is 0, which makes the
$\ad{\f_{7}}$ matrix $$\ad{\f_{7}}=\begin{pmatrix} a_{1}&0&0&0&x_{1}&y_{1} \\
0&2\la_{1}+\la_{2}&0&-y_{4}&x_{2}&y_{2} \\
a_{3}&0&\la_{1}+2\la_{2}&x_{4}&x_{3}&y_{3} \\
0&0&0&\la_{1}+\la_{2}&x_{4} &y_{4} \\ 0&0&0&0&\la_{1}&0 \\
0&0&0&0&0&\la_{2} \end{pmatrix}.$$  Next we conjugate by $\A_{3}$
and find that we can move $a_{3}$ to 0 if $a_{1}\neq
\la_{1}+2\la_{2}$. This gives us another two cases.
\begin{enumerate}
    \item $a_{1} \neq \la_{1}+2\la_{2}$ and we move $a_{3}$ to
    0, or $a_{1}=\la_{1}+2\la_{2}$ and $a_{3}=0$ already.
    \item $a_{1}=\la_{1}+2\la_{2}$, but $a_{3} \neq 0$.
\end{enumerate}

\subsubsection{Subcase 1.1.1:}

Here we have $$\ad{\f_{7}}=\begin{pmatrix} a_{1}&0&0&0&x_{1}&y_{1} \\
0&2\la_{1}+\la_{2}&0&-y_{4}&x_{2}&y_{2} \\
0&0&\la_{1}+2\la_{2}&x_{4}&x_{3}&y_{3} \\
0&0&0&\la_{1}+\la_{2}&x_{4} &y_{4} \\ 0&0&0&0&\la_{1}&0 \\
0&0&0&0&0&\la_{2} \end{pmatrix}.$$ We conjugate by $\A_{10}$ at this
point and find that if $a_{1} \neq \la_{1}$, then we can move the
$x_{1}$ position to 0.  This yields two cases.
\begin{enumerate}
    \item $a_{1} \neq \la_{1}$ and we move $x_{1}$ to 0, or
    $a_{1}=\la_{1}$ and $x_{1}=0$ already.
    \item $a_{1} = \la_{1}$, but $x_{1} \neq 0$.
\end{enumerate}

\subsubsection{Subcase 1.1.1.1:}

In this subcase, we have the $\ad{\f_{7}}$ matrix as $$\ad{\f_{7}}=
\begin{pmatrix} a_{1}&0&0&0&0&y_{1} \\
0&2\la_{1}+\la_{2}&0&-y_{4}&x_{2}&y_{2} \\
0&0&\la_{1}+2\la_{2}&x_{4}&x_{3}&y_{3} \\
0&0&0&\la_{1}+\la_{2}&x_{4} &y_{4} \\ 0&0&0&0&\la_{1}&0 \\
0&0&0&0&0&\la_{2} \end{pmatrix}.$$  We next conjugate by $\A_{13}$.
This allows us to move the $y_{1}$ position to 0, if $a_{1} \neq
\la_{2}$ and hence yields two cases.
\begin{enumerate}
    \item $a_{1} \neq \la_{2}$ and we move $y_{1}$ to 0, or
    $a_{1} = \la_{2}$ and $y_{1}=0$ already.
    \item $a_{1} = \la_{1}$, but $y_{1} \neq 0$.
\end{enumerate}

\subsubsection{Subcase 1.1.1.1.1:}

Here we have the $\ad{\f_{7}}$ matrix $$\ad{\f_{7}}=
\begin{pmatrix} a_{1}&0&0&0&0&0 \\
0&2\la_{1}+\la_{2}&0&-y_{4}&x_{2}&y_{2} \\
0&0&\la_{1}+2\la_{2}&x_{4}&x_{3}&y_{3} \\
0&0&0&\la_{1}+\la_{2}&x_{4} &y_{4} \\ 0&0&0&0&\la_{1}&0 \\
0&0&0&0&0&\la_{2} \end{pmatrix}.$$  We conjugate by $\A_{14}$ and
see that if $\la_{1} \neq 0$, we can move the $y_{2}$ position to 0.
We again have two cases.
\begin{enumerate}
    \item $\la_{1} \neq 0$ and we move $y_{2}$ to 0, or
    $\la_{1}=0$ and $y_{2}=0$ already.
    \item $\la_{1}=0$, but $y_{2} \neq 0$.
\end{enumerate}

\subsubsection{Subcase 1.1.1.1.1.1:}

We have the $\ad{\f_{7}}$ matrix $$\ad{\f_{7}}=
\begin{pmatrix} a_{1}&0&0&0&0&0 \\
0&2\la_{1}+\la_{2}&0&-y_{4}&x_{2}&0 \\
0&0&\la_{1}+2\la_{2}&x_{4}&x_{3}&y_{3} \\
0&0&0&\la_{1}+\la_{2}&x_{4} &y_{4} \\ 0&0&0&0&\la_{1}&0 \\
0&0&0&0&0&\la_{2} \end{pmatrix}.$$  In this section, we conjugate by
$\A_{12}$ and find that we can move the $x_{3}$ position to 0, if
$\la_{2} \neq 0$.  This yields yet another two cases.
\begin{enumerate}
    \item $\la_{2} \neq 0$ and we move $x_{3}$ to 0, or
    $\la_{2}=0$ and $x_{3}=0$ already.
    \item $\la_{2} =0$, but $x_{3} \neq 0$.
\end{enumerate}

\subsubsection{Subcase 1.1.1.1.1.1.1:}

In this section, we have the $\ad{\f_{7}}$ matrix $$\ad{\f_{7}}=
\begin{pmatrix} a_{1}&0&0&0&0&0 \\
0&2\la_{1}+\la_{2}&0&-y_{4}&x_{2}&0 \\
0&0&\la_{1}+2\la_{2}&x_{4}&0&y_{3} \\
0&0&0&\la_{1}+\la_{2}&x_{4} &y_{4} \\ 0&0&0&0&\la_{1}&0 \\
0&0&0&0&0&\la_{2} \end{pmatrix}.$$  We conjugate here by $\A_{11}$.
This allows us to move the $x_{2}$ value to the value of $y_{3}$, if
$\la_{2} \neq -\la_{1}$ and thus yields two cases.
\begin{enumerate}
    \item $\la_{2} \neq -\la_{1}$ and we move the $x_{2}$ position
    to the $y_{3}$ value, or $\la_{2} = -\la_{1}$ and $x_{2}=y_{3}$ already.
    \item $\la_{2} = -\la_{1}$, but $x_{2} \neq y_{3}$.
\end{enumerate}

\subsubsection{Subcase 1.1.1.1.1.1.1.1:}

Here the $\ad{\f_{7}}$ matrix is of the form $$\ad{\f_{7}}=
\begin{pmatrix} a_{1}&0&0&0&0&0 \\
0&2\la_{1}+\la_{2}&0&-y_{4}&y_{3}&0 \\
0&0&\la_{1}+2\la_{2}&x_{4}&0&y_{3} \\
0&0&0&\la_{1}+\la_{2}&x_{4} &y_{4} \\ 0&0&0&0&\la_{1}&0 \\
0&0&0&0&0&\la_{2} \end{pmatrix}.$$  At this point, we make the basis
change $$\f_{i}' = \f_{i} \ \ \text{for} \ \ 1 \leq i \leq 6, \ \
\text{and} \ \ \f_{7}'=\f_{7}-y_{3}\f_{4}-y_{4}\f_{5}+x_{4}\f_{6}.$$
This will yield that $$\ad{\f_{7}'}=
\begin{pmatrix} a_{1}&0&0&0&0&0 \\
0&2\la_{1}+\la_{2}&0&0&0&0 \\
0&0&\la_{1}+2\la_{2}&0&0&0 \\
0&0&0&\la_{1}+\la_{2}&0&0 \\ 0&0&0&0&\la_{1}&0 \\
0&0&0&0&0&\la_{2} \end{pmatrix}.$$  We see here that if $a_{1}=0$,
or if $\la_{1}=\la_{2}=0$, then the algebra decomposes.  In
addition, we can conjugate by $\A_{0}$ with $a=d=0$ and $b=c=1$ to
swap the order of $\la_{1}$ and $\la_{2}$.  We order these two so
that $\frac{\la_{1}}{a_{1}} \geq \frac{\la_{2}}{a_{1}}$.  Finally,
as $a_{1} \neq 0$, we make the change of basis $$\e_{i}=\f_{i}' \ \
\text{for} \ \ 1 \leq i \leq 6, \ \
\e_{7}=-\frac{1}{a_{1}}\f_{7}',$$ and let $a =
\frac{\la_{1}}{a_{1}}$ and $b=\frac{\la_{2}}{a_{1}}$.  This yields
the structure equations
\begin{align}
    \notag [\e_{4},\e_{5}]&=\e_{2}, & [\e_{4},\e_{6}]&= \e_{3},
    &[\e_{5},\e_{6}]&=\e_{4}, & [\e_{1},\e_{7}]&=\e_{1}, &
    [\e_{2},\e_{7}]&= (2a+b)\e_{2}, \\
    \notag [\e_{3},\e_{7}]&=(a+2b)\e_{3}, &
    [\e_{4},\e_{7}]&=(a+b)\e_{4}, & [\e_{5},\e_{7}]&=a\e_{5}, &
    [\e_{6},\e_{7}]&=b\e_{6},
\end{align}
with $a \geq b$ and $a^{2} + b^{2} \neq 0$.  In the table in
Appendix \ref{a2}, this is [7,[6,4],1,1].

\subsubsection{Subcase 1.1.1.1.1.1.1.2:}

In this section, we assume that $x_{2} \neq y_{3}$ and
$\la_{2}=-\la_{1}$, yielding the $\ad{\f_{7}}$ matrix
$$\ad{\f_{7}}=
\begin{pmatrix} a_{1}&0&0&0&0&0 \\
0&\la_{1}&0&-y_{4}&x_{2}&0 \\
0&0&-\la_{1}&x_{4}&0&y_{3} \\
0&0&0&0&x_{4} &y_{4} \\ 0&0&0&0&\la_{1}&0 \\
0&0&0&0&0&-\la_{1} \end{pmatrix}.$$  We make the change of basis
$$\f_{i}' = \f_{i} \ \ \text{for} \ \ 1 \leq i \leq 6, \ \
\text{and} \ \ \f_{7}'=\f_{7}-y_{3}\f_{4}-y_{4}\f_{5}+x_{4}\f_{6}.$$
This will yield the $\ad{\f_{7}'}$ matrix $$\ad{\f_{7}'}=
\begin{pmatrix} a_{1}&0&0&0&0&0 \\
0&\la_{1}&0&&x_{2}-y_{3}&0 \\
0&0&-\la_{1}&0&0&0 \\
0&0&0&0&0 &0\\ 0&0&0&0&\la_{1}&0 \\
0&0&0&0&0&-\la_{1} \end{pmatrix}.$$  Now note that if $a_{1}=0$,
then the algebra decomposes.  Hence $a_{1} \neq 0$.  Then we
conjugate by $\A_{4}$ and as $x_{2} \neq y_{3}$ and $a_{1} \neq 0$,
we can scale $x_{2}-y_{3}$ to $a_{1}$.  Finally we make the change
of basis $$\e_{i}=\f_{i}' \ \ \text{for} \ \ 1 \leq i \leq 6, \ \
\e_{7}=-\frac{1}{a_{1}}\f_{7}',$$ and let $a =
\frac{\la_{1}}{a_{1}}$.  This will yield the structure equations
\begin{align}
    \notag [\e_{4},\e_{5}]&=\e_{2}, & [\e_{4},\e_{6}]&= \e_{3},
    &[\e_{5},\e_{6}]&=\e_{4}, & [\e_{1},\e_{7}]&=\e_{1}, \\
    \notag [\e_{2},\e_{7}]&= a \e_{2}, & [\e_{3},\e_{7}]&=-a\e_{3},
    & [\e_{5},\e_{7}]&=\e_{2}+a\e_{5}, & [\e_{6},\e_{7}]&=-a\e_{6},
\end{align}
with $a \in \R$.  This is [7,[6,4],1,2].

\subsubsection{Subcase 1.1.1.1.1.1.2:}

Here, we assumed that $\la_{2} = 0$, but $x_{3} \neq 0$ giving us
the $\ad{\f_{7}}$ matrix $$\ad{\f_{7}}=
\begin{pmatrix} a_{1}&0&0&0&0&0 \\
0&2\la_{1}&0&-y_{4}&x_{2}&0 \\
0&0&\la_{1}&x_{4}&x_{3}&y_{3} \\
0&0&0&\la_{1}&x_{4} &y_{4} \\ 0&0&0&0&\la_{1}&0 \\
0&0&0&0&0& \end{pmatrix}.$$  We conjugate by $\A_{11}$ and find that
if $\la_{1} \neq 0$, then we can move the $x_{2}$ position to equal
$y_{3}$.  This yields two possibilities.
\begin{enumerate}
    \item $\la_{1} \neq 0$ and we make the $x_{2}$ position
    equal $y_{3}$, or $\la_{1}=0$ and $x_{2}=y_{3}$ already.
    \item $\la_{1} = 0$, but $x_{2} \neq y_{3}$.
\end{enumerate}

\subsubsection{Subcase 1.1.1.1.1.1.2.1:}

We have in this case, that $$\ad{\f_{7}}=
\begin{pmatrix} a_{1}&0&0&0&0&0 \\
0&2\la_{1}&0&-y_{4}&y_{3}&0 \\
0&0&\la_{1}&x_{4}&x_{3}&y_{3} \\
0&0&0&\la_{1}&x_{4} &y_{4} \\ 0&0&0&0&\la_{1}&0 \\
0&0&0&0&0& \end{pmatrix}.$$  We note now that $a_{1} \neq 0$ or the
algebra is decomposable.  This allows us to conjugate by $\A_{4}$
and, as $x_{3}$ is nonzero as well, scale $x_{3}$ to $\pm a_{1}$.
Then we make the basis change $$\e_{i}=\f_{i}, \ \ \text{for} \ \ 1
\leq i \leq 6, \ \ \text{and} \ \ \e_{7}=-\frac{1}{a_{1}}(\f_{7}
-y_{3}\f_{4}-y_{4}\f_{5}+x_{4}\f_{6}),$$ and let $a =
\frac{\la_{1}}{a_{1}}$.  This will yield the structure equations
\begin{align}
    \notag [\e_{4},\e_{5}]&=\e_{2}, & [\e_{4},\e_{6}]&= \e_{3},
    &[\e_{5},\e_{6}]&=\e_{4}, & [\e_{1},\e_{7}]&=\e_{1}, \\
    \notag [\e_{2},\e_{7}]&= 2a \e_{2}, & [\e_{3},\e_{7}]&= a\e_{3},
    & [\e_{4},\e_{7}]&=a\e_{4}, & [\e_{5},\e_{7}]&=\vep\e_{3}+
    a\e_{5},
\end{align}
with $a \in \R$ and $\vep^{2} =1$.  This is [7,[6,4],1,3] in the
table.

\subsubsection{Subcase 1.1.1.1.1.1.2.2:}

We assumed here that $\la_{1}=0$ and $x_{2} \neq y_{3}$.  This
yields the $\ad{\f_{7}}$ matrix $$\ad{\f_{7}}=
\begin{pmatrix} a_{1}&0&0&0&0&0 \\
0&0&0&-y_{4}&x_{2}&0 \\
0&0&0&x_{4}&x_{3}&y_{3} \\
0&0&0&0&x_{4} &y_{4} \\ 0&0&0&0&0&0 \\
0&0&0&0&0& \end{pmatrix}.$$  This implies that $a_{1} \neq 0$ or the
algebra is nilpotent.

We make the change of basis $$\f_{i}' = \f_{i} \ \ \text{for} \ \ 1
\leq i \leq 6, \ \ \text{and} \ \
\f_{7}'=\f_{7}-y_{3}\f_{4}-y_{4}\f_{5}+x_{4}\f_{6}.$$  This will
yield the $\ad{\f_{7}}$ matrix $$\ad{\f_{7}}=
\begin{pmatrix} a_{1}&0&0&0&0&0 \\
0&0&0&0&x_{2}-y_{3}&0 \\
0&0&0&0&x_{3}&0 \\
0&0&0&0&0&0 \\ 0&0&0&0&0&0 \\
0&0&0&0&0& \end{pmatrix}.$$  As $a_{1}$, $x_{3}$, and $x_{2}-y_{3}$
are all nonzero values, we can conjugate by $\A_{4}\A_{7}$ and pick
$s_{4}$ and $s_{7}$ to simultaneously scale $x_{2}-y_{3}$ to $a_{1}$
and $x_{3}$ to $\pm a_{1}$.  Finally, we make the change of basis
$$\e_{i}=\f_{i}' \ \ \text{for} \ \ 1 \leq i \leq 6, \ \
\e_{7}=-\frac{1}{a_{1}}\f_{7}'.$$  The resulting structure equations
are
\begin{align}
    \notag [\e_{4},\e_{5}]&=\e_{2}, & [\e_{4},\e_{6}]&= \e_{3},
    &[\e_{5},\e_{6}]&=\e_{4}, & [\e_{1},\e_{7}]&=\e_{1}, &
    [\e_{5},\e_{7}]&=\e_{2}+\vep \e_{3},
\end{align}
with $\vep^{2} = 1$.  This is [7,[6,4],1,4].

\subsubsection{Subcase 1.1.1.1.1.2:}

In this section, we assumed that $\la_{1}=0$ and $y_{2} \neq 0$.
This yields the $\ad{\f_{7}}$ matrix $$\ad{\f_{7}}=
\begin{pmatrix} a_{1}&0&0&0&0&0 \\
0&\la_{2}&0&-y_{4}&x_{2}&y_{2} \\
0&0&2\la_{2}&x_{4}&x_{3}&y_{3} \\
0&0&0&\la_{2}&x_{4} &y_{4} \\ 0&0&0&0&0&0 \\
0&0&0&0&0&\la_{2} \end{pmatrix}.$$  Next we conjugate by
$\A_{11}\A_{12}$ and find that if $\la_{2} \neq 0$, then we can pick
$s_{11}$ and $s_{12}$ to make the $x_{2}$ position equal $y_{3}$ and
$x_{3}$ equal 0.  Also, if $\la_{2}=0$, then we can conjugate by
$\A_{5}$, without affecting any previous changes, and as $y_{2} \neq
0$, we can still make the $x_{2}$ position equal the $y_{3}$
position. We have two cases then.
\begin{enumerate}
    \item $\la_{2} \neq 0$ and we make $x_{2}$ equal to $y_{3}$ and
    $x_{3}$ equal to 0, or $\la_{2}=0$ and we still make the $x_{2}$ position
    equal to the $y_{3}$ position and the $x_{3}$ position is zero afterwards.
    \item $\la_{2} = 0$ and we still make the $x_{2}$ position
    equal the $y_{3}$ position, but the $x_{3}$ position is nonzero
    afterwards.
\end{enumerate}

\subsubsection{Subcase 1.1.1.1.1.2.1:}

Here we have the $\ad{\f_{7}}$ matrix $$\ad{\f_{7}}=
\begin{pmatrix} a_{1}&0&0&0&0&0 \\
0&\la_{2}&0&-y_{4}&y_{3}&y_{2} \\
0&0&2\la_{2}&x_{4}&0&y_{3} \\
0&0&0&\la_{2}&x_{4} &y_{4} \\ 0&0&0&0&0&0 \\
0&0&0&0&0&\la_{2} \end{pmatrix}.$$  We conjugate by $\A_{0}$ with
$a=d=0$ and $b=c=1$, and we'll end up with the matrix $$\ad{\f_{7}}=
\begin{pmatrix} a_{1}&0&0&0&0&0 \\
0&2\la_{2}&0&x_{4}&-y_{3}&0 \\
0&0&\la_{2}&-y_{4}&-y_{2}&-y_{3} \\
0&0&0&\la_{2}&-y_{4} &-x_{4} \\ 0&0&0&0&0&\la_{2} \\
0&0&0&0&0&0 \end{pmatrix}.$$  If we relabel, $\la_{2}$, $-y_{3}$,
$-x_{4}$, $-y_{4}$, and $-y_{2}$ as $\la_{1}$, $y_{3}$, $y_{4}$,
$x_{4}$, and $x_{3}$ respectively, then this is the same as the
$\ad{\f_{7}}$ matrix in Subcase 1.1.1.1.1.1.2.1.  Thus we don't
obtain anything new from this subcase.

\subsubsection{Subcase 1.1.1.1.1.2.2:}

Here we have that $\la_{2}=0$ and $x_{3} \neq 0$, which yields the
$\ad{\f_{7}}$ matrix $$\ad{\f_{7}}=
\begin{pmatrix} a_{1}&0&0&0&0&0 \\
0&0&0&-y_{4}&x_{2}&y_{2} \\
0&0&0&x_{4}&x_{3}&y_{3} \\
0&0&0&0&x_{4} &y_{4} \\ 0&0&0&0&0&0 \\
0&0&0&0&0&0 \end{pmatrix}.$$  This implies that $a_{1} \neq 0$ or
the algebra is nilpotent.  We then conjugate by $\A_{4}\A_{7}$ and,
as $x_{3}$ and $y_{2}$ are both nonzero as well, scale both $x_{3}$
and $y_{2}$ to $\pm a_{1}$.  Finally, we make the change of basis
$$\e_{i}=\f_{i}, \ \ \text{for} \ \ 1 \leq i \leq 6, \ \ \text{and}
\ \ \e_{7}=-\frac{1}{a_{1}}(\f_{7}
-y_{3}\f_{4}-y_{4}\f_{5}+x_{4}\f_{6}).$$  Then the structure
equations are
\begin{align}
    \notag [\e_{4},\e_{5}]&=\e_{2}, & [\e_{4},\e_{6}]&= \e_{3},
    &[\e_{5},\e_{6}]&=\e_{4}, & [\e_{1},\e_{7}]&=\e_{1}, \\
    \notag [\e_{5},\e_{7}]&=\del \e_{3}, & [\e_{6},\e_{7}]&=\vep \e_{2},
\end{align}
with $\vep^{2} = \del^{2} = 1$.  In Appendix \ref{a2}, this is
[7,[6,4],1,5].

\subsubsection{Subcase 1.1.1.1.2:}

In this section, we have the $\ad{\f_{7}}$ matrix $$\ad{\f_{7}}=
\begin{pmatrix} \la_{2}&0&0&0&0&y_{1} \\
0&2\la_{1}+\la_{2}&0&-y_{4}&x_{2}&y_{2} \\
0&0&\la_{1}+2\la_{2}&x_{4}&x_{3}&y_{3} \\
0&0&0&\la_{1}+\la_{2}&x_{4} &y_{4} \\ 0&0&0&0&\la_{1}&0 \\
0&0&0&0&0&\la_{2} \end{pmatrix},$$ with $y_{1} \neq 0$.  Next, we
conjugate by $\A_{14}$ and find that if $\la_{1} \neq 0$, we can
move the $y_{2}$ position to 0.  However, if $\la_{1}=0$, then we
can conjugate by $\A_{2}$ without affecting any previous changes
and, as $y_{1} \neq 0$, we can still make the $y_{2}$ position equal
0.

We can do something similar for the $y_{3}$ position.  If we
conjugate by $\A_{15}$, then we can change the $y_{3}$ position to
equal the $x_{2}$ position whenever $\la_{2} \neq -\la_{1}$.
However, if $\la_{2} = -\la_{1}$, then we can apply $\A_{3}$ without
affecting any previous changes and, as $y_{1} \neq 0$, we can still
change the $y_{3}$ position to equal $x_{2}$.

Now we conjugate by $\A_{12}$ and we see that we can only make the
$x_{3}$ position equal to 0, if $\la_{2} \neq 0$.  This gives us two
cases.
\begin{enumerate}
    \item $\la_{2} \neq 0$ and we make the $x_{3}$ position
    equal 0, or $\la_{2}=0$ and $x_{3}=0$ already.
    \item $\la_{2}=0$, but $x_{3} \neq 0$.
\end{enumerate}

\subsubsection{Subcase 1.1.1.1.2.1:}

Here we have the $\ad{\f_{7}}$ matrix $$\ad{\f_{7}}=
\begin{pmatrix} \la_{2}&0&0&0&0&y_{1} \\
0&2\la_{1}+\la_{2}&0&-y_{4}&x_{2}&0 \\
0&0&\la_{1}+2\la_{2}&x_{4}&0&x_{2} \\
0&0&0&\la_{1}+\la_{2}&x_{4} &y_{4} \\ 0&0&0&0&\la_{1}&0 \\
0&0&0&0&0&\la_{2} \end{pmatrix}.$$  We have that $\la_{1}$ and
$\la_{2}$ are not simultaneously 0.  We will bifurcate on this here
and consider both of the following cases individually.
\begin{enumerate}
    \item $\la_{1} \neq 0$.
    \item $\la_{1}=0$, which implies that $\la_{2} \neq 0$ or the
    algebra is nilpotent.
\end{enumerate}

\subsubsection{Subcase 1.1.1.1.2.1.1:}

Here, we conjugate by $\A_{1}$ and, as $y_{1}$ and $\la_{1}$ are
both nonzero, scale $y_{1}$ to $\la_{1}$.  Then we make the change
of basis $$\e_{i}=\f_{i}, \ \ \text{for} \ \ 1 \leq i \leq 6, \ \
\text{and} \ \ \e_{7}=-\frac{1}{\la_{1}}(\f_{7}-x_{2}\f_{4}-
y_{4}\f_{5}+x_{4}\f_{6}),$$ and let $a = \frac{\la_{2}}{\la_{1}}$.
We end up with the structure equations
\begin{align}
    \notag [\e_{4},\e_{5}]&=\e_{2}, & [\e_{4},\e_{6}]&= \e_{3},
    &[\e_{5},\e_{6}]&=\e_{4}, & [\e_{1},\e_{7}]&=a\e_{1}, &
    [\e_{2},\e_{7}]&=(a+2)\e_{2}, \\
    \notag [\e_{3},\e_{7}]&=(2a+1)\e_{3}, &
    [\e_{4},\e_{7}]&=(a+1)\e_{4}, & [\e_{5},\e_{7}]&=\e_{5}, &
    [\e_{6},\e_{7}]&=\e_{1}+a\e_{6},
\end{align}
with $a \in \R$.  This is [7,[6,4],1,6].

\subsubsection{Subcase 1.1.1.1.2.1.2:}

In this section, we assumed that $\la_{1}=0$, which yields the
$\ad{\f_{7}}$ matrix $$\ad{\f_{7}}=
\begin{pmatrix} \la_{2}&0&0&0&0&y_{1} \\
0&\la_{2}&0&-y_{4}&x_{2}&0 \\
0&0&2\la_{2}&x_{4}&0&x_{2} \\
0&0&0&\la_{2}&x_{4} &y_{4} \\ 0&0&0&0&0&0 \\
0&0&0&0&0&\la_{2} \end{pmatrix}.$$  Conjugate by $\A_{1}$, and this
time, as $y_{1}$ and $\la_{2}$ are both nonzero, scale $y_{1}$ to
$\la_{2}$.  Then make the change of basis $$\e_{i}=\f_{i}, \ \
\text{for} \ \ 1 \leq i \leq 6, \ \ \text{and} \ \
\e_{7}=-\frac{1}{\la_{2}}(\f_{7}-x_{2}\f_{4}-
y_{4}\f_{5}+x_{4}\f_{6}),$$ which will yield the structure equations
\begin{align}
    \notag [\e_{4},\e_{5}]&=\e_{2}, & [\e_{4},\e_{6}]&= \e_{3},
    &[\e_{5},\e_{6}]&=\e_{4}, & [\e_{1},\e_{7}]&=\e_{1}, \\
    \notag [\e_{2},\e_{7}]&=\e_{2}, &[\e_{3},\e_{7}]&=2\e_{3}, &
    [\e_{4},\e_{7}]&=\e_{4}, & [\e_{6},\e_{7}]&=\e_{1}+\e_{6}.
\end{align}
This is [7,[6,4],1,7].

\subsubsection{Subcase 1.1.1.1.2.2:}

In this section, we assume that $\la_{2} =0$ and $x_{3} \neq 0$.
This yields that the $\ad{\f_{7}}$ matrix is $$\ad{\f_{7}}=
\begin{pmatrix} 0&0&0&0&0&y_{1} \\
0&2\la_{1}&0&-y_{4}&x_{2}&0 \\
0&0&\la_{1}&x_{4}&x_{3}&x_{2} \\
0&0&0&\la_{1}&x_{4} &y_{4} \\ 0&0&0&0&\la_{1}&0 \\
0&0&0&0&0&0 \end{pmatrix}.$$  We next conjugate by $\A_{1}\A_{4}$
and, as $y_{1}$, $x_{3}$, and $\la_{1}$ are all nonzero, scale
$y_{1}$ to $\la_{1}$ and $x_{3}$ to $\pm \la_{1}$.  Finally, we make
the basis change $$\e_{i}=\f_{i}, \ \ \text{for} \ \ 1 \leq i \leq
6, \ \ \text{and} \ \ \e_{7}=-\frac{1}{\la_{1}}(\f_{7}-x_{2}\f_{4}-
y_{4}\f_{5}+x_{4}\f_{6}),$$ and the structure equations become
\begin{align}
    \notag [\e_{4},\e_{5}]&=\e_{2}, & [\e_{4},\e_{6}]&= \e_{3},
    &[\e_{5},\e_{6}]&=\e_{4}, & [\e_{2},\e_{7}]&=2\e_{2}, \\
    \notag [\e_{3},\e_{7}]&=\e_{3}, & [\e_{4},\e_{7}]&=\e_{4}, &
    [\e_{5},\e_{7}]&=\vep \e_{3}+\e_{5}, &[\e_{6},\e_{7}]&= \e_{1},
\end{align}
with $\vep^{2} = 1$.  In the table, this is [7,[6,4],1,8].

\subsubsection{Subcase 1.1.1.2:}

We assumed in this section that $a_{1}=\la_{1}$ but $x_{1} \neq 0$.
This gives us the $\ad{\f_{7}}$ matrix $$\ad{\f_{7}}=\begin{pmatrix} \la_{1}&0&0&0&x_{1}&y_{1} \\
0&2\la_{1}+\la_{2}&0&-y_{4}&x_{2}&y_{2} \\
0&0&\la_{1}+2\la_{2}&x_{4}&x_{3}&y_{3} \\
0&0&0&\la_{1}+\la_{2}&x_{4} &y_{4} \\ 0&0&0&0&\la_{1}&0 \\
0&0&0&0&0&\la_{2} \end{pmatrix}.$$  Here we conjugate by $\A_{13}$
and, if $\la_{2} \neq \la_{1}$, move $y_{1}$ to 0.  However, if
$\la_{2}=\la_{1}$, then we conjugate by $\A_{6}$ instead and, as
$x_{1} \neq 0$, still move $y_{1}$ to 0.  At this point, if we
conjugate by $\A_{0}$ with $a=d=0$ and $b=c=1$, we'll get that
$$\ad{\f_{7}}=\begin{pmatrix} \la_{1}&0&0&0&0&x_{1} \\
0&2\la_{2}+\la_{1}&0&x_{4}&-y_{3}&-x_{3} \\
0&0&\la_{2}+2\la_{1}&-y_{4}&-y_{2}&-x_{2} \\
0&0&0&\la_{2}+\la_{1}&-y_{4} &-x_{4} \\ 0&0&0&0&\la_{2}&0 \\
0&0&0&0&0&\la_{1} \end{pmatrix}.$$  If we relabel $\la_{1}$,
$\la_{2}$, $x_{1}$, $-x_{4}$, $-y_{3}$, $-y_{4}$, $-y_{2}$, and
$-x_{2}$ as $\la_{2}$, $\la_{1}$, $y_{1}$, $y_{4}$, $x_{2}$,
$y_{2}$, $x_{4}$, $x_{3}$, and $y_{3}$ respectively, we'll have the
same $\ad{\f_{7}}$ matrix as in Subcase 1.1.1.1.2.  Thus there are
no new isomorphism classes to be found here.

\subsubsection{Subcase 1.1.2:}

In this section, we assumed that $a_{1}=\la_{1}+2\la_{2}$ and that
$a_{3} \neq 0$.  This will give us the $\ad{\f_{7}}$ matrix
$$\ad{\f_{7}}=\begin{pmatrix} \la_{1}+2\la_{2}&0&0&0&x_{1}&y_{1} \\
0&2\la_{1}+\la_{2}&0&-y_{4}&x_{2}&y_{2} \\
a_{3}&0&\la_{1}+2\la_{2}&x_{4}&x_{3}&y_{3} \\
0&0&0&\la_{1}+\la_{2}&x_{4} &y_{4} \\ 0&0&0&0&\la_{1}&0 \\
0&0&0&0&0&\la_{2} \end{pmatrix}.$$  We conjugate next by $\A_{10}$.
If $\la_{2} \neq 0$, this allows us to move the $x_{1}$ position to
0 and hence gives us two cases.
\begin{enumerate}
    \item $\la_{2} \neq 0$ and we move $x_{1}$ to 0, or $\la_{2}=0$
    and $x_{1}=0$ already.
    \item $\la_{2}=0$, but $x_{1} \neq 0$.
\end{enumerate}

\subsubsection{Subcase 1.1.2.1:}

Here we have the $\ad{\f_{7}}$ matrix $$\ad{\f_{7}}=
\begin{pmatrix} \la_{1}+2\la_{2}&0&0&0&0&y_{1} \\
0&2\la_{1}+\la_{2}&0&-y_{4}&x_{2}&y_{2} \\
a_{3}&0&\la_{1}+2\la_{2}&x_{4}&x_{3}&y_{3} \\
0&0&0&\la_{1}+\la_{2}&x_{4} &y_{4} \\ 0&0&0&0&\la_{1}&0 \\
0&0&0&0&0&\la_{2} \end{pmatrix}.$$  We next conjugate by $\A_{13}$.
If $\la_{2} \neq -\la_{1}$, we move $y_{1}$ to $0$.  This yields
another two cases
\begin{enumerate}
    \item $\la_{2} \neq -\la_{1}$ and we make the $y_{1}$ position
    equal to 0, or $\la_{2} = -\la_{1}$ and $y_{1}=0$ already.
    \item $\la_{2} = -\la_{1}$, but $y_{1} \neq 0$.
\end{enumerate}

\subsubsection{Subcase 1.1.2.1.1:}

In this section, we have $$\ad{\f_{7}}=
\begin{pmatrix} \la_{1}+2\la_{2}&0&0&0&0&0 \\
0&2\la_{1}+\la_{2}&0&-y_{4}&x_{2}&y_{2} \\
a_{3}&0&\la_{1}+2\la_{2}&x_{4}&x_{3}&y_{3} \\
0&0&0&\la_{1}+\la_{2}&x_{4} &y_{4} \\ 0&0&0&0&\la_{1}&0 \\
0&0&0&0&0&\la_{2} \end{pmatrix}.$$  We conjugate by $\A_{12}$ and,
if $\la_{2} \neq 0$, move $x_{3}$ to 0.  If $\la_{2}=0$, then we
conjugate by $\A_{10}$ instead and, as $a_{3} \neq 0$, we still move
$x_{3}$ to 0.

Similarly, we next conjugate by $\A_{15}$ and, if $\la_{2} \neq
-\la_{1}$, then we move $y_{3}$ to $x_{2}$.  If $\la_{2} =
-\la_{1}$, then we conjugate by $\A_{13}$ instead and, as $a_{3}
\neq 0$, we still move $y_{3}$ to $x_{2}$.

Then we conjugate by $\A_{14}$.  If $\la_{1} \neq 0$, then we can
move $y_{2}$ to 0.  This yields two cases.
\begin{enumerate}
    \item $\la_{1} \neq 0$ and we move $y_{2}$ to 0, or $\la_{1}=0$
    and $y_{2}=0$ already.
    \item $\la_{1}=0$, but $y_{2} \neq 0$.
\end{enumerate}

\subsubsection{Subcase 1.1.2.1.1.1:}

We have the $\ad{\f_{7}}$ matrix $$\ad{\f_{7}}=
\begin{pmatrix} \la_{1}+2\la_{2}&0&0&0&0&0 \\
0&2\la_{1}+\la_{2}&0&-y_{4}&x_{2}&0 \\
a_{3}&0&\la_{1}+2\la_{2}&x_{4}&0&x_{2} \\
0&0&0&\la_{1}+\la_{2}&x_{4} &y_{4} \\ 0&0&0&0&\la_{1}&0 \\
0&0&0&0&0&\la_{2} \end{pmatrix}.$$  We know that $\la_{1}$ and
$\la_{2}$ are not simultaneously 0, otherwise the algebra is
nilpotent.  This gives us another two cases.
\begin{enumerate}
    \item $\la_{1} \neq 0$.
    \item $\la_{1}=0$, which implies that $\la_{2} \neq 0$.
\end{enumerate}

\subsubsection{Subcase 1.1.2.1.1.1.1:}

Here we conjugate by $\A_{1}$ and, as $a_{3}$ and $\la_{1}$ are both
nonzero, we scale $a_{3}$ to $\la_{1}$.  Then we make the change of
basis $$\e_{i}=\f_{i}, \ \ \text{for} \ \ 1 \leq i \leq 6, \ \
\text{and} \ \ \e_{7}=-\frac{1}{\la_{1}}(\f_{7}-x_{2}\f_{4}-
y_{4}\f_{5}+x_{4}\f_{6}),$$ and let $a = \frac{\la_{2}}{\la_{1}}$.
This yields the structure equations
\begin{align}
    \notag [\e_{4},\e_{5}]&=\e_{2}, & [\e_{4},\e_{6}]&= \e_{3},
    &[\e_{5},\e_{6}]&=\e_{4}, & [\e_{1},\e_{7}]&=(2a+1)\e_{1}+
    \e_{3},  \\
    \notag [\e_{2},\e_{7}]&=(a+2)\e_{2}, &
    [\e_{3},\e_{7}]&=(2a+1)\e_{3}, & [\e_{4},\e_{7}]&=(a+1)\e_{4}, &
    [\e_{5},\e_{7}]&= \e_{5}, \\
    \notag [\e_{6},\e_{7}]&=a\e_{6},
\end{align}
with $a \in \R$.  This is [7,[6,4],1,9] in the table.

\subsubsection{Subcase 1.1.2.1.1.1.2:}

In this section, we have that $$\ad{\f_{7}}=
\begin{pmatrix} 2\la_{2}&0&0&0&0&0 \\
0&\la_{2}&0&-y_{4}&x_{2}&0 \\
a_{3}&0&2\la_{2}&x_{4}&0&x_{2} \\
0&0&0&\la_{2}&x_{4} &y_{4} \\ 0&0&0&0&0&0 \\
0&0&0&0&0&\la_{2} \end{pmatrix}.$$  Then we conjugate by $\A_{1}$
and, as $a_{3}$ and $\la_{2}$ are both nonzero, we scale $a_{3}$ to
$\la_{2}$.  Finally, we make the change of basis
$$\e_{i}=\f_{i}, \ \ \text{for} \ \ 1 \leq i \leq 6, \ \ \text{and}
\ \ \e_{7}=-\frac{1}{\la_{2}}(\f_{7}-x_{2}\f_{4}-
y_{4}\f_{5}+x_{4}\f_{6}).$$ Then we obtain the structure equations
\begin{align}
    \notag [\e_{4},\e_{5}]&=\e_{2}, & [\e_{4},\e_{6}]&= \e_{3},
    &[\e_{5},\e_{6}]&=\e_{4}, & [\e_{1},\e_{7}]&=2\e_{1}+
    \e_{3},  \\
    \notag [\e_{2},\e_{7}]&=\e_{2}, &
    [\e_{3},\e_{7}]&=2\e_{3}, & [\e_{4},\e_{7}]&=\e_{4}, &
    [\e_{6},\e_{7}]&=\e_{6}.
\end{align}
This is [7,[6,4],1,10].

\subsubsection{Subcase 1.1.2.1.1.2:}

In this section, we assumed that $\la_{1} = 0$ and $y_{2} \neq 0$.
This yields that $$\ad{\f_{7}}=
\begin{pmatrix} 2\la_{2}&0&0&0&0&0 \\
0&\la_{2}&0&-y_{4}&x_{2}&y_{2} \\
a_{3}&0&2\la_{2}&x_{4}&0&x_{2} \\
0&0&0&\la_{2}&x_{4} &y_{4} \\ 0&0&0&0&0&0 \\
0&0&0&0&0&\la_{2} \end{pmatrix}.$$  Note that if $\la_{2}=0$, then
the algebra is nilpotent. Then we conjugate by $\A_{1}\A_{7}$ and,
as $a_{3}$, $y_{2}$, and $\la_{2}$ are all nonzero, scale $a_{3}$ to
$\la_{2}$ and $y_{2}$ to $\pm \la_{2}$. Finally, we make the change
of basis
$$\e_{i}=\f_{i}, \ \ \text{for} \ \ 1 \leq i \leq 6, \ \ \text{and}
\ \ \e_{7}=-\frac{1}{\la_{2}}(\f_{7}-x_{2}\f_{4}-
y_{4}\f_{5}+x_{4}\f_{6}).$$  This yields the structure equations
\begin{align}
    \notag [\e_{4},\e_{5}]&=\e_{2}, & [\e_{4},\e_{6}]&= \e_{3},
    &[\e_{5},\e_{6}]&=\e_{4}, & [\e_{1},\e_{7}]&=2\e_{1}+
    \e_{3},  \\
    \notag [\e_{2},\e_{7}]&=\e_{2}, &
    [\e_{3},\e_{7}]&=2\e_{3}, & [\e_{4},\e_{7}]&=\e_{4}, &
    [\e_{6},\e_{7}]&=\vep \e_{2} + \e_{6},
\end{align}
with $\vep^{2} = 1$.  In Appendix \ref{a2}, this is [7,[6,4],1,11].

\subsubsection{Subcase 1.1.2.1.2:}

Here we had that $\la_{2} = -\la_{1}$ and $y_{1} \neq 0$.  This
yields the $\ad{\f_{7}}$ matrix $$\ad{\f_{7}}=
\begin{pmatrix} -\la_{1}&0&0&0&0&y_{1} \\
0&\la_{1}&0&-y_{4}&x_{2}&y_{2} \\
a_{3}&0&-\la_{1}&x_{4}&x_{3}&y_{3} \\
0&0&0&0&x_{4} &y_{4} \\ 0&0&0&0&\la_{1}&0 \\
0&0&0&0&0&-\la_{1} \end{pmatrix}.$$  Also we see that $\la_{1} \neq
0$ or the algebra is nilpotent.  Then we conjugate by
$\A_{12}\A_{13}\A_{14}$ and, as $\la_{1}$ and $a_{3}$ are both
nonzero, move $y_{2}$ and $x_{3}$ to 0 and $y_{3}$ to $x_{2}$. Next
we conjugate by $\A_{1}\A_{4}$ and, as $a_{3}$, $y_{1}$, and
$\la_{1}$ are all nonzero, scale $a_{3}$ and $y_{1}$ both to
$-\la_{1}$.  Finally, we make the basis change $$\e_{i}=\f_{i}, \ \
\text{for} \ \ 1 \leq i \leq 6, \ \ \text{and} \ \
\e_{7}=\frac{1}{\la_{1}}(\f_{7}-x_{2}\f_{4}-
y_{4}\f_{5}+x_{4}\f_{6}),$$ and obtain the structure equations
\begin{align}
    \notag [\e_{4},\e_{5}]&=\e_{2}, & [\e_{4},\e_{6}]&= \e_{3},
    &[\e_{5},\e_{6}]&=\e_{4}, & [\e_{1},\e_{7}]&=\e_{1}+\e_{3}, \\
    \notag [\e_{2},\e_{7}]&=-\e_{2}, &[\e_{3},\e_{7}]&=\e_{3}, &
    [\e_{5},\e_{7}]&=-\e_{5}, & [\e_{6},\e_{7}]&=\e_{1}+\e_{6}.
\end{align}
This is [7,[6,4],1,12].

\subsubsection{Subcase 1.1.2.2:}

In this section, we have assumed that $\la_{2}=0$ and $x_{1} \neq
0$.  Then $$\ad{\f_{7}}=\begin{pmatrix} \la_{1}&0&0&0&x_{1}&y_{1} \\
0&2\la_{1}&0&-y_{4}&x_{2}&y_{2} \\
a_{3}&0&\la_{1}&x_{4}&x_{3}&y_{3} \\
0&0&0&\la_{1}&x_{4} &y_{4} \\ 0&0&0&0&\la_{1}&0 \\
0&0&0&0&0&0 \end{pmatrix}$$ and we see that $\la_{1} \neq 0$ or the
algebra is nilpotent. We now conjugate by
$\A_{3}\A_{13}\A_{14}\A_{15}$ and, as $\la_{1}$ and $x_{1}$ are both
nonzero, move $y_{1}$, $y_{2}$, and $x_{3}$ to 0 and $y_{3}$ to
$x_{2}$.  Then we conjugate by $\A_{1}\A_{4}$ and, as $x_{1}$,
$a_{3}$, and $\la_{1}$ are all nonzero, scale $a_{3}$ to $\la_{1}$
and $x_{1}$ to $\pm \la_{1}$. Finally, we make the change of basis
$$\e_{i}=\f_{i}, \ \ \text{for} \ \ 1 \leq i \leq 6, \ \ \text{and}
\ \ \e_{7}=\frac{1}{\la_{1}}(\f_{7}-x_{2}\f_{4}-
y_{4}\f_{5}+x_{4}\f_{6}),$$ and arrive at the structure equations
\begin{align}
    \notag [\e_{4},\e_{5}]&=\e_{2}, & [\e_{4},\e_{6}]&= \e_{3},
    &[\e_{5},\e_{6}]&=\e_{4}, & [\e_{1},\e_{7}]&=\e_{1}+\e_{3}, \\
    \notag [\e_{2},\e_{7}]&=2\e_{2}, &[\e_{3},\e_{7}]&=\e_{3}, &
    [\e_{4},\e_{7}]&=\e_{4}, & [\e_{5},\e_{7}]&=\vep\e_{1}+\e_{5},
\end{align}
with $\vep^{2} = 1$.  This is [7,[6,4],1,13] in the table.

\subsubsection{Subcase 1.2:}

We assumed here that $a_{1}=2\la_{1}+\la_{2}$ and $a_{2} \neq 0$.
This will yield that $$\ad{\f_{7}}=\begin{pmatrix}
2\la_{1}+\la_{2}&0&0&0&x_{1}&y_{1} \\
a_{2}&2\la_{1}+\la_{2}&0&-y_{4}&x_{2}&y_{2} \\
a_{3}&0&\la_{1}+2\la_{2}&x_{4}&x_{3}&y_{3} \\
0&0&0&\la_{1}+\la_{2}&x_{4} &y_{4} \\ 0&0&0&0&\la_{1}&0 \\
0&0&0&0&0&\la_{2} \end{pmatrix}.$$  If $\la_{2} \neq \la_{1}$, then
we conjugate by $\A_{3}$ and move $a_{3}$ to 0.  If, on the other
hand, $\la_{2}=\la_{1}$, then we conjugate by $\A_{5}$ and, as
$a_{2} \neq 0$, we still move $a_{3}$ to 0.  Finally, we conjugate
by $\A_{0}$ with $a=d=0$ and $b=c=1$, and arrive at $$\ad{\f_{7}}=
\begin{pmatrix} 2\la_{1}+\la_{2}&0&0&0&y_{1}&x_{1} \\
0&\la_{1}+2\la_{2}&0&x_{4}&-y_{3}&-x_{3} \\
-a_{2}&0&2\la_{1}+\la_{2}&-y_{4}&-y_{2}&-x_{2} \\
0&0&0&\la_{1}+\la_{2}&-y_{4} &-x_{4} \\ 0&0&0&0&\la_{2}&0 \\
0&0&0&0&0&\la_{1} \end{pmatrix}.$$  If we relabel $\la_{1}$,
$\la_{2}$, $-a_{2}$, $y_{1}$, $x_{1}$, $-x_{4}$, $-y_{3}$, $-x_{3}$,
$-y_{4}$, $-y_{2}$, and $-x_{2}$ as $\la_{2}$, $\la_{1}$, $a_{3}$,
$x_{1}$, $y_{1}$, $y_{4}$, $x_{2}$, $y_{2}$, $x_{4}$, $x_{3}$, and
$y_{3}$ respectively, then this is the same $\ad{\f_{7}}$ matrix as
the one in Subcase 1.1.2.

\subsection{Parent Case 2:}

In the second parent case, we classify the algebras with an
$\ad{\f_{7}}$ matrix of the form $$\ad{\f_{7}}=\begin{pmatrix}
a_{1}&0&0&0&x_{1}&y_{1} \\ a_{2}&3\la_{1}&\la_{2}&-y_{4}&x_{2}&y_{2}
\\ a_{3}&-\la_{2}&3\la_{1}&x_{4}&x_{3}&y_{3} \\
0&0&0&2\la_{1}&x_{4}&y_{4} \\ 0&0&0&0&\la_{1}&\la_{2} \\
0&0&0&0&-\la_{2}&\la_{1} \end{pmatrix},$$ with $\la_{2} \neq 0$ and
the eigenvalues ordered and labeled so that $\frac{\la_{1}}{\la_{2}}
\geq 0$.

We first conjugate by $\A_{2}\A_{3}$ and pick $s_{2}$ and $s_{3}$ to
be as follows
\begin{align}
    \notag s_{2}&=\frac{\la_{2}a_{3}-3\la_{1}a_{2}+a_{1}a_{2}}
    {(\la_{2})^{2} + (a_{1}-3\la_{1})^{2}}, & s_{3}&=
    \frac{a_{1}a_{3}-3\la_{1}a_{3}-\la_{2}a_{2}}{(\la_{2})^{2} +
    (a_{1}-3\la_{1})^{2}}.
\end{align}
As $\la_{2} \neq 0$, these denominators are nonzero.  Picking
$s_{2}$ and $s_{3}$ in this manner will move both the $a_{2}$ and
$a_{3}$ positions to 0 simultaneously.  This is the common type of
change that we make when we have a real Jordan block of this type in
the $\add$ matrix we're trying to simplify.

Next conjugate by $\A_{10}\A_{13}$ and, again as $\la_{2} \neq 0$,
we move the $x_{1}$ and $y_{1}$ positions to 0.  If $\la_{1} \neq
0$, then we conjugate by $\A_{11}\A_{12}\A_{14}$, and, as $\la_{2}
\neq 0$, we pick $s_{11}$, $s_{12}$, and $s_{14}$ to simultaneously
move the $y_{2}$ and $x_{3}$ positions to 0 and make the $x_{2}$
position equal to the $y_{3}$ position (call the common value
$y_{3}$). If $\la_{1}=0$, then we conjugate by $\A_{11}\A_{12}$ and,
as $\la_{2} \neq 0$, move $x_{3}$ to 0 and make the $x_{2}$ position
equal to the $y_{3}$ position.  This yields two cases.
\begin{enumerate}
    \item $\la_{1} \neq 0$ and we make all three changes, or
    $\la_{1}=0$ and we move $x_{3}$ to 0 and $x_{2}$ to $y_{3}$
    while $y_{2}=0$ already.
    \item $\la_{1}=0$ and we move $x_{3}$ to 0 and $x_{2}$ to
    $y_{3}$, but $y_{2} \neq 0$.
\end{enumerate}

\subsubsection{Subcase 2.1:}

We have the $\ad{\f_{7}}$ matrix here as
$$\ad{\f_{7}}=\begin{pmatrix} a_{1}&0&0&0&0&0 \\
0&3\la_{1}&\la_{2}&-y_{4}&y_{3}&0
\\ 0&-\la_{2}&3\la_{1}&x_{4}&0&y_{3} \\
0&0&0&2\la_{1}&x_{4}&y_{4} \\ 0&0&0&0&\la_{1}&\la_{2} \\
0&0&0&0&-\la_{2}&\la_{1} \end{pmatrix}.$$  Note that $a_{1}\neq 0$
or the algebra decomposes.  We then make the change of basis
$$\e_{i}=\f_{i}, \ \ \text{for} \ \ 1 \leq i \leq 6, \ \ \text{and} \
\ \e_{7}=-\frac{1}{\la_{2}}(\f_{7}-y_{3}\f_{4}-
y_{4}\f_{5}+x_{4}\f_{6}),$$ and let $a=\frac{a_{1}}{\la_{2}}$ and
$b=\frac{\la_{1}}{\la_{2}}$.  This yields the structure equations
\begin{align}
    \notag [\e_{4},\e_{5}]&=\e_{2}, & [\e_{4},\e_{6}]&= \e_{3},
    &[\e_{5},\e_{6}]&=\e_{4}, & [\e_{1},\e_{7}]&=a\e_{1}, &
    [\e_{2},\e_{7}]&=3b\e_{2}-\e_{3}, \\
    \notag [\e_{3},\e_{7}]&=\e_{2}+3b\e_{3}, &
    [\e_{4},\e_{7}]&=2b\e_{4}, & [\e_{5},\e_{7}]&=b\e_{5}-\e_{6}, &
    [\e_{6},\e_{7}]&=\e_{5}+b\e_{6},
\end{align}
with $a \neq 0$ and $b \geq 0$.  In the table in Appendix \ref{a2},
this is [7,[6,4],2,1].

\subsubsection{Subcase 2.2:}

Here we assumed that $\la_{1}=0$ and, while we moved $x_{3}$ to 0
and $x_{2}$ to $y_{3}$, that $y_{2} \neq 0$.  This gives us the
$\ad{\f_{7}}$ matrix $$\ad{\f_{7}}=\begin{pmatrix} a_{1}&0&0&0&0&0 \\
0&0&\la_{2}&-y_{4}&y_{3}&y_{2}
\\ 0&-\la_{2}&0&x_{4}&0&y_{3} \\
0&0&0&2\la_{1}&x_{4}&y_{4} \\ 0&0&0&0&0&\la_{2} \\
0&0&0&0&-\la_{2}&0 \end{pmatrix}.$$  Again, note that if $a_{1}=0$,
then the algebra is decomposable.  Then conjugate by $\A_{4}\A_{7}$
and let $s_{7}=s_{4}$.  As $\la_{2}$ and $y_{2}$ are both nonzero,
this will allow us to pick $s_{4}$ to scale $y_{2}$ to $\pm
\la_{2}$.  Then make the change of basis $$\e_{i}=\f_{i}, \ \
\text{for} \ \ 1 \leq i \leq 6, \ \ \text{and} \ \
\e_{7}=-\frac{1}{\la_{2}}(\f_{7}-y_{3}\f_{4}-
y_{4}\f_{5}+x_{4}\f_{6}),$$ and let $a=\frac{a_{1}}{\la_{2}}$.  This
yields the structure equations
\begin{align}
    \notag [\e_{4},\e_{5}]&=\e_{2}, & [\e_{4},\e_{6}]&= \e_{3},
    &[\e_{5},\e_{6}]&=\e_{4}, & [\e_{1},\e_{7}]&=a\e_{1},  \\
    \notag [\e_{2},\e_{7}]&=-\e_{3}, & [\e_{3},\e_{7}]&=\e_{2},
    & [\e_{5},\e_{7}]&=-\e_{6}, &
    [\e_{6},\e_{7}]&=\vep \e_{2}+\e_{5},
\end{align}
with $a \neq 0$ and $\vep^{2}=1$.  This is [7,[6,4],2,2].

\subsection{Parent Case 3:}

In this final section, we classify those algebras whose
$\ad{\f_{7}}$ matrix is of the form $$\ad{\f_{7}}=
\begin{pmatrix} a_{1}&0&0&0&x_{1}&y_{1} \\
a_{2}&3\la&1&-y_{4}&x_{2}&y_{2} \\
a_{3}&0&3\la&x_{4}&x_{3}&y_{3} \\
0&0&0&2\la&x_{4}&y_{4} \\ 0&0&0&0&\la&1 \\
0&0&0&0&0&\la \end{pmatrix}.$$  First, if $a_{1}\neq 3\la$, we
conjugate by $\A_{2}\A_{3}$ and pick $s_{2}$ and $s_{3}$ as follows
\begin{align}
    \notag s_{2}&=\frac{a_{3}}{a_{1}-3\la}, &
    s_{3}&=\frac{a_{1}a_{2}-3\la a_{2}+a_{3}}{(a_{1}-3\la)^{2}}.
\end{align}
This will move both the $a_{2}$ and $a_{3}$ positions to 0
simultaneously.  If, on the other hand, $a_{1}=3\la$, then we
conjugate by $\A_{3}$ and let $s_{3}=-a_{2}$.  This will move the
$a_{2}$ position to 0.  This is usually how we'll proceed if a real
Jordan block of this type is in the $\add$ matrix we're trying to
simplify. We have two cases.
\begin{enumerate}
    \item $a_{1} \neq 3\la$ and we move $a_{2}$ and $a_{3}$ to 0, or
    $a_{1}=3\la$ and we still move $a_{2}$ to 0, while $a_{3}=0$
    already.
    \item $a_{1}=3\la$ and we still move $a_{2}$ to 0, but
    $a_{3}\neq 0$.
\end{enumerate}

\subsubsection{Subcase 3.1:}

Here we have the $\ad{\f_{7}}$ matrix $$\ad{\f_{7}}=
\begin{pmatrix} a_{1}&0&0&0&x_{1}&y_{1} \\
0&3\la&1&-y_{4}&x_{2}&y_{2} \\
0&0&3\la&x_{4}&x_{3}&y_{3} \\
0&0&0&2\la&x_{4}&y_{4} \\ 0&0&0&0&\la&1 \\
0&0&0&0&0&\la \end{pmatrix}.$$  Now, if $a_{1} \neq \la$, then we
conjugate by $\A_{10}\A_{13}$ and move $x_{1}$ and $y_{1}$ to 0.  If
$a_{1}=\la$, then we conjugate by just $\A_{10}$ instead and pick
$s_{10}$ to move $y_{1}$ to 0.  This gives us two cases.
\begin{enumerate}
    \item $a_{1} \neq \la$ and we move $x_{1}$ and $y_{1}$ to 0, or
    $a_{1}=\la$ and we move $y_{1}$ to 0 and $x_{1}=0$
    already.
    \item $a_{1}=\la$ and we still move $y_{1}$ to 0, but $x_{1}
    \neq 0$.
\end{enumerate}

\subsubsection{Subcase 3.1.1:}

In this section, we consider the case when the $\ad{\f_{7}}$ matrix
is of the form
$$\ad{\f_{7}}=
\begin{pmatrix} a_{1}&0&0&0&0&0 \\
0&3\la&1&-y_{4}&x_{2}&y_{2} \\
0&0&3\la&x_{4}&x_{3}&y_{3} \\
0&0&0&2\la&x_{4}&y_{4} \\ 0&0&0&0&\la&1 \\
0&0&0&0&0&\la \end{pmatrix}.$$  If $\la \neq 0$, then we conjugate
by $\A_{11}\A_{12}\A_{14}$ and pick $s_{11}$, $s_{12}$, and $s_{14}$
to move the $y_{2}$ and $x_{3}$ positions to 0 and make the $x_{2}$
position equal to the $y_{3}$ position (label the common value
$y_{3}$).  If $\la=0$, then we conjugate by $\A_{11}\A_{12}$ and
pick $s_{11}$ and $s_{12}$ to move $y_{2}$ to 0 and $x_{2}$ to
$y_{3}$.  This gives us two cases.
\begin{enumerate}
    \item $\la \neq 0$ and we move $x_{3}$ and $y_{2}$ to 0 and
    $x_{2}$ to $y_{3}$, or $\la=0$ and we still move $y_{2}$ to 0
    and $x_{2}$ to $y_{3}$, while $x_{3}=0$ already.
    \item $\la =0$ and while we still move $y_{2}$ to 0 and $x_{2}$
    to $y_{3}$, $x_{3} \neq 0$.
\end{enumerate}

\subsubsection{Subcase 3.1.1.1:}

Here we end up with $$\ad{\f_{7}}=
\begin{pmatrix} a_{1}&0&0&0&0&0 \\
0&3\la&1&-y_{4}&y_{3}&0 \\
0&0&3\la&x_{4}&0&y_{3} \\
0&0&0&2\la&x_{4}&y_{4} \\ 0&0&0&0&\la&1 \\
0&0&0&0&0&\la \end{pmatrix}.$$  Note that $a_{1} \neq 0$ or the
algebra decomposes.  This allows us to make the change of basis
$$\e_{1}=\f_{1}, \ \ \e_{3}=(a_{1})^{2}\f_{3}, \ \
\e_{5}=\f_{5}, \ \ \e_{i}=a_{1}\f_{i}, \ \ \text{for} \ \ i \in
\{2,4,6\}, \ \ \text{and} \ \ \e_{7}=-\frac{1}{a_{1}}(\f_{7}-
y_{3}\f_{4}- y_{4}\f_{5}+x_{4}\f_{6}),$$ and let $a =
\frac{\la}{a_{1}}$.  Then we have the structure equations
\begin{align}
    \notag [\e_{4},\e_{5}]&=\e_{2}, & [\e_{4},\e_{6}]&= \e_{3},
    &[\e_{5},\e_{6}]&=\e_{4}, & [\e_{1},\e_{7}]&=\e_{1}, &
    [\e_{2},\e_{7}]&=3a\e_{2}, \\
    \notag [\e_{3},\e_{7}]&=\e_{2}+3a\e_{3}, &
    [\e_{4},\e_{7}]&=2a\e_{4}, & [\e_{5},\e_{7}]&=a\e_{5}, &
    [\e_{6},\e_{7}]&=\e_{5}+a\e_{6},
\end{align}
with $a \in \R$.  In the table in Appendix \ref{a2}, this is
[7,[6,4],3,1].

\subsubsection{Subcase 3.1.1.2:}

In this section, we have that $$\ad{\f_{7}}=
\begin{pmatrix} a_{1}&0&0&0&0&0 \\
0&0&1&-y_{4}&y_{3}&0 \\
0&0&0&x_{4}&x_{3}&y_{3} \\
0&0&0&0&x_{4}&y_{4} \\ 0&0&0&0&0&1 \\
0&0&0&0&0&0 \end{pmatrix},$$ with $x_{3} \neq 0$.  Note also that
$a_{1} \neq 0$ or the algebra is nilpotent and decomposable. Now we
conjugate by $\A_{4}\A_{7}$ and let $s_{7}=s_{4}$.  As $a_{1}$ and
$x_{3}$ are both nonzero, this allows us to scale $x_{3}$ to $\pm
(a_{1})^{3}$.  Then we make the change of basis $$\e_{1}=\f_{1}, \ \
\e_{3}=(a_{1})^{2}\f_{3}, \ \ \e_{5}=\f_{5}, \ \ \e_{i}=a_{1}\f_{i},
\ \ \text{for} \ \ i \in \{2,4,6\}, \ \ \text{and} \ \
\e_{7}=-\frac{1}{a_{1}}(\f_{7}- y_{3}\f_{4}-
y_{4}\f_{5}+x_{4}\f_{6}),$$ which will yield the structure equations
\begin{align}
    \notag [\e_{4},\e_{5}]&=\e_{2}, & [\e_{4},\e_{6}]&= \e_{3},
    &[\e_{5},\e_{6}]&=\e_{4}, & [\e_{1},\e_{7}]&=\e_{1}, \\
    \notag [\e_{3},\e_{7}]&=\e_{2}, & [\e_{5},\e_{7}]&=\vep\e_{3}, &
    [\e_{6},\e_{7}]&=\e_{5},
\end{align}
with $\vep^{2} = 1$.  This is [7,[6,4],3,2].

\subsubsection{Subcase 3.1.2:}

We assumed here that $a_{1}=\la$ and that while we moved $y_{1}$ to
0, $x_{1} \neq 0$. $$\ad{\f_{7}}=
\begin{pmatrix} \la&0&0&0&x_{1}&0 \\
0&3\la&1&-y_{4}&x_{2}&y_{2} \\
0&0&3\la&x_{4}&x_{3}&y_{3} \\
0&0&0&2\la&x_{4}&y_{4} \\ 0&0&0&0&\la&1 \\
0&0&0&0&0&\la \end{pmatrix}.$$  This implies that $\la \neq 0$ or
the algebra is nilpotent.  Then we conjugate by
$\A_{11}\A_{12}\A_{14}$ and, as $\la \neq 0$, move $y_{2}$ and
$x_{3}$ to 0 and make the $x_{2}$ position equal the $y_{3}$
position.  Next we conjugate by $\A_{1}$ and, as $\la$ and $x_{1}$
are both nonzero, scale $x_{1}$ to $\la$.  Finally, we make the
basis change $$\e_{1}=\f_{1}, \ \ \e_{3}=\la^{2}\f_{3}, \ \
\e_{5}=\f_{5}, \ \ \e_{i}=\la\f_{i}, \ \ \text{for} \ \ i \in
\{2,4,6\}, \ \ \text{and} \ \ \e_{7}=-\frac{1}{\la}(\f_{7}-
y_{3}\f_{4}- y_{4}\f_{5}+x_{4}\f_{6}),$$ and we arrive at the
structure equations
\begin{align}
    \notag [\e_{4},\e_{5}]&=\e_{2}, & [\e_{4},\e_{6}]&= \e_{3},
    &[\e_{5},\e_{6}]&=\e_{4}, & [\e_{1},\e_{7}]&=\e_{1}, &
    [\e_{2},\e_{7}]&=3\e_{2}, \\
    \notag [\e_{3},\e_{7}]&=\e_{2}+3\e_{3}, &
    [\e_{4},\e_{7}]&=2\e_{4}, & [\e_{5},\e_{7}]&=\e_{1}+\e_{5}, &
    [\e_{6},\e_{7}]&=\e_{5}+\e_{6}.
\end{align}
In the table, this is [7,[6,4],3,3].

\subsubsection{Subcase 3.2:}

In this section, we assumed that $a_{1}=3\la$ and while we moved
$a_{2}$ to 0, $a_{3} \neq 0$.  This yields that $\ad{\f_{7}}$ matrix
$$\ad{\f_{7}}=
\begin{pmatrix} 3\la&0&0&0&x_{1}&y_{1} \\
0&3\la&1&-y_{4}&x_{2}&y_{2} \\
a_{3}&0&3\la&x_{4}&x_{3}&y_{3} \\
0&0&0&2\la&x_{4}&y_{4} \\ 0&0&0&0&\la&1 \\
0&0&0&0&0&\la \end{pmatrix}.$$  This implies that $\la \neq 0$ or
the algebra is nilpotent.  This allows us to conjugate by
$\A_{10}\A_{11}\A_{12}\A_{13}\A_{14}$ and pick $s_{10}, \ldots,
s_{14}$ to move $x_{1}$, $y_{1}$, $y_{2}$, and $x_{3}$ to 0 and make
the $x_{2}$ position equal to the $y_{3}$ position.  Then we
conjugate by $\A_{1}$ and, as $a_{3}$ and $\la$ are both nonzero,
scale $a_{3}$ to $\la^{3}$.  Finally, we make the basis change
$$\e_{1}=\f_{1}, \ \ \e_{3}=\la^{2}\f_{3}, \ \ \e_{5}=\f_{5}, \ \
\e_{i}=\la\f_{i}, \ \ \text{for} \ \ i \in \{2,4,6\}, \ \ \text{and}
\ \ \e_{7}=-\frac{1}{\la}(\f_{7}- y_{3}\f_{4}-
y_{4}\f_{5}+x_{4}\f_{6}),$$ and arrive at the structure equations
\begin{align}
    \notag [\e_{4},\e_{5}]&=\e_{2}, & [\e_{4},\e_{6}]&= \e_{3},
    &[\e_{5},\e_{6}]&=\e_{4}, & [\e_{1},\e_{7}]&=3\e_{1}+\e_{3}, &
    [\e_{2},\e_{7}]&=3\e_{2}, \\
    \notag [\e_{3},\e_{7}]&=\e_{2}+3\e_{3}, &
    [\e_{4},\e_{7}]&=2\e_{4}, & [\e_{5},\e_{7}]&=\e_{5}, &
    [\e_{6},\e_{7}]&=\e_{5}+\e_{6}.
\end{align}
This is [7,[6,4],3,4]. \vspace{.1 in}

This completes the classification of seven dimensional algebras with
this nilradical.

\section{A Derivation Algebra with Semisimple Part Isomorphic to
$\si(4,\R)$}

In this section, we classify those seven dimensional Lie algebras,
$\g$, whose nilradical, $NR(\g)$, is isomorphic to the six
dimensional nilpotent algebra with structure equations
\begin{align}
    \notag [\f_{3},\f_{5}]&=\f_{2}, &[\f_{4},\f_{6}]&=\f_{2}.
\end{align}
This is Nilradical 5 from the six dimensional nilradicals listed in
Appendix \ref{a1}.  We again assume that we have a basis for $\g$
such that the first six vectors, $\f_{1}, \ldots, \f_{6}$, form a
basis for $NR(\g)$ and have the structure equations given above. Let
$\f_{7}$ be any vector not in $NR(\g)$.  Then by the Jacobi
property, $\ad{\f_{7}}$ must be of the form $$\ad{\f_{7}} =
\begin{pmatrix} a_{1}&0&c_{1}&d_{1}&x_{1}&y_{1} \\ a_{2}&
b_{2}&c_{2}&d_{2}&x_{2}&y_{2} \\
0&0&c_{3}&d_{3}&x_{3}&x_{4} \\
0&0&c_{4}&d_{4}&x_{4}&y_{4} \\ 0&0&c_{5}&c_{6}&b_{2}-c_{3}&-c_{4}
\\ 0&0&c_{6}&d_{6}&-d_{3}&b_{2}-d_{4} \end{pmatrix}.$$  Let
$m=\frac{b_{2}}{2}$, $n=c_{3}-\frac{b_{2}}{2}$, and
$p=d_{4}-\frac{b_{2}}{2}$.  Then we can rewrite $\ad{\f_{7}}$ as
$$\ad{\f_{7}} =
\begin{pmatrix} a_{1}&0&c_{1}&d_{1}&x_{1}&y_{1} \\ a_{2}&
2m&c_{2}&d_{2}&x_{2}&y_{2} \\
0&0&m+n&d_{3}&x_{3}&x_{4} \\
0&0&c_{4}&m+p&x_{4}&y_{4} \\ 0&0&c_{5}&c_{6}&m-n&-c_{4}
\\ 0&0&c_{6}&d_{6}&-d_{3}&m-p \end{pmatrix}.$$  Then the lower right
hand $4 \times 4$ block is of the form $mI_{4}+a$ where $I_{4}$
denotes the $4 \times 4$ identity and $a$ is of the form $$a=\begin{pmatrix} A&B \\
C&-A^{t} \end{pmatrix},$$ where $A,B,C$ are all $2 \times 2$
matrices, $A$ is arbitrary, and $B$ and $C$ are symmetric.  Then $a$
is an arbitrary element in the symplectic Lie algebra $\si(4,\R)$
\cite{sattinger}.  This will become especially important in a
moment.

Next we look at the other nonzero $\add$ matrices.  They are
\begin{align}
    \notag \ad{\f_{3}}&= \begin{pmatrix} 0&0&0&0&0&0 \\ 0&0&0&0&1&0
    \\ 0&0&0&0&0&0 \\ 0&0&0&0&0&0 \\ 0&0&0&0&0&0 \\ 0&0&0&0&0&0
    \end{pmatrix}, &\ad{\f_{4}}&=\begin{pmatrix} 0&0&0&0&0&0 \\
    0&0&0&0&0&1 \\ 0&0&0&0&0&0 \\ 0&0&0&0&0&0 \\ 0&0&0&0&0&0 \\
    0&0&0&0&0&0 \end{pmatrix}, \displaybreak[0] \\
    \notag \ad{\f_{5}}&= \begin{pmatrix}
    0&0&0&0&0&0 \\ 0&0&-1&0&0&0 \\ 0&0&0&0&0&0 \\ 0&0&0&0&0&0 \\
    0&0&0&0&0&0 \\ 0&0&0&0&0&0 \end{pmatrix}, & \ad{\f_{6}}&=
    \begin{pmatrix} 0&0&0&0&0&0 \\ 0&0&0&-1&0&0
    \\ 0&0&0&0&0&0 \\ 0&0&0&0&0&0 \\ 0&0&0&0&0&0 \\ 0&0&0&0&0&0
    \end{pmatrix}.
\end{align}
This will effectively allow us to move the $c_{2}$, $d_{2}$,
$x_{2}$, and $y_{2}$ positions to 0 by perturbing $\f_{7}$.

To deal with the rest of the parameters in $\ad{\f_{7}}$, we must
consider the automorphisms of $NR(\g)$.  Using Maple to compute a
basis for the derivation algebra and computing its Levi
decomposition, we find that a basis for the semisimple part is
formed by the following matrices.
\begin{align}
    \notag \begin{pmatrix} 0&0&0&0&0&0 \\ 0&0&0&0&0&0 \\ 0&0&1&0&0&0
    \\0&0&0&0&0&0 \\ 0&0&0&0&-1&0 \\ 0&0&0&0&0&0 \end{pmatrix} & &
    \begin{pmatrix} 0&0&0&0&0&0 \\ 0&0&0&0&0&0 \\ 0&0&0&0&0&0 \\
    0&0&1&0&0&0 \\ 0&0&0&0&0&-1 \\ 0&0&0&0&0&0 \end{pmatrix}, & &
    \begin{pmatrix} 0&0&0&0&0&0 \\ 0&0&0&0&0&0 \\ 0&0&0&1&0&0 \\
    0&0&0&0&0&0 \\ 0&0&0&0&0&0 \\ 0&0&0&0&-1&0 \end{pmatrix}, \displaybreak[0]\\
    \notag \begin{pmatrix} 0&0&0&0&0&0 \\ 0&0&0&0&0&0 \\ 0&0&0&0&0&0
    \\0&0&0&1&0&0 \\ 0&0&0&0&0&0 \\ 0&0&0&0&0&-1 \end{pmatrix}, & &
    \begin{pmatrix} 0&0&0&0&0&0 \\ 0&0&0&0&0&0 \\ 0&0&0&0&0&0
    \\0&0&0&0&0&0 \\ 0&0&1&0&0&0 \\ 0&0&0&0&0&0 \end{pmatrix}, & &
    \begin{pmatrix} 0&0&0&0&0&0 \\ 0&0&0&0&0&0 \\ 0&0&0&0&0&0
    \\0&0&0&0&0&0 \\ 0&0&0&1&0&0 \\ 0&0&1&0&0&0 \end{pmatrix}, \displaybreak[0]\\
    \notag \begin{pmatrix} 0&0&0&0&0&0 \\ 0&0&0&0&0&0 \\ 0&0&0&0&0&0
    \\0&0&0&0&0&0 \\ 0&0&0&0&0&0 \\ 0&0&0&1&0&0 \end{pmatrix}, & &
    \begin{pmatrix} 0&0&0&0&0&0 \\ 0&0&0&0&0&0 \\ 0&0&0&0&1&0
    \\0&0&0&0&0&0 \\ 0&0&0&0&0&0 \\ 0&0&0&0&0&0 \end{pmatrix}, & &
    \begin{pmatrix} 0&0&0&0&0&0 \\ 0&0&0&0&0&0 \\ 0&0&0&0&0&1
    \\0&0&0&0&1&0 \\ 0&0&0&0&0&0 \\ 0&0&0&0&0&0 \end{pmatrix}, \displaybreak[0] \\
    \notag \begin{pmatrix} 0&0&0&0&0&0 \\ 0&0&0&0&0&0 \\ 0&0&0&0&0&0
    \\0&0&0&0&0&1 \\ 0&0&0&0&0&0 \\ 0&0&0&0&0&0 \end{pmatrix},
\end{align}
which we'll call $D_{1}, \ldots, D_{10}$ respectively.  The lower
right hand $4 \times 4$ submatrices of $D_{1}, \ldots D_{10}$ are
all linearly independent and in $\si(4,\R)$ as they are of the form
discussed above.  As $\dim \si(4,\R)=10$, then this set of
submatrices form a basis for $\si(4,\R)$, and clearly the semisimple
part of our derivation algebra is isomorphic to $\si(4,\R)$.  Then
there exists an automorphism of $NR(\g)$ of the form
$$\begin{pmatrix} I_{2}&0 \\ 0& S \end{pmatrix},$$ where $I_{2}$ is
the $2 \times 2$ identity and $S$ is an arbitrary element in the
symplectic group $Sp(4,\R)$.  Call this automorphism $\A_{0}$.

If we use Maple to compute a full basis of derivations and
exponentiate them, we find that the automorphism group of $NR(\g)$
is generated by the following one parameter groups of
transformations
\begin{align}
    \notag \A_{1}&=\begin{pmatrix} s_{1}&0&0&0&0&0 \\ 0&1&0&0&0&0 \\
    0&0&1&0&0&0 \\ 0&0&0&1&0&0 \\ 0&0&0&0&1&0 \\ 0&0&0&0&0&1
    \end{pmatrix}, & \A_{2}&=\begin{pmatrix} 1&0&0&0&0&0 \\ s_{2}&1&0&0&0&0 \\
    0&0&1&0&0&0 \\ 0&0&0&1&0&0 \\ 0&0&0&0&1&0 \\ 0&0&0&0&0&1
    \end{pmatrix}, \displaybreak[0] \\
    \notag \A_{3}&=\begin{pmatrix} 1&0&0&0&0&0 \\ 0&1&0&0&0&0 \\
    0&0&s_{3}&0&0&0 \\ 0&0&0&1&0&0 \\ 0&0&0&0&1&0 \\ 0&0&0&0&0&1
    \end{pmatrix}, & \A_{4}&=\begin{pmatrix} 1&0&s_{4}&0&0&0 \\ 0&1&0&0&0&0 \\
    0&0&1&0&0&0 \\ 0&0&0&1&0&0 \\ 0&0&0&0&1&0 \\ 0&0&0&0&0&1
    \end{pmatrix}, \displaybreak[0]\\
    \notag \A_{5}&=\begin{pmatrix} 1&0&0&0&0&0 \\ 0&1&s_{5}&0&0&0 \\
    0&0&1&0&0&0 \\ 0&0&0&1&0&0 \\ 0&0&0&0&1&0 \\ 0&0&0&0&0&1
    \end{pmatrix}, & \A_{6}&=\begin{pmatrix} 1&0&0&0&0&0 \\ 0&1&0&0&0&0 \\
    0&0&s_{6}&0&0&0 \\ 0&0&0&1&0&0 \\ 0&0&0&0&\frac{1}{s_{6}}&0 \\ 0&0&0&0&0&1
    \end{pmatrix}, \displaybreak[0] \\
    \notag \A_{7}&=\begin{pmatrix} 1&0&0&0&0&0 \\ 0&1&0&0&0&0 \\
    0&0&1&0&0&0 \\ 0&0&s_{7}&1&0&0 \\ 0&0&0&0&1&-s_{7} \\ 0&0&0&0&0&1
    \end{pmatrix}, & \A_{8}&=\begin{pmatrix} 1&0&0&0&0&0 \\ 0&1&0&0&0&0 \\
    0&0&1&0&0&0 \\ 0&0&0&1&0&0 \\ 0&0&s_{8}&0&1&0 \\ 0&0&0&0&0&1
    \end{pmatrix}, \displaybreak[0] \\
    \notag \A_{9}&=\begin{pmatrix} 1&0&0&0&0&0 \\ 0&1&0&0&0&0 \\
    0&0&1&0&0&0 \\ 0&0&0&1&0&0 \\ 0&0&0&s_{9}&1&0 \\ 0&0&s_{9}&0&0&1
    \end{pmatrix}, & \A_{10}&=\begin{pmatrix} 1&0&0&s_{10}&0&0 \\ 0&1&0&0&0&0 \\
    0&0&1&0&0&0 \\ 0&0&0&1&0&0 \\ 0&0&0&0&1&0 \\ 0&0&0&0&0&1
    \end{pmatrix}, \displaybreak[0] \\
    \notag \A_{11}&=\begin{pmatrix} 1&0&0&0&0&0 \\ 0&1&0&s_{11}&0&0 \\
    0&0&1&0&0&0 \\ 0&0&0&1&0&0 \\ 0&0&0&0&1&0 \\ 0&0&0&0&0&1
    \end{pmatrix}, & \A_{12}&=\begin{pmatrix} 1&0&0&0&0&0 \\ 0&1&0&0&0&0 \\
    0&0&1&s_{12}&0&0 \\ 0&0&0&1&0&0 \\ 0&0&0&0&1&0 \\ 0&0&0&0&-s_{12}&1
    \end{pmatrix}, \displaybreak[0] \\
    \notag \A_{13}&=\begin{pmatrix} 1&0&0&0&0&0 \\ 0&1&0&0&0&0 \\
    0&0&1&0&0&0 \\ 0&0&0&s_{13}&0&0 \\ 0&0&0&0&1&0 \\
    0&0&0&0&0&\frac{1}{s_{13}}
    \end{pmatrix}, & \A_{14}&=\begin{pmatrix} 1&0&0&0&0&0 \\ 0&1&0&0&0&0 \\
    0&0&1&0&0&0 \\ 0&0&0&1&0&0 \\ 0&0&0&0&1&0 \\ 0&0&0&s_{14}&0&1
    \end{pmatrix}, \displaybreak[0] \\
    \notag \A_{15}&=\begin{pmatrix} 1&0&0&0&s_{15}&0 \\ 0&1&0&0&0&0 \\
    0&0&1&0&0&0 \\ 0&0&0&1&0&0 \\ 0&0&0&0&1&0 \\ 0&0&0&0&0&1
    \end{pmatrix}, & \A_{16}&=\begin{pmatrix} 1&0&0&0&0&0 \\ 0&1&0&0&s_{16}&0 \\
    0&0&1&0&0&0 \\ 0&0&0&1&0&0 \\ 0&0&0&0&1&0 \\ 0&0&0&0&0&1
    \end{pmatrix}, \displaybreak[0] \\
    \notag \A_{17}&=\begin{pmatrix} 1&0&0&0&0&0 \\ 0&1&0&0&0&0 \\
    0&0&1&0&s_{17}&0 \\ 0&0&0&1&0&0 \\ 0&0&0&0&1&0 \\ 0&0&0&0&0&1
    \end{pmatrix}, & \A_{18}&=\begin{pmatrix} 1&0&0&0&0&0 \\ 0&1&0&0&0&0 \\
    0&0&1&0&0&s_{18} \\ 0&0&0&1&s_{18}&0 \\ 0&0&0&0&1&0 \\ 0&0&0&0&0&1
    \end{pmatrix}, \displaybreak[0] \\
    \notag \A_{19}&=\begin{pmatrix} 1&0&0&0&0&s_{19} \\ 0&1&0&0&0&0 \\
    0&0&1&0&0&0 \\ 0&0&0&1&0&0 \\ 0&0&0&0&1&0 \\ 0&0&0&0&0&1
    \end{pmatrix}, & \A_{20}&=\begin{pmatrix} 1&0&0&0&0&0 \\ 0&1&0&0&0&s_{20} \\
    0&0&1&0&0&0 \\ 0&0&0&1&0&0 \\ 0&0&0&0&1&0 \\ 0&0&0&0&0&1
    \end{pmatrix}, \displaybreak[0] \\
    \notag \A_{21}&=\begin{pmatrix} 1&0&0&0&0&0 \\ 0&1&0&0&0&0 \\
    0&0&1&0&0&s_{21} \\ 0&0&0&1&0&0 \\ 0&0&0&0&1&0 \\ 0&0&0&0&0&1
    \end{pmatrix}.
\end{align}
Note that conjugating $\ad{\f_{7}}$ by $\A_{5}$, $\A_{11}$,
$\A_{16}$, or $\A_{20}$ will affect most significantly those entries
in $\ad{\f_{7}}$ that we're going to move to zero by perturbing
$\f_{7}$; consequently, these will be less useful in simplifying
$\ad{\f_{7}}$ that the others.

We consider now conjugating $\ad{\f_{7}}$ by $\A_{0}$.  By block
multiplication, the symplectic submatrix of $\A_{0}$, $S$, will only
affect the lower right hand $4 \times 4$ piece of $\ad{\f_{7}}$. And
as that part of $\ad{\f_{7}}$ could be written as $mI_{4}+a$, we
have that
\begin{align}
    \notag S^{-1}(mI_{4}+a)S &= S^{-1}(mI_{4})S + S^{-1}aS =
    mI_{4}+S^{-1}aS,
\end{align}
As $S \in Sp(4,\R)$ and $a \in \si(4,\R)$ arbitrarily.  Then, by the
section on the real symplectic canonical form in Chapter \ref{c2},
conjugation by $\A_{0}$ can be used to put $a$ into real symplectic
canonical form. This yields ten parent cases.
\begin{align}
    \notag 1. \ \ \ad{\f_{7}}&=\begin{pmatrix}
    a_{1}&0&c_{1}&d_{1}&x_{1}&y_{1} \\ a_{2}&
    2m&c_{2}&d_{2}&x_{2}&y_{2} \\ 0&0&m+\la&0&0&0 \\ 0&0&0&m+\mu&0&0
    \\ 0&0&0&0&m-\la&0 \\ 0&0&0&0&0&m-\mu \end{pmatrix}.
    \displaybreak[0] \\
    \notag 2. \ \ \ad{\f_{7}}&=\begin{pmatrix}
    a_{1}&0&c_{1}&d_{1}&x_{1}&y_{1} \\ a_{2}&
    2m&c_{2}&d_{2}&x_{2}&y_{2} \\ 0&0&m+\la&0&0&0 \\ 0&0&0&m&0&\vep
    \\ 0&0&0&0&m-\la&0 \\ 0&0&0&0&0&m \end{pmatrix}, \ \text{with
    $\vep^{2} = 1$.} \displaybreak[0] \\
    \notag 3. \ \ \ad{\f_{7}}&=\begin{pmatrix}
    a_{1}&0&c_{1}&d_{1}&x_{1}&y_{1} \\ a_{2}&
    2m&c_{2}&d_{2}&x_{2}&y_{2} \\ 0&0&m+\la&1&0&0 \\ 0&0&0&m+\la&0&0
    \\ 0&0&0&0&m-\la&0 \\ 0&0&0&0&-1&m-\la \end{pmatrix}.
    \displaybreak[0] \\
    \notag 4. \ \ \ad{\f_{7}}&=\begin{pmatrix}
    a_{1}&0&c_{1}&d_{1}&x_{1}&y_{1} \\ a_{2}&
    2m&c_{2}&d_{2}&x_{2}&y_{2} \\ 0&0&m&0&\vep&0 \\ 0&0&0&m&0&\vep
    \\ 0&0&0&0&m&0 \\ 0&0&0&0&0&m \end{pmatrix}, \ \text{with
    $\vep^{2}=1$.} \displaybreak[0] \\
    \notag 5. \ \ \ad{\f_{7}}&=\begin{pmatrix}
    a_{1}&0&c_{1}&d_{1}&x_{1}&y_{1} \\ a_{2}&
    2m&c_{2}&d_{2}&x_{2}&y_{2} \\ 0&0&m&1&0&0 \\ 0&0&0&m&0&\vep
    \\ 0&0&0&0&m&0 \\ 0&0&0&0&-1&m \end{pmatrix}, \ \text{with
    $\vep^{2}=1$.} \displaybreak[0] \\
    \notag 6. \ \ \ad{\f_{7}}&=\begin{pmatrix}
    a_{1}&0&c_{1}&d_{1}&x_{1}&y_{1} \\ a_{2}&
    2m&c_{2}&d_{2}&x_{2}&y_{2} \\ 0&0&m+\la&0&0&0 \\ 0&0&0&m&0&\vep
    \mu\\ 0&0&0&0&m-\la&0 \\ 0&0&0&-\vep \mu&0&m \end{pmatrix}, \
    \text{with $\mu \neq 0$ and $\vep^{2}=1$.} \displaybreak[0] \\
    \notag 7. \ \ \ad{\f_{7}}&=\begin{pmatrix}
    a_{1}&0&c_{1}&d_{1}&x_{1}&y_{1} \\ a_{2}&
    2m&c_{2}&d_{2}&x_{2}&y_{2} \\ 0&0&m&0&\vep&0 \\ 0&0&0&m&0&\del
    \mu\\ 0&0&0&0&m&0 \\ 0&0&0&-\del\mu&0&m \end{pmatrix},
    \text{with $\mu \neq 0$ and $\vep^{2}=\del^{2}=1$.} \displaybreak[0] \\
    \notag 8. \ \ \ad{\f_{7}}&=\begin{pmatrix}
    a_{1}&0&c_{1}&d_{1}&x_{1}&y_{1} \\ a_{2}&
    2m&c_{2}&d_{2}&x_{2}&y_{2} \\ 0&0&m+\la&\mu&0&0 \\ 0&0&-\mu&m+\la&0&0
    \\ 0&0&0&0&m-\la&\mu \\ 0&0&0&0&-\mu&m-\la \end{pmatrix}, \
    \begin{array}{l} \text{with $\mu \neq 0$ and the eigenvalues} \\
    \text{ordered so that $\frac{\la}{\mu} \geq 0$.} \end{array} \displaybreak[0] \\
    \notag 9. \ \ \ad{\f_{7}}&=\begin{pmatrix}
    a_{1}&0&c_{1}&d_{1}&x_{1}&y_{1} \\ a_{2}&
    2m&c_{2}&d_{2}&x_{2}&y_{2} \\ 0&0&m&0&\vep \mu&0 \\
    0&0&0&m&0&\vep \eta \\ 0&0&-\vep \mu &0&m&0 \\ 0&0&0&0-\vep \eta&0&m
    \end{pmatrix}, \ \text{with $\eta \neq -\mu$, $\eta,\mu \neq 0$,
    and $\vep^{2}=1$.} \displaybreak[0] \\
    \notag 10. \ \ \ad{\f_{7}}&=\begin{pmatrix}
    a_{1}&0&c_{1}&d_{1}&x_{1}&y_{1} \\ a_{2}&
    2m&c_{2}&d_{2}&x_{2}&y_{2} \\ 0&0&m&\mu&1&0 \\ 0&0&-\mu&m&0&1
    \\ 0&0&0&0&m&\mu \\ 0&0&0&0&-\mu&m \end{pmatrix}, \
    \begin{array}{l} \text{with $\mu \neq 0$ and the eigenvalues} \\
    \text{ordered so that $\frac{m}{\mu} \geq 0$.} \end{array}
\end{align}
From here, the classification follows in a similar manner as the
previous sections and so we omit the remainder of the classification
proof.  The resultant algebras, however, are in the table in
Appendix \ref{a2}.

\section{A Derivation Algebra with Semisimple Part Isomorphic to
$\so(3,1,\R)$.}

In this section, we classify those seven dimensional Lie algebras,
$\g$, whose nilradical, $NR(\g)$, is isomorphic to the six
dimensional nilpotent algebra with structure equations
\begin{align}
    \notag [\f_{3},\f_{5}]&=\f_{2}, &[\f_{3},\f_{6}]&=\f_{1},
    &[\f_{4},\f_{5}]&=-\f_{1}, &[\f_{4},\f_{6}]&=\f_{2}.
\end{align}
This is Nilradical 9 from the six dimensional nilradicals listed in
Appendix \ref{a1}.  We again assume that we have a basis for $\g$
such that the first six vectors, $\f_{1}, \ldots, \f_{6}$, form a
basis for $NR(\g)$ and have the structure equations given above. Let
$\f_{7}$ be any vector not in $NR(\g)$.  Then by the Jacobi
property, $\ad{\f_{7}}$ must be of the form $$\ad{\f_{7}}=
\begin{pmatrix} c_{3}+x_{5}& -c_{4}+x_{6}& c_{1}&d_{1}&x_{1}&y_{1} \\
c_{4}-x_{6} & c_{3}+x_{5} &c_{2}&d_{2}&x_{2}&y_{2} \\
0&0&c_{3}&-c_{4}&x_{3}&x_{4} \\ 0&0&c_{4}&c_{3}&x_{4}&-x_{3} \\
0&0&c_{5}&c_{6}&x_{5}&-x_{6} \\ 0&0&c_{6}&-c_{5}&x_{6}&-x_{5}
\end{pmatrix}.$$ Let $m=\frac{1}{2}(c_{3}+x_{5})$, $n=\frac{1}{2}
(-c_{4}+x_{6})$, $p=\frac{1}{2} (c_{3}-x_{5})$, and $q=\frac{1}{2}
(c_{4}+x_{6})$.  Then we can rewrite $\ad{\f_{7}}$ as $$\ad{\f_{7}}=
\begin{pmatrix} 2m&2n&c_{1}&d_{1}&x_{1}&y_{1} \\
-2n&2m&c_{2}&d_{2}&x_{2}&y_{2} \\ 0&0&m+p&n-q&x_{3}&x_{4} \\
0&0&-n+q&m+p&x_{4}&-x_{3} \\ 0&0&c_{5}&c_{6}&m-p&-n-q \\
0&0&c_{6}&-c_{5}&n+q&m-p \end{pmatrix}.$$  Let $I_{4}$ denote the $4
\times 4$ identity and let $$K=\begin{pmatrix} 0&1&0&0& \\ -1&0&0&0
\\ 0&0&0&-1 \\ 0&0&1&0 \end{pmatrix}.$$  Then the lower right hand
$4 \times 4$ submatrix of $\ad{\f_{7}}$ is of the form $mI_{4} + nK
+ a$ where $a$ is of the form $$a=\begin{pmatrix} A&B \\ C&-A^{t}
\end{pmatrix},$$ with $A,B,C$ all $2 \times2$ matrices such that
$A=\bigl{(}\begin{smallmatrix} p&-q \\ q&p
\end{smallmatrix}\bigr{)}$ and $B$ and $C$ are trace-free symmetric.
Then $a$ is an arbitrary element of $\h(J_{2})$, the nonstandard
representation of $\so(3,1,\R)$ that we studied in Chapter \ref{c2}.
Again, this will be important later.

Next, we look at the other nonzero $\add$ matrices.  They are
\begin{align}
    \notag \ad{\f_{3}} &= \begin{pmatrix} 0&0&0&0&0&1 \\ 0&0&0&0&1&0
    \\ 0&0&0&0&0&0 \\ 0&0&0&0&0&0 \\ 0&0&0&0&0&0 \\ 0&0&0&0&0&0
    \end{pmatrix}, & \ad{\f_{4}}&= \begin{pmatrix} 0&0&0&0&-1&0 \\
    0&0&0&0&0&1 \\ 0&0&0&0&0&0 \\ 0&0&0&0&0&0 \\0&0&0&0&0&0 \\
    0&0&0&0&0&0 \end{pmatrix}, \\
    \notag \ad{\f_{5}} &= \begin{pmatrix} 0&0&0&1&0&0 \\
    0&0&-1&0&0&0 \\ 0&0&0&0&0&0 \\ 0&0&0&0&0&0 \\ 0&0&0&0&0&0
    \\0&0&0&0&0&0 \end{pmatrix}, & \ad{\f_{6}}&=\begin{pmatrix} 0&0&-1&0&0&0 \\
    0&0&0&-1&0&0 \\ 0&0&0&0&0&0 \\ 0&0&0&0&0&0 \\ 0&0&0&0&0&0 \\
    0&0&0&0&0&0 \end{pmatrix}.
\end{align}
If we can move $y_{1}$ to $x_{2}$, $y_{2}$ to $-x_{1}$, $d_{1}$ to
$-c_{2}$, and $d_{2}$ to $c_{1}$ in $\ad{\f_{7}}$, then we can
perturb $\f_{7}$ to annihilate the entire upper left hand $2 \times
4$ block.

In order to do this and to take care of the remaining parameters, we
need to compute the derivation algebra.  We use Maple to compute a
basis for the Lie algebra of derivations of $NR(\g)$ and its Levi
decomposition.  After doing this, we'll find that a basis for the
semisimple part of the derivation algebra is given by
\begin{align}
    \notag \begin{pmatrix} 0&0&0&0&0&0 \\ 0&0&0&0&0&0 \\ 0&0&1&0&0&0
    \\ 0&0&0&1&0&0 \\ 0&0&0&0&-1&0 \\ 0&0&0&0&0&-1 \end{pmatrix},
    & & \begin{pmatrix} 0&0&0&0&0&0 \\ 0&0&0&0&0&0 \\ 0&0&0&-1&0&0
    \\ 0&0&1&0&0&0 \\ 0&0&0&0&0&-1 \\ 0&0&0&0&1&0 \end{pmatrix}, \displaybreak[0] \\
    \notag \begin{pmatrix} 0&0&0&0&0&0 \\ 0&0&0&0&0&0 \\ 0&0&0&0&0&0
    \\ 0&0&0&0&0&0 \\ 0&0&1&0&0&0 \\ 0&0&0&-1&0&0  \end{pmatrix}, &
    & \begin{pmatrix} 0&0&0&0&0&0 \\ 0&0&0&0&0&0 \\ 0&0&0&0&0&0 \\
    0&0&0&0&0&0 \\ 0&0&0&1&0&0 \\ 0&0&1&0&0&0 \end{pmatrix}, \displaybreak[0] \\
    \notag \begin{pmatrix} 0&0&0&0&0&0 \\ 0&0&0&0&0&0 \\ 0&0&0&0&1&0
    \\ 0&0&0&0&0&-1 \\ 0&0&0&0&0&0 \\ 0&0&0&0&0&0 \end{pmatrix}, & &
    \begin{pmatrix} 0&0&0&0&0&0 \\ 0&0&0&0&0&0 \\ 0&0&0&0&0&1 \\
    0&0&0&0&1&0 \\ 0&0&0&0&0&0 \\ 0&0&0&0&0&0 \end{pmatrix},
\end{align}
which we'll call $D_{1},\ldots, D_{6}$ respectively.  Note that the
lower right hand $4 \times 4$ submatrices of $D_{1},\ldots, D_{6}$
form a basis $\h(J_{2})$. This implies that the semisimple part of
the derivation algebra is naturally isomorphic to $\h(J_{2})$.  In
addition, this implies that the matrix
$$A_{0}=\begin{pmatrix} I_{2} & 0 \\ 0& S \end{pmatrix},$$ with $S$
an arbitrary member of $H(J_{2})$ is an automorphism of $NR(\g)$.

At this point, we use Maple to compute a complete basis of the Lie
algebra of derivations of $NR(\g)$.  We then exponentiate them to
find that the automorphism group of $NR(\g)$ is generated by the
following one parameter groups of transformations
\begin{align}
    \notag \A_{1}&= \begin{pmatrix} s_{1}&0&0&0&0&0 \\ 0&s_{1}&0&0&0&0 \\
    0&0&1&0&0&0 \\ 0&0&0&1&0&0 \\ 0&0&0&0&s_{1}&0 \\ 0&0&0&0&0&s_{1}
    \end{pmatrix}, & \A_{2}&= \begin{pmatrix} \cos(s_{2})&-\sin(s_{2})&0&0&0&0 \\ \sin(s_{2})&\cos(s_{2})&0&0&0&0 \\
    0&0&1&0&0&0 \\ 0&0&0&1&0&0 \\ 0&0&0&0&\cos(s_{2})&\sin(s_{2}) \\ 0&0&0&0&-\sin(s_{2})&\cos(s_{2})
    \end{pmatrix}, \displaybreak[0] \\
    \notag \A_{3}&= \begin{pmatrix} 1&0&s_{3}&0&0&0 \\ 0&1&0&0&0&0 \\
    0&0&1&0&0&0 \\ 0&0&0&1&0&0 \\ 0&0&0&0&1&0 \\ 0&0&0&0&0&1
    \end{pmatrix}, & \A_{4}&= \begin{pmatrix} 1&0&0&0&0&0 \\ 0&1&s_{4}&0&0&0 \\
    0&0&1&0&0&0 \\ 0&0&0&1&0&0 \\ 0&0&0&0&1&0 \\ 0&0&0&0&0&1
    \end{pmatrix}, \displaybreak[0] \\
    \notag \A_{5}&= \begin{pmatrix} 1&0&0&0&0&0 \\ 0&1&0&0&0&0 \\
    0&0&s_{5}&0&0&0 \\ 0&0&0&s_{5}&0&0 \\ 0&0&0&0&\frac{1}{s_{5}}&0 \\
    0&0&0&0&0&\frac{1}{s_{5}}
    \end{pmatrix}, & \A_{6}&= \begin{pmatrix} 1&0&0&0&0&0 \\ 0&1&0&0&0&0 \\
    0&0&\cos(s_{6})&-\sin(s_{6})&0&0 \\ 0&0&\sin(s_{6})&\cos(s_{6})&0&0 \\ 0&0&0&0&\cos(s_{6})&-\sin(s_{6})
    \\ 0&0&0&0&\sin(s_{6})&\cos(s_{6})
    \end{pmatrix}, \displaybreak[0] \\
    \notag \A_{7}&= \begin{pmatrix} 1&0&0&0&0&0 \\ 0&1&0&0&0&0 \\
    0&0&1&0&0&0 \\ 0&0&0&1&0&0 \\ 0&0&s_{7}&0&1&0 \\ 0&0&0&-s_{7}&0&1
    \end{pmatrix}, & \A_{8}&= \begin{pmatrix} 1&0&0&0&0&0 \\ 0&1&0&0&0&0 \\
    0&0&1&0&0&0 \\ 0&0&0&1&0&0 \\ 0&0&0&s_{8}&1&0 \\ 0&0&s_{8}&0&0&1
    \end{pmatrix}, \displaybreak[0] \\
    \notag \A_{9}&= \begin{pmatrix} 1&0&0&s_{9}&0&0 \\ 0&1&0&0&0&0 \\
    0&0&1&0&0&0 \\ 0&0&0&1&0&0 \\ 0&0&0&0&1&0 \\ 0&0&0&0&0&1
    \end{pmatrix}, & \A_{10}&= \begin{pmatrix} 1&0&0&0&0&0 \\ 0&1&0&s_{10}&0&0 \\
    0&0&1&0&0&0 \\ 0&0&0&1&0&0 \\ 0&0&0&0&1&0 \\ 0&0&0&0&0&1
    \end{pmatrix}, \displaybreak[0] \\
    \notag \A_{11}&= \begin{pmatrix} 1&0&0&0&s_{11}&0 \\ 0&1&0&0&0&0 \\
    0&0&1&0&0&0 \\ 0&0&0&1&0&0 \\ 0&0&0&0&1&0 \\ 0&0&0&0&0&1
    \end{pmatrix}, & \A_{12}&= \begin{pmatrix} 1&0&0&0&0&0 \\ 0&1&0&0&s_{12}&0 \\
    0&0&1&0&0&0 \\ 0&0&0&1&0&0 \\ 0&0&0&0&1&0 \\ 0&0&0&0&0&1
    \end{pmatrix}, \displaybreak[0] \\
    \notag \A_{13}&= \begin{pmatrix} 1&0&0&0&0&0 \\ 0&1&0&0&0&0 \\
    0&0&1&0&s_{13}&0 \\ 0&0&0&1&0&-s_{13} \\ 0&0&0&0&1&0 \\ 0&0&0&0&0&1
    \end{pmatrix}, & \A_{14}&= \begin{pmatrix} 1&0&0&0&0&0 \\ 0&1&0&0&0&0 \\
    0&0&1&0&0&s_{14} \\ 0&0&0&1&s_{14}&0 \\ 0&0&0&0&1&0 \\ 0&0&0&0&0&1
    \end{pmatrix}, \displaybreak[0] \\
    \notag \A_{15}&= \begin{pmatrix} 1&0&0&0&0&s_{15} \\ 0&1&0&0&0&0 \\
    0&0&1&0&0&0 \\ 0&0&0&1&0&0 \\ 0&0&0&0&1&0 \\ 0&0&0&0&0&1
    \end{pmatrix}, & \A_{16}&= \begin{pmatrix} 1&0&0&0&0&0 \\ 0&1&0&0&0&s_{16} \\
    0&0&1&0&0&0 \\ 0&0&0&1&0&0 \\ 0&0&0&0&1&0 \\ 0&0&0&0&0&1
    \end{pmatrix}.
\end{align}

Now note that if $S \in H(J_{2})$, then $S^{t}J_{i}S = J_{i}$ for
all $i \in \{1,2\}$.  This implies that $J_{i}S=(S^{t})^{-1}J_{i}$
and $S^{t}J_{i}=J_{i}S^{-1}$. Then, as $(J_{i})^{-1}=-J_{i}$, we
have
$$S^{-1}J_{i} = (-J_{i}S)^{-1} = (-(S^{t})^{-1}J_{i})^{-1} =J_{i}S^{t}$$
for all $i \in \{1,2\}$.  Hence we have that
$$S^{-1}J_{1}J_{2}S = J_{1}S^{t}J_{2}S = J_{1}J_{2}S^{-1}S = J_{1}J_{2}.$$ Moreover,
note that $J_{1}J_{2} = K$ and hence $S^{-1}KS = K$.

By block multiplication, we have that the lower right hand $4 \times
4$ block of $(A_{0})^{-1} \ad{\f_{7}} A_{0}$ will be $S^{-1} (m
I_{4} + nK + a)S$.  This yields
\begin{align}
    \notag S^{-1}(m I_{4} + nK + a) S &= m (S^{-1} I_{4} S) + n
    (S^{-1} K S) + S^{-1}aS = m I_{4} + nK + S^{-1}aS.
\end{align}
As $S$ is an arbitrary member of the group $H(J_{2})$ and $a$ is an
arbitrary member of $\h(J_{2})$, then we pick $S$ to put $a$ into
real $\h$-symplectic canonical form. This yields three parent cases
\begin{align}
    \notag 1. \ \ \ad{\f_{7}}&=\begin{pmatrix} 2m&2n&c_{1}&d_{1}&x_{1}&y_{1} \\
    -2n&2m&c_{2}&d_{2}&x_{2}&y_{2} \\ 0&0&m+\la&n&0&0 \\
    0&0&-n&m+\la&0&0 \\ 0&0&0&0&m-\la&-n \\ 0&0&0&0&n&m-\la
    \end{pmatrix}. \\
    \notag 2. \ \ \ad{\f_{7}}&=\begin{pmatrix} 2m&2n&c_{1}&d_{1}&x_{1}&y_{1} \\
    -2n&2m&c_{2}&d_{2}&x_{2}&y_{2} \\ 0&0&m&n&\vep&0 \\
    0&0&-n&m&0&-\vep \\ 0&0&0&0&m&-n \\ 0&0&0&0&n&m
    \end{pmatrix}, \ \text{with $\vep^{2}=1$.} \\
    \notag 3. \ \ \ad{\f_{7}}&=\begin{pmatrix} 2m&2n&c_{1}&d_{1}&x_{1}&y_{1} \\
    -2n&2m&c_{2}&d_{2}&x_{2}&y_{2} \\ 0&0&m+\la&n-\vep\mu&0&0 \\
    0&0&-n+\vep\mu&m+\la&0&0 \\ 0&0&0&0&m-\la&-n-\vep\mu \\ 0&0&0&0&n+\vep\mu&m-\la
    \end{pmatrix}, \ \text{with $\mu > 0$ and $\vep^{2}=1$.}
\end{align}
From here, the classification runs similar to the previous section
and so we omit the remainder of the classification proof.  However,
the algebras are listed in Appendix \ref{a2}.

This completes our sample of the classification of the seven
dimensional algebras and the text of this paper.  The multiplication
tables of all the isomorphism classes of indecomposable solvable Lie
algebras of dimension two through dimension seven with codimension
one nilradicals can be found in Appendix \ref{a2}.

%
\setlength{\baselineskip}{11pt}


\newpage
\setlength{\baselineskip}{22pt}
\bibliographystyle{siam}
\bibliography{References}


%
%
%
%
%
%
%
%
%
%
%
%
%
%
%
%
%
%
\newpage

\addcontentsline{toc}{chapter}{APPENDICES}

\begin{center}

\null \vskip 260pt APPENDICES \\

\end{center}

\vfill \eject

\newpage

\appendix
%
\setlength{\baselineskip}{22pt}

%
%

\chapter{A LIST OF NILPOTENT ALGEBRAS USED AS NILRADICALS}\label{a1}

In this appendix, we simply list and number the multiplication
tables of the nilpotent algebras from dimension three through
dimension six (for dimension one and two the only nilpotent algebra
is the abelian algebra). These are what we used as nilradicals for
our classification.

The lists were compiled by other authors. Dimensions three through
five are from Winternitz's list; wherever possible, we've indicated
the appropriate reference in Mubarakzyanov's as well
\cite{winternitz-1,mubar-1}.  In the classification, Winternitz and
Mubarakzyanov's names are abbreviated to Win and Mubar respectively.
Dimension six is from Gong's classification \cite{gong}. In using
Gong's classification, we have made a change of basis for every
algebra. If $[\e_{1},\ldots,\e_{6}]$ is the basis that Gong used,
then we used the basis $[-\e_{6},\ldots,-\e_{1}]$ for every six
dimensional nilpotent algebra except number five; for that algebra,
we instead used the basis $[-\e_{6},-\e_{5},-\e_{4},-\e_{2},
-\e_{3},-\e_{1}]$. This was done for consistency in our
classification.

In addition, we have imposed the numbering and will use it in the
numbering scheme of Appendix \ref{a2}. Also we don't list the
abelian algebra of each dimension, $m$, but we will refer to it as
$[m,0]$.

\twocolumn \setlength{\parindent}{0 pt}

\section{Dimension Three}

\vbox{ {\bf[3, 1]}
\medskip

\hbox{

$
\begin{tabular}[t]{L| C  C  C } \ &\ \e_1  &\e_2   &\e_3 \\ \hline
\, \e_1 & . & . & .\\
\, \e_2 &  & . & \, \e_1\\
\, \e_3 &  &  & .\\
\end{tabular}$

\qquad
\begin{tabular}[t] {K K} \\[-10pt]
Win: &[3, 1]\\
\end{tabular}
} }

\bigskip

\section{Dimension Four}

\vbox{ {\bf[4, 1]}
\medskip

\hbox{

$\begin{tabular}[t]{L| C  C  C  C } \ &\ \e_1  &\e_2   &\e_3   &\e_4 \\
\hline
\, \e_1 & . & . & . & .\\
\, \e_2 &  & . & \, \e_1 & .\\
\, \e_3 &  &  & . & .\\
\, \e_4 &  &  &  & .\\
\end{tabular}$

\qquad
\begin{tabular}[t] {K K} \\[-10pt]
\end{tabular}
} }

\bigskip
\vbox{ {\bf[4, 2]}
\medskip

\hbox{

$\begin{tabular}[t]{L| C  C  C  C } \ &\ \e_1  &\e_2   &\e_3   &\e_4 \\
\hline
\, \e_1 & . & . & . & .\\
\, \e_2 &  & . & . & \, \e_1\\
\, \e_3 &  &  & . & \, \e_2\\
\, \e_4 &  &  &  & .\\
\end{tabular}$

\qquad
\begin{tabular}[t] {K K} \\[-10pt]
Win: &[4, 1]\\
\end{tabular}
} }

\bigskip

\section{Dimension Five}

\vbox{ {\bf[5, 1]}
\medskip

\hbox{

$\begin{tabular}[t]{L| C  C  C  C  C } \ &\ \e_1  &\e_2   &\e_3
&\e_4 &\e_5 \\ \hline
\, \e_1 & . & . & . & . & .\\
\, \e_2 &  & . & \, \e_1 & . & .\\
\, \e_3 &  &  & . & . & .\\
\, \e_4 &  &  &  & . & .\\
\, \e_5 &  &  &  &  & .\\
\end{tabular}$

\qquad
\begin{tabular}[t] {K K} \\[-10pt]
\end{tabular}
} }

\bigskip
\vbox{ {\bf[5, 2]}
\medskip

\hbox{

$\begin{tabular}[t]{L| C  C  C  C  C } \ &\ \e_1  &\e_2   &\e_3
&\e_4 &\e_5 \\ \hline
\, \e_1 & . & . & . & . & .\\
\, \e_2 &  & . & . & \, \e_1 & .\\
\, \e_3 &  &  & . & \, \e_2 & .\\
\, \e_4 &  &  &  & . & .\\
\, \e_5 &  &  &  &  & .\\
\end{tabular}$

\qquad
\begin{tabular}[t] {K K} \\[-10pt]
\end{tabular}
} }

\bigskip
\vbox{ {\bf[5, 3]}
\medskip

\hbox{

$\begin{tabular}[t]{L| C  C  C  C  C } \ &\ \e_1  &\e_2   &\e_3
&\e_4 &\e_5 \\ \hline
\, \e_1 & . & . & . & . & .\\
\, \e_2 &  & . & . & . & .\\
\, \e_3 &  &  & . & . & \, \e_1\\
\, \e_4 &  &  &  & . & \, \e_2\\
\, \e_5 &  &  &  &  & .\\
\end{tabular}$

\qquad
\begin{tabular}[t] {K K} \\[-10pt]
Win: &[5, 1]\\
Mubar: &[5, 1]\\
\end{tabular}
} }

\bigskip
\vbox{ {\bf[5, 4]}
\medskip

\hbox{

$\begin{tabular}[t]{L| C  C  C  C  C } \ &\ \e_1  &\e_2   &\e_3
&\e_4 &\e_5 \\ \hline
\, \e_1 & . & . & . & . & .\\
\, \e_2 &  & . & . & . & \, \e_1\\
\, \e_3 &  &  & . & . & \, \e_2\\
\, \e_4 &  &  &  & . & \, \e_3\\
\, \e_5 &  &  &  &  & .\\
\end{tabular}$

\qquad
\begin{tabular}[t] {K K} \\[-10pt]
Win: &[5, 2]\\
Mubar: &[5, 2]\\
\end{tabular}
} }

\bigskip
\vbox{ {\bf[5, 5]}
\medskip

\hbox{

$\begin{tabular}[t]{L| C  C  C  C  C } \ &\ \e_1  &\e_2   &\e_3
&\e_4 &\e_5 \\ \hline
\, \e_1 & . & . & . & . & .\\
\, \e_2 &  & . & . & \, \e_3 & \, \e_1\\
\, \e_3 &  &  & . & . & .\\
\, \e_4 &  &  &  & . & \, \e_2\\
\, \e_5 &  &  &  &  & .\\
\end{tabular}$

\qquad
\begin{tabular}[t] {K K} \\[-10pt]
Win: &[5, 3]\\
Mubar: &[5, 3]\\
\end{tabular}
} }

\bigskip
\vbox{ {\bf[5, 6]}
\medskip

\hbox{

$\begin{tabular}[t]{L| C  C  C  C  C } \ &\ \e_1  &\e_2   &\e_3
&\e_4 &\e_5 \\ \hline
\, \e_1 & . & . & . & . & .\\
\, \e_2 &  & . & . & \, \e_1 & .\\
\, \e_3 &  &  & . & . & \, \e_1\\
\, \e_4 &  &  &  & . & .\\
\, \e_5 &  &  &  &  & .\\
\end{tabular}$

\qquad
\begin{tabular}[t] {K K} \\[-10pt]
Win: &[5, 4]\\
Mubar: &[5, 4]\\
\end{tabular}
} }

\bigskip
\vbox{ {\bf[5, 7]}
\medskip

\hbox{

$\begin{tabular}[t]{L| C  C  C  C  C } \ &\ \e_1  &\e_2   &\e_3
&\e_4 &\e_5 \\ \hline
\, \e_1 & . & . & . & . & .\\
\, \e_2 &  & . & . & . & \, \e_1\\
\, \e_3 &  &  & . & \, \e_1 & \, \e_2\\
\, \e_4 &  &  &  & . & .\\
\, \e_5 &  &  &  &  & .\\
\end{tabular}$

\qquad
\begin{tabular}[t] {K K} \\[-10pt]
Win: &[5, 5]\\
Mubar: &[5, 5]\\
\end{tabular}
} }

\bigskip
\vbox{ {\bf[5, 8]}
\medskip

\hbox{

$\begin{tabular}[t]{L| C  C  C  C  C } \ &\ \e_1  &\e_2   &\e_3
&\e_4 &\e_5 \\ \hline
\, \e_1 & . & . & . & . & .\\
\, \e_2 &  & . & . & . & \, \e_1\\
\, \e_3 &  &  & . & \, \e_1 & \, \e_2\\
\, \e_4 &  &  &  & . & \, \e_3\\
\, \e_5 &  &  &  &  & .\\
\end{tabular}$

\qquad
\begin{tabular}[t] {K K} \\[-10pt]
Win: &[5, 6]\\
Mubar: &[5, 6]\\
\end{tabular}
} }

\bigskip

\section{Dimension Six}

\vbox{ {\bf[6, 1]}
\medskip

\hbox{

$\begin{tabular}[t]{L| C  C  C  C  C  C } \ &\ \e_1  &\e_2   &\e_3
&\e_4   &\e_5   &\e_6 \\ \hline
\, \e_1 & . & . & . & . & . & .\\
\, \e_2 &  & . & . & . & . & .\\
\, \e_3 &  &  & . & . & . & \, \e_2\\
\, \e_4 &  &  &  & . & \, \e_2 & \, \e_3\\
\, \e_5 &  &  &  &  & . & \, \e_4\\
\, \e_6 &  &  &  &  &  & .\\
\end{tabular}$

\qquad
\begin{tabular}[t] {K K} \\[-10pt]
Gong: &[6, N6211]\\
\end{tabular}
} }

\bigskip
\vbox{ {\bf[6, 2]}
\medskip

\hbox{

$\begin{tabular}[t]{L| C  C  C  C  C  C } \ &\ \e_1  &\e_2   &\e_3
&\e_4   &\e_5   &\e_6 \\ \hline
\, \e_1 & . & . & . & . & . & .\\
\, \e_2 &  & . & . & . & . & .\\
\, \e_3 &  &  & . & . & . & \, \e_2\\
\, \e_4 &  &  &  & . & . & \, \e_3\\
\, \e_5 &  &  &  &  & . & \, \e_4\\
\, \e_6 &  &  &  &  &  & .\\
\end{tabular}$

\qquad
\begin{tabular}[t] {K K} \\[-10pt]
Gong: &[6, N637]\\
\end{tabular}
} }

\bigskip
\vbox{ {\bf[6, 3]}
\medskip

\hbox{

$\begin{tabular}[t]{L| C  C  C  C  C  C } \ &\ \e_1  &\e_2   &\e_3
&\e_4   &\e_5   &\e_6 \\ \hline
\, \e_1 & . & . & . & . & . & .\\
\, \e_2 &  & . & . & . & . & .\\
\, \e_3 &  &  & . & . & . & \, \e_2\\
\, \e_4 &  &  &  & . & \, \e_2 & .\\
\, \e_5 &  &  &  &  & . & \, \e_3\\
\, \e_6 &  &  &  &  &  & .\\
\end{tabular}$

\qquad
\begin{tabular}[t] {K K} \\[-10pt]
Gong: &[6, N638]\\
\end{tabular}
} }

\bigskip
\vbox{ {\bf[6, 4]}
\medskip

\hbox{

$\begin{tabular}[t]{L| C  C  C  C  C  C } \ &\ \e_1  &\e_2   &\e_3
&\e_4   &\e_5   &\e_6 \\ \hline
\, \e_1 & . & . & . & . & . & .\\
\, \e_2 &  & . & . & . & . & .\\
\, \e_3 &  &  & . & . & . & .\\
\, \e_4 &  &  &  & . & \, \e_2 & \, \e_3\\
\, \e_5 &  &  &  &  & . & \, \e_4\\
\, \e_6 &  &  &  &  &  & .\\
\end{tabular}$

\qquad
\begin{tabular}[t] {K K} \\[-10pt]
Gong: &[6, N639]\\
\end{tabular}
} }

\bigskip
\vbox{ {\bf[6, 5]}
\medskip

\hbox{

$\begin{tabular}[t]{L| C  C  C  C  C  C } \ &\ \e_1  &\e_2   &\e_3
&\e_4   &\e_5   &\e_6 \\ \hline
\, \e_1 & . & . & . & . & . & .\\
\, \e_2 &  & . & . & . & . & .\\
\, \e_3 &  &  & . & . & \, \e_2 & .\\
\, \e_4 &  &  &  & . & . & \, \e_2\\
\, \e_5 &  &  &  &  & . & .\\
\, \e_6 &  &  &  &  &  & .\\
\end{tabular}$

\qquad
\begin{tabular}[t] {K K} \\[-10pt]
Gong: &[6, N641]\\
\end{tabular}
} }

\bigskip
\vbox{ {\bf[6, 6]}
\medskip

\hbox{

$\begin{tabular}[t]{L| C  C  C  C  C  C } \ &\ \e_1  &\e_2   &\e_3
&\e_4   &\e_5   &\e_6 \\ \hline
\, \e_1 & . & . & . & . & . & .\\
\, \e_2 &  & . & . & . & . & .\\
\, \e_3 &  &  & . & . & . & .\\
\, \e_4 &  &  &  & . & . & \, \e_2\\
\, \e_5 &  &  &  &  & . & \, \e_3\\
\, \e_6 &  &  &  &  &  & .\\
\end{tabular}$

\qquad
\begin{tabular}[t] {K K} \\[-10pt]
Gong: &[6, N642]\\
\end{tabular}
} }

\bigskip
\vbox{ {\bf[6, 7]}
\medskip

\hbox{

$\begin{tabular}[t]{L| C  C  C  C  C  C } \ &\ \e_1  &\e_2   &\e_3
&\e_4   &\e_5   &\e_6 \\ \hline
\, \e_1 & . & . & . & . & . & .\\
\, \e_2 &  & . & . & . & . & .\\
\, \e_3 &  &  & . & . & . & .\\
\, \e_4 &  &  &  & . & . & \, \e_3\\
\, \e_5 &  &  &  &  & . & \, \e_4\\
\, \e_6 &  &  &  &  &  & .\\
\end{tabular}$

\qquad
\begin{tabular}[t] {K K} \\[-10pt]
Gong: &[6, N643]\\
\end{tabular}
} }

\bigskip
\vbox{ {\bf[6, 8]}
\medskip

\hbox{

$\begin{tabular}[t]{L| C  C  C  C  C  C } \ &\ \e_1  &\e_2   &\e_3
&\e_4   &\e_5   &\e_6 \\ \hline
\, \e_1 & . & . & . & . & . & .\\
\, \e_2 &  & . & . & . & . & .\\
\, \e_3 &  &  & . & \, \e_1 & . & .\\
\, \e_4 &  &  &  & . & . & .\\
\, \e_5 &  &  &  &  & . & \, \e_2\\
\, \e_6 &  &  &  &  &  & .\\
\end{tabular}$

\qquad
\begin{tabular}[t] {K K} \\[-10pt]
Gong: &[6, N644]\\
\end{tabular}
} }

\bigskip
\vbox{ {\bf[6, 9]}
\medskip

\hbox{

$\begin{tabular}[t]{L| C  C  C  C  C  C } \ &\ \e_1  &\e_2   &\e_3
&\e_4   &\e_5   &\e_6 \\ \hline
\, \e_1 & . & . & . & . & . & .\\
\, \e_2 &  & . & . & . & . & .\\
\, \e_3 &  &  & . & . & \, \e_2 & \, \e_1\\
\, \e_4 &  &  &  & . & -\, \e_1 & \, \e_2\\
\, \e_5 &  &  &  &  & . & .\\
\, \e_6 &  &  &  &  &  & .\\
\end{tabular}$

\qquad
\begin{tabular}[t] {K K} \\[-10pt]
Gong: &[6, N644a]\\
\end{tabular}
} }

\bigskip
\vbox{ {\bf[6, 10]}
\medskip

\hbox{

$\begin{tabular}[t]{L| C  C  C  C  C  C } \ &\ \e_1  &\e_2   &\e_3
&\e_4   &\e_5   &\e_6 \\ \hline
\, \e_1 & . & . & . & . & . & .\\
\, \e_2 &  & . & . & . & . & .\\
\, \e_3 &  &  & . & . & . & .\\
\, \e_4 &  &  &  & . & . & .\\
\, \e_5 &  &  &  &  & . & \, \e_4\\
\, \e_6 &  &  &  &  &  & .\\
\end{tabular}$

\qquad
\begin{tabular}[t] {K K} \\[-10pt]
Gong: &[6, N65]\\
\end{tabular}
} }

\bigskip
\vbox{ {\bf[6, 11]}
\medskip

\hbox{

$\begin{tabular}[t]{L| C  C  C  C  C  C } \ &\ \e_1  &\e_2   &\e_3
&\e_4   &\e_5   &\e_6 \\ \hline
\, \e_1 & . & . & . & . & . & .\\
\, \e_2 &  & . & . & . & . & \, \e_1\\
\, \e_3 &  &  & . & . & \, \e_1 & \, \e_2\\
\, \e_4 &  &  &  & . & \, \e_2 & \, \e_3\\
\, \e_5 &  &  &  &  & . & \, \e_4\\
\, \e_6 &  &  &  &  &  & .\\
\end{tabular}$

\qquad
\begin{tabular}[t] {K K} \\[-10pt]
Gong: &[6, N611]\\
\end{tabular}
} }

\bigskip
\vbox{ {\bf[6, 12]}
\medskip

\hbox{

$\begin{tabular}[t]{L| C  C  C  C  C  C } \ &\ \e_1  &\e_2   &\e_3
&\e_4   &\e_5   &\e_6 \\ \hline
\, \e_1 & . & . & . & . & . & .\\
\, \e_2 &  & . & . & . & \, \e_1 & .\\
\, \e_3 &  &  & . & -\, \e_1 & . & \, \e_2\\
\, \e_4 &  &  &  & . & \, \e_2 & \, \e_3\\
\, \e_5 &  &  &  &  & . & \, \e_4\\
\, \e_6 &  &  &  &  &  & .\\
\end{tabular}$

\qquad
\begin{tabular}[t] {K K} \\[-10pt]
Gong: &[6, N612]\\
\end{tabular}
} }

\bigskip
\vbox{ {\bf[6, 13]}
\medskip

\hbox{

$\begin{tabular}[t]{L| C  C  C  C  C  C } \ &\ \e_1  &\e_2   &\e_3
&\e_4   &\e_5   &\e_6 \\ \hline
\, \e_1 & . & . & . & . & . & .\\
\, \e_2 &  & . & . & . & . & \, \e_1\\
\, \e_3 &  &  & . & . & . & \, \e_2\\
\, \e_4 &  &  &  & . & \, \e_1 & \, \e_3\\
\, \e_5 &  &  &  &  & . & \, \e_4\\
\, \e_6 &  &  &  &  &  & .\\
\end{tabular}$

\qquad
\begin{tabular}[t] {K K} \\[-10pt]
Gong: &[6, N613]\\
\end{tabular}
} }

\bigskip
\vbox{ {\bf[6, 14]}
\medskip

\hbox{

$\begin{tabular}[t]{L| C  C  C  C  C  C } \ &\ \e_1  &\e_2   &\e_3
&\e_4   &\e_5   &\e_6 \\ \hline
\, \e_1 & . & . & . & . & . & .\\
\, \e_2 &  & . & . & . & \, \e_1 & .\\
\, \e_3 &  &  & . & . & . & \, \e_1\\
\, \e_4 &  &  &  & . & \, \e_1 & \, \e_3\\
\, \e_5 &  &  &  &  & . & \, \e_4\\
\, \e_6 &  &  &  &  &  & .\\
\end{tabular}$

\qquad
\begin{tabular}[t] {K K} \\[-10pt]
Gong: &[6, N614]\\
\end{tabular}
} }

\bigskip
\vbox{ {\bf[6, 15]}
\medskip

\hbox{

$\begin{tabular}[t]{L| C  C  C  C  C  C } \ &\ \e_1  &\e_2   &\e_3
&\e_4   &\e_5   &\e_6 \\ \hline
\, \e_1 & . & . & . & . & . & .\\
\, \e_2 &  & . & . & . & . & \, \e_1\\
\, \e_3 &  &  & . & . & . & \, \e_2\\
\, \e_4 &  &  &  & . & . & \, \e_3\\
\, \e_5 &  &  &  &  & . & \, \e_4\\
\, \e_6 &  &  &  &  &  & .\\
\end{tabular}$

\qquad
\begin{tabular}[t] {K K} \\[-10pt]
Gong: &[6, N621]\\
\end{tabular}
} }

\bigskip
\vbox{ {\bf[6, 16]}
\medskip

\hbox{

$\begin{tabular}[t]{L| C  C  C  C  C  C } \ &\ \e_1  &\e_2   &\e_3
&\e_4   &\e_5   &\e_6 \\ \hline
\, \e_1 & . & . & . & . & . & .\\
\, \e_2 &  & . & . & . & \, \e_1 & .\\
\, \e_3 &  &  & . & -\, \e_1 & . & \, \e_2\\
\, \e_4 &  &  &  & . & . & \, \e_3\\
\, \e_5 &  &  &  &  & . & \, \e_4\\
\, \e_6 &  &  &  &  &  & .\\
\end{tabular}$

\qquad
\begin{tabular}[t] {K K} \\[-10pt]
Gong: &[6, N622]\\
\end{tabular}
} }

\bigskip
\vbox{ {\bf[6, 17]}
\medskip

\hbox{

$\begin{tabular}[t]{L| C  C  C  C  C  C } \ &\ \e_1  &\e_2   &\e_3
&\e_4   &\e_5   &\e_6 \\ \hline
\, \e_1 & . & . & . & . & . & .\\
\, \e_2 &  & . & . & . & . & \, \e_1\\
\, \e_3 &  &  & . & -\, \e_1 & . & \, \e_2\\
\, \e_4 &  &  &  & . & \, \e_2 & .\\
\, \e_5 &  &  &  &  & . & \, \e_3\\
\, \e_6 &  &  &  &  &  & .\\
\end{tabular}$

\qquad
\begin{tabular}[t] {K K} \\[-10pt]
Gong: &[6, N623]\\
\end{tabular}
} }

\bigskip
\vbox{ {\bf[6, 18]}
\medskip

\hbox{

$\begin{tabular}[t]{L| C  C  C  C  C  C } \ &\ \e_1  &\e_2   &\e_3
&\e_4   &\e_5   &\e_6 \\ \hline
\, \e_1 & . & . & . & . & . & .\\
\, \e_2 &  & . & . & . & \, \e_1 & .\\
\, \e_3 &  &  & . & . & . & \, \e_1\\
\, \e_4 &  &  &  & . & . & \, \e_3\\
\, \e_5 &  &  &  &  & . & \, \e_4\\
\, \e_6 &  &  &  &  &  & .\\
\end{tabular}$

\qquad
\begin{tabular}[t] {K K} \\[-10pt]
Gong: &[6, N624]\\
\end{tabular}
} }

\bigskip
\vbox{ {\bf[6, 19]}
\medskip

\hbox{

$\begin{tabular}[t]{L| C  C  C  C  C  C } \ &\ \e_1  &\e_2   &\e_3
&\e_4   &\e_5   &\e_6 \\ \hline
\, \e_1 & . & . & . & . & . & .\\
\, \e_2 &  & . & . & . & . & \, \e_1\\
\, \e_3 &  &  & . & . & \, \e_1 & .\\
\, \e_4 &  &  &  & . & \, \e_2 & \, \e_3\\
\, \e_5 &  &  &  &  & . & \, \e_4\\
\, \e_6 &  &  &  &  &  & .\\
\end{tabular}$

\qquad
\begin{tabular}[t] {K K} \\[-10pt]
Gong: &[6, N625]\\
\end{tabular}
} }

\bigskip
\vbox{ {\bf[6, 20]}
\medskip

\hbox{

$\begin{tabular}[t]{L| C  C  C  C  C  C } \ &\ \e_1  &\e_2   &\e_3
&\e_4   &\e_5   &\e_6 \\ \hline
\, \e_1 & . & . & . & . & . & .\\
\, \e_2 &  & . & . & . & -\, \e_1 & .\\
\, \e_3 &  &  & . & . & . & -\, \e_1\\
\, \e_4 &  &  &  & . & \, \e_2 & \, \e_3\\
\, \e_5 &  &  &  &  & . & \, \e_4\\
\, \e_6 &  &  &  &  &  & .\\
\end{tabular}$

\qquad
\begin{tabular}[t] {K K} \\[-10pt]
Gong: &[6, N625a]\\
\end{tabular}
} }

\bigskip
\vbox{ {\bf[6, 21]}
\medskip

\hbox{

$\begin{tabular}[t]{L| C  C  C  C  C  C } \ &\ \e_1  &\e_2   &\e_3
&\e_4   &\e_5   &\e_6 \\ \hline
\, \e_1 & . & . & . & . & . & .\\
\, \e_2 &  & . & . & \, \e_1 & . & .\\
\, \e_3 &  &  & . & . & . & \, \e_1\\
\, \e_4 &  &  &  & . & . & \, \e_2\\
\, \e_5 &  &  &  &  & . & \, \e_3\\
\, \e_6 &  &  &  &  &  & .\\
\end{tabular}$

\qquad
\begin{tabular}[t] {K K} \\[-10pt]
Gong: &[6, N626]\\
\end{tabular}
} }

\bigskip
\vbox{ {\bf[6, 22]}
\medskip

\hbox{

$\begin{tabular}[t]{L| C  C  C  C  C  C } \ &\ \e_1  &\e_2   &\e_3
&\e_4   &\e_5   &\e_6 \\ \hline
\, \e_1 & . & . & . & . & . & .\\
\, \e_2 &  & . & . & . & . & .\\
\, \e_3 &  &  & . & . & . & \, \e_2\\
\, \e_4 &  &  &  & . & \, \e_1 & \, \e_3\\
\, \e_5 &  &  &  &  & . & \, \e_4\\
\, \e_6 &  &  &  &  &  & .\\
\end{tabular}$

\qquad
\begin{tabular}[t] {K K} \\[-10pt]
Gong: &[6, N627]\\
\end{tabular}
} }

\bigskip
\vbox{ {\bf[6, 23]}
\medskip

\hbox{

$\begin{tabular}[t]{L| C  C  C  C  C  C } \ &\ \e_1  &\e_2   &\e_3
&\e_4   &\e_5   &\e_6 \\ \hline
\, \e_1 & . & . & . & . & . & .\\
\, \e_2 &  & . & . & . & . & .\\
\, \e_3 &  &  & . & . & \, \e_2 & \, \e_1\\
\, \e_4 &  &  &  & . & . & \, \e_2\\
\, \e_5 &  &  &  &  & . & \, \e_4\\
\, \e_6 &  &  &  &  &  & .\\
\end{tabular}$

\qquad
\begin{tabular}[t] {K K} \\[-10pt]
Gong: &[6, N628]\\
\end{tabular}
} }

\bigskip
\vbox{ {\bf[6, 24]}
\medskip

\hbox{

$\begin{tabular}[t]{L| C  C  C  C  C  C } \ &\ \e_1  &\e_2   &\e_3
&\e_4   &\e_5   &\e_6 \\ \hline
\, \e_1 & . & . & . & . & . & .\\
\, \e_2 &  & . & . & . & . & .\\
\, \e_3 &  &  & . & . & \, \e_1 & .\\
\, \e_4 &  &  &  & . & \, \e_1 & \, \e_2\\
\, \e_5 &  &  &  &  & . & \, \e_4\\
\, \e_6 &  &  &  &  &  & .\\
\end{tabular}$

\qquad
\begin{tabular}[t] {K K} \\[-10pt]
Gong: &[6, N629]\\
\end{tabular}
} }

\bigskip
\vbox{ {\bf[6, 25]}
\medskip

\hbox{

$\begin{tabular}[t]{L| C  C  C  C  C  C } \ &\ \e_1  &\e_2   &\e_3
&\e_4   &\e_5   &\e_6 \\ \hline
\, \e_1 & . & . & . & . & . & .\\
\, \e_2 &  & . & . & . & . & .\\
\, \e_3 &  &  & . & . & \, \e_2 & \, \e_1\\
\, \e_4 &  &  &  & . & -\, \e_1 & \, \e_2\\
\, \e_5 &  &  &  &  & . & \, \e_4\\
\, \e_6 &  &  &  &  &  & .\\
\end{tabular}$

\qquad
\begin{tabular}[t] {K K} \\[-10pt]
Gong: &[6, N629a]\\
\end{tabular}
} }

\bigskip
\vbox{ {\bf[6, 26]}
\medskip

\hbox{

$\begin{tabular}[t]{L| C  C  C  C  C  C } \ &\ \e_1  &\e_2   &\e_3
&\e_4   &\e_5   &\e_6 \\ \hline
\, \e_1 & . & . & . & . & . & .\\
\, \e_2 &  & . & . & . & . & .\\
\, \e_3 &  &  & . & . & \, \e_2 & .\\
\, \e_4 &  &  &  & . & \, \e_1 & \, \e_2\\
\, \e_5 &  &  &  &  & . & \, \e_4\\
\, \e_6 &  &  &  &  &  & .\\
\end{tabular}$

\qquad
\begin{tabular}[t] {K K} \\[-10pt]
Gong: &[6, N6210]\\
\end{tabular}
} }

\bigskip
\vbox{ {\bf[6, 27]}
\medskip

\hbox{

$\begin{tabular}[t]{L| C  C  C  C  C  C } \ &\ \e_1  &\e_2   &\e_3
&\e_4   &\e_5   &\e_6 \\ \hline
\, \e_1 & . & . & . & . & . & .\\
\, \e_2 &  & . & . & . & \, \e_1 & .\\
\, \e_3 &  &  & . & \, \e_1 & . & .\\
\, \e_4 &  &  &  & . & . & \, \e_2\\
\, \e_5 &  &  &  &  & . & \, \e_3\\
\, \e_6 &  &  &  &  &  & .\\
\end{tabular}$

\qquad
\begin{tabular}[t] {K K} \\[-10pt]
Gong: &[6, N631]\\
\end{tabular}
} }

\bigskip
\vbox{ {\bf[6, 28]}
\medskip

\hbox{

$\begin{tabular}[t]{L| C  C  C  C  C  C } \ &\ \e_1  &\e_2   &\e_3
&\e_4   &\e_5   &\e_6 \\ \hline
\, \e_1 & . & . & . & . & . & .\\
\, \e_2 &  & . & . & \, \e_1 & . & .\\
\, \e_3 &  &  & . & . & \, \e_1 & .\\
\, \e_4 &  &  &  & . & . & \, \e_2\\
\, \e_5 &  &  &  &  & . & \, \e_3\\
\, \e_6 &  &  &  &  &  & .\\
\end{tabular}$

\qquad
\begin{tabular}[t] {K K} \\[-10pt]
Gong: &[6, N631a]\\
\end{tabular}
} }

\bigskip
\vbox{ {\bf[6, 29]}
\medskip

\hbox{

$\begin{tabular}[t]{L| C  C  C  C  C  C } \ &\ \e_1  &\e_2   &\e_3
&\e_4   &\e_5   &\e_6 \\ \hline
\, \e_1 & . & . & . & . & . & .\\
\, \e_2 &  & . & \, \e_1 & . & . & .\\
\, \e_3 &  &  & . & . & . & .\\
\, \e_4 &  &  &  & . & . & \, \e_1\\
\, \e_5 &  &  &  &  & . & \, \e_4\\
\, \e_6 &  &  &  &  &  & .\\
\end{tabular}$

\qquad
\begin{tabular}[t] {K K} \\[-10pt]
Gong: &[6, N632]\\
\end{tabular}
} }

\bigskip
\vbox{ {\bf[6, 30]}
\medskip

\hbox{

$\begin{tabular}[t]{L| C  C  C  C  C  C } \ &\ \e_1  &\e_2   &\e_3
&\e_4   &\e_5   &\e_6 \\ \hline
\, \e_1 & . & . & . & . & . & .\\
\, \e_2 &  & . & . & . & . & .\\
\, \e_3 &  &  & . & . & . & \, \e_1\\
\, \e_4 &  &  &  & . & \, \e_2 & .\\
\, \e_5 &  &  &  &  & . & \, \e_4\\
\, \e_6 &  &  &  &  &  & .\\
\end{tabular}$

\qquad
\begin{tabular}[t] {K K} \\[-10pt]
Gong: &[6, N633]\\
\end{tabular}
} }

\bigskip
\vbox{ {\bf[6, 31]}
\medskip

\hbox{

$\begin{tabular}[t]{L| C  C  C  C  C  C } \ &\ \e_1  &\e_2   &\e_3
&\e_4   &\e_5   &\e_6 \\ \hline
\, \e_1 & . & . & . & . & . & .\\
\, \e_2 &  & . & . & . & . & .\\
\, \e_3 &  &  & . & . & \, \e_1 & .\\
\, \e_4 &  &  &  & . & \, \e_2 & .\\
\, \e_5 &  &  &  &  & . & \, \e_4\\
\, \e_6 &  &  &  &  &  & .\\
\end{tabular}$

\qquad
\begin{tabular}[t] {K K} \\[-10pt]
Gong: &[6, N634]\\
\end{tabular}
} }

\bigskip
\vbox{ {\bf[6, 32]}
\medskip

\hbox{

$\begin{tabular}[t]{L| C  C  C  C  C  C } \ &\ \e_1  &\e_2   &\e_3
&\e_4   &\e_5   &\e_6 \\ \hline
\, \e_1 & . & . & . & . & . & .\\
\, \e_2 &  & . & . & . & . & .\\
\, \e_3 &  &  & . & . & . & \, \e_1\\
\, \e_4 &  &  &  & . & \, \e_1 & .\\
\, \e_5 &  &  &  &  & . & \, \e_2\\
\, \e_6 &  &  &  &  &  & .\\
\end{tabular}$

\qquad
\begin{tabular}[t] {K K} \\[-10pt]
Gong: &[6, N635]\\
\end{tabular}
} }

\bigskip
\vbox{ {\bf[6, 33]}
\medskip

\hbox{

$\begin{tabular}[t]{L| C  C  C  C  C  C } \ &\ \e_1  &\e_2   &\e_3
&\e_4   &\e_5   &\e_6 \\ \hline
\, \e_1 & . & . & . & . & . & .\\
\, \e_2 &  & . & . & . & . & .\\
\, \e_3 &  &  & . & . & . & .\\
\, \e_4 &  &  &  & . & \, \e_1 & \, \e_2\\
\, \e_5 &  &  &  &  & . & \, \e_3\\
\, \e_6 &  &  &  &  &  & .\\
\end{tabular}$

\qquad
\begin{tabular}[t] {K K} \\[-10pt]
Gong: &[6, N636it]\\
\end{tabular}
} }

\bigskip

%
\setlength{\baselineskip}{11pt}

%
\setlength{\baselineskip}{22pt}

%
%
\onecolumn \setlength{\parindent}{40 pt}
\chapter{MULTIPLICATION TABLES}\label{a2}

In this appendix, we will list all of the isomorphism classes of two
through seven dimensional real solvable indecomposable Lie algebras
with codimension one nilradicals. Before we list them, however, we
will describe the numbering system to the reader.  Any given Lie
algebra on the list will be given a number sequence.  For example,
$$[5,[4,1],3,2].$$  The first number in the list corresponds to the
dimension of the Lie algebra; our example is a five dimensional
algebra.  The second list corresponds to the algebra's nilradical.
This first number is the dimension of the nilradical and the second
number corresponds to the numbering of the nilradicals given in
Appendix \ref{a1} with the convention that as the abelian nilradical
was given no number, we will number it with a 0. Our example has the
first four dimensional non-abelian nilradical.  The third number
corresponds to the number of the parent case.  For instance, if
there was a $2 \times 2$ block that was moved into real Jordan
canonical form, it would create three ``parent'' cases, one for each
possible real Jordan canonical form. Our example is in the third
parent case. Finally, the fourth number is the number of the algebra
produced in that parent case. To put it all together, our example is
the second algebra that came from the third parent case of the first
non-abelian nilradical of a five dimensional Lie algebra.

Wherever possible, we have indicated the appropriate reference to
the classification lists of both Winternitz and Mubarakzyanov
\cite{winternitz-1, mubar-1, mubar-2}.  As in Appendix \ref{a1},
Winternitz and Mubarakzyanov's names are abbreviated to Win and
Mubar respectively.

\twocolumn \setlength{\parindent}{0 pt}

\section{Dimension Two}

\vbox{ {\bf[2, [1, 0], 1, 1]}
\medskip

\hbox{

$
\begin{tabular}[t]{L| C  C } \ &\ \e_1  &\e_2 \\ \hline
\, \e_1 & . & \, \e_1\\
\, \e_2 &  & .\\
\end{tabular}$

\qquad
\begin{tabular}[t] {K K} \\[-10pt]
\end{tabular}
} }

\bigskip

\section{Dimension Three}

\vbox{ {\bf[3, [2, 0], 1, 1]}
\medskip

\hbox{

$\begin{tabular}[t]{L| C  C  C } \ &\ \e_1  &\e_2   &\e_3 \\ \hline
\, \e_1 & . & . & \, \e_1\\
\, \e_2 &  & . & a\, \e_2\\
\, \e_3 &  &  & .\\
\end{tabular}$

\qquad
\begin{tabular}[t] {K K} \\[-10pt]
Win: &[3, 3]\\
Win: &[3, 4]\\
Win: &[3, 5]\\
\end{tabular}
}
\begin{tabular}[t] {K K}
parameters:  $ [a]$ \\
$[[-1 \leq a, a \leq 1, a \neq 0]]$\\
\end{tabular}
\medskip
}

\bigskip
\vbox{ {\bf[3, [2, 0], 2, 1]}
\medskip

\hbox{

$\begin{tabular}[t]{L| C  C  C } \ &\ \e_1  &\e_2   &\e_3 \\ \hline
\, \e_1 & . & . & a\, \e_1-\, \e_2\\
\, \e_2 &  & . & \, \e_1+a\, \e_2\\
\, \e_3 &  &  & .\\
\end{tabular}$

\qquad
\begin{tabular}[t] {K K} \\[-10pt]
Win: &[3, 6]\\
Win: &[3, 7]\\
\end{tabular}
}
\begin{tabular}[t] {K K}
parameters:  $ [a]$ \\
$[[0 \leq a]]$\\
\end{tabular}
\medskip
}

\bigskip
\vbox{ {\bf[3, [2, 0], 3, 1]}
\medskip

\hbox{

$\begin{tabular}[t]{L| C  C  C } \ &\ \e_1  &\e_2   &\e_3 \\ \hline
\, \e_1 & . & . & \, \e_1\\
\, \e_2 &  & . & \, \e_1+\, \e_2\\
\, \e_3 &  &  & .\\
\end{tabular}$

\qquad
\begin{tabular}[t] {K K} \\[-10pt]
Win: &[3, 2]\\
\end{tabular}
} }

\bigskip

\section{Dimension Four}

\vbox{ {\bf[4, [3, 0], 1, 1]}
\medskip

\hbox{

$\begin{tabular}[t]{L| C  C  C  C } \ &\ \e_1  &\e_2   &\e_3   &\e_4 \\
\hline
\, \e_1 & . & . & . & \, \e_1\\
\, \e_2 &  & . & . & a\, \e_2\\
\, \e_3 &  &  & . & b\, \e_3\\
\, \e_4 &  &  &  & .\\
\end{tabular}$

\qquad
\begin{tabular}[t] {K K} \\[-10pt]
Win: &[4, 5]\\
\end{tabular}
}
\begin{tabular}[t] {K K}
parameters:  $ [a,b]$ \\
$[[b \leq a, -1 \leq b, a \leq 1, a \neq 0, b \neq 0]]$\\
\end{tabular}
\medskip
}

\bigskip
\vbox{ {\bf[4, [3, 0], 2, 1]}
\medskip

\hbox{

$\begin{tabular}[t]{L| C  C  C  C } \ &\ \e_1  &\e_2   &\e_3   &\e_4 \\
\hline
\, \e_1 & . & . & . & a\, \e_1-\, \e_2\\
\, \e_2 &  & . & . & \, \e_1+a\, \e_2\\
\, \e_3 &  &  & . & b\, \e_3\\
\, \e_4 &  &  &  & .\\
\end{tabular}$

\qquad
\begin{tabular}[t] {K K} \\[-10pt]
Win: &[4, 6]\\
\end{tabular}
}
\begin{tabular}[t] {K K}
parameters:  $ [a,b]$ \\
$[[0 \leq a, b \neq 0]]$\\
\end{tabular}
\medskip
}

\bigskip
\vbox{ {\bf[4, [3, 0], 3, 1]}
\medskip

\hbox{

$\begin{tabular}[t]{L| C  C  C  C } \ &\ \e_1  &\e_2   &\e_3   &\e_4 \\
\hline
\, \e_1 & . & . & . & a\, \e_1\\
\, \e_2 &  & . & . & \, \e_1+a\, \e_2\\
\, \e_3 &  &  & . & \, \e_3\\
\, \e_4 &  &  &  & .\\
\end{tabular}$

\qquad
\begin{tabular}[t] {K K} \\[-10pt]
Win: &[4, 2]\\
Win: &[4, 3]\\
\end{tabular}
}
\begin{tabular}[t] {K K}
parameters:  $ [a]$ \\
$[[]]$\\
\end{tabular}
\medskip
}

\bigskip
\vbox{ {\bf[4, [3, 0], 4, 1]}
\medskip

\hbox{

$\begin{tabular}[t]{L| C  C  C  C } \ &\ \e_1  &\e_2   &\e_3   &\e_4 \\
\hline
\, \e_1 & . & . & . & \, \e_1\\
\, \e_2 &  & . & . & \, \e_1+\, \e_2\\
\, \e_3 &  &  & . & \, \e_2+\, \e_3\\
\, \e_4 &  &  &  & .\\
\end{tabular}$

\qquad
\begin{tabular}[t] {K K} \\[-10pt]
Win: &[4, 4]\\
\end{tabular}
} }

\bigskip
\vbox{ {\bf[4, [3, 1], 1, 1]}
\medskip

\hbox{

$\begin{tabular}[t]{L| C  C  C  C } \ &\ \e_1  &\e_2   &\e_3   &\e_4 \\
\hline
\, \e_1 & . & . & . & (a+1)\, \e_1\\
\, \e_2 &  & . & \, \e_1 & \, \e_2\\
\, \e_3 &  &  & . & a\, \e_3\\
\, \e_4 &  &  &  & .\\
\end{tabular}$

\qquad
\begin{tabular}[t] {K K} \\[-10pt]
Win: &[4, 8]\\
Win: &[4, 9]\\
\end{tabular}
}
\begin{tabular}[t] {K K}
parameters:  $ [a]$ \\
$[[-1 \leq a, a \leq 1]]$\\
\end{tabular}
\medskip
}

\bigskip
\vbox{ {\bf[4, [3, 1], 2, 1]}
\medskip

\hbox{

$\begin{tabular}[t]{L| C  C  C  C } \ &\ \e_1  &\e_2   &\e_3   &\e_4 \\
\hline
\, \e_1 & . & . & . & 2\,a\, \e_1\\
\, \e_2 &  & . & \, \e_1 & a\, \e_2-\, \e_3\\
\, \e_3 &  &  & . & \, \e_2+a\, \e_3\\
\, \e_4 &  &  &  & .\\
\end{tabular}$

\qquad
\begin{tabular}[t] {K K} \\[-10pt]
Win: &[4, 10]\\
Win: &[4, 11]\\
\end{tabular}
}
\begin{tabular}[t] {K K}
parameters:  $ [a]$ \\
$[[0 \leq a]]$\\
\end{tabular}
\medskip
}

\bigskip
\vbox{ {\bf[4, [3, 1], 3, 1]}
\medskip

\hbox{

$\begin{tabular}[t]{L| C  C  C  C } \ &\ \e_1  &\e_2   &\e_3   &\e_4 \\
\hline
\, \e_1 & . & . & . & 2\, \e_1\\
\, \e_2 &  & . & \, \e_1 & \, \e_2\\
\, \e_3 &  &  & . & \, \e_2+\, \e_3\\
\, \e_4 &  &  &  & .\\
\end{tabular}$

\qquad
\begin{tabular}[t] {K K} \\[-10pt]
Win: &[4, 7]\\
\end{tabular}
} }

\bigskip

\section{Dimension Five}

\vbox{ {\bf[5, [4, 0], 1, 1]}
\medskip

\hbox{

$\begin{tabular}[t]{L| C  C  C  C  C } \ &\ \e_1  &\e_2   &\e_3
&\e_4 &\e_5 \\ \hline
\, \e_1 & . & . & . & . & \, \e_1\\
\, \e_2 &  & . & . & . & a\, \e_2\\
\, \e_3 &  &  & . & . & b\, \e_3\\
\, \e_4 &  &  &  & . & c\, \e_4\\
\, \e_5 &  &  &  &  & .\\
\end{tabular}$

\qquad
\begin{tabular}[t] {K K} \\[-10pt]
Win: &[5, 7]\\
Win: &[5, 8]\\
Mubar: &[5, 7]\\
\end{tabular}
}
\begin{tabular}[t] {K K}
parameters:  $ [a,b,c]$ \\
$[[-1 \leq c, c \leq b, b \leq a, a \leq 1, a \neq 0, b \neq 0, c \neq 0]]$\\
\end{tabular}
\medskip
}

\bigskip
\vbox{ {\bf[5, [4, 0], 2, 1]}
\medskip

\hbox{

$\begin{tabular}[t]{L| C  C  C  C  C } \ &\ \e_1  &\e_2   &\e_3
&\e_4 &\e_5 \\ \hline
\, \e_1 & . & . & . & . & a\, \e_1-\, \e_2\\
\, \e_2 &  & . & . & . & \, \e_1+a\, \e_2\\
\, \e_3 &  &  & . & . & b\, \e_3\\
\, \e_4 &  &  &  & . & c\, \e_4\\
\, \e_5 &  &  &  &  & .\\
\end{tabular}$

\qquad
\begin{tabular}[t] {K K} \\[-10pt]
Win: &[5, 13]\\
Win: &[5, 14]\\
Mubar: &[5, 13]\\
Mubar: &[5, 14]\\
\end{tabular}
}
\begin{tabular}[t] {K K}
parameters:  $ [a,b,c]$ \\
$[[0 \leq a, c \leq b, b \neq 0, c \neq 0]]$\\
\end{tabular}
\medskip
}

\bigskip
\vbox{ {\bf[5, [4, 0], 3, 1]}
\medskip

\hbox{

$\begin{tabular}[t]{L| C  C  C  C  C } \ &\ \e_1  &\e_2   &\e_3
&\e_4 &\e_5 \\ \hline
\, \e_1 & . & . & . & . & a\, \e_1\\
\, \e_2 &  & . & . & . & \, \e_1+a\, \e_2\\
\, \e_3 &  &  & . & . & \, \e_3\\
\, \e_4 &  &  &  & . & b\, \e_4\\
\, \e_5 &  &  &  &  & .\\
\end{tabular}$

\qquad
\begin{tabular}[t] {K K} \\[-10pt]
Win: &[5, 8]\\
Win: &[5, 9]\\
Mubar: &[5, 8]\\
Mubar: &[5, 9]\\
\end{tabular}
}
\begin{tabular}[t] {K K}
parameters:  $ [a,b]$ \\
$[[-1 \leq b, b \leq 1, b \neq 0]]$\\
\end{tabular}
\medskip
}

\bigskip
\vbox{ {\bf[5, [4, 0], 4, 1]}
\medskip

\hbox{

$\begin{tabular}[t]{L| C  C  C  C  C } \ &\ \e_1  &\e_2   &\e_3
&\e_4 &\e_5 \\ \hline
\, \e_1 & . & . & . & . & a\, \e_1-\, \e_2\\
\, \e_2 &  & . & . & . & \, \e_1+a\, \e_2\\
\, \e_3 &  &  & . & . & b\, \e_3-c\, \e_4\\
\, \e_4 &  &  &  & . & c\, \e_3+b\, \e_4\\
\, \e_5 &  &  &  &  & .\\
\end{tabular}$

\qquad
\begin{tabular}[t] {K K} \\[-10pt]
Win: &[5, 17]\\
Mubar: &[5, 17]\\
\end{tabular}
}
\begin{tabular}[t] {K K}
parameters:  $ [a,b,c]$ \\
$[[0 \leq a, 0 < c, c \leq 1]]$\\
\end{tabular}
\medskip
}

\bigskip
\vbox{ {\bf[5, [4, 0], 5, 1]}
\medskip

\hbox{

$\begin{tabular}[t]{L| C  C  C  C  C } \ &\ \e_1  &\e_2   &\e_3
&\e_4 &\e_5 \\ \hline
\, \e_1 & . & . & . & . & a\, \e_1\\
\, \e_2 &  & . & . & . & \, \e_1+a\, \e_2\\
\, \e_3 &  &  & . & . & b\, \e_3-\, \e_4\\
\, \e_4 &  &  &  & . & \, \e_3+b\, \e_4\\
\, \e_5 &  &  &  &  & .\\
\end{tabular}$

\qquad
\begin{tabular}[t] {K K} \\[-10pt]
Win: &[5, 14]\\
Win: &[5, 16]\\
Mubar: &[5, 14]\\
Mubar: &[5, 16]\\
\end{tabular}
}
\begin{tabular}[t] {K K}
parameters:  $ [a,b]$ \\
$[[0 \leq b]]$\\
\end{tabular}
\medskip
}

\bigskip
\vbox{ {\bf[5, [4, 0], 6, 1]}
\medskip

\hbox{

$\begin{tabular}[t]{L| C  C  C  C  C } \ &\ \e_1  &\e_2   &\e_3
&\e_4 &\e_5 \\ \hline
\, \e_1 & . & . & . & . & \, \e_1\\
\, \e_2 &  & . & . & . & \, \e_1+\, \e_2\\
\, \e_3 &  &  & . & . & a\, \e_3\\
\, \e_4 &  &  &  & . & \, \e_3+a\, \e_4\\
\, \e_5 &  &  &  &  & .\\
\end{tabular}$

\qquad
\begin{tabular}[t] {K K} \\[-10pt]
Win: &[5, 15]\\
Mubar: &[5, 15]\\
\end{tabular}
}
\begin{tabular}[t] {K K}
parameters:  $ [a]$ \\
$[[-1 \leq a, a \leq 1]]$\\
\end{tabular}
\medskip
}

\bigskip
\vbox{ {\bf[5, [4, 0], 7, 1]}
\medskip

\hbox{

$\begin{tabular}[t]{L| C  C  C  C  C } \ &\ \e_1  &\e_2   &\e_3
&\e_4 &\e_5 \\ \hline
\, \e_1 & . & . & . & . & a\, \e_1\\
\, \e_2 &  & . & . & . & \, \e_1+a\, \e_2\\
\, \e_3 &  &  & . & . & \, \e_2+a\, \e_3\\
\, \e_4 &  &  &  & . & \, \e_4\\
\, \e_5 &  &  &  &  & .\\
\end{tabular}$

\qquad
\begin{tabular}[t] {K K} \\[-10pt]
Win: &[5, 10]\\
Win: &[5, 11]\\
Mubar: &[5, 10]\\
Mubar: &[5, 11]\\
\end{tabular}
}
\begin{tabular}[t] {K K}
parameters:  $ [a]$ \\
$[[]]$\\
\end{tabular}
\medskip
}

\bigskip
\vbox{ {\bf[5, [4, 0], 8, 1]}
\medskip

\hbox{

$\begin{tabular}[t]{L| C  C  C  C  C } \ &\ \e_1  &\e_2   &\e_3
&\e_4 &\e_5 \\ \hline
\, \e_1 & . & . & . & . & a\, \e_1-\, \e_2\\
\, \e_2 &  & . & . & . & \, \e_1+a\, \e_2\\
\, \e_3 &  &  & . & . & \, \e_1+a\, \e_3-\, \e_4\\
\, \e_4 &  &  &  & . & \, \e_2+\, \e_3+a\, \e_4\\
\, \e_5 &  &  &  &  & .\\
\end{tabular}$

\qquad
\begin{tabular}[t] {K K} \\[-10pt]
Win: &[5, 18]\\
Mubar: &[5, 18]\\
\end{tabular}
}
\begin{tabular}[t] {K K}
parameters:  $ [a]$ \\
$[[0 \leq a]]$\\
\end{tabular}
\medskip
}

\bigskip
\vbox{ {\bf[5, [4, 0], 9, 1]}
\medskip

\hbox{

$\begin{tabular}[t]{L| C  C  C  C  C } \ &\ \e_1  &\e_2   &\e_3
&\e_4 &\e_5 \\ \hline
\, \e_1 & . & . & . & . & \, \e_1\\
\, \e_2 &  & . & . & . & \, \e_1+\, \e_2\\
\, \e_3 &  &  & . & . & \, \e_2+\, \e_3\\
\, \e_4 &  &  &  & . & \, \e_3+\, \e_4\\
\, \e_5 &  &  &  &  & .\\
\end{tabular}$

\qquad
\begin{tabular}[t] {K K} \\[-10pt]
Win: &[5, 12]\\
Mubar: &[5, 12]\\
\end{tabular}
} }

\bigskip
\vbox{ {\bf[5, [4, 1], 1, 1]}
\medskip

\hbox{

$\begin{tabular}[t]{L| C  C  C  C  C } \ &\ \e_1  &\e_2   &\e_3
&\e_4 &\e_5 \\ \hline
\, \e_1 & . & . & . & . & (a+1)\, \e_1\\
\, \e_2 &  & . & \, \e_1 & . & \, \e_2\\
\, \e_3 &  &  & . & . & a\, \e_3\\
\, \e_4 &  &  &  & . & b\, \e_4\\
\, \e_5 &  &  &  &  & .\\
\end{tabular}$

\qquad
\begin{tabular}[t] {K K} \\[-10pt]
Win: &[5, 19]\\
Mubar: &[5, 19]\\
\end{tabular}
}
\begin{tabular}[t] {K K}
parameters:  $ [a,b]$ \\
$[[-1 \leq a, a \leq 1, b \neq 0]]$\\
\end{tabular}
\medskip
}

\bigskip
\vbox{ {\bf[5, [4, 1], 1, 2]}
\medskip

\hbox{

$\begin{tabular}[t]{L| C  C  C  C  C } \ &\ \e_1  &\e_2   &\e_3
&\e_4 &\e_5 \\ \hline
\, \e_1 & . & . & . & . & (a+1)\, \e_1\\
\, \e_2 &  & . & \, \e_1 & . & \, \e_2+\, \e_4\\
\, \e_3 &  &  & . & . & a\, \e_3\\
\, \e_4 &  &  &  & . & \, \e_4\\
\, \e_5 &  &  &  &  & .\\
\end{tabular}$

\qquad
\begin{tabular}[t] {K K} \\[-10pt]
Win: &[5, 28]\\
Mubar: &[5, 28]\\
\end{tabular}
}
\begin{tabular}[t] {K K}
parameters:  $ [a]$ \\
$[[]]$\\
\end{tabular}
\medskip
}

\bigskip
\vbox{ {\bf[5, [4, 1], 1, 3]}
\medskip

\hbox{

$\begin{tabular}[t]{L| C  C  C  C  C } \ &\ \e_1  &\e_2   &\e_3
&\e_4 &\e_5 \\ \hline
\, \e_1 & . & . & . & . & \, \e_1\\
\, \e_2 &  & . & \, \e_1 & . & \, \e_4\\
\, \e_3 &  &  & . & . & \, \e_3\\
\, \e_4 &  &  &  & . & .\\
\, \e_5 &  &  &  &  & .\\
\end{tabular}$

\qquad
\begin{tabular}[t] {K K} \\[-10pt]
Win: &[5, 29]\\
Mubar: &[5, 29]\\
\end{tabular}
} }

\bigskip
\vbox{ {\bf[5, [4, 1], 1, 4]}
\medskip

\hbox{

$\begin{tabular}[t]{L| C  C  C  C  C } \ &\ \e_1  &\e_2   &\e_3
&\e_4 &\e_5 \\ \hline
\, \e_1 & . & . & . & . & (a+1)\, \e_1\\
\, \e_2 &  & . & \, \e_1 & . & \, \e_2\\
\, \e_3 &  &  & . & . & a\, \e_3\\
\, \e_4 &  &  &  & . & \, \e_1+(a+1)\, \e_4\\
\, \e_5 &  &  &  &  & .\\
\end{tabular}$

\qquad
\begin{tabular}[t] {K K} \\[-10pt]
Win: &[5, 20]\\
Mubar: &[5, 20]\\
\end{tabular}
}
\begin{tabular}[t] {K K}
parameters:  $ [a]$ \\
$[[-1 \leq a, a \leq 1]]$\\
\end{tabular}
\medskip
}

\bigskip
\vbox{ {\bf[5, [4, 1], 1, 5]}
\medskip

\hbox{

$\begin{tabular}[t]{L| C  C  C  C  C } \ &\ \e_1  &\e_2   &\e_3
&\e_4 &\e_5 \\ \hline
\, \e_1 & . & . & . & . & \, \e_1\\
\, \e_2 &  & . & \, \e_1 & . & \, \e_2+\, \e_4\\
\, \e_3 &  &  & . & . & .\\
\, \e_4 &  &  &  & . & \, \e_1+\, \e_4\\
\, \e_5 &  &  &  &  & .\\
\end{tabular}$

\qquad
\begin{tabular}[t] {K K} \\[-10pt]
Win: &[5, 27]\\
Mubar: &[5, 27]\\
\end{tabular}
} }

\bigskip
\vbox{ {\bf[5, [4, 1], 2, 1]}
\medskip

\hbox{

$\begin{tabular}[t]{L| C  C  C  C  C } \ &\ \e_1  &\e_2   &\e_3
&\e_4 &\e_5 \\ \hline
\, \e_1 & . & . & . & . & 2\,a\, \e_1\\
\, \e_2 &  & . & \, \e_1 & . & a\, \e_2-\, \e_3\\
\, \e_3 &  &  & . & . & \, \e_2+a\, \e_3\\
\, \e_4 &  &  &  & . & b\, \e_4\\
\, \e_5 &  &  &  &  & .\\
\end{tabular}$

\qquad
\begin{tabular}[t] {K K} \\[-10pt]
Win: &[5, 25]\\
Mubar: &[5, 25]\\
\end{tabular}
}
\begin{tabular}[t] {K K}
parameters:  $ [a,b]$ \\
$[[0 \leq a, b \neq 0]]$\\
\end{tabular}
\medskip
}

\bigskip
\vbox{ {\bf[5, [4, 1], 2, 2]}
\medskip

\hbox{

$\begin{tabular}[t]{L| C  C  C  C  C } \ &\ \e_1  &\e_2   &\e_3
&\e_4 &\e_5 \\ \hline
\, \e_1 & . & . & . & . & 2\,a\, \e_1\\
\, \e_2 &  & . & \, \e_1 & . & a\, \e_2-\, \e_3\\
\, \e_3 &  &  & . & . & \, \e_2+a\, \e_3\\
\, \e_4 &  &  &  & . & \, \e_1+2\,a\, \e_4\\
\, \e_5 &  &  &  &  & .\\
\end{tabular}$

\qquad
\begin{tabular}[t] {K K} \\[-10pt]
Win: &[5, 26]\\
Mubar: &[5, 26]\\
\end{tabular}
}
\begin{tabular}[t] {K K}
parameters:  $ [a]$ \\
$[[0 \leq a]]$\\
\end{tabular}
\medskip
}

\bigskip
\vbox{ {\bf[5, [4, 1], 3, 1]}
\medskip

\hbox{

$\begin{tabular}[t]{L| C  C  C  C  C } \ &\ \e_1  &\e_2   &\e_3
&\e_4 &\e_5 \\ \hline
\, \e_1 & . & . & . & . & 2\,a\, \e_1\\
\, \e_2 &  & . & \, \e_1 & . & a\, \e_2\\
\, \e_3 &  &  & . & . & \, \e_2+a\, \e_3\\
\, \e_4 &  &  &  & . & \, \e_4\\
\, \e_5 &  &  &  &  & .\\
\end{tabular}$

\qquad
\begin{tabular}[t] {K K} \\[-10pt]
Win: &[5, 22]\\
Win: &[5, 23]\\
Mubar: &[5, 22]\\
Mubar: &[5, 23]\\
\end{tabular}
}
\begin{tabular}[t] {K K}
parameters:  $ [a]$ \\
$[[]]$\\
\end{tabular}
\medskip
}

\bigskip
\vbox{ {\bf[5, [4, 1], 3, 2]}
\medskip

\hbox{

$\begin{tabular}[t]{L| C  C  C  C  C } \ &\ \e_1  &\e_2   &\e_3
&\e_4 &\e_5 \\ \hline
\, \e_1 & . & . & . & . & 2\, \e_1\\
\, \e_2 &  & . & \, \e_1 & . & \, \e_2\\
\, \e_3 &  &  & . & . & \, \e_2+\, \e_3\\
\, \e_4 &  &  &  & . & \, \e_1+2\, \e_4\\
\, \e_5 &  &  &  &  & .\\
\end{tabular}$

\qquad
\begin{tabular}[t] {K K} \\[-10pt]
Win: &[5, 24]\\
Mubar: &[5, 24]\\
\end{tabular}
} }

\bigskip
\vbox{ {\bf[5, [4, 1], 3, 3]}
\medskip

\hbox{

$\begin{tabular}[t]{L| C  C  C  C  C } \ &\ \e_1  &\e_2   &\e_3
&\e_4 &\e_5 \\ \hline
\, \e_1 & . & . & . & . & 2\, \e_1\\
\, \e_2 &  & . & \, \e_1 & . & \, \e_2+\, \e_4\\
\, \e_3 &  &  & . & . & \, \e_2+\, \e_3\\
\, \e_4 &  &  &  & . & \, \e_4\\
\, \e_5 &  &  &  &  & .\\
\end{tabular}$

\qquad
\begin{tabular}[t] {K K} \\[-10pt]
Win: &[5, 21]\\
Mubar: &[5, 21]\\
\end{tabular}
} }

\bigskip
\vbox{ {\bf[5, [4, 2], 1, 1]}
\medskip

\hbox{

$\begin{tabular}[t]{L| C  C  C  C  C } \ &\ \e_1  &\e_2   &\e_3
&\e_4 &\e_5 \\ \hline
\, \e_1 & . & . & . & . & (a+2)\, \e_1\\
\, \e_2 &  & . & . & \, \e_1 & (a+1)\, \e_2\\
\, \e_3 &  &  & . & \, \e_2 & a\, \e_3\\
\, \e_4 &  &  &  & . & \, \e_4\\
\, \e_5 &  &  &  &  & .\\
\end{tabular}$

\qquad
\begin{tabular}[t] {K K} \\[-10pt]
Win: &[5, 30]\\
Mubar: &[5, 30]\\
\end{tabular}
}
\begin{tabular}[t] {K K}
parameters:  $ [a]$ \\
$[[]]$\\
\end{tabular}
\medskip
}

\bigskip
\vbox{ {\bf[5, [4, 2], 1, 2]}
\medskip

\hbox{

$\begin{tabular}[t]{L| C  C  C  C  C } \ &\ \e_1  &\e_2   &\e_3
&\e_4 &\e_5 \\ \hline
\, \e_1 & . & . & . & . & \, \e_1\\
\, \e_2 &  & . & . & \, \e_1 & \, \e_2\\
\, \e_3 &  &  & . & \, \e_2 & \, \e_3\\
\, \e_4 &  &  &  & . & .\\
\, \e_5 &  &  &  &  & .\\
\end{tabular}$

\qquad
\begin{tabular}[t] {K K} \\[-10pt]
\end{tabular}
} }

\bigskip
\vbox{ {\bf[5, [4, 2], 1, 3]}
\medskip

\hbox{

$\begin{tabular}[t]{L| C  C  C  C  C } \ &\ \e_1  &\e_2   &\e_3
&\e_4 &\e_5 \\ \hline
\, \e_1 & . & . & . & . & 3\, \e_1\\
\, \e_2 &  & . & . & \, \e_1 & 2\, \e_2\\
\, \e_3 &  &  & . & \, \e_2 & \, \e_3\\
\, \e_4 &  &  &  & . & \, \e_3+\, \e_4\\
\, \e_5 &  &  &  &  & .\\
\end{tabular}$

\qquad
\begin{tabular}[t] {K K} \\[-10pt]
Win: &[5, 31]\\
Mubar: &[5, 31]\\
\end{tabular}
} }

\bigskip
\vbox{ {\bf[5, [4, 2], 1, 4]}
\medskip

\hbox{

$\begin{tabular}[t]{L| C  C  C  C  C } \ &\ \e_1  &\e_2   &\e_3
&\e_4 &\e_5 \\ \hline
\, \e_1 & . & . & . & . & \, \e_1\\
\, \e_2 &  & . & . & \, \e_1 & \, \e_2\\
\, \e_3 &  &  & . & \, \e_2 & \vep\, \e_1+\, \e_3\\
\, \e_4 &  &  &  & . & .\\
\, \e_5 &  &  &  &  & .\\
\end{tabular}$

\qquad
\begin{tabular}[t] {K K} \\[-10pt]
Win: &[5, 32]\\
Mubar: &[5, 32]\\
\end{tabular}
}
\begin{tabular}[t] {K K}
parameters:  $ [\vep]$ \\
$[[\vep^2 = 1]]$\\
\end{tabular}
\medskip
}

\bigskip

\section{Dimension Six}

\vbox{ {\bf[6, [5, 0], 1, 1]}
\medskip

\hbox{

$\begin{tabular}[t]{L| C  C  C  C  C  C } \ &\ \e_1  &\e_2   &\e_3
&\e_4   &\e_5   &\e_6 \\ \hline
\, \e_1 & . & . & . & . & . & \, \e_1\\
\, \e_2 &  & . & . & . & . & a\, \e_2\\
\, \e_3 &  &  & . & . & . & b\, \e_3\\
\, \e_4 &  &  &  & . & . & c\, \e_4\\
\, \e_5 &  &  &  &  & . & d\, \e_5\\
\, \e_6 &  &  &  &  &  & .\\
\end{tabular}$

\qquad
\begin{tabular}[t] {K K} \\[-10pt]
Mubar: &[6, 1]\\
\end{tabular}
}
\begin{tabular}[t] {K K}
parameters:  $ [a,b,c,d]$ \\
$[[-1 \leq d, d \leq c, c \leq b, b \leq a, a \leq 1, a \neq 0, b \neq 0,$ \\
 $c \neq 0, d \neq 0]]$\\
\end{tabular}
\medskip
}

\bigskip
\vbox{ {\bf[6, [5, 0], 2, 1]}
\medskip

\hbox{

$\begin{tabular}[t]{L| C  C  C  C  C  C } \ &\ \e_1  &\e_2   &\e_3
&\e_4   &\e_5   &\e_6 \\ \hline
\, \e_1 & . & . & . & . & . & a\, \e_1-\, \e_2\\
\, \e_2 &  & . & . & . & . & \, \e_1+a\, \e_2\\
\, \e_3 &  &  & . & . & . & b\, \e_3\\
\, \e_4 &  &  &  & . & . & c\, \e_4\\
\, \e_5 &  &  &  &  & . & d\, \e_5\\
\, \e_6 &  &  &  &  &  & .\\
\end{tabular}$

\qquad
\begin{tabular}[t] {K K} \\[-10pt]
Mubar: &[6, 8]\\
\end{tabular}
}
\begin{tabular}[t] {K K}
parameters:  $ [a,b,c,d]$ \\
$[[d \leq c, c \leq b, 0 \leq a, b \neq 0, c \neq 0, d \neq 0]]$\\
\end{tabular}
\medskip
}

\bigskip
\vbox{ {\bf[6, [5, 0], 3, 1]}
\medskip

\hbox{

$\begin{tabular}[t]{L| C  C  C  C  C  C } \ &\ \e_1  &\e_2   &\e_3
&\e_4   &\e_5   &\e_6 \\ \hline
\, \e_1 & . & . & . & . & . & a\, \e_1\\
\, \e_2 &  & . & . & . & . & \, \e_1+a\, \e_2\\
\, \e_3 &  &  & . & . & . & \, \e_3\\
\, \e_4 &  &  &  & . & . & b\, \e_4\\
\, \e_5 &  &  &  &  & . & c\, \e_5\\
\, \e_6 &  &  &  &  &  & .\\
\end{tabular}$

\qquad
\begin{tabular}[t] {K K} \\[-10pt]
Mubar: &[6, 2]\\
\end{tabular}
}
\begin{tabular}[t] {K K}
parameters:  $ [a,b,c]$ \\
$[[-1 \leq c, c \leq b, b \leq 1, b \neq 0, c \neq 0]]$\\
\end{tabular}
\medskip
}

\bigskip
\vbox{ {\bf[6, [5, 0], 4, 1]}
\medskip

\hbox{

$\begin{tabular}[t]{L| C  C  C  C  C  C } \ &\ \e_1  &\e_2   &\e_3
&\e_4   &\e_5   &\e_6 \\ \hline
\, \e_1 & . & . & . & . & . & (a-1)\, \e_1\\
\, \e_2 &  & . & . & . & . & \, \e_1+a\, \e_2\\
\, \e_3 &  &  & . & . & . & b\, \e_3-c\, \e_4\\
\, \e_4 &  &  &  & . & . & c\, \e_3+b\, \e_4\\
\, \e_5 &  &  &  &  & . & d\, \e_5\\
\, \e_6 &  &  &  &  &  & .\\
\end{tabular}$

\qquad
\begin{tabular}[t] {K K} \\[-10pt]
Mubar: &[6, 11]\\
\end{tabular}
}
\begin{tabular}[t] {K K}
parameters:  $ [a,b,c,d]$ \\
$[[0 \leq a, 0 < c, c \leq 1, d \neq 0]]$\\
\end{tabular}
\medskip
}

\bigskip
\vbox{ {\bf[6, [5, 0], 5, 1]}
\medskip

\hbox{

$\begin{tabular}[t]{L| C  C  C  C  C  C } \ &\ \e_1  &\e_2   &\e_3
&\e_4   &\e_5   &\e_6 \\ \hline
\, \e_1 & . & . & . & . & . & a\, \e_1\\
\, \e_2 &  & . & . & . & . & \, \e_1+a\, \e_2\\
\, \e_3 &  &  & . & . & . & b\, \e_3-\, \e_4\\
\, \e_4 &  &  &  & . & . & \, \e_3+b\, \e_4\\
\, \e_5 &  &  &  &  & . & c\, \e_5\\
\, \e_6 &  &  &  &  &  & .\\
\end{tabular}$

\qquad
\begin{tabular}[t] {K K} \\[-10pt]
Mubar: &[6, 9]\\
\end{tabular}
}
\begin{tabular}[t] {K K}
parameters:  $ [a,b,c]$ \\
$[[c \neq 0, 0 \leq b]]$\\
\end{tabular}
\medskip
}

\bigskip
\vbox{ {\bf[6, [5, 0], 6, 1]}
\medskip

\hbox{

$\begin{tabular}[t]{L| C  C  C  C  C  C } \ &\ \e_1  &\e_2   &\e_3
&\e_4   &\e_5   &\e_6 \\ \hline
\, \e_1 & . & . & . & . & . & a\, \e_1\\
\, \e_2 &  & . & . & . & . & \, \e_1+a\, \e_2\\
\, \e_3 &  &  & . & . & . & b\, \e_3\\
\, \e_4 &  &  &  & . & . & \, \e_3+b\, \e_4\\
\, \e_5 &  &  &  &  & . & \, \e_5\\
\, \e_6 &  &  &  &  &  & .\\
\end{tabular}$

\qquad
\begin{tabular}[t] {K K} \\[-10pt]
Mubar: &[6, 6]\\
\end{tabular}
}
\begin{tabular}[t] {K K}
parameters:  $ [a,b]$ \\
$[[b \leq a]]$\\
\end{tabular}
\medskip
}

\bigskip
\vbox{ {\bf[6, [5, 0], 7, 1]}
\medskip

\hbox{

$\begin{tabular}[t]{L| C  C  C  C  C  C } \ &\ \e_1  &\e_2   &\e_3
&\e_4   &\e_5   &\e_6 \\ \hline
\, \e_1 & . & . & . & . & . & a\, \e_1\\
\, \e_2 &  & . & . & . & . & \, \e_1+a\, \e_2\\
\, \e_3 &  &  & . & . & . & \, \e_2+a\, \e_3\\
\, \e_4 &  &  &  & . & . & \, \e_4\\
\, \e_5 &  &  &  &  & . & b\, \e_5\\
\, \e_6 &  &  &  &  &  & .\\
\end{tabular}$

\qquad
\begin{tabular}[t] {K K} \\[-10pt]
Mubar: &[6, 3]\\
\end{tabular}
}
\begin{tabular}[t] {K K}
parameters:  $ [a,b]$ \\
$[[-1 \leq b, b \leq 1, b \neq 0]]$\\
\end{tabular}
\medskip
}

\bigskip
\vbox{ {\bf[6, [5, 0], 8, 1]}
\medskip

\hbox{

$\begin{tabular}[t]{L| C  C  C  C  C  C } \ &\ \e_1  &\e_2   &\e_3
&\e_4   &\e_5   &\e_6 \\ \hline
\, \e_1 & . & . & . & . & . & a\, \e_1\\
\, \e_2 &  & . & . & . & . & \, \e_1+a\, \e_2\\
\, \e_3 &  &  & . & . & . & \, \e_2+a\, \e_3\\
\, \e_4 &  &  &  & . & . & b\, \e_4-\, \e_5\\
\, \e_5 &  &  &  &  & . & \, \e_4+b\, \e_5\\
\, \e_6 &  &  &  &  &  & .\\
\end{tabular}$

\qquad
\begin{tabular}[t] {K K} \\[-10pt]
Mubar: &[6, 10]\\
\end{tabular}
}
\begin{tabular}[t] {K K}
parameters:  $ [a,b]$ \\
$[[0 \leq b]]$\\
\end{tabular}
\medskip
}

\bigskip
\vbox{ {\bf[6, [5, 0], 9, 1]}
\medskip

\hbox{

$\begin{tabular}[t]{L| C  C  C  C  C  C } \ &\ \e_1  &\e_2   &\e_3
&\e_4   &\e_5   &\e_6 \\ \hline
\, \e_1 & . & . & . & . & . & \, \e_1\\
\, \e_2 &  & . & . & . & . & \, \e_1+\, \e_2\\
\, \e_3 &  &  & . & . & . & \, \e_2+\, \e_3\\
\, \e_4 &  &  &  & . & . & a\, \e_4\\
\, \e_5 &  &  &  &  & . & \, \e_4+a\, \e_5\\
\, \e_6 &  &  &  &  &  & .\\
\end{tabular}$

\qquad
\begin{tabular}[t] {K K} \\[-10pt]
Mubar: &[6, 10]\\
\end{tabular}
}
\begin{tabular}[t] {K K}
parameters:  $ [a]$ \\
$[[]]$\\
\end{tabular}
\medskip
}

\bigskip
\vbox{ {\bf[6, [5, 0], 9, 2]}
\medskip

\hbox{

$\begin{tabular}[t]{L| C  C  C  C  C  C } \ &\ \e_1  &\e_2   &\e_3
&\e_4   &\e_5   &\e_6 \\ \hline
\, \e_1 & . & . & . & . & . & .\\
\, \e_2 &  & . & . & . & . & \, \e_1\\
\, \e_3 &  &  & . & . & . & \, \e_2\\
\, \e_4 &  &  &  & . & . & \, \e_4\\
\, \e_5 &  &  &  &  & . & \, \e_4+\, \e_5\\
\, \e_6 &  &  &  &  &  & .\\
\end{tabular}$

\qquad
\begin{tabular}[t] {K K} \\[-10pt]
Mubar: &[6, 10]\\
\end{tabular}
} }

\bigskip
\vbox{ {\bf[6, [5, 0], 10, 1]}
\medskip

\hbox{

$\begin{tabular}[t]{L| C  C  C  C  C  C } \ &\ \e_1  &\e_2   &\e_3
&\e_4   &\e_5   &\e_6 \\ \hline
\, \e_1 & . & . & . & . & . & a\, \e_1-\, \e_2\\
\, \e_2 &  & . & . & . & . & \, \e_1+a\, \e_2\\
\, \e_3 &  &  & . & . & . & \, \e_1+a\, \e_3-\, \e_4\\
\, \e_4 &  &  &  & . & . & \, \e_2+\, \e_3+a\, \e_4\\
\, \e_5 &  &  &  &  & . & b\, \e_5\\
\, \e_6 &  &  &  &  &  & .\\
\end{tabular}$

\qquad
\begin{tabular}[t] {K K} \\[-10pt]
Mubar: &[6, 12]\\
\end{tabular}
}
\begin{tabular}[t] {K K}
parameters:  $ [a,b]$ \\
$[[0 \leq a, b \neq 0]]$\\
\end{tabular}
\medskip
}

\bigskip
\vbox{ {\bf[6, [5, 0], 11, 1]}
\medskip

\hbox{

$\begin{tabular}[t]{L| C  C  C  C  C  C } \ &\ \e_1  &\e_2   &\e_3
&\e_4   &\e_5   &\e_6 \\ \hline
\, \e_1 & . & . & . & . & . & a\, \e_1\\
\, \e_2 &  & . & . & . & . & \, \e_1+a\, \e_2\\
\, \e_3 &  &  & . & . & . & \, \e_2+a\, \e_3\\
\, \e_4 &  &  &  & . & . & \, \e_3+a\, \e_4\\
\, \e_5 &  &  &  &  & . & \, \e_5\\
\, \e_6 &  &  &  &  &  & .\\
\end{tabular}$

\qquad
\begin{tabular}[t] {K K} \\[-10pt]
Mubar: &[6, 4]\\
\end{tabular}
}
\begin{tabular}[t] {K K}
parameters:  $ [a]$ \\
$[[]]$\\
\end{tabular}
\medskip
}

\bigskip
\vbox{ {\bf[6, [5, 0], 12, 1]}
\medskip

\hbox{

$\begin{tabular}[t]{L| C  C  C  C  C  C } \ &\ \e_1  &\e_2   &\e_3
&\e_4   &\e_5   &\e_6 \\ \hline
\, \e_1 & . & . & . & . & . & \, \e_1\\
\, \e_2 &  & . & . & . & . & \, \e_1+\, \e_2\\
\, \e_3 &  &  & . & . & . & \, \e_2+\, \e_3\\
\, \e_4 &  &  &  & . & . & \, \e_3+\, \e_4\\
\, \e_5 &  &  &  &  & . & \, \e_4+\, \e_5\\
\, \e_6 &  &  &  &  &  & .\\
\end{tabular}$

\qquad
\begin{tabular}[t] {K K} \\[-10pt]
Mubar: &[6, 5]\\
\end{tabular}
} }

\bigskip
\vbox{ {\bf[6, [5, 1], 1, 1]}
\medskip

\hbox{

$\begin{tabular}[t]{L| C  C  C  C  C  C } \ &\ \e_1  &\e_2   &\e_3
&\e_4   &\e_5   &\e_6 \\ \hline
\, \e_1 & . & . & . & . & . & (a+b)\, \e_1\\
\, \e_2 &  & . & \, \e_1 & . & . & a\, \e_2\\
\, \e_3 &  &  & . & . & . & b\, \e_3\\
\, \e_4 &  &  &  & . & . & \, \e_4\\
\, \e_5 &  &  &  &  & . & c\, \e_5\\
\, \e_6 &  &  &  &  &  & .\\
\end{tabular}$

\qquad
\begin{tabular}[t] {K K} \\[-10pt]
Mubar: &[6, 13]\\
\end{tabular}
}
\begin{tabular}[t] {K K}
parameters:  $ [a,b,c]$ \\
$[[b \leq a, a^2+b^2 \neq 0, -1 \leq c, c \leq 1, c \neq 0]]$\\
\end{tabular}
\medskip
}

\bigskip
\vbox{ {\bf[6, [5, 1], 1, 2]}
\medskip

\hbox{

$\begin{tabular}[t]{L| C  C  C  C  C  C } \ &\ \e_1  &\e_2   &\e_3
&\e_4   &\e_5   &\e_6 \\ \hline
\, \e_1 & . & . & . & . & . & (a+b)\, \e_1\\
\, \e_2 &  & . & \, \e_1 & . & . & a\, \e_2\\
\, \e_3 &  &  & . & . & . & b\, \e_3\\
\, \e_4 &  &  &  & . & . & \, \e_4\\
\, \e_5 &  &  &  &  & . & \, \e_1+(a+b)\, \e_5\\
\, \e_6 &  &  &  &  &  & .\\
\end{tabular}$

\qquad
\begin{tabular}[t] {K K} \\[-10pt]
Mubar: &[6, 14]\\
\end{tabular}
}
\begin{tabular}[t] {K K}
parameters:  $ [a,b]$ \\
$[[b \leq a]]$\\
\end{tabular}
\medskip
}

\bigskip
\vbox{ {\bf[6, [5, 1], 1, 3]}
\medskip

\hbox{

$\begin{tabular}[t]{L| C  C  C  C  C  C } \ &\ \e_1  &\e_2   &\e_3
&\e_4   &\e_5   &\e_6 \\ \hline
\, \e_1 & . & . & . & . & . & (a+b)\, \e_1\\
\, \e_2 &  & . & \, \e_1 & . & . & a\, \e_2+\, \e_4\\
\, \e_3 &  &  & . & . & . & b\, \e_3\\
\, \e_4 &  &  &  & . & . & a\, \e_4\\
\, \e_5 &  &  &  &  & . & \, \e_5\\
\, \e_6 &  &  &  &  &  & .\\
\end{tabular}$

\qquad
\begin{tabular}[t] {K K} \\[-10pt]
Mubar: &[6, 17]\\
Mubar: &[6, 18]\\
\end{tabular}
}
\begin{tabular}[t] {K K}
parameters:  $ [a,b]$ \\
$[[]]$\\
\end{tabular}
\medskip
}

\bigskip
\vbox{ {\bf[6, [5, 1], 1, 4]}
\medskip

\hbox{

$\begin{tabular}[t]{L| C  C  C  C  C  C } \ &\ \e_1  &\e_2   &\e_3
&\e_4   &\e_5   &\e_6 \\ \hline
\, \e_1 & . & . & . & . & . & (a+1)\, \e_1\\
\, \e_2 &  & . & \, \e_1 & . & . & \, \e_2+\, \e_4\\
\, \e_3 &  &  & . & . & . & a\, \e_3\\
\, \e_4 &  &  &  & . & . & \, \e_4\\
\, \e_5 &  &  &  &  & . & \, \e_1+(a+1)\, \e_5\\
\, \e_6 &  &  &  &  &  & .\\
\end{tabular}$

\qquad
\begin{tabular}[t] {K K} \\[-10pt]
Mubar: &[6, 19]\\
\end{tabular}
}
\begin{tabular}[t] {K K}
parameters:  $ [a]$ \\
$[[]]$\\
\end{tabular}
\medskip
}

\bigskip
\vbox{ {\bf[6, [5, 1], 1, 5]}
\medskip

\hbox{

$\begin{tabular}[t]{L| C  C  C  C  C  C } \ &\ \e_1  &\e_2   &\e_3
&\e_4   &\e_5   &\e_6 \\ \hline
\, \e_1 & . & . & . & . & . & \, \e_1\\
\, \e_2 &  & . & \, \e_1 & . & . & \, \e_4\\
\, \e_3 &  &  & . & . & . & \, \e_3\\
\, \e_4 &  &  &  & . & . & .\\
\, \e_5 &  &  &  &  & . & \, \e_1+\, \e_5\\
\, \e_6 &  &  &  &  &  & .\\
\end{tabular}$

\qquad
\begin{tabular}[t] {K K} \\[-10pt]
\end{tabular}
} }

\bigskip
\vbox{ {\bf[6, [5, 1], 1, 6]}
\medskip

\hbox{

$\begin{tabular}[t]{L| C  C  C  C  C  C } \ &\ \e_1  &\e_2   &\e_3
&\e_4   &\e_5   &\e_6 \\ \hline
\, \e_1 & . & . & . & . & . & a\, \e_1\\
\, \e_2 &  & . & \, \e_1 & . & . & a\, \e_2+\, \e_4\\
\, \e_3 &  &  & . & . & . & .\\
\, \e_4 &  &  &  & . & . & \, \e_1+a\, \e_4\\
\, \e_5 &  &  &  &  & . & \, \e_5\\
\, \e_6 &  &  &  &  &  & .\\
\end{tabular}$

\qquad
\begin{tabular}[t] {K K} \\[-10pt]
Mubar: &[6, 17]\\
Mubar: &[6, 20]\\
\end{tabular}
}
\begin{tabular}[t] {K K}
parameters:  $ [a]$ \\
$[[]]$\\
\end{tabular}
\medskip
}

\bigskip
\vbox{ {\bf[6, [5, 1], 1, 7]}
\medskip

\hbox{

$\begin{tabular}[t]{L| C  C  C  C  C  C } \ &\ \e_1  &\e_2   &\e_3
&\e_4   &\e_5   &\e_6 \\ \hline
\, \e_1 & . & . & . & . & . & (a+1)\, \e_1\\
\, \e_2 &  & . & \, \e_1 & . & . & \, \e_2+\, \e_4\\
\, \e_3 &  &  & . & . & . & a\, \e_3+\, \e_5\\
\, \e_4 &  &  &  & . & . & \, \e_4\\
\, \e_5 &  &  &  &  & . & a\, \e_5\\
\, \e_6 &  &  &  &  &  & .\\
\end{tabular}$

\qquad
\begin{tabular}[t] {K K} \\[-10pt]
Mubar: &[6, 15]\\
\end{tabular}%
} }

\bigskip

\section{Dimension Seven}

\vbox{ {\bf[7, [6, 0], 1, 1]}
\medskip

\hbox{

$%
} }

\bigskip

%
\setlength{\baselineskip}{11pt}


\end{document}